\documentclass[12pt,authoryear]{article}
\newcommand{\blind}{1}
\usepackage{amssymb,amsthm}
\usepackage{amsmath}
\usepackage{natbib}
\usepackage{bibunits}
\usepackage{array}
\usepackage{graphicx}
\newcolumntype{C}[1]{>{\centering\arraybackslash}p{#1}}
\defaultbibliographystyle{plainnat}  % 或elsarticle-harv
\defaultbibliography{bibfile_main} % 默认文献库

\usepackage{graphicx,enumitem,subcaption,physics}
\usepackage{tikz}
\usetikzlibrary{calc}
\usetikzlibrary{positioning} 
\usetikzlibrary{arrows.meta}
\usetikzlibrary{shadows}
\usepackage{algorithm}   
\usepackage{algpseudocode} 

\usepackage{booktabs}
\usepackage{multirow}
\usepackage{setspace}
\newtheorem{theorem}{Theorem}[section]
\newtheorem{remark}{Remark}[section]
\newtheorem{definition}{Definition}[section]
\newtheorem{proposition}{Proposition}[section]
\newtheorem{assumption}{Assumption}[section]
\newtheorem{example}{Example}[section]
\newtheorem{lemma}{Lemma}[section]
\newtheorem{corollary}{Corollary}[section]

\usepackage{chngcntr}
\counterwithin{table}{section}
\counterwithin{algorithm}{section}
\counterwithin{figure}{section}
\numberwithin{equation}{section}
\newcommand{\sy}[1]{{\color{black}{#1}}}

\usepackage[colorlinks=true, linkcolor=blue, urlcolor=black, citecolor=blue]{hyperref}
\usepackage{xcolor} 

\newcommand{\magenta}[1]{\textcolor{black}{#1}}

\newcommand{\indep}{\perp\!\!\!\perp}

\allowdisplaybreaks[4]
\begin{document}		
	\if1\blind
	{
		\title{\vspace{-0.2cm} Double Machine Learning of Continuous Treatment Effects with \sy{Additive} Instrumental Variables}
		\author{Shuyuan Chen, Peng Zhang, Yifan Cui\thanks{The authors were partially supported by the National Key R\&D Program of China (2024YFA1015600), the National Natural Science Foundation of China (12471266 and U23A2064). Correspondence to \href{cuiyf@zju.edu.cn}{cuiyf@zju.edu.cn}.}}
		\date{Zhejiang University    \vspace{-1cm}} 
		\maketitle
	} \fi

	\if0\blind
	{
		\bigskip
		\bigskip
		\bigskip
		\begin{center}
			{\LARGE\bf Double Machine Learning of Continuous Treatment Effects with \sy{Additive} Instrumental Variables}
		\end{center}
		\medskip
	} \fi

	\begin{abstract}
		
		Estimating causal effects of continuous treatments is a common problem in practice, for example, in studying average dose-response functions. Classical analyses typically assume that all confounders are fully observed, whereas in real-world applications, unmeasured confounding often persists.
		In this article, we propose a novel framework for the identification of average dose-response functions using instrumental variables, thereby mitigating bias induced by unobserved confounders.
		We introduce the concept of a uniform regular weighting function and consider covering the treatment space with a finite collection of open sets. On each of these sets, such a weighting function exists, allowing us to identify the average dose-response function locally within the corresponding region.
		For estimation, we propose an augmented inverse probability weighted score for continuous treatments with instrumental variables under a debiased machine learning framework, and provide practical guidance for adaptively constructing regular weighting functions from the data. \sy{We also propose a falsification test for the additive instrumental variable condition and develop a testing procedure for assessing the validity of regular weighting functions.}
		We further establish the asymptotic properties of the resulting estimators based on kernel regression or empirical risk minimization \sy{as well as the theoretical validity of the proposed testing procedures}.
		Finally, we conduct both simulation and empirical studies to assess the finite-sample performance of the proposed methods.
	\end{abstract}
	
	\noindent
	{\it Keywords:}
	Average dose-response function,
	Continuous treatment,
	Debiased machine learning,
	Finite open cover,
	Instrumental variable,
	Kernel regression,
	Regular weighting function.
	
	\begin{bibunit} 
		\section{Introduction}
		Estimating causal effects of continuous treatments is a common problem in practice, for example, in studying dose-response relationships. This problem has been extensively investigated in the literature. For instance, \citet{hirano2004propensity,ImaiVanDyk2004,GalvaoWang2015} adapted generalized propensity score-based approaches to handle continuous treatments. \citet{DiazVanderLaan2013} considered global modeling of the average dose-response function (ADRF) and developed a doubly robust substitution estimator of the risk.
		
		More recently, \citet{kennedy2017non} propose a doubly robust influence function for estimating the ADRF using kernel regression techniques, whereas \citet{SemenovaChernozhukov2021,Colangelo16072025} develop doubly robust estimators for ADRFs within the debiased machine learning (DML) framework \citep{Chernozhukov2018}. 
		Meanwhile, \citet{bonvini2022fastconvergenceratesdoseresponse} summarize and generalize these approaches, allowing full exploitation of the smoothness of the target curve and achieving the oracle convergence rate even when the nuisance functions must be estimated; \sy{\citet{luedtke2024one} regard the ADRF as a parameter taking values in a Hilbert space, and construct a regularized EIF by expanding the parameter along an orthonormal basis in a manner analogous to sieve approximation; \citet{lim2025estimation} estimate the ADRF by constructing an orthogonal loss function that satisfies Neyman orthogonality.}
		Moreover, \citet{FongHazlettImai2018,HulingGreiferChen2023} study covariate balancing and independence weights for estimating the ADRF, \citet{ai2022estimation} estimate the counterfactual distribution and quantile functions in continuous treatment models, \citet{huang2025estimation} estimate the ratios of the quantile and tail mean at different treatment statuses, and \citet{huang2023nonparametric} estimate the continuous treatment effect with measurement error.

		All the above discussions on identifying continuous treatment effects rely on the assumption of no unmeasured confounders (NUC). Instrumental variable (IV) methods offer an effective approach to address unmeasured confounding. These methods rely on exogenous treatment variation generated by instruments that affect treatment assignment while remaining conditionally unrelated to latent confounding. Under the monotonicity assumption, IV methods recover causal effects for the subpopulation of compliers \citep{imbens1994identification,angrist1995two}. Building on these works, \citet{vytlacil2002independence,kennedy2019robust} propose estimators for the local IV effect curve, which captures treatment effects among individuals who comply when the instrument crosses a given threshold. \citet{newey2003instrumental,chen2012estimation} propose nonparametric IV models to identify the structural mean function under conditional moment restrictions.

		On the other hand, \citet{wang2018bounded,Hartwig2023} identify the average treatment effect (ATE) by imposing no-interaction assumptions between the IV and latent confounders in the treatment model, in settings where both the IV and the treatment are binary. Building on this framework, \citet{tchetgen2018marginal,michael2024instrumental} extend the approach to longitudinal data, while \citet{richardson2017nonparametric,cui2023instrumental,wang2023instrumental} apply this IV framework to identify and estimate causal Cox hazard ratios. 
		Recently, motivated by nonparametric instrumental variable methods, \citet{chen2025identificationdebiasedlearningcausal} propose an additive IV condition to identify mean potential outcomes with multi-categorical treatments and general IVs, and further identify the average treatment effect under a no unmeasured common effect modifier assumption \citep{cui2021semiparametric}. Independently and concurrently, similar results are obtained by \citet{dong2025marginalcausaleffectestimation}. 
		Additionally, \citet{liu2025multiplicativeinstrumentalvariablemodel} propose an alternative identification condition that rules out any multiplicative interaction between the IV and latent confounders in the treatment model,
		\citet{syrgkanis2019machine} estimate heterogeneous treatment effects using IVs within the DML framework, 
		and \citet{ye2023instrumented} estimate the ATE under a difference-in-difference IV framework.

		Based on the existing literature, we propose a general IV framework for identifying the ADRF with continuous treatments. Specifically, we propose the definitions of the regular weighting function (RWF) and additive IV, and further show that a finite collection of RWFs suffices for identification over any compact subset of the treatment space. We employ semiparametric theory \citep{vanderLaan2003, tsiatis2006semiparametric} to derive an augmented inverse probability weighting (AIPW) score function with the mixed-bias property \citep{rotnitzky2021characterization}, whose conditional expectation given the treatment exactly equals the ADRF. Moreover, inspired by \citet{SemenovaChernozhukov2021, bonvini2022fastconvergenceratesdoseresponse}, we develop a general cross-fitting algorithm under the DML framework to compute the AIPW score. Leveraging local linear kernel regression (LLKR) methods \citep{fan1996local,wasserman2006all,tsybakov2009introduction} and empirical risk minimization~\citep{BartlettBousquetMendelson2005,bonvini2022fastconvergenceratesdoseresponse,foster2023orthogonalstatisticallearning}, we estimate the ADRF nonparametrically and locally with the assistance of an additive IV. Overall, our approach accommodates unmeasured confounding, recovers the dose-response relationship nonparametrically, and extends the influence function of the ADRF in the NUC setting to the IV framework.

		We outline the structure of this article as follows. 
		Section~\ref{sec: model framework} introduces the foundational assumptions under the IV framework for continuous treatments, including a high-level IV relevance condition based on the $\chi^2$-divergence, the concepts of an RWF with local stability, a uniform RWF, and an additive IV condition for continuous treatments. 
		Section~\ref{sec: methodology} derives an AIPW score for the ADRF, describes a cross-fitting procedure for computing the AIPW score, and employs LLKR methods for nonparametric and local estimation of the ADRF. \sy{A falsification test is established for the additive IV condition, and a hypothesis testing procedure is developed for assessing the validity of the RWFs.} In Section~\ref{sec: practical guidance}, we provide practical guidance on how to construct a cover of a set such that each subset is associated with a corresponding uniform RWF.
		In Section~\ref{sec: asymptotic theory}, we establish the convergence rate and asymptotic normality of the proposed estimator. Section~\ref{sec: simulations} reports simulation results assessing the finite-sample performance of our method. Finally, Section~\ref{sec: empirical illustration} applies the proposed approach to investigate the effect of years of education on total annual earnings using data from the Job Training Partnership Act study.
		All scripts and data required to replicate the numerical results in our paper are available in \href{https://github.com/chensy123-sys/Continuous_treatment_IV}{https://github.com/chensy123-sys/Continuous\_treatment\_IV}.

		% IV relevance -> weakness of binary IV -> RWF definition -> existence of RWF -> URWF definition -> non-existence of global RWF -> finite open cover
		\section{The model framework}\label{sec: model framework}
		\subsection{Assumptions and notation}\label{subsec: assumptions and notation}
		We begin by introducing the basic notation used throughout the paper. \magenta{Let $d_Z$, $d_L$, and $d_U$ denote three positive integers. Let 
			$Y \in \mathcal{Y} \subseteq \mathbb{R}$ be the outcome of interest, 
			$A \in \mathcal{A} \subseteq \mathbb{R}$ a continuous treatment, 
			$Z \in \mathcal{Z} \subseteq \mathbb{R}^{d_Z}$ an IV, 
			which may be discrete or continuous, 
			$L \in \mathcal{L} \subseteq \mathbb{R}^{d_L}$ the observed confounders, 
			and $U \in \mathcal{U} \subseteq \mathbb{R}^{d_U}$ the unobserved confounders.}
		We assume that $\mathcal{A}$, the support of $A$, 
		is a closed subset of $\mathbb{R}$.
		
		Let $O := [L, Z, A, Y]$ denote the observed data. 
		The observations consist of independent and identically distributed samples 
		$\{O_i\}_{i=1}^n$. 
		Under the potential outcomes framework, 
		let $Y(a)$ denote the potential outcome that would be observed if $A = a$. 
		Let $p_{A \mid Z, U, L}(a \mid z, u, l)$ denote the conditional probability density function 
		of $A$ given $Z, U, L$, $p_A(a)$ its marginal density, and $P_A(a)$ its marginal distribution. 
		We use similar notation for other conditional and marginal distributions as well. 
		In addition, define $\mathcal{B}(a,h) := (a-h, a+h)$. For any subset $\mathcal{N}\subseteq\mathcal{A}$, let $\mathring{\mathcal{N}}$ denote the interior of $\mathcal{N}$, $\overline{\mathcal{N}}$ denote the closure of $\mathcal{N}$, and $f'(a)$ denote the derivative of $f(a)$.
		
		Throughout this paper, we assume that $p_{A|L}(a|L)>0$ for all $a\in\mathring{\mathcal{A}}$ almost surely, so that it shares the same support as the marginal distribution $P_A(a)$. 
		Next, we introduce several fundamental assumptions for the continuous treatment setting with IVs.
		\begin{assumption}[Consistency]\label{as: consistency}
			$Y=Y(A)$.
		\end{assumption}
		
		\begin{assumption}[Latent ignorability]\label{as: latent ignorability}
			For any $a\in\mathring{\mathcal{A}}$, $Y(a)\indep \{A,Z\}\mid U,L$.
		\end{assumption}
		
		\begin{assumption}[IV independence]\label{as: IV independence}
			$Z\indep U\mid L$.
		\end{assumption}
		
		\begin{assumption}[Continuity and uniform boundedness]\label{as: continuity}
			For all $a\in\mathring{\mathcal{A}}$, $p_{A| Z,L}(a| Z,L)$ is continuous with respect to $a$ almost surely. In addition, there exists $M(a)$, such that $p_{A| Z,L}(a| Z,L)\leq M(a)$ almost surely.
		\end{assumption}
		% \begin{assumption}[Continuity and uniform boundedness]\label{as: continuity}
			% 	For all $z\in\mathcal{Z}$ and $l\in\mathcal{L}$, $p_{A| Z,L}(a|z,l)$ is continuous for all $a\in\mathring{\mathcal{A}}$. For all $a\in\mathring{\mathcal{A}}$, there exists $M(a)$, such that $p_{A| Z,L}(a| Z,L)\leq M(a)$ and $\chi^2[p_{Z|A,L}(\cdot|a,L)\| p_{Z|L}(\cdot|L)]\leq M(a)$ almost surely.
			% \end{assumption}
		
		Assumption~\ref{as: consistency} is a standard consistency condition in causal inference, stating that the observed outcome equals the potential outcome corresponding to the received treatment. 
		Assumption~\ref{as: latent ignorability} states that, conditional on observed covariates and unmeasured confounders, $Y(a)$ is independent of both the treatment and the IVs, implying that $Z$ affects $Y$ only through $A$. 
		Assumption~\ref{as: IV independence} requires that the IV is independent of the unmeasured confounders. 
		Finally, Assumption~\ref{as: continuity} assumes that the conditional  densities of the treatment variable associated with the IVs and measured confounders do not vary sharply within the $\mathring{\mathcal{A}}$.

		Moreover, identification of the ADRF also requires that the treatment and instrument exhibit sufficient variability across the covariate space. 
		The following two assumptions ensure that the treatment is sufficiently positive across covariates and that the IV is relevant and informative for the treatment. 
		\begin{assumption}[Positivity]\label{as: positivity}
			For any $a\in\mathring{\mathcal{A}}$, there exists an $\epsilon_1(a)>0$ such that 
			$p_{A\mid L}(a\mid L)\geq \epsilon_1(a) \text{ almost surely.}$
		\end{assumption}
		\begin{assumption}[IV relevance]\label{as: IV relevance}
			For two probability distributions with densities $q$ and $p$, define the $\chi^2$-divergence between them as $\chi^2[q(\cdot)\| p(\cdot)]:=\int (\frac{q(z)}{p(z)} -1)^2 p(z)\mathrm{d} z$.
			For any $a\in\mathring{\mathcal{A}}$, there exists an $\epsilon_2(a)>0$ such that 
			$\chi^2[p_{Z|A,L}(\cdot|a,L)\| p_{Z|L}(\cdot|L)]\geq \epsilon_2(a)$ almost surely.
		\end{assumption}
		\begin{remark}
			For each fixed $a$, we impose positivity and IV relevance conditions uniformly over $L$.
			For identification, it suffices that for all $a \in \mathring{\mathcal{A}}$, 
			both $\chi^2\!\left[p_{Z\mid A,L}(\cdot\mid a,L)\,\|\, p_{Z\mid L}(\cdot\mid L)\right]$ 
			and $p_{A\mid L}(a\mid L)$ are almost surely nonzero.
		\end{remark}
		
		Assumption~\ref{as: positivity} states that $p_{A\mid L}(a\mid L)$ is uniformly bounded away from zero almost surely, ensuring the stability of the proposed estimator.  
		Assumption~\ref{as: IV relevance} can be viewed as an IV relevance condition, which requires that the instrument exerts a non-negligible influence on the treatment at any $a\in\mathring{\mathcal{A}}$.
		For notational convenience, whenever no confusion arises, we write $\epsilon_0(a):=\min\{\epsilon_1(a),\epsilon_2(a)\}$.
		We next provide an equivalent characterization of Assumption~\ref{as: IV relevance}, which also facilitates the graphical interpretation of Proposition~\ref{prop: weakness of binary IV}.
		\begin{proposition}\label{prop: positivity equivalence}
			Under Assumptions~\ref{as: continuity} and \ref{as: positivity}, Assumption~\ref{as: IV relevance} holds if and only if for any $a\in\mathring{\mathcal{A}}$, there exists an $\epsilon_3(a)>0$ such that 
			$\mathrm{Var}[p_{A|Z,L}(a\mid Z,L) \mid L] \geq \epsilon_3(a) \text{ almost surely.}$
		\end{proposition}
		%The statement in Proposition~\ref{prop: positivity equivalence} represents a strong form of IV relevance, corresponding to the positivity condition discussed in \citet{chen2025identificationdebiasedlearningcausal}.
		\magenta{Interestingly, in the continuous treatment setting, binary IVs generally fail to satisfy Assumption~\ref{as: IV relevance}, which is fundamentally different from the discrete treatment case \citep{chen2025identificationdebiasedlearningcausal}. This limitation is formalized in the following proposition.}
		\begin{proposition}[Weakness of a binary IV]\label{prop: weakness of binary IV}
			Under Assumptions~\ref{as: continuity} and \ref{as: positivity}, let $Z$ be a binary IV and $\mathcal{A}$ a closed interval. 
			Further assume that, for any $z\in\mathcal{Z}$ and $l\in\mathcal{L}$, the density functions $p_{A\mid Z,L}(a\mid z,l)$ share the same support $\mathcal{A}$. 
			Then, for any $l \in \mathcal{L}$ such that $p_{A\mid Z,L}(a\mid z,l)$ is continuous in $a$, there exists $a_0 \in \mathring{\mathcal{A}}$ such that $p_{A\mid Z,L}(a_0 \mid 0, l) = p_{A\mid Z,L}(a_0 \mid 1, l),$ 
            therefore $\chi^2[p_{Z|A,L}(\cdot|a_0,l)\| p_{Z|L}(\cdot|l)]=0$.
		\end{proposition}

		\begin{figure}
			\centering
			\begin{subfigure}{0.33\textwidth}
				\centering
				\includegraphics[width=\linewidth]{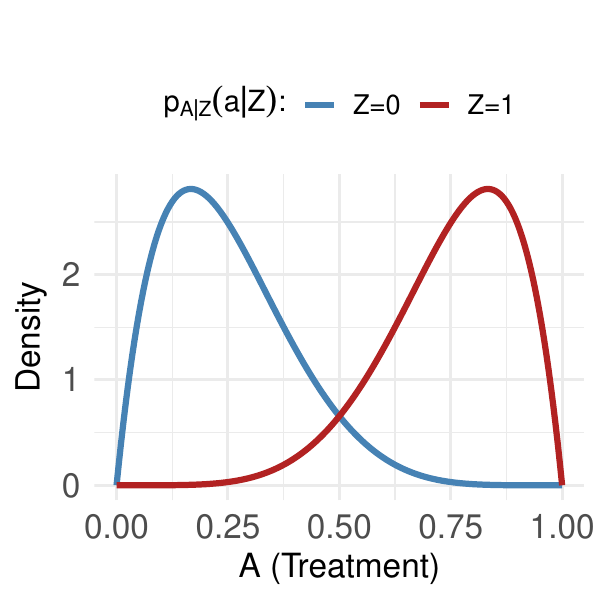}
				\caption{}
			\end{subfigure}
			\begin{subfigure}{0.33\textwidth}
				\centering
				\includegraphics[width=\linewidth]{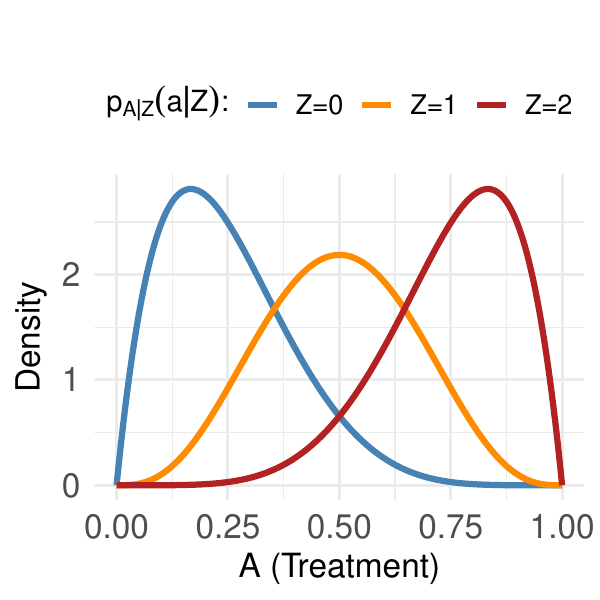}
				\caption{}
			\end{subfigure}
			
			\caption{
				(a) Density curves of $p_{A| Z}(a|0),p_{A| Z}(a|1)$ when $Z$ is binary.
				(b) Density curves of $p_{A| Z}(a|0),p_{A| Z}(a|1),p_{A| Z}(a|2)$ when $Z$ has three classes. 
			}
			\label{fig: binary IV}
		\end{figure}
		This proposition shows that, under the conditions of Proposition~\ref{prop: weakness of binary IV}, a binary instrument $Z$ cannot satisfy the IV relevance condition in Assumption~\ref{as: IV relevance} uniformly over $\mathring{\mathcal{A}}$.  
		For illustration, consider the case $L = \emptyset$. In Figure~\ref{fig: binary IV}(a), where $\mathcal{A} = [0,1]$, the two conditional densities $p_{A\mid Z}(a\mid 0)$ and $p_{A\mid Z}(a\mid 1)$ must intersect at least once within $\mathring{\mathcal{A}} = (0,1)$, leading to a violation of Assumption~\ref{as: IV relevance}.  
		In contrast, when $Z$ takes three categories, as shown in Figure~\ref{fig: binary IV}(b), the intersection points need not coincide, making Assumption~\ref{as: IV relevance} more likely to hold.

		\subsection{Regular weighting functions}\label{subsec: regular weighting function}
		In this subsection, we define a concept central to ADRF identification, closely tied to the IV relevance condition in Assumption~\ref{as: IV relevance}.
		For clarity, all results in this subsection rely only on Assumptions~\ref{as: continuity}--\ref{as: IV relevance}. 
		First, for any measurable function $\pi(Z,L)$, define
		$\kappa_{\pi}^o(A,L) := \mathbb{E}[\pi(Z,L) \mid A,L] - \mathbb{E}[\pi(Z,L) \mid L].$
		\magenta{
			We view $\kappa_{\pi}^o$ as a functional mapping $\pi$ to $\kappa_{\pi}^o(A,L)$, capturing how effectively $\pi(Z,L)$ leverages information in $Z$ to predict $A$. Values near zero indicate that $\pi$ uses little information from $Z$, while larger magnitudes suggest that $\pi$ is more informative. As this quantity later appears in the denominator of the identification formula, it must remain bounded away from zero. We thus use $\kappa_{\pi}^o$ to define a key concept for identifying the ADRF.}
		\begin{definition}\label{defn: RWF}
			A measurable function $\pi(Z, L): \mathcal{Z} \times \mathcal{L} \to\mathbb{R}$ is said to be a regular weighting function (RWF) for $A = a\in\mathring{\mathcal{A}}$ if it is uniformly bounded and there exists a constant $\epsilon_\pi(a) > 0$ such that
			$|\kappa_{\pi}^o(a,L)| \geq \epsilon_\pi(a)$ almost surely.
		\end{definition}
		\magenta{
			Intuitively, the existence of an RWF for $A = a$ implies that variation in the IVs has a non-negligible effect on the treatment at that level. Thus, an RWF can also be interpreted as a \textit{relevant} weighting function.} % This relevance is precisely why the definition of an RWF is significant when establishing identification.}
	It is therefore natural to ask whether such an RWF exists at a single interior point.  We present a necessary and sufficient condition for its existence.
	\begin{proposition}[Existence of RWFs]\label{prop: existence of RWFs}
		Under Assumption~\ref{as: continuity}, for any $a \in \mathring{\mathcal{A}}$, an RWF for $A=a$ exists if and only if 
		$\kappa_{\pi_a}^o(a,L)=\mathrm{Var}[\pi_a(Z,L)\mid L]$ is uniformly bounded above and below almost surely, 
		where $\pi_a(Z,L) := p_{A \mid Z,L}(a \mid Z,L) / p_{A \mid L}(a \mid L)$. 
		Moreover, if there exists an RWF for $A=a$, then $\pi_a(Z, L)$ must be an RWF for $A=a$.
	\end{proposition}
	
	Actually, one can verify that
	$\chi^2[p_{Z|A,L}(\cdot|a,L)\| p_{Z|L}(\cdot|L)]=\mathrm{Var}\left[\pi_a(Z,L)\middle| L\right].$
	Proposition~\ref{prop: existence of RWFs} indicates that the IV relevance condition in Assumption~\ref{as: IV relevance} is equivalent to the statement that, for each $a\in\mathring{\mathcal{A}}$, there exists a corresponding RWF. Notably, the subscript $a$ indicates that $\pi_a$ is ``optimal'' for a fixed $a\in\mathring{\mathcal{A}}$.
	That is, if $\pi_a(Z, L)$ is not an RWF for $A = a$, then there is no other RWF at this point. This observation simplifies the search for a suitable RWF in practice. 
	
	In fact, there exists a family of ``optimal'' weighting functions that perform well when used as RWFs specifically for $A=a$. 
	This family can be generated by multiplying $\pi_a(Z, L)$ by any function $f(L)$ that is uniformly bounded above and below. 
	For instance, under Assumption~\ref{as: positivity}, one may simply take 
	$p_{A\mid Z,L}(a\mid Z,L)$ as an ``optimal'' weighting function, 
	since it belongs to the same family of ``optimal'' functions as $\pi_a(Z, L)$.
	
	\subsection{Uniform regular weighting function}\label{subsec: Uniform RWF}
	\magenta{
		So far, we have derived a family of RWFs for each $a \in \mathring{\mathcal{A}}$. 
		However, to identify the ADRF uniformly over a subset $\mathcal{N} \subseteq \mathring{\mathcal{A}}$, it is impractical to specify separate RWFs for each $a\in\mathcal{N}$. 
		Instead, it is desirable that all $a \in \mathcal{N}$ share a common RWF. 
		Hence, we adopt a uniform RWF defined on $\mathcal{N}$ instead of pointwise RWFs for every $a \in \mathcal{N}$.
	}
	
	\begin{definition}\label{defn: URWF}
		A measurable function $\pi(Z, L): \mathcal{Z} \times \mathcal{L} \to\mathbb{R}$ is said to be a uniform RWF (URWF) for a subset $\mathcal{N} \subseteq \mathring{\mathcal{A}}$ if it is uniformly bounded and there exists a constant $\epsilon_\pi(\mathcal{N}) > 0$ such that, for all $a \in \mathcal{N}$,
		$|\kappa_{\pi}^o(a,L)|\geq \epsilon_\pi(\mathcal{N})$ almost surely.
	\end{definition}
	% It is crucial to distinguish between an RWF at a single point $A = a$ and a URWF defined over a set $\mathcal{N} \subseteq \mathring{\mathcal{A}}$.
	\magenta{
		The definition of a URWF is necessary when identifying the ADRF over a set $\mathcal{N}\subseteq\mathring{\mathcal{A}}$ that contains an uncountable number of points.}
	Once a URWF is specified, it is important to establish the conditions under which it exists.  
	Intuitively, if $\pi(Z, L)$ is an RWF for $A = a$, it can be extended to a URWF in a neighborhood of $a$.  
	\begin{proposition}[Local stability of RWFs]
		\label{prop: local stability of RWF}
		Under Assumption~\ref{as: continuity}, let $\pi(Z,L)$ be an RWF for $A = a_0$ with $a_0 \in \mathring{\mathcal{A}}$, and assume that $\kappa_{\pi}^o(a,l)$ is equicontinuous at $a_0$, uniformly over $l\in\mathcal{L}$, i.e., for any $\epsilon(a_0)>0$, 
		there exists a positive $r_0$ such that, for any $l \in \mathcal{L},\;a\in \mathcal{B}(a_0, r_0)$,
		$|\kappa_{\pi}^o(a,l) - \kappa_{\pi}^o(a_0,l)| \leq \epsilon(a_0).$ Then, 
		there exists a small $h > 0$ such that $\pi(Z,L)$ is a URWF for $\mathcal{B}(a_0,h)$.
	\end{proposition}
	Proposition~\ref{prop: local stability of RWF} implies that if $\pi(Z, L)$ serves as an RWF at $A = a_0 \in \mathring{\mathcal{A}}$, then it can also serve as a URWF for a neighborhood around $a_0$. % provided that $\kappa_{\pi}^o(a,l)$ is equicontinuous at $A = a_0$ over $L$  
	This demonstrates the local stability of an RWF when the corresponding $\kappa_{\pi}^o(a,l)$ is not too sharp or irregular.
	
	So far, Proposition~\ref{prop: local stability of RWF} has established that a URWF can be constructed within a sufficiently small neighborhood under Assumptions~\ref{as: continuity} and \ref{as: IV relevance}.  
	A natural theoretical question is whether there exists a global URWF over the entire interior $\mathring{\mathcal{A}}$. Unfortunately, such a global URWF does not exist when $\mathcal{A}$ is a closed interval.
	
	\begin{proposition}[Nonexistence of a URWF for $\mathring{\mathcal{A}}$]
		\label{prop: non-existence of a URWF}
		Suppose Assumption~\ref{as: continuity} holds for every $L=l$, and let $\mathcal{A}$ be a closed interval and $\pi(Z,L)$ a measurable function such that $\kappa_{\pi}^o(a,l)$ is continuous in $a$.  
		Then, for any $l_0 \in \mathcal{L}$, there exists $a(l_0) \in \mathring{\mathcal{A}}$ such that 
		$\kappa_{\pi}^o(a(l_0), l_0) = 0.$
		Moreover, if $l_0 \in \mathring{\mathcal{L}}$ and the directional derivative of $\kappa_{\pi}^o(a,l)$ with respect to $l$ exists and is nonzero at $(a(l_0), l_0)$, then there exists a neighborhood $\mathcal{B}(a(l_0),h) \subseteq\mathring{\mathcal{A}}$ such that $\pi(Z,L)$ is not an RWF for any $A = a \in \mathcal{B}(a(l_0), h)$.
	\end{proposition}
	\begin{remark}\label{re: comparison to chen2025}
		\magenta{In the discrete case, the definitions of RWF and URWF do not differ in any essential way. The limitation described in Proposition~\ref{prop: non-existence of a URWF} arises only for continuous treatments and does not occur in the discrete setting, a distinction that stems from the completeness of the real line.}
	\end{remark}
	\begin{figure}[t]
		\centering
		\begin{subfigure}{0.33\textwidth}
			\centering
			\includegraphics[width=\linewidth]{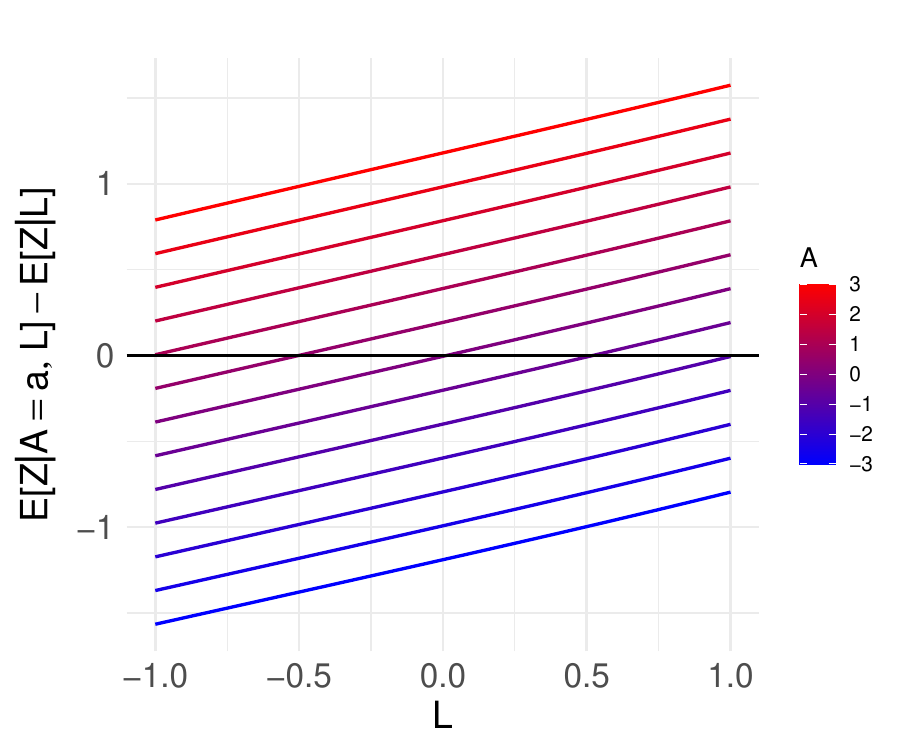}
			\caption{}
		\end{subfigure}
		\begin{subfigure}{0.33\textwidth}
			\centering
			\includegraphics[width=\linewidth]{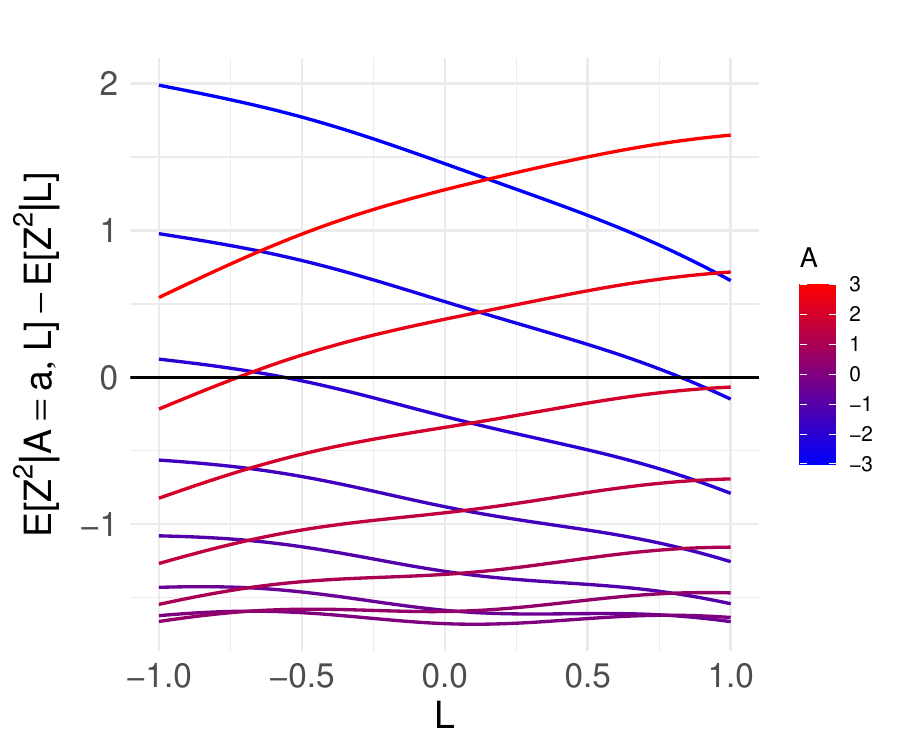}
			\caption{}
		\end{subfigure}
		\caption{
			(a) Illustration of $\kappa_{\pi}^o(A,L)$ for different treatment values $A=a$ with $\pi(Z,L)=Z$.
			(b) Illustration of $\kappa_{\pi}^o(A,L)$ for different treatment values $A=a$ with $\pi(Z,L)=Z^2$.
		}
		\label{fig: critical_problem}
	\end{figure}

	Proposition~\ref{prop: non-existence of a URWF} implies that different regions of $a \in \mathring{\mathcal{A}}$ actually require different choices of RWFs. For instance, 
	Figures~\ref{fig: critical_problem}(a) and (b) graphically illustrate the phenomena described in Propositions~\ref{prop: local stability of RWF} and \ref{prop: non-existence of a URWF}. 
	In both figures, each curve represents the value of $\kappa_{\pi}^o(a, L)$ for a fixed $a$, and each color corresponds to a specific treatment level $a \in \mathring{\mathcal{A}}$. The concrete data-generating process is described in Appendix~\ref{append: generating figures in section 2}. 
	
	In Figure~\ref{fig: critical_problem}(a), we set $\pi(Z, L) = Z$.  
	As shown, $\pi(Z, L) = Z$ acts as a URWF for $[-3, -2] \cup [2, 3]$, since all corresponding curves of $\kappa_{\pi}^o(a, L)$ remain strictly above or below the $x$-axis.  
	However, it fails to be an RWF for all $A = a \in [-0.2, 0.2]$, where the curves intersect the $x$-axis.  
	In contrast, in Figure~\ref{fig: critical_problem}(b), we set $\pi(Z, L) = Z^2$.  
	In this case, $\pi(Z, L) = Z^2$ serves as a URWF for $[-0.5, 0.5]$, but not for $A = a \in [-2.5, -2] \cup [2, 2.5]$.
	
	Since a global URWF over $\mathring{\mathcal{A}}$ may fail to exist, it is natural to instead cover $\mathring{\mathcal{A}}$ (or its subsets) by finitely many regions, each admitting a URWF. Within each region, local identification can then be achieved using a single URWF.
	In fact, based on the local existence and stability of URWFs, it is feasible to cover any compact subset of $\mathring{\mathcal{A}}$ by a finite collection of neighborhoods, each admitting a URWF.  
	
	\begin{proposition}[Finite covering number]\label{prop: finite covering number}
		Under Assumptions~\ref{as: continuity}--\ref{as: IV relevance},  
		for every $a \in \mathring{\mathcal{A}}$ with $\pi_a(Z,L) := p_{A\mid Z,L}(a \mid Z,L)/p_{A\mid L}(a\mid L)$,  
		suppose that $\kappa_{\pi_a}^o(a,l)$ is equicontinuous at $a$, uniformly over $l\in\mathcal{L}$.  
		Then, for any compact subset $\mathcal{A}_c \subseteq \mathring{\mathcal{A}}$,  
		there exists a finite collection of open intervals $\{\mathcal{B}(a_{m}, h_m)\}_{m=1}^M$, such that $\mathcal{A}_c \subseteq \bigcup_{m=1}^M \mathcal{B}(a_{m}, h_m)$.  
		Moreover, for each $m = 1, \ldots, M$, the function $\pi_{a_{m}}(Z,L)$ serves as a URWF for $\mathcal{B}(a_{m}, h_m)$.
	\end{proposition}
	\magenta{This proposition essentially illustrates the main motivation behind our practical guidance for constructing URWFs. See the discussion in Section~\ref{sec: practical guidance} for further details.}
	
	\subsection{Additive instrumental variables}\label{subsec: additive instrumental variable}
	\magenta{While we have established the IV relevance condition and its intrinsic link to the proposed RWF and URWF framework, these components alone do not guarantee the identification of the ADRF.}
	In this subsection, we introduce a key condition that plays a central role in the identification of the ADRF $\theta(a):=\mathbb{E}[Y(a)]$.
	
	\begin{definition}[Additive IV]\label{defn: AIV}
		For each fixed $a \in \mathring{\mathcal{A}}$, we say that $Z$ is an additive IV (AIV) for $A = a$ if there exist functions $b_a(U,L)$ and $c_a(Z,L)$ such that
		$p_{A\mid Z,U,L}(a \mid Z, U, L) = b_a(U, L) + c_a(Z, L).$ 
		Moreover, $Z$ is said to be an AIV for $A$ if, for every $a \in \mathring{\mathcal{A}}$, $Z$ is an AIV for $A=a$.
	\end{definition}
	
	The AIV definition is closely related to the no-interaction condition between $Z$ and $U$ in the treatment model \citep{wang2018bounded, Hartwig2023, wang2023instrumental, michael2024instrumental}, and similar formulations appear in the invalid IV literature \citep{Tchetgen_Tchetgen_2021, sun2022selective, ye2024geniusmawiirobustmendelianrandomization}.  
	\sy{Specifically, under the NUC framework, i.e., when $U=\emptyset$, the additive IV condition is automatically satisfied by taking $b_a(U,L)=0$ and $c_a(Z,L)=p_{A|Z,L}(a|Z,L)$.}
	% Since these works focus on discrete treatments, they do not directly extend to continuous treatments. 
	To clarify the proposed AIV condition, we derive an equivalent characterization.
	
	\begin{proposition}[Equivalent characterization of AIV]\label{prop: AIV characterization}
		Under Assumptions~\ref{as: IV independence}--\ref{as: positivity}, for each fixed $a\in \mathring{\mathcal{A}}$ and any $\pi(Z,L)$, define a weighting function
		\begin{align}\label{eq: omega}
			\omega_{a,\pi}(U,L) := \frac{\mathrm{Cov}\{ p_{A|Z,U,L}(a\mid Z,U,L),\, \pi(Z,L) \mid U,L \}}
			{\mathrm{Cov}\{ p_{A|Z,U,L}(a\mid Z,U,L),\, \pi(Z,L) \mid L \}}.
		\end{align}
		Then, if $\pi(Z,L)$ is an RWF for $A=a$, $\mathbb{E}[\omega_{a,\pi}(U,L)\mid L] \equiv 1$. Moreover, $Z$ is an AIV for $A = a$ if and only if, for any $\pi(Z,L)$, $\omega_{a,\pi}(U,L) \equiv 1.$
	\end{proposition}
    % \begin{align*}
			% 	\omega_{a,\pi}(U,L) := \frac{p_{A|U,L}(a\mid U,L)}{p_{A|L}(a\mid L)} \times 
			% 	\dfrac{\mathbb{E}[\pi(Z,L)\mid A=a,U,L] - \mathbb{E}[\pi(Z,L)\mid L]}{\mathbb{E}[\pi(Z,L)\mid A=a,L] - \mathbb{E}[\pi(Z,L)\mid L]}.
			% \end{align*}
	In Proposition~\ref{prop: AIV characterization}, $\omega_{a,\pi}(U,L)$ can be interpreted as a weighting function whose expectation equals one. This weighting function plays a central role in the identification of the ADRF.
	\begin{remark}\label{re: w function}
		When the treatment variable $A$ is discrete, we modify the definition of the weighting function $ \omega_{a,\pi}(U, L) $ by replacing $ p_{A\mid Z, U, L}(a\mid Z, U, L)$ with $\Pr(A = a \mid Z, U, L) $, and the weighting function reduces to that in \cite{chen2025identificationdebiasedlearningcausal}:
		\[\omega_{a,\pi}(U,L)= \frac{\mathrm{Cov}\{ I(A = a),\, \pi(Z,L) \mid U,L \}}
		{\mathrm{Cov}\{ I(A = a),\, \pi(Z,L) \mid L \}}.\]
	\end{remark}
	
	\begin{remark}[Generating an AIV from a mixture model]\label{re: AIV DGP}
		\magenta{We illustrate a setting satisfying the AIV assumption: with probability $p_0(L)$, $A \sim \tilde{p}_{A|Z,L}(a|Z,L)$ independent of $U$; otherwise, $A \sim \tilde{p}_{A|U,L}(a|U,L)$ independent of $Z$. By the law of total probability, 
			$p_{A|Z,U,L}(a\mid Z,U,L) = p_0(L)\, \tilde{p}_{A|Z,L}(a\mid Z,L) + (1-p_0(L))\, 
			\tilde{p}_{A|U,L}(a\mid U,L)$,
			so $Z$ satisfies the AIV structure for any $A=a$.}
	\end{remark}
	
	\section{Methodology}\label{sec: methodology}
	\subsection{Identification}\label{subsec: identification}
	In this subsection, we identify the ADRF $\theta(a):=\mathbb{E}[Y(a)]$ with the help of an AIV for each $a\in\mathring{\mathcal{A}}$. First, for a weighting function $\pi(Z,L)$, define $Z_{\pi} := \pi(Z,L)$, and the nuisance function
	$$\mu_{\pi}^o(a,l) := \frac{\mathbb{E}[Y Z_{\pi} \mid A=a,L=l] - \mathbb{E}[Y \mid A=a,L=l]\,\mathbb{E}[Z_{\pi} \mid L=l]}{\mathbb{E}[Z_{\pi} \mid A=a,L=l] - \mathbb{E}[Z_{\pi} \mid L=l]}.$$
	We formulate our main identification theorem as follows.
	\begin{theorem}[Identification]\label{thm: identification}
		Under Assumptions~\ref{as: consistency}--\ref{as: positivity}, for a fixed $ a\in\mathring{\mathcal{A}} $, 
		suppose that $\pi(Z,L)$ is an RWF for $A=a$. Then
		\begin{align}
			\mathbb{E}\left[ \mu_{\pi}^o(a,L)\right]
			=& \mathbb{E}[\mathbb{E}\{Y(a)\mid U,L\}\,\omega_{a,\pi}(U,L)],
			\label{eq: identification E[Y(a)]}\\
			\magenta{\mathbb{E}\left[ \mu_{\pi}^o(a,L)\right] - \mathbb{E}[Y(a)]=}& \magenta{\mathbb{E}\left[\mathrm{Cov}\{\mathbb{E}\{Y(a)\mid U,L\},\,\omega_{a,\pi}(U,L)\mid L\}\right].}
			\label{eq: identification E[Y(a)] bias}
		\end{align}
		Furthermore, if $Z$ is an AIV for $A=a$, \magenta{we have $\theta(a):=\mathbb{E}[Y(a)] = \mathbb{E}[\mu_{\pi}^o(a,L)]$.}
	\end{theorem}
	\begin{remark}\label{re: GRDF}\magenta{
			For any function $\phi$, Theorem~\ref{thm: identification} still holds if $Y$, $Y(a)$ are replaced by $\phi(Y)$, $\phi(Y(a))$ respectively. For example, setting $\phi(Y) = I(Y \leq t)$ allows us to identify the distribution function of $Y(a)$, namely $\Pr(Y(a) \leq t)$. This demonstrates that the identification framework may be extended to more general dose-response functions, such as quantile dose-response functions \citep{ai2026datadrivenuniforminferencegeneral}.}
	\end{remark}
	\sy{Theorem~\ref{thm: identification} shows that the ADRF is identifiable under a valid AIV condition. In Appendix~\ref{append: NPIV}, we provide an alternative nonparametric IV interpretation of this identification result, which further explains why $\mathbb{E}[\mu_{\pi}^o(a,L)]$ serves as the identification functional for the ADRF.}
	%\begin{remark}\label{re: IV relevance violation}
	Note that by Proposition~\ref{prop: existence of RWFs}, the existence of an RWF at $A=a$ implies that $\mathrm{Var}[\pi_a(Z,L)\mid L]$ is bounded away from zero at $A=a$ almost surely. Thus, in Theorem~\ref{thm: identification}, Assumption~\ref{as: IV relevance} is implicitly required for $A=a$ rather than uniformly over $a\in \mathring{\mathcal{A}}$. Similarly, Assumption~\ref{as: positivity} only needs to hold at $A=a$ as well.
	%\end{remark}

	\sy{In practice, the AIV assumption can be violated in certain applications, e.g., in a linear treatment assignment model of the form $A= Z + U + L + \epsilon.$ Below, we provide several remarks to discuss the interpretation of the functional $\mathbb{E}[\mu_{\pi}^o(a, L)]$ when the AIV assumption does not hold.}\magenta{
		\begin{remark}[Contrast identification]\label{re: derivative identification}
			Under Assumptions~\ref{as: consistency}--\ref{as: positivity}, for any $a_1,a_2\in\mathring{\mathcal{A}}$, assume that $\pi_1$, $\pi_2$ are two RWFs for $A=a_1,a_2$, respectively. Then, from Equation~\eqref{eq: identification E[Y(a)] bias},
			\begin{align*}
				&\mathbb{E}[\mu_{\pi_1}^o(a_1,L)] -\mathbb{E}[\mu_{\pi_2}^o(a_2,L)] 
				- \{\mathbb{E}[Y(a_1)]-\mathbb{E}[Y(a_2)]\}\\
				=&\mathbb{E}\left[\mathrm{Cov}\left\{\mathbb{E}[Y(a_1)-Y(a_2)|U,L], \omega_{a_1,\pi_1}(U,L) \mid L\right\}\right]\\
				&-\mathbb{E}\left[\mathrm{Cov}\left\{\mathbb{E}[Y(a_2)|U,L], \omega_{a_2,\pi_2}(U,L) -\omega_{a_1,\pi_1}(U,L)\mid L\right\}\right].
			\end{align*}
			If $\mathbb{E}[Y(a_1)-Y(a_2)|U,L]$ and $\omega_{a_2,\pi_2}(U,L)-\omega_{a_1,\pi_1}(U,L)$ do not depend on $U$, $\mathbb{E}[Y(a_1)-Y(a_2)]$ remains identifiable.  This can be viewed as a natural generalization of the no unmeasured common effect modifier assumption in \citet{cui2021semiparametric}. See Appendix~\ref{append: weak AIV} for further discussion on a novel weak AIV (WAIV) assumption, i.e., $\omega_{a,\pi}(U,L)$ does not depend on  $(a,\pi)$, which implies $\omega_{a_2,\pi_2}(U,L)-\omega_{a_1,\pi_1}(U,L)$ does not depend on $U$. Moreover, as shown in Corollary~\ref{cor: derivative identification}, we note that this idea can further be extended to identify the derivative of the ADRF, $\nabla_a \mathbb{E}[Y(a)]$. 
		\end{remark}
		\begin{remark}[Covariate shift]\label{re: covariate shift}
			The right-hand side of Equation~\eqref{eq: identification E[Y(a)]} admits a weighted-average interpretation of $\mathbb{E}[Y(a)\mid U,L]$ when $\omega_{a,\pi}(U,L)\ge 0$, linking our framework to the assumption-lean perspective of \citet{vansteelandt2022assumption,vansteelandt2024assumption}. 
			In particular, by Proposition~\ref{prop: AIV characterization}, 
			$
			p_{a,\pi}(U\mid L) := \omega_{a,\pi}(U,L)\, p_{U\mid L}(U\mid L)
			$
			defines a shifted conditional distribution since $\mathbb{E}[\omega_{a,\pi}(U,L)\mid L]=1$. Let $\mathbb{E}_{a,\pi}$ denote expectation under this distribution. Then
			$
			\mathbb{E}_{a,\pi}[Y(a)]
			= \mathbb{E}_{a,\pi}\!\left[\mathbb{E}\{Y(a)\mid U,L\}\right]
			= \mathbb{E}\!\left[\mathbb{E}\{Y(a)\mid U,L\}\,\omega_{a,\pi}(U,L)\right]
			= \mathbb{E}\!\left[\mu_{\pi}^o(a,L)\right].
			$
			Notably, $\omega_{a,\pi}(U,L)$ is likely to remain positive as long as the influence of $U$ on the numerator of Equation~\eqref{eq: omega} is not too large. 
			% As we can see from Equation~\eqref{eq: omega}, the weighting function $\omega_{a,\pi}(U,L)$ is positive as long as the influence of unobserved confounders $U$ on the $\mathrm{Cov}\{ p_{A\mid Z,U,L}(a\mid Z,U,L),\, \pi(Z,L) \mid U,L \}$ is not too large relative to the variation captured by $L$. 
		\end{remark}
		\begin{remark}[Partial identification]\label{re: partial identification}
			The right-hand side of Equation~\eqref{eq: identification E[Y(a)] bias} characterizes the gap between the estimand $\mathbb{E}[\mu_{\pi}^o(a,L)]$ and the true ADRF $\mathbb{E}[Y(a)]$. Even without the AIV assumption,
			\begin{align*}
				\left|\mathbb{E}[\mu_{\pi}^o(a,L)] - \mathbb{E}[Y(a)]\right|
				\leq\mathbb{E}\!\left[
				\mathrm{Var}\{\mathbb{E}[Y(a)\mid U,L]\mid L\}\right]^{1/2}\cdot
				\mathbb{E}\!\left[
				\mathrm{Var}\{\omega_{a,\pi}(U,L)\mid L\}\right]^{1/2}.
			\end{align*}
			Hence, partial identification remains possible provided these variance terms are controlled well. 
	\end{remark}}
	
	\sy{In summary, when the AIV condition holds, the pointwise ADRF is identified by $\mathbb{E}[\mu_{\pi}^o(a,L)].$ When only the WAIV condition holds, together with the additional restriction that $\mathbb{E}[Y(a_1)-Y(a_2)\mid U,L]$ does not depend on $U$, causal contrasts between treatment levels can still be identified. When neither AIV nor WAIV holds, the functional $\mathbb{E}[\mu_{\pi}^o(a,L)]$ admits a covariate-shift interpretation.}
	Next, we show two degenerate examples to illustrate Theorem~\ref{thm: identification}.
	\begin{example}[NUC]
		When $U=\emptyset$, set $Z\equiv A$ and $\pi(Z,L)=f(A)$. Equation~\eqref{eq: identification E[Y(a)]} degenerates:
		\begin{align*}
			\mathbb{E}[Y(a)] 
			= &\mathbb{E}\left[ 
			\dfrac{ \{f(a)-\mathbb{E}[f(A)\mid L]\}\mathbb{E}[Y\mid A=a,L]}{ f(a) - \mathbb{E}[f(A)\mid L]} 
			\right]
			= \mathbb{E}\left[\mathbb{E}[Y\mid A=a,L]\right].
		\end{align*}
		This demonstrates that the identification formula degenerates when there are no unmeasured confounders.
	\end{example}
	
	\begin{example}[Binary IV]
		Assume that $ Z $ is binary, and $\pi(Z,L)=Z$ is an RWF for $ A = a $. Then, 
		if $ Z $ is an AIV for $ A $,
		Theorem~\ref{thm: identification} implies that 
		\begin{align*}
			\mathbb{E}[Y(a)] =
			\mathbb{E}\left[\dfrac{\mathbb{E}[YZ \mid A=a,L] - \mathbb{E}[Y \mid A=a,L]\mathbb{E}[Z \mid L]}
			{\mathbb{E}[Z \mid A=a,L] - \mathbb{E}[Z \mid L]}\right].
		\end{align*}
	\end{example}
	From Proposition~\ref{prop: local stability of RWF}, when identifying the ADRF over a small subset $\mathcal{N}\subseteq\mathring{\mathcal{A}}$, one may rely on a single fixed $\pi(Z,L)$ as the URWF for $\mathcal{N}$. However, Proposition~\ref{prop: non-existence of a URWF} cautions against using a single fixed $\pi(Z,L)$ for the entire $\mathring{\mathcal{A}}$, as doing so may lead to critical errors.
	As a compromise, by Proposition~\ref{prop: finite covering number}, it is preferable to select different RWFs for different values of $A$. 
	
	For a graphical illustration, consider Figure~\ref{fig: critical_problem} again. When $A$ is close to zero, Figure~\ref{fig: critical_problem}(b) suggests that $\pi(Z,L)=Z^2$ is an appropriate RWF for identification. In contrast, when $A$ is far from zero, Figure~\ref{fig: critical_problem}(a) indicates that $\pi(Z,L)=Z$ can be used as an RWF.
	Consequently, in practice, we should fix an RWF $\pi(Z,L)$ for $A = a_0$, to identify the ADRF only within a neighborhood of $a_0 \in \mathring{\mathcal{A}}$.  
	% To conclude this subsection, 
	% Then, by Theorem~\ref{thm: identification}, we have $\theta(a)= \mathbb{E}[\mu_{\pi}^o(a,L)]$ whenever $Z$ is an AIV for $A=a$ and $\pi(Z,L)$ is an RWF for $A=a$.

	\subsection{Semiparametric theory}\label{subsec: semiparametric theory}

	In practice, the true form of $\mu_{\pi}^o(A,L)$ is unknown. Moreover, in the absence of parametric assumptions, the functional $\mathbb{E}[\mu_{\pi}^o(a,L)]$ lacks pathwise differentiability, which implies that root-$n$ consistent estimators for it are not readily available. To proceed, we follow the approach of \citet{kennedy2017non} and leverage semiparametric theory to establish an AIPW score function $\varphi_{\pi}(O)$ that satisfies $\mathbb{E}\left[\varphi_{\pi}(O) \middle| A=a\right] = \theta(a)$ when $Z$ is an AIV for $A=a$.
	
	First, we define an estimand that represents the average outcome under an intervention that randomly assigns treatment based on the marginal density $q(a)p_A(a)$:
	\begin{align}\label{eq: psi defn}
		\psi_{\pi,q}^o := \int_{\mathcal{A}} \int_{\mathcal{L}} \mu_{\pi}^o(a, l) \mathrm{d}P_L(l) q(a)\mathrm{d}P_A(a).
	\end{align}
	The quantity $\psi_{\pi,q}^o$ can be interpreted as a covariate-shifted causal estimand, corresponding to the scenario where the conditional treatment distribution $p_{A|L}(a \mid l)$ is replaced by a rescaled version $q(a)p_A(a)$.
	Before giving the semiparametric theory for $\psi_{\pi,q}^o$, we define nuisance functions as 
	$\rho_{\pi}^o(L):=\mathbb{E}[Z_{\pi}\mid L]$, 
	$\kappa_{\pi}^o(A,L) := \mathbb{E}[Z_{\pi}\mid A,L] - \mathbb{E}[Z_{\pi}\mid L]$, 
	$\eta^o(A,L):= \mathbb{E}[Y\mid A,L]$, 
	$\delta^o(A,L) :=  p_{A}(A)/p_{A|L}(A\mid L).$
	% \begin{align*}
		% 	\begin{array}{ll}
			% 		\rho_{\pi}^o(L):=\mathbb{E}[Z_{\pi}\mid L], &\kappa_{\pi}^o(A,L) := \mathbb{E}[Z_{\pi}\mid A,L] - \mathbb{E}[Z_{\pi}\mid L], \\
			% 		\eta^o(A,L):= \mathbb{E}[Y\mid A,L], &\delta^o(A,L) :=  p_{A}(A)/p_{A|L}(A\mid L).
			% 	\end{array}
		% \end{align*}
	We unify the nuisance functions into a nuisance vector
	$\alpha_{\pi}^o(A,L) = [\mu_{\pi}^o(A,L),\rho_{\pi}^o(L), \kappa_{\pi}^o(A,L), \eta^o(A,L), \delta^o(A,L)].$
	% \begin{align*}
		% 	\alpha_{\pi}(A,L) = [\mu_{\pi}(A,L),\rho_{\pi}(L), \kappa_{\pi}(A,L), \eta(A,L), \delta(A,L)].
		% \end{align*}
	Next, we derive the efficient influence function (EIF) for $\psi_{\pi,q}^o$ in the following theorem.
	\begin{theorem}[EIF of $\psi_{\pi,q}^o$]\label{thm: psi_q}
		Suppose that $\pi(Z,L)$ is a URWF over the support of $q(a)$. Then, the EIF for $\psi_{\pi,q}^o$ in the nonparametric model is given by
		\begin{align*}
			&\delta^o(A,L)(Z_{\pi}-\rho_{\pi}^o(L))
			\left\{
			\dfrac{Y-\mu_{\pi}^o(A,L)}{\kappa_{\pi}^o(A,L)}q(A)
			-\displaystyle\int \dfrac{\eta^o(a,L) - \mu_{\pi}^o(a,L)}
			{\kappa_{\pi}^o(a,L)} q(a)\mathrm{d} P_A(a)\right\}\\
			&+q(A)\int \mu_{\pi}^o(A,l)\mathrm{d} P_L(l)
			+\int \mu_{\pi}^o(a,L) q(a)\mathrm{d} P_A(a) - 2\psi_{\pi,q}^o.
		\end{align*}
	\end{theorem}
	\begin{remark}\label{re: q(a) defn}
		In \citet{kennedy2017non}, they simply set $q(A) \equiv 1$. Here, we derive the EIF for a general $q(a)$ because, as shown in Proposition~\ref{prop: non-existence of a URWF}, a URWF may not exist over $\mathring{\mathcal{A}}$. In such cases, introducing $q(a)$ effectively restricts the range of values, so that $\pi(Z,L)$ only needs to be a URWF over the support of $q(a)$, rather than over the entire $\mathring{\mathcal{A}}$. 
	\end{remark}
	Define $z_{\pi}:=\pi(z,l)$. Motivated by the EIF of $\psi_{\pi,q}^o$, we establish an AIPW score function as
	\begin{equation}
		\begin{aligned}
			\varphi_{\pi}(O;\alpha_{\pi},P_O) := &\delta(A,L)
			\dfrac{(Z_{\pi}- \rho_{\pi}(L))}{\kappa_{\pi}(A,L)}\{Y-\mu_{\pi}(A,L)\}\\
			&+\int \mu_{\pi}(A,l) - (z_{\pi}-\rho_{\pi}(l))\dfrac{\eta(A,l)-\mu_{\pi}(A,l)}{\kappa_{\pi}(A,l)}\mathrm{d} P_{O}(o).
			\label{eq: varphi defn}
		\end{aligned}
	\end{equation}
	\sy{Akin to \citet{kennedy2017non}, as shown in Lemma~\ref{lem: mixed bias}, $\varphi_{\pi}(O;\alpha_{\pi},P_O)$ admits a mixed bias property and $\mathbb{E}[\varphi_{\pi}(O; \alpha_{\pi}^o, P_O)\mid A=a]=\theta(a)$.
		Therefore, $\varphi_{\pi}(O;\alpha_{\pi},P_O)$ can be viewed as an AIPW score for the pointwise target $\theta(a)$.
	}

	Notably, when $A$ is discrete, provided that $\delta^o(a,L)$ is redefined as $\Pr(A=a)/\Pr(A=a \mid L)$, we still have $\mathbb{E}[\varphi_{\pi}(O;\alpha_{\pi}^o, P_O) \mid A=a] = \mathbb{E}[Y(a)]$.  
	This indicates that, although the proposed AIPW score is designed for continuous treatments, it also applies when $A$ is discrete. 
	A detailed comparison with \citet{chen2025identificationdebiasedlearningcausal} is provided in Appendix~\ref{append: multi-categorical treatments}. 
	Moreover, the construction naturally accommodates the degenerate case where $L$ is empty. 
	In this setting, the AIPW score differs slightly from Equation~\eqref{eq: varphi defn}; see Appendix~\ref{append: degenerate AIPW score} for more details.

	\subsection{General cross-fitting procedures}\label{subsec: estimation procedure}
	In this subsection, we construct the AIPW score vector for $\theta(a)$ using a general cross-fitting procedure. This procedure has previously appeared in the estimation of heterogeneous treatment effects \citep{kennedy2023optimaldoublyrobustestimation} and ADRFs \citep{bonvini2022fastconvergenceratesdoseresponse}.
	Concretely, we randomly partition the samples into $K$ folds, denoted by $\{I_k\}_{k=1}^K$, and let $I_{-k} := \bigcup_{j \neq k} I_j$ denote the complement of $I_k$. When training the nuisance functions $\hat\alpha_{\pi}^{(n,-k)}$ and $\hat{P}_O^{(n,-k)}$ for the $k$-th fold, we further divide $I_{-k}$ into $I_{-k}^{(1)}$ and $I_{-k}^{(2)}$. In particular, the sizes of $I_{-k}^{(1)}$ and $I_{-k}^{(2)}$ can be freely adjusted. 
	
	For notational simplicity, let $\mathbb{E}_{nk}[O] := \sum_{i \in I_k} O_i/|I_k|$ denote the sample average within the $k$-th fold. 
	We further define conditional expectations based on the folds: 
	$\mathbb{E}_k[O] := \mathbb{E}[O \mid O_i, i \notin I_k]$, 
	$\mathbb{E}_{-k}[O] := \mathbb{E}[O \mid O_i, i \notin I_{-k}]$, and 
	$\mathbb{E}_{-k}^{(l)}[O] := \mathbb{E}[O \mid O_i, i \notin I_{-k}^{(l)}]$ for $l = 1, 2$.
	
	The samples in $I_{-k}^{(1)}$ are employed to train $\hat\alpha_{\pi}^{(n,-k)}$, whereas those in $I_{-k}^{(2)}$ are used to train $\hat{P}_O^{(n,-k)}$. 
	This design guarantees that the training sets for $\hat\alpha_{\pi}^{(n,-k)}$ and $\hat{P}_O^{(n,-k)}$ are mutually independent. 
	Subsequently, we compute the AIPW score vector $[\varphi_{\pi}(O_i; \hat\alpha_{\pi}^{(n,-k_i)}, \hat{P}_O^{(n,-k_i)})]_{i=1}^n$, where $k_i$ is the unique fold index such that $O_i\in I_{k_i}$, and $I_{-k_i}^{(1)}$ and $I_{-k_i}^{(2)}$ are the two subsamples obtained by splitting $I_{-k_i}$. 
	The entire cross-fitting procedure is summarized in Algorithm~\ref{alg: cross-fitting}. 
	Following the same logic, an alternative cross-fitting strategy is provided in Appendix~\ref{append: alternative cross-fitting procedure}.

	\begin{algorithm}[t]
		\caption{Cross-fitting procedure}
		\label{alg: cross-fitting}
		\begin{algorithmic}[1]
			\State \textbf{Input:} Number of folds $K$ and $\{O_i\}_{i=1}^n$.
			\State Randomly split the sample into $K$ folds $\{I_k\}_{k=1}^K$.
			\For{$k = 1, \ldots, K$}
			\State Split $I_{-k}$ into two subsamples $I_{-k}^{(1)}$ and $I_{-k}^{(2)}$.
			\State Train nuisance estimators on $I_{-k}^{(1)}$: 
			$\hat{\alpha}_{\pi}^{(n,-k)} := 
			\big[\hat\mu_{\pi}^{(n,-k)}, \hat\rho_{\pi}^{(n,-k)}, 
			\hat\kappa_{\pi}^{(n,-k)}, \hat\eta^{(n,-k)}, 
			\hat\delta^{(n,-k)}\big].$
			\State Estimate the empirical distribution $\hat{P}_O^{(n,-k)}$ using $I_{-k}^{(2)}$.
			\EndFor
			\State \textbf{Output:} AIPW score vector $[\varphi_{\pi}(O_i; \hat\alpha_{\pi}^{(n,-k_i)}, \hat{P}_O^{(n,-k_i)}) ]_{i=1}^n$, where $O_i\in I_{k_i}$ for each $i$.
		\end{algorithmic}
	\end{algorithm}

	\subsection{Local linear kernel regression} \label{subsec: LLKR}
	In this subsection, we estimate the ADRF locally using LLKR by regressing the AIPW scores on the treatment variable. 
	Specifically, we select a kernel weighting function $K(a)$, bandwidth $h$, and a target point $a \in \mathring{\mathcal{A}}$. Let $K_h(a):=K(a/h)/h$, $g(a):=[1,a]^{\top} $, and $e_1:=[1,0]^{\top}$. We then calculate the AIPW estimator $\hat{\theta}_{\pi,h}^{(n)}(a)$ by solving:
	\begin{align}\label{eq: LLKR}
		e_1^{\top}  \underset{\beta \in \mathbb{R}^2}{\arg\min}
		\frac{1}{n}\sum_{k=1}^{K} 
		\sum_{i\in I_k}
		K_h(A_i - a)
		\left\{\varphi_{\pi}(O_i; \hat\alpha_{\pi}^{(n,-k)}, \hat{P}_O^{(n,-k)})
		- g\left(\dfrac{A_i-a}{h}\right)^{\top} \beta\right\}^2.
	\end{align}
	Notably, the kernel $K_h(A - a)$ restricts attention to a local neighborhood around $a$. 
	As a result, the outputs of Algorithm~\ref{alg: cross-fitting} are required to be bounded and regular only for values of $A$ within this neighborhood, reflecting the inherently local nature of the estimation.
	
	For bandwidth selection, classical criteria such as leave-one-out cross-validation (LOOCV), 
	generalized cross-validation (GCV), and the $\text{C}_{\text{p}}$-criterion must also be localized, since a URWF may not exist over $\mathring{\mathcal{A}}$ (Proposition~\ref{prop: non-existence of a URWF}). 
	Concretely, define a localized LOOCV for bandwidth selection as:
	\begin{align*}
		\hat{h}_{\text{opt}}^{(n)}(\mathcal{N})=\underset{h>0}{\arg\min} \dfrac{1}{n}\sum_{k=1}^K\sum_{i\in I_k}\left\{
		\dfrac{\varphi_{\pi}(O_i;\hat\alpha_{\pi}^{(n,-k)}, \hat{P}_O^{(n,-k)})-\hat{\theta}_{\pi,h}^{(n)}(A_i)}
		{1-w_{ni}(A_i,h)}
		\right\}^2I\{A_i\in \mathcal{N}\},
	\end{align*}
	where $w_{ni}(A_i,h)$, defined in Appendix~\ref{append: additional lemma and notation}, is the $i$-th diagonal element in the smoothing matrix. 
	This criterion minimizes the localized residual mean squared error (RMSE). 
	The localized GCV and $\text{C}_{\text{p}}$-criterion can be formulated similarly, just by focusing only on the samples in $\mathcal{N}$.

	\sy{
		\subsection{A falsification test}\label{subsec: falsification test}
		In practice, the validity of the AIV condition is generally unknown. To address this issue, we note that under Assumptions~\ref{as: consistency}--\ref{as: positivity}, if $Z$ is an AIV for $A$, $\mathbb{E}[\mu_{\pi}(a,L)]$ must be invariant across different choices of RWFs. Hence, a significant discrepancy between estimates obtained from two RWFs provides evidence against the AIV condition. 
		
		Based on this observation, we provide a Hansen J-type statistic \citep{hansen1982large} for testing the AIV condition pointwise. In particular, for a positive integer $J$, define $\varphi_{\Pi}(O;\alpha_{\Pi},P_O):=[\varphi_{\pi_j}(O;\alpha_{\pi_j},P_O)]_{j=1}^J\in\mathbb{R}^J$, and $\mathbf{1}_J:=[1,\ldots, 1]^{\top}\in\mathbb{R}^J$. Let ``$\otimes$'' denote the Kronecker product, $\mathbf{I}_{2J}\in\mathbb{R}^{2J\times 2J}$  an identity matrix, $Q\in\mathbb{R}^{2J\times 2J}$ a weighting matrix, and $\|v\|_Q^2:=v^{\top} Q v$ for all $v\in\mathbb{R}^{2J}$. Based on this, for any $\beta\in\mathbb{R}^2$, define a loss function
		\begin{align*}
			&\mathcal{L}_{\Pi,h,a}^{(n)}(\beta,Q):= 
			\left\|\frac{1}{K}\sum_{k=1}^{K}\mathbb{E}_{nk}\left[\Psi_{\Pi,h,a}(O;\beta,\hat\alpha_{\Pi}^{(n,-k)}, \hat{P}_O^{(n,-k)})\right]
			\right\|_{Q}^2,\\
			&\Psi_{\Pi,h,a}(O;\beta,\alpha_{\Pi},P_O):= K_h(A - a)g\left(\dfrac{A-a}{h}\right)\otimes
			\left\{\varphi_{\Pi}(O; \alpha_{\Pi}, P_O)
			- \mathbf{1}_J g\left(\dfrac{A-a}{h}\right)^{\top} \beta\right\}.
		\end{align*}
		Define an initial estimator $\hat{\beta}_{\Pi,h}^{(n,\mathrm{init})}(a):=\underset{\beta\in\mathbb{R}^2}{\arg\min}\{\mathcal{L}_{\Pi,h,a}^{(n)}(\beta,\mathbf{I}_{2J})\}$, and a debiased generalized method of moments (GMM) estimator $\hat{\beta}_{\Pi,h}^{(n)}(a) :=\underset{\beta\in\mathbb{R}^2}{\arg\min}\{\mathcal{L}_{\Pi,h,a}^{(n)}(\beta,\hat{Q}_{\Pi,h}^{(n)}(a)^{-1})\},$ where
		\begin{align*}
			\hat{Q}_{\Pi,h}^{(n)}(a)=&\frac{h}{K}\sum_{k=1}^{K}\mathbb{E}_{nk}\left[\Psi_{\Pi,h,a}(O;\hat{\beta}_{\Pi,h}^{(n,\mathrm{init})}(a),\hat\alpha_{\Pi}^{(n,-k)}, \hat{P}_O^{(n,-k)})\Psi_{\Pi,h,a}(O;\hat{\beta}_{\Pi,h}^{(n,\mathrm{init})}(a),\hat\alpha_{\Pi}^{(n,-k)}, \hat{P}_O^{(n,-k)})^{\top}\right].
		\end{align*}
		
		Define $m_j(a):=\mathbb{E}[\varphi_{\pi_j}(O;\alpha_{\pi_j}^o,P_O)| A=a]$. We consider the null hypothesis $\widetilde{\mathcal{H}}_0(a):$ for any $j_1,j_2=1,\ldots,J,\; m_{j_1}(a)=m_{j_2}(a)$ and $m_{j_1}'(a)=m_{j_2}'(a)$, which requires all URWF-induced conditional means to share the same value and first derivative at $a$. The alternative hypothesis is $\widetilde{\mathcal{H}}_1(a): \text{ there exist } j_1\neq j_2 \text{ such that } m_{j_1}(a)\neq m_{j_2}(a)$ or $m_{j_1}'(a)\neq m_{j_2}'(a)$. By Theorem~\ref{thm: identification}, under Assumptions~\ref{as: consistency}--\ref{as: positivity}, if $Z$ is an AIV for $A$, then $\widetilde{\mathcal{H}}_0(a)$ must hold. Hence, rejection of $\widetilde{\mathcal{H}}_0(a)$ provides evidence that at least one of the assumptions, including the AIV condition, is violated.
		
		In particular, under $\widetilde{\mathcal{H}}_0(a)$, the Hansen J-type statistic satisfies $nh\mathcal{L}_{\Pi,h,a}^{(n)}(\hat{\beta}_{\Pi,h}^{(n)}(a),\hat{Q}_{\Pi,h}^{(n)}(a)^{-1}) \overset{d}{\to}\chi^2_{2J-2}$ (see Theorem~\ref{thm: falsify AIV condition} for the corresponding theoretical justification). Therefore, $\widetilde{\mathcal{H}}_0(a)$ is rejected when the observed value of this statistic is sufficiently large. The corresponding $p$-value is $\hat{p}_{a,\Pi}^{(n)} :=1-F_{\chi^2_{2J-2}}( nh\mathcal{L}_{\Pi,h,a}^{(n)}(\hat{\beta}_{\Pi,h}^{(n)}(a), \hat{Q}_{\Pi,h}^{(n)}(a)^{-1}) )$, where $F_{\chi^2_{2J-2}}$ denotes the cumulative distribution function of the $\chi^2_{2J-2}$ distribution. Thus, the proposed procedure provides a falsification test for the necessary condition $\widetilde{\mathcal{H}}_0(a)$ implied by the AIV condition.}

	\subsection{RWF testing procedure}\label{subsec: RWF testing procedure}
	\sy{All analyses above assume that RWFs are known.}
	In practice, once a weighting function $\pi(Z, L)$ is prespecified, it is natural to ask whether it serves as an RWF for any treatment level $A=a$. This motivates the development of a hypothesis testing procedure for assessing the RWF condition. To address this question, we first consider a hypothesis test that captures the key aspects of RWF when $\mathcal{L}:=\{1,\ldots, |\mathcal{L}|\}$ is discrete. Concretely, we define the estimand 
	$\zeta_{\pi}^o(a,l) := p_{A\mid L}(a\mid l)\kappa_{\pi}^o(a,l).$ 
	We then test the null hypothesis $\mathcal{H}_0(a)$: $\zeta_{\pi}^o(a,l_0) = 0$ for some $l_0\in\mathcal{L}$ versus $ \mathcal{H}_1(a): \zeta_{\pi}^o(a,l) \neq 0$ for all $l\in\mathcal{L}$.
	
	When $\mathcal{L}$ is a compact subset of $\mathbb{R}^{k}$, 
	this testing procedure jointly assesses the RWF condition and the positivity condition in Assumption~\ref{as: positivity}, because a violation of either condition may lead to a failure to reject $\mathcal{H}_0(a)$. An advantage of this formulation is that the resulting $p$-value is directly related to the stability of the subsequent estimator.
	
	\begin{algorithm}[t]\sy{
			\caption{The $p$-value for testing $ \mathcal{H}_0(a):$ $\zeta_{\pi}^o(a,l_0) = 0$ for some $l_0\in\mathcal{L}$ (discrete $L$).}
			\label{alg: test RWF discrete}
			\begin{algorithmic}[1]
				\State \textbf{Input:}
				A prespecified weighting function $\pi$; a target point $a$; number of cross-fitting folds $K$; testing samples $\{O_i\}_{i=1}^n$.
				\State Set the bandwidth $h_A = \, n^{-1/4}\hat{\sigma}_{A}$, where $\hat{\sigma}_{A}$ is the sample standard deviation of $\{A_i\}_{i=1}^n$.
				\State Randomly partition the samples $\{O_i\}_{i=1}^n$ into $K$ folds $\{I_k\}_{k=1}^K$.
				\For{$k = 1, \ldots, K$}
				\State Using the samples in $I_{-k}$, fit nuisance functions
				$$\hat{\gamma}_{h_A,a}^{(n,-k)}(L):=\hat{\mathbb{E}}^{(n,-k)}[K_{h_A}(A-a)\mid L],\,\hat{\rho}_{\pi}^{(n,-k)}(L):=\hat{\mathbb{E}}^{(n,-k)}[\pi(Z,L)\mid L].$$
				\State For each $i\in I_k$, calculate $\hat{S}_{\pi,h_A}(a,O_i)=(K_{h_A}(A_i-a) - \hat{\gamma}_{h_A,a}^{(n,-k)}(L_i))(\pi(Z_i,L_i) - \hat{\rho}_{\pi}^{(n,-k)}(L_i)).$
				\EndFor
				\State For each $l\in\mathcal{L}$, calculate $n_L=\sum_{i=1}^n I(L=l)$, 
				$\tilde{\zeta}_{\pi,h_A}^{(n)}(a,l)= \frac{1}{n_L}\sum_{i=1}^n\hat{S}_{\pi,h_A}(a,O_i)I(L_i=l),\;\text{ and }$
				$\tilde{\sigma}_{\pi,\zeta,h_A}^{(n)}(a,l)^2:= nh_A\sum_{i=1}^n (\hat{S}_{\pi,h_A}(a,O_i)-\tilde{\zeta}_{\pi,h_A}^{(n)}(a,l))^2 I(L_i=l)/n_L^2.$
				\State \textbf{Output:} A $p$-value $\tilde{p}_{a,\pi}^{(n)} := \max\limits_{l\in\mathcal{L}}\{2 - 2\Phi(\sqrt{nh_A}|\tilde{\zeta}_{\pi,h_A}^{(n)}(a,l)|/\tilde{\sigma}_{\pi,\zeta,h_A}^{(n)}(a,l))\}$.
				%				, where $n_l:=\sum_{i=1}^nI(L_i=l)$.
		\end{algorithmic}}
	\end{algorithm}

	\begin{algorithm}[t]\sy{
			\caption{The $p$-value for testing $\mathcal{H}_0(a): \zeta_{\pi}^o(a,l_0) = 0$ for some $l_0\in\mathcal{L}$ (continuous $L$).}
			\label{alg: test RWF continuous}
			\begin{algorithmic}[1]
				\State \textbf{Input:}
				A prespecified weighting function $\pi$; a target point $a$; the number of cross-fitting folds $K$; a testing sample $\{O_i\}_{i=1}^n$; and the number of bootstrap replications $B$.
				\State Set the bandwidths $h:=[h_A,h_L]:=[n^{-1/5}\hat{\sigma}_{A},n^{-1/5}\hat{\sigma}_{L}]$, where $\hat{\sigma}_{A}$, $\hat{\sigma}_{L}$ are the sample standard deviations of $\{A_i\}_{i=1}^n$ and $\{L_i\}_{i=1}^n$.
				%				\State Randomly partition the samples into $K$ folds $\{I_k\}_{k=1}^K$.
				%				\State For all $k$, using the samples in $I_{-k}$, fit nuisance functions
				%				$$\hat{\gamma}_{h_A,a}^{(n,-k)}(L):=\hat{\mathbb{E}}^{(n,-k)}[K_{h_A}(A-a)\mid L],\;\hat{\rho}_{\pi}^{(n,-k)}(L):=\hat{\mathbb{E}}^{(n,-k)}[\pi(Z,L)\mid L].$$
				%				\State For all $k$ and $i\in I_k$, calculate $\hat{S}_{\pi,h_A}(a,O_i)=\{K_{h_A}(A_i-a) - \hat{\gamma}_{h_A,a}^{(n,-k)}(L_i)\}\{\pi(Z_i,L_i) - \hat{\rho}_{\pi}^{(n,-k)}(L_i)\}.$
				\State Calculate $\hat{S}_{\pi,h_A}(a,O_i)$ as in Algorithm~\ref{alg: test RWF discrete}.
				\State Select $T:=\lfloor \sqrt{n}\rfloor$ points $\{l_t\}_{t=1}^T$ evenly from $\mathcal{L}$.
				\State For each $t=1,\ldots,T$, calculate 
				$\hat{\zeta}_{\pi,h}^{(n)}(a,l_t)=e_1^{\top}\hat\beta_{\pi,h}^{(n)}(a,l_t)$, where $g(l)=[1,l]^{\top}$, and
				\begin{align*}
					\hat\beta_{\pi,h}^{(n)}(a,l)=& \underset{\beta \in \mathbb{R}^2}{\arg\min}
					\frac{1}{n}\sum_{i=1}^{n} 
					K_{h_L}(L_i - l)
					\left\{\hat{S}_{\pi,h_A}(a,O_i) -  g\left(\frac{L_i-l}{h_L}\right)^{\top}\beta\right\}^2.
				\end{align*}
				\State Calculate the estimated correlation matrix  $\hat{R}_{\pi,h}^{(n)}(a):=\left\{\frac{\hat{\Gamma}_{\pi,h}^{(n)}(a,l_{t_1},l_{t_2})}{\sqrt{\hat{\Gamma}_{\pi,h}^{(n)}(a,l_{t_1},l_{t_1})\hat{\Gamma}_{\pi,h}^{(n)}(a,l_{t_2},l_{t_2})}}\right\}_{t_1,t_2=1}^T$ and 
				$\hat{\sigma}_{\pi,\zeta,h}^{(n)}(a,l) := \sqrt{\hat{\Gamma}_{\pi,h}^{(n)}(a,l,l)}$, 
				where $\hat{\Gamma}_{\pi,h}^{(n)}(a,l,l') := \frac{h_Ah_L}{n}\sum_{i=1}^n\hat{W}_{\pi,h,i}(a,l)\hat{W}_{\pi,h,i}(a,l'),$ and
				\begin{align*}
					\hat{W}_{\pi,h,i}(a,l) :=& e_1^{\top}\left\{\frac{1}{n}\sum_{j=1}^nK_{h_L}(L_j - l)g\left(\dfrac{L_j-l}{h_L}\right)g\left(\dfrac{L_j-l}{h_L}\right)^{\top}\right\}^{-1} \\
					& \times K_{h_L}(L_i - l)g\left(\dfrac{L_i-l}{h_L}\right) \{\hat{S}_{\pi,h_A}(a,O_i)-g\left(\dfrac{L_i-l}{h_L}\right)^{\top}\hat\beta_{\pi,h}^{(n)}(a,l)\}.
				\end{align*}
				\State Generate $B$ i.i.d. multivariate normal samples $\{\hat{D}_b\}_{b=1}^B$, where $\hat{D}_b=[\hat{d}_{b,l_t}]_{t=1}^T$ with mean zero and covariance matrix
				$\hat{R}_{\pi,h}^{(n)}(a)$. Calculate
				$\hat{D}_{b}^*=\max\limits_{t=1,\ldots, T}\{|\hat{d}_{b,l_t}|\}$.
				\State Output $\hat{p}_{a,\pi}^{(n)} = \dfrac{1}{B}\sum\limits_{b=1}^B I\left\{\hat{D}_b^*\geq \min\limits_{t=1,\ldots,T} 
				\left|\sqrt{nh_Ah_L}\hat{\zeta}_{\pi,h}^{(n)}(a,l_t)/\hat{\sigma}_{\pi,\zeta,h}^{(n)}(a,l_t)\right|\right\}$.
		\end{algorithmic}}
	\end{algorithm}

	In particular, when the measured confounder $L$ is discrete, the testing procedure is presented in Algorithm~\ref{alg: test RWF discrete}, where $\Phi(\cdot)$ denotes the cumulative distribution function of the standard normal distribution and
	$\gamma_{h_A,a}^o(L):=\mathbb{E}\!\left[K_{h_A}(A-a)\mid L\right].$
	The procedure is motivated by the intersection-union testing principle \citep{roger1996bioequivalence}. Intuitively, we employ undersmoothing by setting the bandwidth to $h_A=n^{-1/4}\hat{\sigma}_A$. 
	\sy{The debiased estimator $\tilde{\zeta}_{\pi,h_A}^{(n)}(a,l)$ defined in Algorithm~\ref{alg: test RWF discrete} is motivated by the covariance structure that $\zeta_{\pi}^o(a,l)=\mathrm{Cov}\{p_{A|Z,L}(a| Z,L),\pi(Z,L)| L=l\}$.}
	The proposed $\tilde{\zeta}_{\pi,h_A}^{(n)}(a,l)$ consistently estimates $\zeta_{\pi}^o(a,l)$ as the bandwidth $h_A\to 0$ and $nh_A\to\infty$, while its asymptotic normality ensures the asymptotic validity of the resulting $p$-value $\tilde{p}_{a,\pi}^{(n)}$. 
	
	\sy{By contrast, when $L$ is continuous, Algorithm~\ref{alg: test RWF discrete} is no longer applicable because the sample cannot be partitioned into infinitely many groups according to the values of $L$. Motivated by \citet{takatsu2023debiasedinferencecovariateadjustedregression}, we therefore propose a bootstrap procedure for computing the $p$-value associated with the null hypothesis $ \mathcal{H}_0(a)$ in Algorithm~\ref{alg: test RWF continuous}. 
		Specifically, in Algorithm~\ref{alg: test RWF continuous}, an LLKR estimator $\hat{\zeta}_{\pi,h}^{(n)}(a,l)$ is proposed to estimate $\zeta_{\pi}^{o}(a,l)$ at a finite collection of evaluation points $\{l_t\}_{t=1}^T$. The number of evaluation points, denoted by $T$, is allowed to increase with the sample size, thereby providing an increasingly fine approximation over $\mathcal{L}$. Finally, a Gaussian parametric bootstrap is employed to approximate the distribution of the supremum of the standardized estimation error over these evaluation points, which is then used to construct the $p$-value for testing $\mathcal{H}_0(a)$.
		The validity of $p$-values proposed in Algorithms~\ref{alg: test RWF discrete},  \ref{alg: test RWF continuous} and the asymptotic statistical properties of $\tilde{\zeta}_{\pi,h_A}^{(n)}(a,l)$ and  $\hat{\zeta}_{\pi,h}^{(n)}(a,l)$ will be established in Section~\ref{subsec: RWF theorectical results}.}

	\section{Practical guidance}\label{sec: practical guidance}
	\magenta{In practice, to estimate the ADRF over a compact set $\mathcal{A}_c \subseteq \mathring{\mathcal{A}}$, one must first prespecify a collection of URWFs $\{\pi_{[m]}\}_{m=1}^M$ for subsets $\{\mathcal{N}_m\}_{m=1}^M$ covering $\mathcal{A}_c$. 
		Motivated by Proposition~\ref{prop: existence of RWFs}, one practical approach is to select points $\{a_{m}\}_{m=1}^M$ uniformly over $\mathcal{A}_c$, and estimate $\pi_{[m]}(Z,L) := \hat{\pi}_{a_{m}}(Z,L)$, serving as the RWF at $A=a_{m}$. 
		In practice, $\pi_a(Z, L):= p_{A|Z, L}(a | Z, L) / p_{A|L}(a | L)$ can be estimated by fitting smoothing spline regressions of $K_{h_A}(A-a)$ on $(Z, L)$ and $L$ respectively to obtain estimates of the two conditional densities, and then computing their ratio to obtain $\hat{\pi}_a(Z, L)$.}
	
	\magenta{Next, we introduce a graphical method for allocating a set $\mathcal{N}_m$ to each weighting function $\pi_{[m]}$. 
		First, select a sequence of points $\{a_j^*\}_{j=1}^J$ evenly spaced over $\mathcal{A}_c$, and apply the RWF testing procedure in Algorithm~\ref{alg: test RWF discrete} (or Algorithm~\ref{alg: test RWF continuous}) for each prespecified $\pi_{[m]}$ and $a_j^*$. 
		Then, plot $M$ curves, each displaying the $p$-values $\{\tilde{p}_{a_j^*,\pi_{[m]}}^{(n)}\}_{j=1}^J$ (or $\{\hat{p}_{a_j^*,\pi_{[m]}}^{(n)}\}_{j=1}^J$) for a fixed $\pi_{[m]}$ (see, e.g., Figures~\ref{fig: JTPA test RWF},~\ref{fig: RWF testing}). 
		Based on the resulting $p$-value plot, one may treat any $\{\pi_{[m]}, \mathcal{N}_m\}_{m=1}^M$ as a candidate construction, provided that 
		(1) $\bigcup_{m=1}^M \mathcal{N}_m$ covers $\mathcal{A}_c$, and 
		(2) for each $m=1,\ldots,M$, the $p$-values $\tilde{p}_{a_j^*,\pi_{[m]}}^{(n)}$ (or $\hat{p}_{a_j^*,\pi_{[m]}}^{(n)}$) are below the prespecified significance level $\alpha$ if $a_j^*\in \mathcal{N}_m$.}
	
	\magenta{In principle, $\{\pi_{[m]}, \mathcal{N}_m\}_{m=1}^M$ can be constructed from the $p$-value plot, provided the required conditions are satisfied. For convenience, we also provide a simple but useful method to facilitate this construction. Specifically, when $\pi_{a_m}$ is approximated well, by Proposition~\ref{prop: local stability of RWF}, $\pi_{[m]}$ can serve as the URWF over $\mathcal{N}_m= \overline{\mathcal{B}(a_m,r_0)}$ for sufficiently small $r_0$, thus satisfying the second condition. On the other hand, $\bigcup_{m=1}^M\mathcal{N}_m$ can cover $\mathcal{A}_c$ when  $\{a_m\}_{m=1}^M$ is sufficiently dense in $\mathcal{A}_c$.} 
	Finally, it should be emphasized that the $p$-value plot serves only as a data-guided diagnostic for the RWF condition.

	\magenta{
		With $\{\pi_{[m]},\mathcal{N}_m\}_{m=1}^M$ treated as fixed, Algorithm~\ref{alg: cross-fitting} yields the AIPW score vector, based on which the LLKR estimator in Equation~\eqref{eq: LLKR} is constructed. Finally, following \citet{bonvini2022fastconvergenceratesdoseresponse}, we note that the AIPW score can also be used to estimate the ADRF via empirical risk minimization (ERM) methods, rather than solely through LLKR. See Appendix~\ref{append: oracle bound for ERM} for further discussion.
	}

	\section{Asymptotic theory}\label{sec: asymptotic theory}
	\subsection{The ADRF estimator}
	We first establish the asymptotic theory for the LLKR estimator $\hat{\theta}_{\pi,h}^{(n)}(a)$. For any measurable function $g(L)$, $f(A,L)$, and subset $\mathcal{N}\subseteq\mathring{\mathcal{A}}$, define the norm 
	$\|g\|_2^2:= \int_{\mathcal{L}}|g(l)|^2\mathrm{d} P_L(l)$, 
	$\|f\|_{\mathcal{N},2}^2:=\int_{\mathcal{L}} \sup_{a\in \mathcal{N}}|f(a,l)|^2\mathrm{d} P_L(l),$
	where $P_L(l)$ is the marginal distribution of $L$. Define the norm for the nuisance vector $\alpha_{\pi}$ as
	$\|\alpha_{\pi}\|_{\mathcal{N},2}^2:=
	\|\rho_{\pi}\|_2^2
	+\|\mu_{\pi}\|_{\mathcal{N},2}^2
	+\|\kappa_{\pi}\|_{\mathcal{N},2}^2
	+\|\eta\|_{\mathcal{N},2}^2
	+\|\delta\|_{\mathcal{N},2}^2.$
	Next, we introduce several common regularity conditions in nonparametric analysis. 
	\begin{assumption}\label{as: regular condition} 
		Let $\mathcal{N}\subseteq \mathring{\mathcal{A}}$ be an interval. Assume the following conditions hold: 
		\begin{enumerate}[label=(\alph*), itemsep=0pt, topsep=0pt,leftmargin=0.75cm]
			\item $\mathcal{Z}$ and $\mathcal{Y}$ are bounded subsets of the Euclidean spaces. \label{cd: ZY boundedness}
			
			\item The bandwidth $h=h_n$ satisfies $h \to 0$ and $n h \to +\infty$ as $n \to +\infty$. \label{cd: bandwidth rate}
			
			\item $K(s)$ is a positive, continuous, symmetric function with support $[-1,1]$, satisfying 
			$\int_{\mathbb{R}} K(s)\, \mathrm{d}s = 1$, 
			$\int_{\mathbb{R}} K(s) s\, \mathrm{d}s = 0$, 
			and $\int_{\mathbb{R}} K(s) s^2\, \mathrm{d}s > 0$. \label{cd: kernel condition}
			
			\item $\theta(a)$ is twice continuously differentiable for all $a \in \mathcal{N}$. \label{cd: theta condition}
			
			\item The conditional variance $\sup_{a\in\mathcal{N}}\{\mathrm{Var}[\varphi_{\pi}(O;\alpha_{\pi}^o,P_O) \mid A=a]\}<\infty$. \label{cd: varphi}
			
			\item For each $k=1,\ldots,K$ in Algorithm~\ref{alg: cross-fitting}, $|I_{-k}^{(2)}| \to \infty$ as $n \to \infty$. \label{cd: plnk rate}
			
			\item For each $k=1,\ldots,K$, there exists a constant $\epsilon_0>0$ such that for all $n$, 
			$|\hat\kappa_{\pi}^{(n,-k)}(a,l)|\ge \epsilon_0$ for all $(a,l) \in \mathcal{N}\times \mathcal{L}$ almost surely. Moreover, all functions in $\hat\alpha_{\pi}^{(n,-k)}(A,L)$ are uniformly bounded across all $n,k$ pairs almost surely. \label{cd: alpha uniform boundedness}
			
			\item For each $k=1,\ldots,K$ and $a \in \mathcal{N}$, $\mathbb{E}\big[\|\hat\alpha_{\pi}^{(n,-k)} - \alpha_{\pi}^o\|_{\mathcal{B}(a,h),2}^2\big] = o(1)$. Moreover, there exists a sequence $c(n)$ such that
			$\|\hat\rho_{\pi}^{(n,-k)}-\rho_{\pi}^o\|_2\|\hat\eta^{(n,-k)}-\eta^o\|_{\mathcal{B}(a,h),2}+$
			$\|\hat\rho_{\pi}^{(n,-k)}-\rho_{\pi}^o\|_2\|\hat\delta^{(n,-k)}-\delta^o\|_{\mathcal{B}(a,h),2}+$
			$\|\hat\mu_{\pi}^{(n,-k)}-\mu_{\pi}^o\|_{\mathcal{B}(a,h),2}  \|\hat\kappa_{\pi}^{(n,-k)}-\kappa_{\pi}^o\|_{\mathcal{B}(a,h),2}+$
			$\|\hat\mu_{\pi}^{(n,-k)}-\mu_{\pi}^o\|_{\mathcal{B}(a,h),2}  \|\hat\delta^{(n,-k)}-\delta^o\|_{\mathcal{B}(a,h),2}=o_p(c(n))$.\label{cd: alpha rate}
		\end{enumerate}
	\end{assumption}
	We provide a brief overview of the regularity conditions in Assumption~\ref{as: regular condition}. Conditions~\ref{cd: ZY boundedness}--\ref{cd: theta condition} are standard in the kernel regression literature \citep{tsybakov2009introduction}. Condition~\ref{cd: bandwidth rate} encodes the classical bias-variance trade-off: $h \to 0$ ensures asymptotic unbiasedness, while $n h \to \infty$ ensures vanishing variance of order $1/(n h)$. Condition~\ref{cd: theta condition} imposes smoothness on $\theta(a)$, controlling the estimation bias.
	Condition~\ref{cd: varphi} controls the variance of $\hat{\theta}_{\pi,h}^{(n)}(a)$, which is a relatively mild requirement. 
	% From Equation~\eqref{eq: varphi defn}, under Assumptions~\ref{as: continuity}--\ref{as: IV relevance}, if $\mathbb{E}[Y^2 \mid A=a,L]$ and $\mathbb{E}[Z^2 \mid A=a,L]$ are uniformly bounded in $L$, and $\pi(Z,L)$ is a URWF for $\mathcal{N}$, then Condition~\ref{cd: varphi} is automatically satisfied. 
	Condition~\ref{cd: plnk rate} requires that the sample size used for training $\hat{P}_O^{(n,-k)}$ grows sufficiently fast to ensure consistency. 
	Condition~\ref{cd: alpha uniform boundedness} imposes stability and regularity on the nuisance functions in $\hat{\alpha}_{\pi}^{(n,-k)}$. 
	Finally, Condition~\ref{cd: alpha rate} specifies the convergence rate of $\hat{\alpha}_{\pi}^{(n,-k)}$ towards $\alpha_{\pi}^o$.
	Next, we clarify the convergence rates of $\hat{\theta}_{\pi,h}^{(n)}(a)$ as follows.
	\begin{theorem}[Convergence rate]\label{thm: convergence rate}
		Under Assumptions~\ref{as: consistency}--\ref{as: positivity}, suppose that $\pi(Z,L)$ is a URWF for an interval $\mathcal{N}$ satisfying Assumption~\ref{as: regular condition}. 
		Then, if $Z$ is an AIV for all $A=a \in \mathcal{N}$, the estimator $\hat{\theta}_{\pi,h}^{(n)}(a)$ obtained from Algorithm~\ref{alg: cross-fitting} satisfies for all $a\in\mathcal{N}$,
		\begin{align}\label{eq: convergence rate}
			\hat{\theta}_{\pi,h}^{(n)}(a) - \theta(a) = O_p\Big(\frac{1}{\sqrt{nh}}\Big) + O_p(h^2) + o_p(c(n)).
		\end{align}
	\end{theorem}
	We explain each term in Equation~\eqref{eq: convergence rate} as follows.
	The first term corresponds to the variance term in LLKR; the second term is the bias term in LLKR; the third term represents the mixed bias caused by plugging in $\hat\alpha_{\pi}^{(n,-k)}$ into the estimation.
	Intuitively, by taking $h = n^{-1/5}$ and $c(n) = O(n^{-2/5})$, one obtains $|\hat{\theta}_{\pi,h}^{(n)}(a) - \theta(a)| = O_p(n^{-2/5}).$ 
	In fact, this convergence rate achieves the oracle minimax lower bound for LLKR.
	In particular, a sufficient condition for $c(n)=O(n^{-2/5})$ is that $\|\hat\alpha_{\pi}^{(n,-k)}-\alpha_{\pi}^o\|_{\mathcal{B}(a,h),2}=o_p(n^{-1/5})$, which is a slightly weaker condition than the convergence rate $o_p(n^{-1/4})$ in the DML literature \citep{Chernozhukov2018}.

	\begin{theorem}[Asymptotic normality of $\hat{\theta}_{\pi,h}^{(n)}(a)$]\label{thm: clt}
		Under the conditions of Theorem~\ref{thm: convergence rate}, if $c(n) = O(1/\sqrt{nh})$, then for any $a \in \mathcal{N}$ that satisfies Assumption~\ref{as: regular condition}, we have
		\begin{equation}\label{eq: CLT}
			\begin{aligned}
				& \sqrt{nh} \left\{\hat{\theta}_{\pi,h}^{(n)}(a) - \theta(a) - \mathrm{bias}(a) \right\}
				\overset{d}{\to} N(0, \sigma_{\pi,\theta}^o(a)^2), \quad \text{where} \\
				& \sigma_{\pi,\theta}^o(a)^2 := \frac{\int K(s)^2 \, \mathrm{d}s}{p_A(a)} \, \mathrm{Var}\left[\varphi_{\pi}(O;\alpha_{\pi}^o, P_O) \mid A=a\right].
			\end{aligned}
		\end{equation}
		In Equation~\eqref{eq: CLT}, the bias term $\mathrm{bias}(a)$, defined in Appendix~\ref{append: proof of theorems}, satisfies
		\begin{align*}
			\mathrm{bias}(a) = \frac{h^2}{2} \theta''(a) \int K(s) s^2 \, \mathrm{d}s + o_p(h^2).
		\end{align*}
	\end{theorem}
	This theorem provides a valid variance representation $\sigma_{\pi,\theta}^o(a)^2$ for $\hat{\theta}_{\pi,h}^{(n)}(a)$. 
	Specifically, when there exist consistent estimators for $p_A(a)$ and $\mathrm{Var}[\varphi_{\pi}(O;\alpha_{\pi}^o, P_O)\mid A=a]$, one can construct a consistent variance estimator $\hat{\sigma}_{\pi,\theta}^{(n)}(a)$ for $\sigma_{\pi,\theta}^o(a)$. 
	If $h = o(n^{-1/5})$ and $h^{-1} = o(n)$, the bias term in Equation~\eqref{eq: CLT} can be neglected. In this setting, a $95\%$ confidence interval for $\theta(a)$ can be formulated as $\hat{\theta}_{\pi,h}^{(n)}(a) \pm 1.96 \hat{\sigma}_{\pi,\theta}^{(n)}(a)/\sqrt{nh}.$
	\magenta{Unfortunately, this guarantee does not hold when the bandwidth is selected adaptively via LOOCV. Relatedly, \citet{takatsu2023debiasedinferencecovariateadjustedregression} propose a debiased inference approach to mitigate the $\mathrm{bias}(a)$, which is beyond the scope of this article.}
	\sy{
		Next, we show a generalization of Theorem~\ref{thm: clt}, which uses multiple RWFs to construct a debiased GMM estimator for the ADRFs.
		\begin{theorem}[Asymptotic normality of $\hat{\beta}_{\Pi,h}^{(n)}(a)$]\label{thm: falsify AIV condition}
			Under Assumptions~\ref{as: consistency}--\ref{as: positivity}, let $\{\pi_j\}_{j=1}^J$ be $J$ different URWFs for an interval $\mathcal{N}\subseteq \mathring{\mathcal{A}}$ satisfying Assumption~\ref{as: regular condition}, and $Z$ be an AIV for $A$. 
			Suppose that $\mathrm{Var}\left\{\varphi_{\Pi}(O; \alpha_{\Pi}^o, P_O)| A=a\right\}$ is positive definite for all $a\in\mathcal{N}$.
			Then, for all $a\in\mathcal{N}$, if $c(n)=o(1/\sqrt{nh})$, $h=o(n^{-1/5})$, we have $\sqrt{nh}\{\hat{\beta}_{\Pi,h}^{(n)}(a) - \beta_{h}^o(a)\}\overset{d}{\to} N(0,\Sigma_{\Pi,\beta}^o(a)^2)$, where $\beta_{h}^o(a)=[\theta(a),h\theta'(a)]^{\top}$, and 
			\begin{align*}
				\Sigma_{\Pi,\beta}^o(a)^2:=\frac{\mathrm{diag}\left\{
					\int K\left(s\right)^2 \mathrm{d}s,\; \int K\left(s\right)^2 s^2\mathrm{d}s/(\int K(s)s^2\mathrm{d}s)^2\right\}}{p_A(a)\mathbf{1}_J^{\top}
					\mathrm{Var}\left\{\varphi_{\Pi}(O; \alpha_{\Pi}^o, P_O)| A=a\right\}^{-1}\mathbf{1}_J}\in\mathbb{R}^{2\times 2}.
			\end{align*}
			In addition, for $J>1$, $nh\mathcal{L}_{\Pi,h,a}^{(n)}(\hat{\beta}_{\Pi,h}^{(n)}(a),\hat{Q}_{\Pi,h}^{(n)}(a)^{-1})\overset{d}{\to}\chi^2_{2J-2}$.
		\end{theorem}
		Specifically, when $\Pi=\{\pi\}$, the first component of the asymptotic variance $\Sigma_{\Pi,\beta}^o(a)^2$ reduces to $\sigma_{\pi,\theta}^o(a)^2$, thereby recovering the result in Theorem~\ref{thm: clt}. Moreover, when $J>1$, Theorem~\ref{thm: falsify AIV condition} provides the theoretical justification for the falsification test introduced in Section~\ref{subsec: falsification test}.
	}

	\subsection{The RWF testing}\label{subsec: RWF theorectical results}
	In this subsection, we establish the validity of the $p$-values produced by the RWF testing procedures described in Algorithms~\ref{alg: test RWF discrete} and~\ref{alg: test RWF continuous}.
	First, we provide a theoretical analysis of the estimator $\tilde{\zeta}_{\pi,h_A}^{(n)}(a,l)$ when $\mathcal{L}$ is discrete.
	\begin{theorem}[Asymptotic normality of $\tilde{\zeta}_{\pi,h_A}^{(n)}(a,l)$]
		\label{thm: RWF testing 1}
		Under Assumptions~\ref{as: regular condition}(a), (c), assume that $\nabla_a^2 p_{A|Z,L}(a|z,l)$ exists at $a$ for any $(z,l)\in\mathcal{Z}\times\mathcal{L}$, where $L$ is discrete. 
		Moreover, assume that $nh_A\to\infty$ , $h_A=o(1)$, 
		$\mathbb{E}[\|\hat{\gamma}_{h_A,a}^{(n,-k)}-\gamma_{h_A,a}^o\|_2^2]=o(h_A^{-1/2})$, 
		$\mathbb{E}[\|\hat{\rho}_{\pi}^{(n,-k)}-\rho_{\pi}^o\|_2^2] = o(1)$, and 
		$\|\hat{\gamma}_{h_A,a}^{(n,-k)}-\gamma_{h_A,a}^o\|_2 \, \|\hat{\rho}_{\pi}^{(n,-k)}-\rho_{\pi}^o\|_2 = o_p(1/\sqrt{nh_A})$.
		Then, we have 
		$\tilde{\zeta}_{\pi,h_A}^{(n)}(a,l) = \zeta_{\pi}^o(a,l) + O_p(1/\sqrt{nh_A} + h_A^2).$
		Define $\mathrm{bias}_{\zeta}(a,l) := \frac{h_A^2}{2} \, \mathbb{E}[p_{A|Z,L}''(a| Z,L)(\pi(Z,L)-\rho_{\pi}^o(L))\mid L=l] \int K(x) x^2 \, \mathrm{d}x.$ If $h_A = O(n^{-1/5})$,  then
		\begin{align*}
			&\sqrt{nh_A} \big\{\tilde{\zeta}_{\pi,h_A}^{(n)}(a,l) - \zeta_{\pi}^o(a,l) - \mathrm{bias}_{\zeta}(a,l)\big\} \overset{d}{\to} N(0,\sigma_{\pi,\zeta}^*(a,l)^2), 
			\quad\text{where}\\
			&\sigma_{\pi,\zeta}^*(a,l)^2 := \dfrac{p_{A|L}(a\mid l)}{p_L(l)}\int K(x)^2\mathrm{d}x\cdot  \mathbb{E}[\{\pi(Z,L)-\rho_{\pi}^o(L)\}^2|A=a,L=l].
		\end{align*}
	\end{theorem}
	In Theorem~\ref{thm: RWF testing 1}, if $h_A=o(n^{-1/5})$, the bias term $\mathrm{bias}_{\zeta}(a,l)$ is negligible. Therefore, we adopt the convenient undersmoothing choice $h_A=n^{-1/4}\hat{\sigma}_A$ in Algorithm~\ref{alg: test RWF discrete}. We estimate $\sigma_{\pi,\zeta}^*(a,l)$ using $\tilde{\sigma}_{\pi,\zeta,h_A}^{(n)}(a,l)$, the conditional sample standard deviation of the scores $\{\hat{S}_{\pi,h_A}(a,O_i)\}_{i=1}^n$ in Algorithm~\ref{alg: test RWF discrete}.
	Because $L$ has a finite discrete support, the pointwise asymptotic theory is sufficient to justify the validity of the intersection–union $p$-value $\tilde{p}_{a,\pi}^{(n)}$ proposed in Algorithm~\ref{alg: test RWF discrete}.
	
	\sy{
		By contrast, when $L$ is continuous, we provide theoretical guarantees for the validity of the $p$-value $\hat{p}_{a,\pi}^{(n)}$ proposed in Algorithm~\ref{alg: test RWF continuous}. The analysis is more involved than that in the discrete case because it requires a uniform approximation to the entire estimated process over $l\in\mathcal{L}$.
		
		\begin{theorem}[Asymptotic normality of $\hat{\zeta}_{\pi,h}^{(n)}(a,l)$]
			\label{thm: RWF testing 2}
			Under Assumption~\ref{as: regular condition}(a), (c) with $h=[h_A,h_L]$, assume that $\nabla_a^2 p_{A|Z,L}(a|z,l)$ exists at $a$ for any $(z,l)\in\mathcal{Z}\times\mathcal{L}$, where $\mathcal{L}\subseteq\mathbb{R}$ is a compact subset, and that $\zeta_{\pi}^o(a,l)$ is twice continuously differentiable with respect to $(a,l)\in\mathcal{A}\times\mathcal{L}$. 
			Also, assume that $nh_Ah_L\to\infty$, $h_A,h_L=o(1)$,
			$\mathbb{E}[\|\hat{\gamma}_{h_A,a}^{(n,-k)}-\gamma_{h_A,a}^o\|_\infty^2]=o(h_A^{-1/2})$, 
			$\mathbb{E}[\|\hat{\rho}_{\pi}^{(n,-k)}-\rho_{\pi}^o\|_\infty^2] = o(1)$, and 
			$\|\hat{\gamma}_{h_A,a}^{(n,-k)}-\gamma_{h_A,a}^o\|_\infty \, \|\hat{\rho}_{\pi}^{(n,-k)}-\rho_{\pi}^o\|_\infty = o_p(1/\sqrt{n h_Ah_L})$.
			Then, we have 
			$\hat{\zeta}_{\pi,h}^{(n)}(a,l) = \zeta_{\pi}^o(a,l) + O_p(1/\sqrt{nh_Ah_L} + h_A^2 + h_L^2).$ 
			If $h_A,h_L = o(n^{-1/6})$, then for all fixed  $a\in\mathring{\mathcal{A}}$ and $l\in\mathring{\mathcal{L}}$,
			$\sqrt{n h_Ah_L} \big\{\hat{\zeta}_{\pi,h}^{(n)}(a,l) - \zeta_{\pi}^o(a,l) \big\} \overset{d}{\to} N(0,\sigma_{\pi,\zeta}^o(a,l)^2),$ where
			\begin{align*}
				&\sigma_{\pi,\zeta}^o(a,l)^2 := \frac{p_{A|L}(a\mid l)}{p_L(l)} \left(\int K(x)^2\mathrm{d} x\right)^2
				\mathbb{E}[\{\pi(Z,L) - \rho_{\pi}^o(L)\}^2\mid A=a,L=l].
			\end{align*}
		\end{theorem}
		Theorem~\ref{thm: RWF testing 2} establishes the convergence rate and pointwise asymptotic normality of the estimator $\hat{\zeta}_{\pi,h}^{(n)}(a,l)$ when $L\in\mathbb{R}$ is continuous. The result can also be extended to the case when $\mathcal{L}\subseteq\mathbb{R}^{d_L}$. In particular, the asymptotic normality result established in Theorem~\ref{thm: RWF testing 2} continues to hold, provided that $h_A,h_L=o\left(n^{-1/(5+d_L)}\right).$
		
		However, because $L$ is continuous, pointwise asymptotic normality alone is insufficient to establish the asymptotic validity of the $p$-value proposed in Algorithm~\ref{alg: test RWF continuous}. Instead, we need to establish the validity of a simultaneous confidence band for $\{\zeta_{\pi}^o(a,l):l\in\mathcal{L}\}$, which provides the theoretical foundation for constructing the proposed $p$-value.
		
		\begin{theorem}[Simultaneous confidence band of $\hat{\zeta}_{\pi,h}^{(n)}(a,l)$]\label{thm: RWF testing 3}
			Under the conditions of Theorem~\ref{thm: RWF testing 2}, assume that 
			(a) $\mathbb{E}[\|\hat{\gamma}_{h_A,a}^{(n,-k)} - \gamma_{h_A,a}^o\|_\infty^2] = o(\{h_A\}^{-1/2}/\log(n)^3)$, 
			$\mathbb{E}[\|\hat{\rho}_{\pi}^{(n,-k)} -\rho_{\pi}^o\|_\infty^2] = o(\log (n)^{-3})$;
			(b) $\|\hat{\gamma}_{h_A,a}^{(n,-k)}-\gamma_{h_A,a}^o\|_\infty \, \|\hat{\rho}_{\pi}^{(n,-k)}-\rho_{\pi}^o\|_\infty = o_p(\{n h_Ah_L\log (n)\}^{-1/2})$;
			(c) $h_A,h_L = O(n^{-c})$ for some $c>1/6$,
			$T^{-1}=o(h_L^{1+p})$, and $T=O(n^q)$ for some $p,q>0$;
			(d) $\sup_{l,l'\in \mathcal{L} }|\hat{\Gamma}_{\pi,h}^{(n)}(a,l,l') - \Gamma_{\pi,h}^o(a,l,l')|=o_p(\log(n)^{-2})$, where $ S_{\pi,h_A}^o(a,O_i)=(K_{h_A}(A_i-a) - \gamma_{h_A,a}^o(L_i))(\pi(Z_i,L_i) - \rho_{\pi}^o(L_i))$, and
			\begin{align*}
				\Gamma_{\pi,h}^o(a,l,l'):=\dfrac{h_Ah_L}{p_L(l)p_L(l')}\mathbb{E}\left[K_{h_L}(L-l)K_{h_L}(L-l')\left\{S_{\pi,h_A}^o(a,O)-\mathbb{E}[S_{\pi,h_A}^o(a,O)\mid L]\right\}^2\right].
			\end{align*}
			Then, $\sup_{l\in\mathcal{L}}|\hat{\zeta}_{\pi,h}^{(n)}(a,l) - \zeta_{\pi}^o(a,l)| = O_p(\sqrt{\frac{\log(n)}{nh_Ah_L}})$, and
			\begin{align*}
				\sup_{x\in\mathbb{R}}\left|
				\Pr (\sup_{l\in\mathcal{L}}\sqrt{nh_Ah_L}\left|
				\frac{\hat{\zeta}_{\pi,h}^{(n)}(a,l)- \zeta_{\pi}^o(a,l)}{\hat{\sigma}_{\pi,\zeta,h}^{(n)}(a,l)}
				\right| \leq x)-
				\Pr\left(\max\limits_{t=1,\ldots,T}|\hat{d}_{b,l_t}|\leq x\middle| \{O_i\}_{i=1}^n\right)
				\right|=o_p(1).
			\end{align*}
		\end{theorem}
		
		The proof of Theorem~\ref{thm: RWF testing 3} relies primarily on the Gaussian approximation techniques developed by \citet{chernozhukov2014gaussian}.
		Intuitively, Theorem~\ref{thm: RWF testing 3} requires a trade-off in the choice of the number of evaluation points $T$. 
		Specifically, $T$ cannot be too large, since taking the maximum over an increasingly fine grid inflates the stochastic variability of the resulting process. 
		On the other hand, $T$ cannot be too small, as a coarse grid approximation may introduce a non-negligible discretization bias. 
		
		Theorem~\ref{thm: RWF testing 3} further implies that a simultaneous $1-\alpha$ confidence band for $\zeta_{\pi}^o(a,l)$ can be constructed as $[\hat{\zeta}_{\pi,h}^{(n)}(a,l) - \frac{\hat{c}_{1-\alpha}}{\sqrt{nh_Ah_L}}\hat\sigma_{\pi,\zeta,h}^{(n)}(a,l), \hat{\zeta}_{\pi,h}^{(n)}(a,l) + \frac{\hat{c}_{1-\alpha}}{\sqrt{nh_Ah_L}} \hat\sigma_{\pi,\zeta,h}^{(n)}(a,l)]$, where $\hat{c}_{1-\alpha}$ is the $(1-\alpha)$-th quantile of $\{\hat{D}_{b}^*\}_{b=1}^B$, and $\hat{D}_{b}^*=\max\limits_{t=1,\ldots, T}\{|\hat{d}_{b,l_t}|\}$.
		Next, based upon the validity of this simultaneous confidence band, we establish the validity of the $p$-value $\hat{p}_{a,\pi}^{(n)}$ proposed in Algorithm~\ref{alg: test RWF continuous} as follows. 
		
		\begin{theorem}[Validity of $\hat{p}_{a,\pi}^{(n)}$]
			\label{thm: RWF testing 4}
			Under the conditions of Theorem~\ref{thm: RWF testing 3}, suppose that $B\to\infty$ as $n\to\infty$. Then, for any distribution $P_0$ satisfying $\mathcal{H}_0(a)$,   $\underset{n\to+\infty}{\overline{\lim}}\Pr_{P_0}\left(\hat{p}_{a,\pi}^{(n)}\leq \alpha\right)\leq \alpha$; for any distribution $P_1$ satisfying $\mathcal{H}_1(a)$, $\underset{n\to+\infty}{\lim}\Pr_{P_1}\left(
			\hat{p}_{a,\pi}^{(n)}\leq \alpha\right)=1$.
		\end{theorem}
		Theorem~\ref{thm: RWF testing 4} shows that $\hat{p}_{a,\pi}^{(n)}$ is a pointwise asymptotically valid $p$-value. Moreover, the power of the resulting test converges to one under every fixed alternative distribution $P\in \mathcal{H}_1(a)$.}

	%Concretely, the distribution of $A\mid Z,U,L$ satisfies that
	%\begin{align*}
	%	p_{A|Z,U,L}(a|Z,U,L)=&0.5 p_{A|Z,U,L,\varepsilon_A}(a|Z,U,L,1) 
	%	+ 0.5 p_{A|Z,U,L,\varepsilon_A}(a|Z,U,L,0) \\
	%	=&0.5 p_{A|Z,L,\varepsilon_A}(a|Z,U,L,1) 
	%	+ 0.5 p_{A|U,L,\varepsilon_A}(a|Z,U,L,0).
	%\end{align*}
	\section{Simulations}\label{sec: simulations}
	We conduct simulation studies to investigate the finite-sample performance of the proposed methods. 
	Specifically, we simulate independent random variables 
	$\varepsilon_L, \varepsilon_U, \varepsilon_Z\sim \mathrm{Unif}(0,1)$, $\varepsilon_{AZ},\varepsilon_{AU}\sim N(0,1)$, and $\varepsilon_A \sim \mathrm{Ber}(0.7)$. 
	We then generate the observed and latent variables according to  
	$L = \varepsilon_L - 0.5$,  
	$U = 3(\varepsilon_U - 0.5)$, 
	$Z = -0.5L + 3(\varepsilon_Z - 0.5)$, 
	$\mathrm{logit}(p_z) = 2Z+\varepsilon_{AZ}$, 
	$\mathrm{logit}(p_u) = -2U+\varepsilon_{AU}$, 
	$A = 2\varepsilon_A p_z+ 2(1-\varepsilon_A)p_u - 1$, 
	$Y = A + U - 0.5 L.$ 
	
	By construction, $\mathcal{A} = [-1,1]$, and $Z$ serves as an AIV
	for each treatment level $a \in \mathcal{A}$. We focus on a target subinterval 
	$\mathcal{A}_c := [-0.75, 0.75]$. 
	\magenta{As described in Section~\ref{sec: practical guidance}, we construct a cover of $\mathcal{A}_c$ by defining $\mathcal{N}_m:= [0.5m - 1.25, 0.5m - 0.75]$ for $m = 1,2,3$. For each $m$, we generate 10{,}000 additional samples to construct $\pi_{[m]} = \hat{\pi}_{a_{m}}(Z,L)$ at $a_{m} = 0.5m - 1$. %, with the intention that $\pi_{[m]}$ serves as the URWF for $\mathcal{N}_m$. 
		We present an analysis of the RWF testing procedure (Algorithm~\ref{alg: test RWF continuous}) in Appendix~\ref{append: additional simulation results}, demonstrating that $\pi_{[m]}$ provides a suitable choice of URWF for $\mathcal{N}_m$, for $m = 1, 2, 3$.}

	\begin{figure}[t]
		\centering
		\begin{subfigure}{0.48\textwidth}
			\centering
			\includegraphics[width=\linewidth]{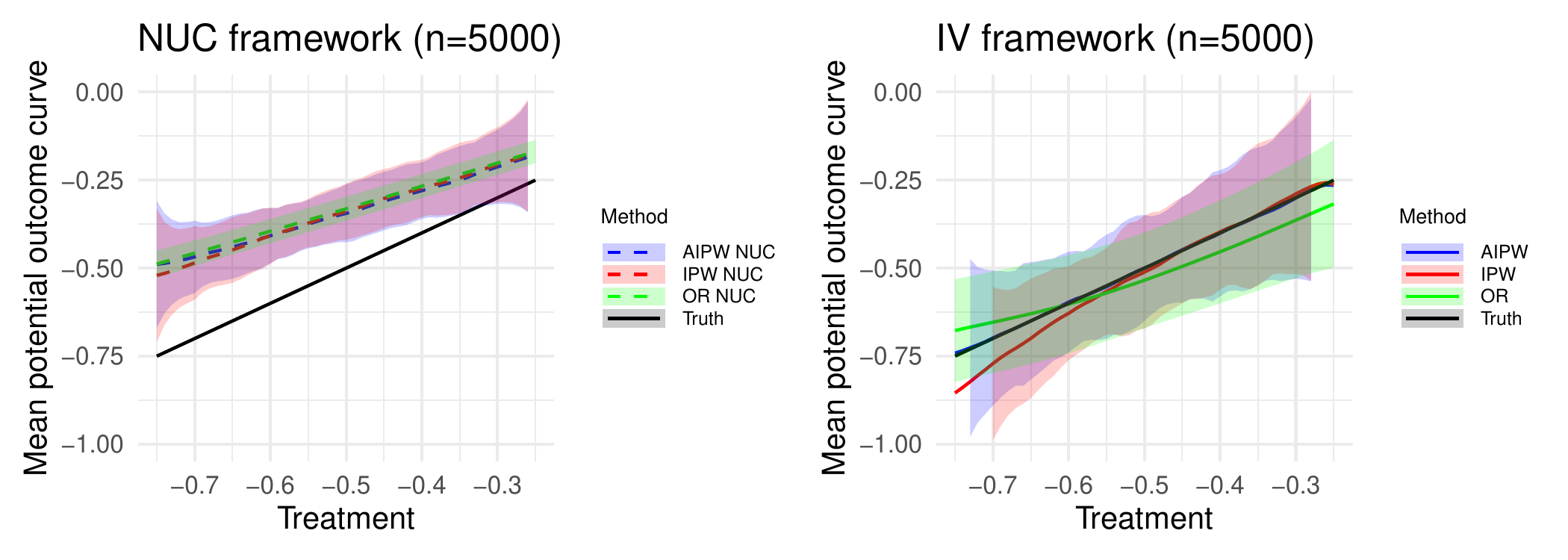}
			\caption{$\mathcal{N}_1 = [-0.75,-0.25],\,n=5000$}
		\end{subfigure}
		\hfill
		\begin{subfigure}{0.48\textwidth}
			\centering
			\includegraphics[width=\linewidth]{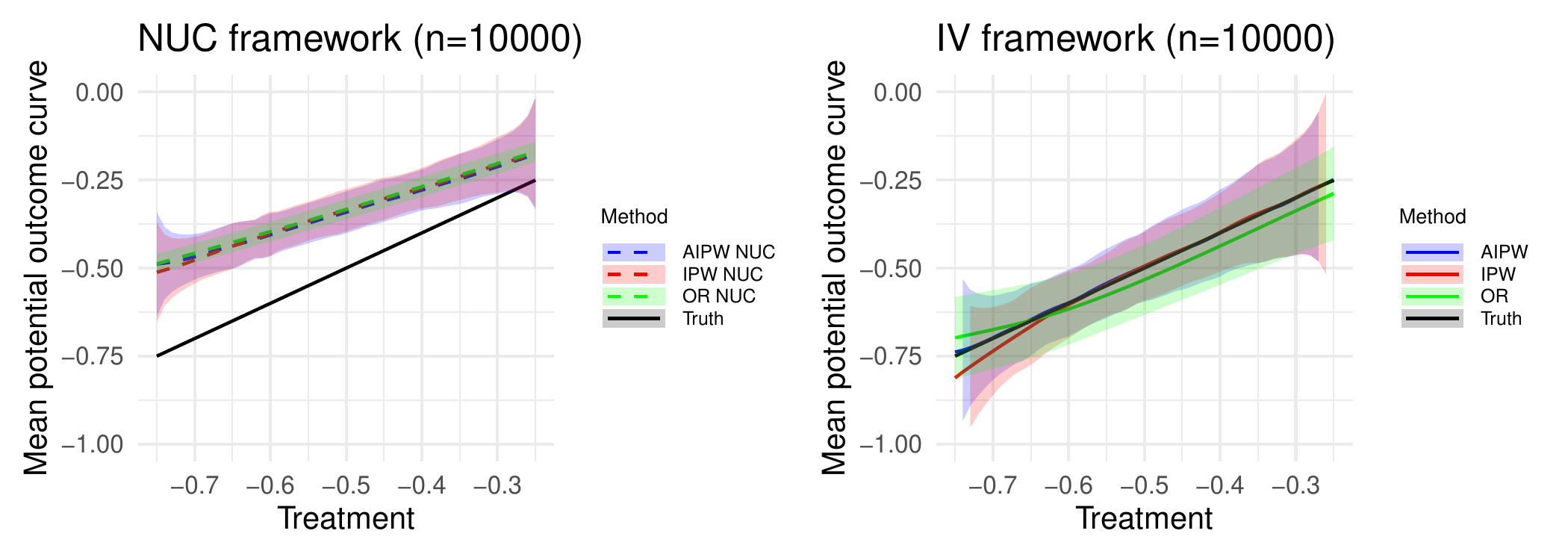}
			\caption{$\mathcal{N}_1 = [-0.75,-0.25],\,n=10000$}
		\end{subfigure}
		\begin{subfigure}{0.48\textwidth}
			\centering
			\includegraphics[width=\linewidth]{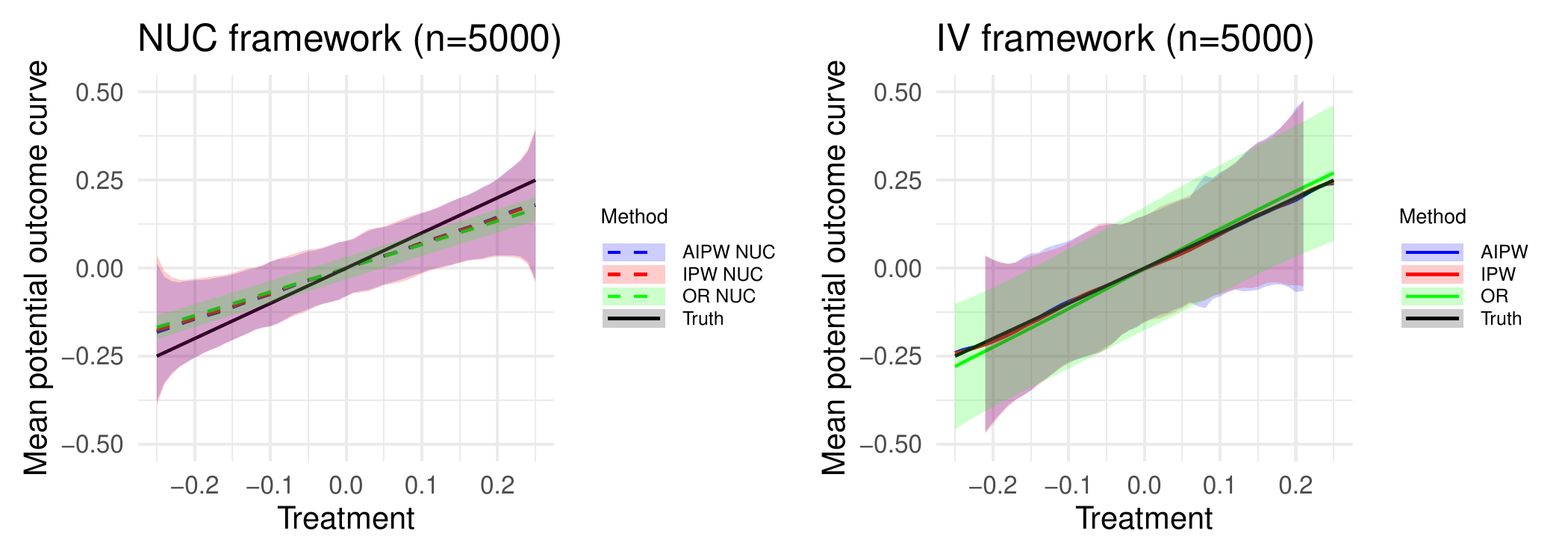}
			\caption{$\mathcal{N}_2 = [-0.25,0.25],\,n=5000$}
		\end{subfigure}
		\hfill
		\begin{subfigure}{0.48\textwidth}
			\centering
			\includegraphics[width=\linewidth]{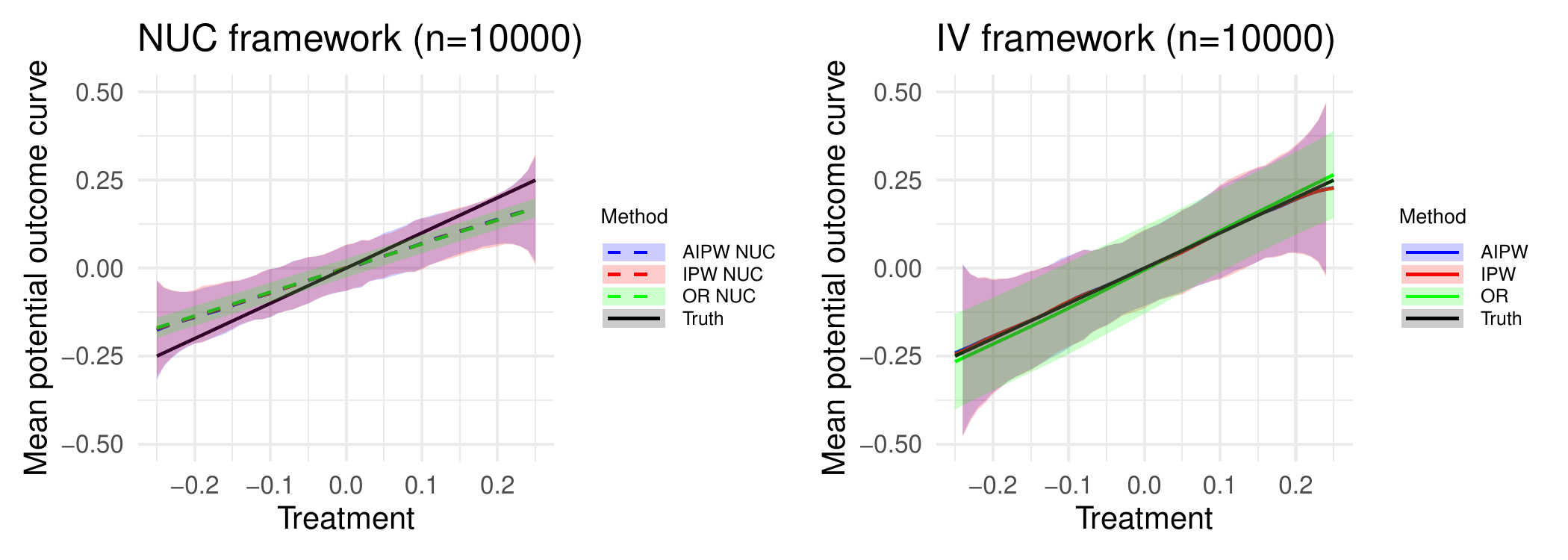}
			\caption{$\mathcal{N}_2 = [-0.25,0.25],\,n=10000$}
		\end{subfigure}
		\hfill
		\begin{subfigure}{0.48\textwidth}
			\centering
			\includegraphics[width=\linewidth]{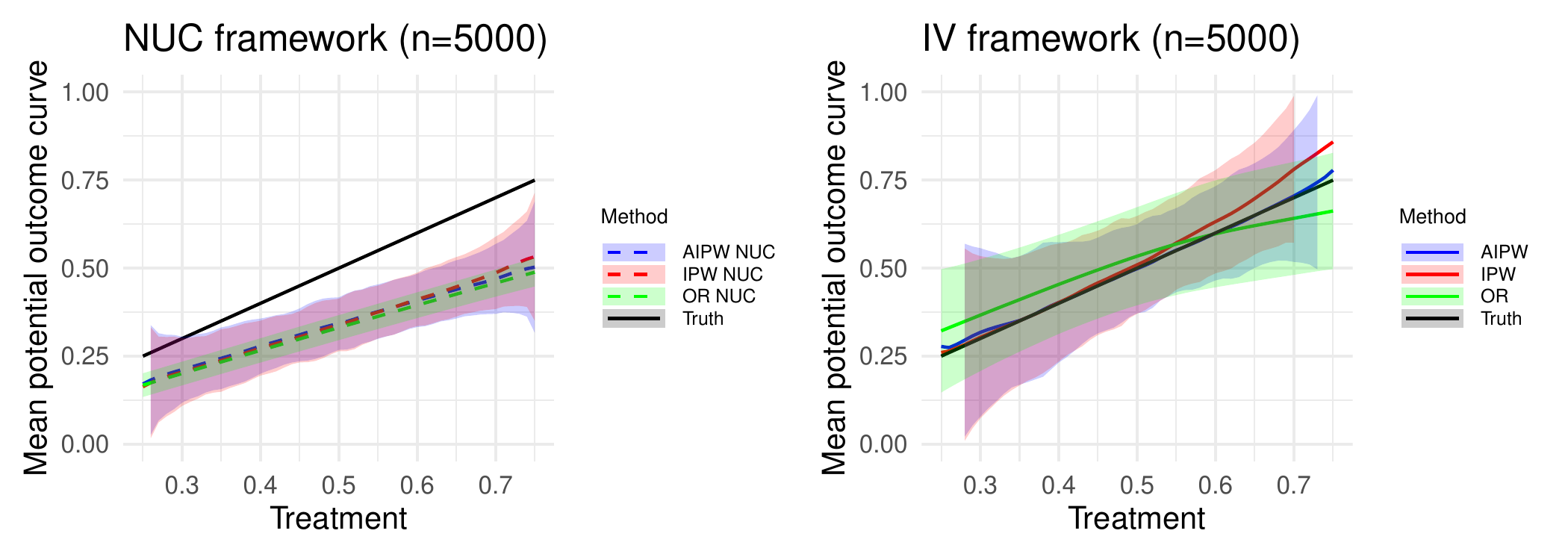}
			\caption{$\mathcal{N}_3 = [0.25,0.75],\,n=5000$}
		\end{subfigure}
		\hfill
		\begin{subfigure}{0.48\textwidth}
			\centering
			\includegraphics[width=\linewidth]{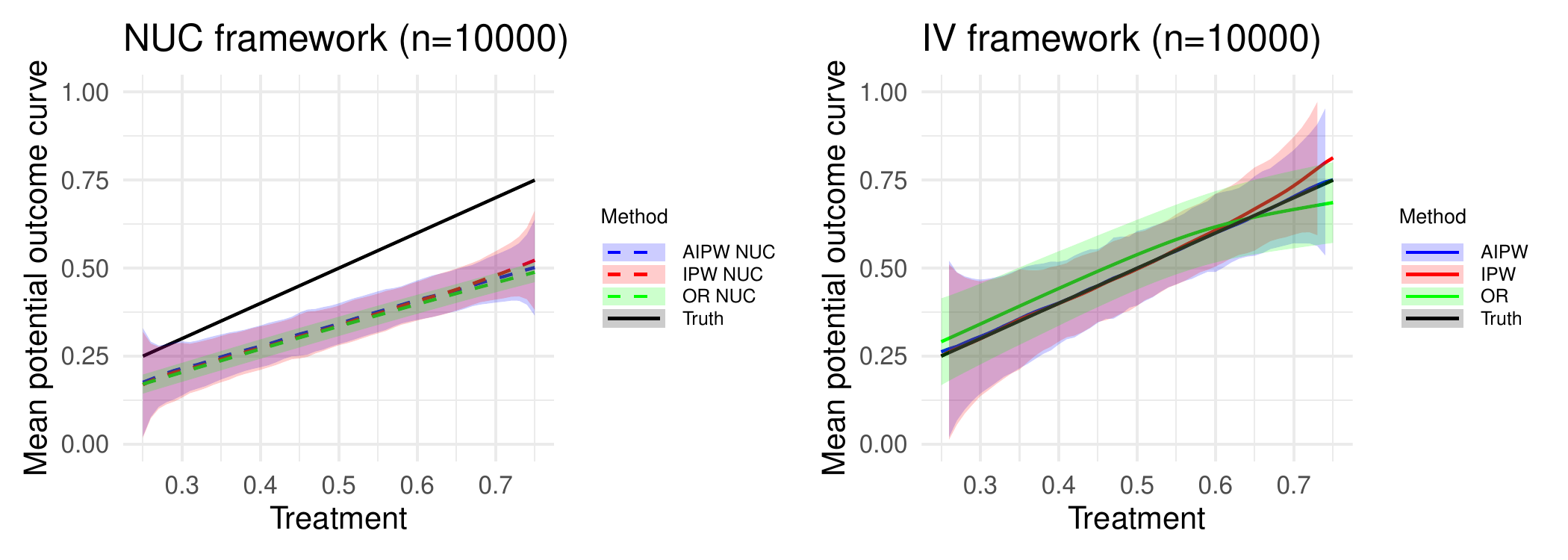}
			\caption{$\mathcal{N}_3 = [0.25,0.75],\,n=10000$}
		\end{subfigure}
		\caption{Empirical mean and pointwise $95\%$ \magenta{uncertainty bands} of estimated ADRF for $n=5000,\,10000$ with 400 replications.}
		\label{fig: simulation results}
	\end{figure}
	
	To demonstrate the effectiveness of our proposed IV-based framework, we compare our estimators with those introduced in \citet{kennedy2017non}. 
	For completeness, the formal definitions of the inverse probability weighting (IPW) and outcome regression (OR) estimators in the IV setting and their counterparts under the NUC framework are provided in Appendix~\ref{append: naive alternative scores}. 
	We set the sample size to be 5,000 or 10,000. The bandwidth is adaptively selected using the localized LOOCV criterion. 
	\magenta{We use $K=5$ folds for cross-fitting, a common choice that balances computational cost and the accuracy of nuisance parameter estimation \citep{Chernozhukov2018}.} 
	We use the Epanechnikov kernel $K(a)=0.75\max\{1 - a^2,0\}$. We train the nuisance functions using spline models fitted with \texttt{gam} from the \texttt{R} package \texttt{mgcv} and the kernel density estimators \texttt{npcdens}, \texttt{npcdensbw} in the \texttt{R} package \texttt{np}.
	
	For graphical illustration, Figure~\ref{fig: simulation results} shows plots of curves estimated from the six methods, displaying the average estimates and the pointwise $95\%$ uncertainty bands of the ADRFs obtained from 400 repetitions.
	\magenta{Notably, we report 95\% uncertainty bands rather than confidence bands, as adaptively selecting the bandwidth via the LOOCV procedure may introduce non-negligible bias.}
	When adopting the NUC framework (left panel of Figure~\ref{fig: simulation results}(a)--(f)), noticeable deviations from the true curve are observed. In contrast, the IV framework (right panel of Figure~\ref{fig: simulation results}(a)--(f)) effectively reduces the estimation bias, albeit with a slight increase in variance compared with the NUC framework. 
	See Appendix~\ref{append: additional simulation results} for further discussion.

	\section{Real data application}\label{sec: empirical illustration}
	We use the dataset provided in the supplementary material of \citet{han2024optimal} for empirical illustration.
	It combines the Job Training Partnership Act (JTPA) dataset with data from the US
	Census and the National Center for Education Statistics (NCES).
	This dataset contains the number of high schools per square mile as an IV, the years of education as the treatment, pre-program earnings as the outcome of interest, measured in dollars and representing total annual income, and sex as a confounder between $A$ and $Y$.
	We are interested in estimating the effect of years of education on pre-program earnings.
	
	Intuitively, local variation in the number of high schools per square mile affects educational level and thus the treatment, and it is unlikely to directly influence pre-program earnings. Hence, it constitutes a reasonable IV in this setting.
	\sy{The AIV condition and its mixture-model interpretation may provide a reasonable approximation in this application. For some individuals, educational attainment could be influenced by local school availability, whereas for others it may depend more strongly on family background, ability, preferences, or other unmeasured socioeconomic factors. This suggests that treatment assignment might be viewed as reflecting both an instrument-related component and a confounded component, which is broadly in line with the motivation underlying the AIV formulation.}

	The dataset consists of 9,223 individuals. Since there are only approximately 500 individuals with fewer than 8 or more than 15 years of education, we grouped those with less than 8 years into the ``8 years'' category and those with more than 15 years into the ``15 years'' category.
	Notably, although years of education are observed discretely, they can be treated as realizations of a continuous treatment variable defined on the real line. See Section~\ref{subsec: semiparametric theory} and Appendix~\ref{append: multi-categorical treatments} for further discussion.

\begin{figure}[t]
		\centering
		\begin{subfigure}{0.32\textwidth}
			\centering
			\includegraphics[width=\linewidth]{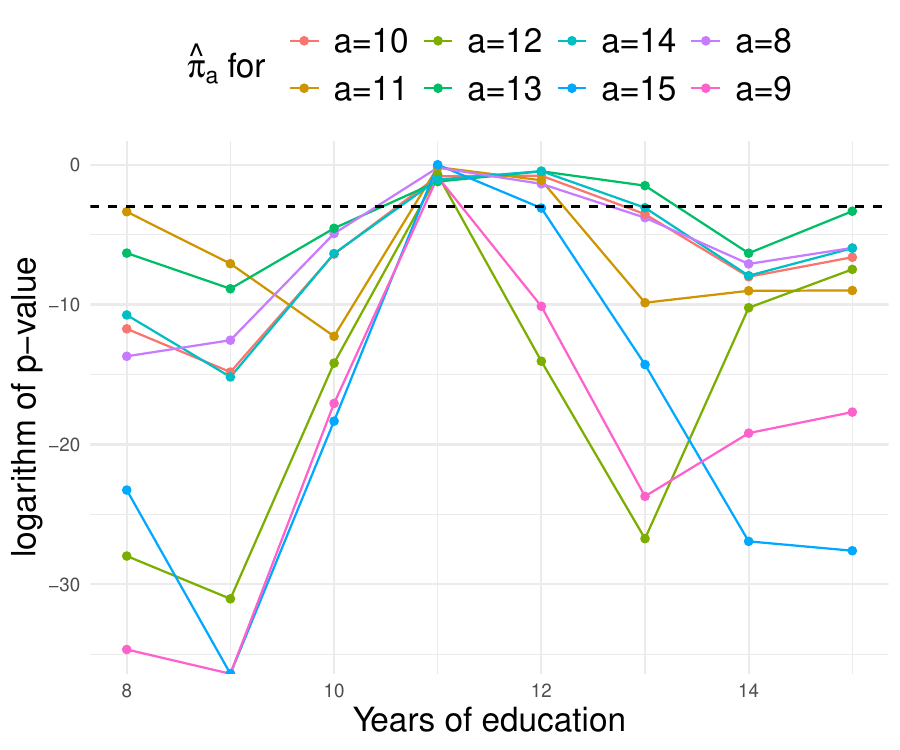}
			\caption{Adjusted for sex.}
		\end{subfigure}
		% \hfill
		\begin{subfigure}{0.32\textwidth}
			\centering
			\includegraphics[width=\linewidth]{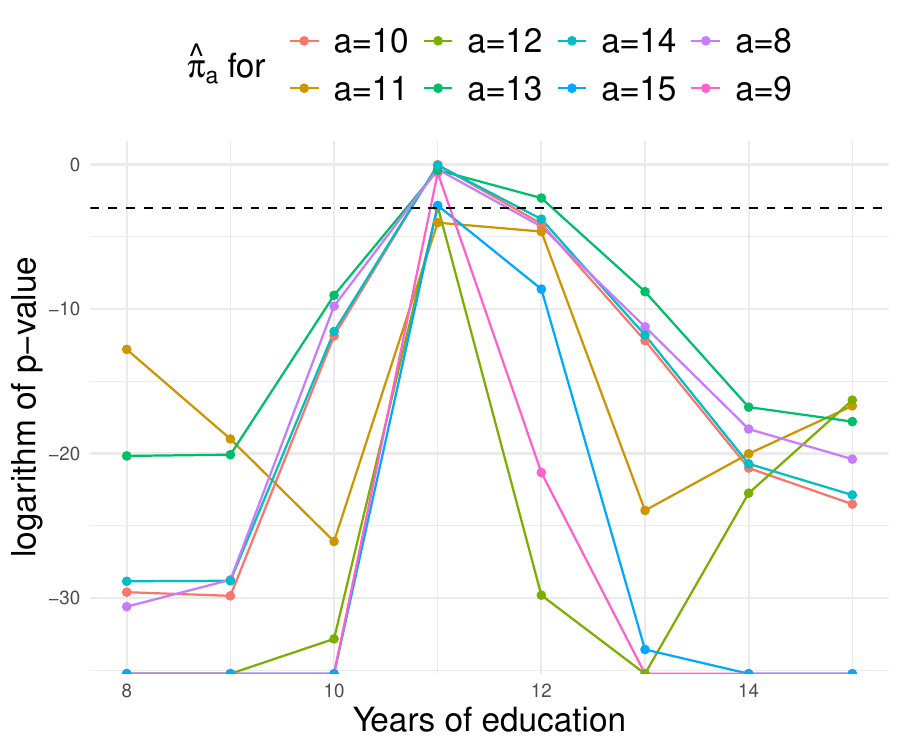}
			\caption{No covariate adjustment.}
		\end{subfigure}
		\begin{subfigure}{0.32\textwidth}
			\centering
			\includegraphics[width=\linewidth]{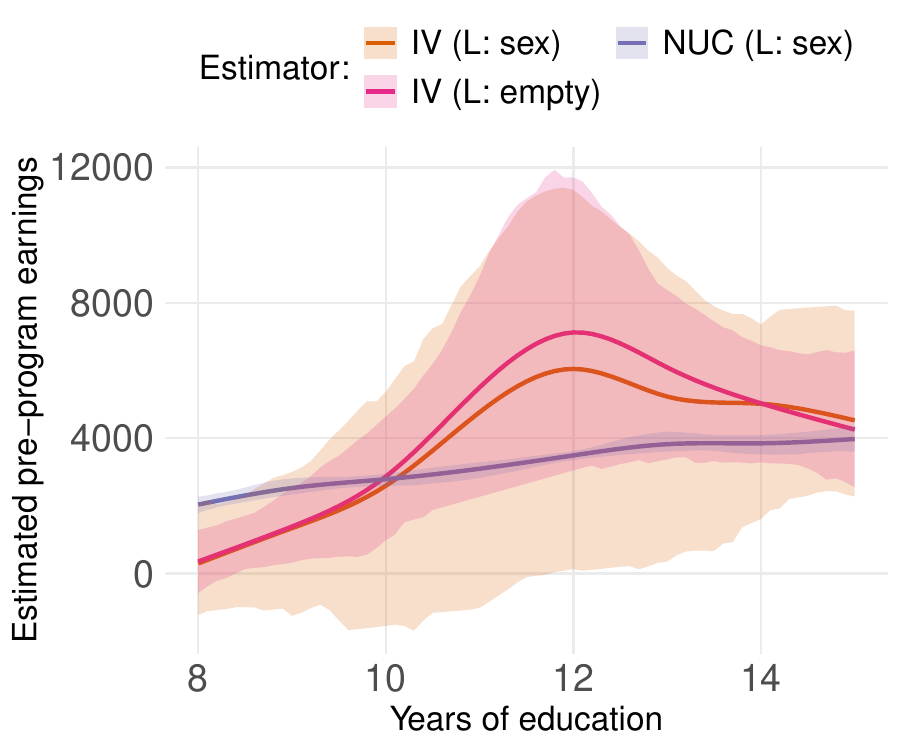}
			\caption{Final results.}
		\end{subfigure}
		\caption{The first two plots show the RWF testing results for eight prespecified weighting functions $\{\hat{\pi}_a\}_{a=8}^{15}$, with the dashed line indicating the significance threshold at $p = 0.05$. The third plot shows the estimated curve (solid line) and the 95\% uncertainty band based on 400 bootstrap replications.}
		\label{fig: JTPA test RWF}
	\end{figure}

	\magenta{For illustration purposes, we conduct our analyses under two specifications: one that sets  $L = \emptyset$ and one that includes sex as $L$, treating $Z$ as a fully exogenous IV in both cases. In particular, when $L = \emptyset$, sex can be regarded as part of the unmeasured confounder $U$.}
	Before presenting the estimated ADRFs, we first estimate $\hat{\pi}_{a}(Z,L)$ for $a$ ranging from 8 to 15, where the density functions are fitted using the \texttt{R} package \texttt{mgcv}.  
	We conduct the RWF testing procedure in Algorithm~\ref{alg: test RWF discrete} for all eight prespecified weighting functions. Figures~\ref{fig: JTPA test RWF}(a) and (b) present the logarithm of the $p$-values across years of education. Figure~\ref{fig: JTPA test RWF}~(a) adjusts for sex as a confounder, while in Figure~\ref{fig: JTPA test RWF}~(b), $L$ is empty. 
	At each $A = a$, smaller p-values provide stronger evidence against the null hypothesis $\mathcal{H}_0(a)$ that the weighting function fails the relevance condition at the specified treatment level. \magenta{Overall, $\hat{\pi}_9$ is a reasonable candidate for $A = 8,9,10$, $\hat{\pi}_{12}$ for $A = 12$, and $\hat{\pi}_{15}$ for $A = 13,14,15$. Accordingly, we select $\hat{\pi}_9$,  $\hat{\pi}_{12}$, and $\hat{\pi}_{15}$ as the RWFs for the corresponding values of $A$.} 
	
	\magenta{
		\sy{For the 1,167 subjects with $A=11$, the corresponding $p$-values in Figure~\ref{fig: JTPA test RWF} are close to or exceed $0.05$, providing insufficient evidence that the relevance condition is satisfied at this treatment level. We therefore restrict the target treatment domain to the observed levels excluding $A=11$ when computing the LLKR estimator, where the proposed estimator may suffer from weak identification and instability.
			Here, we note that observations with $A=11$ are still retained for estimating the nuisance parameter $\alpha_{\pi}^{o}$. In addition, this exclusion does not affect the consistency of the ADRF estimators $\hat{\theta}_{\pi,h}^{(n)}$.
	}}

	\magenta{To estimate the ADRF globally, we merge the AIPW scores corresponding to the three selected RWFs. Specifically, for each treatment level $A$, we use the AIPW score estimated with the corresponding RWF for that level (e.g., at $A=8,9,10$, we use the AIPW score of $\hat{\pi}_9$). Then, we fit a smoothing spline regression with the smoothing parameter adaptively selected via the generalized cross-validation criterion with crossing folds $K=5$. 
		The use of spline methods is theoretically justified, as their properties underpin the error bounds of the ERM framework discussed in Appendix~\ref{append: oracle bound for ERM}.}
	
	\magenta{We conduct 400 bootstrap replications to calculate the standard deviation and construct 95\% uncertainty bands for all treatment levels. To assess the potential impact of unmeasured confounders, we also report the corresponding estimates under the NUC framework, where sex is included as a confounder.
		The final results are shown in Figure~\ref{fig: JTPA test RWF}(c). 
		First, the results consistently highlight the positive effect of education, and those obtained under the IV framework further amplify this pattern for $A\leq 12$, while the NUC estimates are much more stable than those obtained using the IV method, providing a conservative depiction of the overall trend.
		Second, the IV method shows that income slightly decreases when $A \geq 12$, indicating that further increases in education beyond a certain level may no longer have a positive effect on wages. In contrast, this pattern is not reflected in the results obtained using NUC.
		Third, the estimated pre-program earnings are less stable when $L=\{\text{sex}\}$ than when $L$ is empty, although the overall pattern remains similar. This reflects the increased variability introduced by adjusting for additional confounders.}
	
	\sy{We note that there may be other factors that cause the discrepancy between IV and NUC curves, including finite-sample variability due to additional covariate adjustment, instability in nuisance function estimation, weak RWF behavior, or potential deviations from the AIV assumption. 
		Therefore, the differences between the IV and NUC estimates should not be solely attributed to unmeasured confounding, but rather viewed as a descriptive sensitivity analysis under different identification assumptions.}

	\section{Discussion}\label{sec: discussion}
	In this article, we propose a novel IV framework for the identification of ADRFs. 
	We demonstrate that it is not only feasible but also necessary to employ the idea of a finite open cover, whereby any compact subset of the treatment domain can be covered by a collection of open sets. For each of these open sets, we can construct a corresponding URWF, enabling the identification of the ADRF within each region. \magenta{We propose an RWF testing procedure and provide practical guidance for constructing open covers, each of which admits a URWF.}
	We propose the corresponding AIPW score functions and implement a general cross-fitting procedure within the DML framework to compute them. 
	Moreover, we establish the convergence rate of the proposed estimator with kernel regression or empirical risk minimization.

	Several directions for future research remain to be explored.
	\magenta{
		First, developing uniformly valid confidence bands for the estimated ADRF remains an interesting direction for future research. Recent work, such as \citet{takatsu2023debiasedinferencecovariateadjustedregression, ai2026datadrivenuniforminferencegeneral,doss2026doubly}, may provide useful guidance.
		Second, it would be of interest to develop constancy tests under our setting following recent developments of \citet{westling2022nonparametric,doss2024nonparametric}.
		Third, following \citet{branson2023causal}, future work could focus on developing methods that are robust to violations of positivity and IV relevance conditions in the continuous treatments setting.}
	Fourth, our framework also sheds light on developing personalized dose-finding strategies \citep{Chen01102016, kallus2018policyevaluationoptimizationcontinuous, ai2024datadrivenpolicylearningcontinuous} within the proposed IV framework. 
	\sy{Finally, a separate strand of the literature studies the estimation and inference of the ADRF under the NUC setting using alternative methodological frameworks \citep{luedtke2024one,lim2025estimation}. We believe that analogous extensions to the IV setting may also be possible and represent a promising direction for future research.}

	% \bibliographystyle{plainnat}
	% \bibliography{bibfile_main}
	
	{\setlength{\baselineskip}{0.75\baselineskip}
		\putbib[bibfile_main]
	}
\end{bibunit}

\clearpage
\begin{bibunit}
	\appendix
	\section*{Appendix}
	\magenta{Appendix~\ref{append: oracle bound for ERM} presents the methodology and theoretical properties of empirical risk minimization (ERM).
		Appendix~\ref{append: weak AIV} introduces a new concept, termed weak AIV, which is useful for identifying the ADRF contrast.
		Appendix~\ref{append: NPIV} provides an alternative interpretation of the identification of the ADRF under the AIV condition from the perspective of nonparametric instrumental variables and proximal causal inference. Appendix~\ref{append: comparison to existing literature} provides a comparison between our AIV framework and the existing distributional IV literature.
		Appendix~\ref{append: additional data analysis} provides additional data analyses related to Figure~\ref{fig: critical_problem} and Sections~\ref{sec: simulations} and \ref{sec: empirical illustration}.
		Appendix~\ref{append: additional discussion} offers further discussion on the AIPW score and the cross-fitting procedure in Algorithm~\ref{alg: cross-fitting}.
		Appendices~\ref{append: proof of propositions}--\ref{append: proof of lemmas} contain the proofs of the propositions, theorems, and lemmas presented in this article.}

	\section{Oracle bounds for ERM}\label{append: oracle bound for ERM}
	In this section, we follow the framework of \citet{bonvini2022fastconvergenceratesdoseresponse,foster2023orthogonalstatisticallearning} to accommodate the ERM problem under the IV framework.
	Let $\mathcal{F}$ denote a function space, containing the functions that map from $\mathcal{A}$ to $\mathbb{R}$. Let $f$ be any function in $\mathcal{F}$. For a fixed subset $\mathcal{N}\subseteq \mathring{\mathcal{A}}$ and $f:\mathcal{A}\to\mathbb{R}$, define 
	$\|f\|_{\mathcal{N}}^2:=\mathbb{E}[f(A)^2I(A\in\mathcal{N})]/\Pr(A\in\mathcal{N}).$
	For simplicity, we assume that there are no regularization terms.
	
	Then, we solve the localized ERM problem
	\begin{align}
		\hat{f}_{\pi,\mathcal{N}}^{(n)}:= \underset{f\in\mathcal{F}}{\arg\min}\dfrac{1}{n}\sum_{k=1}^K\sum_{i\in I_k}
		\left\{\varphi_{\pi}(O_i;\hat{\alpha}_{\pi}^{(n,-k)},\hat{P}_O^{(n,-k)}) - f(A_i)\right\}^2I\{A_i\in \mathcal{N}\}.
		\label{exp: 30}
	\end{align}
	Let $f_{\mathcal{N}}^*:=\arg\min_{f\in\mathcal{F}}\|f-\theta\|_{\mathcal{N}}$, the function that has the smallest distance in $\mathcal{F}$ from the true $\theta$, and $\mathcal{F}_{\mathcal{N}}^* =  \{f-f_{\mathcal{N}}^*: f \in\mathcal{F}\}$.
	Define the local Rademacher complexity:
	\begin{align*}
		\mathcal{R}_n(\mathcal{F}_{\mathcal{N}}^*,\delta)=\mathbb{E}\left[
		\sup_{f\in\mathcal{F}_{\mathcal{N}}^*:\|f\|_{\mathcal{N}}\leq\delta}\left\{
		\left|\dfrac{1}{n}\sum_{i=1}^n\epsilon_i f(A_i)\right|
		\right\}\middle| A_i\in\mathcal{N}\text{ for all $i$}
		\right],
	\end{align*}
	where $\epsilon_1,\ldots, \epsilon_n$ are i.i.d. Rademacher random variables, independent of the sample.
	
	\begin{theorem}[Oracle bounds for ERM]\label{thm: oracle bounds for ERM}
		Suppose $\mathcal{F}_{\mathcal{N}}^*$ is star-shaped \magenta{and uniformly bounded over $\mathcal{N}$}, and $\theta(a)$ is uniformly bounded. Let $\xi_n$ be any solution to 
		$\mathcal{R}_n(\mathcal{F}_{\mathcal{N}}^*,\xi)\leq \xi^2$ that satisfies
		$\xi_n^2\gtrsim \max\left\{\frac{1}{2n}, \log\log(n)/n\right\},$ which represents 
		the critical radius of $\mathcal{F}_{\mathcal{N}}^*$ \citep{BartlettBousquetMendelson2005}.
		Then, under Assumptions~\ref{as: consistency}--\ref{as: positivity} and Conditions~\ref{as: regular condition}\ref{cd: ZY boundedness}, \ref{cd: alpha uniform boundedness}
		in Assumption~\ref{as: regular condition}, suppose that $\pi(Z,L)$ is a URWF for $\mathcal{N}$, $Z$ is an AIV for $A$, and that $|I_{-k}^{(2)}|\to+\infty$, and for some sequence $d(n)\to 0$,
		\begin{align*}
			\mathbb{E}\left[\begin{array}{l}
				\|(\hat\rho_{\pi}^{(n,-k)}-\rho_{\pi}^o)(\hat\eta^{(n,-k)}-\eta^o)\|_{\mathcal{N}}^2
				+\|(\hat\delta^{(n,-k)}-\delta^o)(\hat\mu_{\pi}^{(n,-k)}-\mu_{\pi}^o)\|_{\mathcal{N}}^2\\
				+\|(\hat\kappa_{\pi}^{(n,-k)}-\kappa_{\pi}^o)(\hat\mu_{\pi}^{(n,-k)}-\mu_{\pi}^o)\|_{\mathcal{N}}^2
				+\|(\hat\rho_{\pi}^{(n,-k)}-\rho_{\pi}^o)(\hat\delta^{(n,-k)}-\delta^o)\|_{\mathcal{N}}^2
			\end{array}\right]
			=o(d(n)).
		\end{align*}
		Then 
		$\mathbb{E}\left[\|\hat{f}_{\pi,\mathcal{N}}^{(n)}-\theta\|_{\mathcal{N}}^2\right]
		\lesssim
		\|f_{\mathcal{N}}^*-\theta\|_{\mathcal{N}}^2+\xi_n^2+o(d(n))
		+O(\sum_{k=1}^K1/|I_{-k}^{(2)}|).$
	\end{theorem}
	In this theorem, $\|f_{\mathcal{N}}^* - \theta\|_{\mathcal{N}}^2$ represents the approximation error, whereas $\xi_n^2$ corresponds to the estimation error \citep{BartlettBousquetMendelson2005,wainwright2019high,foster2023orthogonalstatisticallearning}.
	The third term captures the bias introduced by substituting $\hat{\alpha}_{\pi}^{(n,-k)}$ during estimation, while the fourth term reflects the additional variance arising from the replacement of $P_O$ with $\hat{P}_O^{(n,-k)}$ in Equation~\eqref{exp: 30}.
	
	In fact, when nuisance components are estimated with sufficiently high accuracy, the excess risk of the ERM converges to the oracle rate.
	In particular, if the estimation errors of $\hat{\alpha}_{\pi}^{(n,-k)}$ and $\hat{P}_O^{(n,-k)}$ are asymptotically negligible relative to the main estimation term $\xi_n^2$, the bias and variance terms vanish in the limit.
	Consequently, the overall risk of the ERM estimator matches that of the oracle estimator that has access to the true nuisance functions, achieving the same convergence bound.

	% \magenta{
		\section{Weak AIV}\label{append: weak AIV}
		This section extends Remark~\ref{re: derivative identification} by introducing a weaker version of the AIV condition. We begin with the formal definition.        
		\begin{definition}[Weak AIV]\label{defn: WAIV}
			We say that $Z$ is a weak AIV (WAIV) for $A$ if there exist functions $b_a(U,L)$, $c_a(Z,L)$, and $d(U,L)$ such that for all $a \in \mathring{\mathcal{A}}$,
			\[
			p_{A\mid Z,U,L}(a \mid Z, U, L) = b_a(U, L) + d(U,L)c_a(Z, L).
			\]
		\end{definition}
		As the name suggests, WAIV is a strictly weaker condition than AIV, allowing for a more flexible dependence structure between $(Z,U)$ in the treatment mechanism. The following result provides an equivalent characterization.
		\begin{proposition}[Equivalent characterization of WAIV]\label{prop: equivalent characterization of WAIV}
			Under Assumptions~\ref{as: IV independence}--\ref{as: positivity}, $Z$ is a WAIV for $A$ if and only if $\omega_{a,\pi}(U,L)$ does not depend on $(a,\pi)$. 
		\end{proposition}
		Building on this characterization, we further show that the AIV (WAIV) condition is invariant under monotone transformations of the treatment.
		\begin{proposition}[AIV (WAIV) transformation]\label{prop: AIV transformation}
			If $Z$ is an AIV (WAIV) for $A$, then for any strictly monotone and continuously differentiable function $q(\cdot)$, $Z$ is also an AIV (WAIV) for $q(A)$.
		\end{proposition}
		The next example illustrates how the WAIV condition can arise from a simple data-generating process (DGP).
		\begin{example}[Generating a WAIV from a confounded mixture model]\label{re: WAIV DGP}
			Consider the following DGP: with probability $p_0(U,L)$, 
			$
			A \sim \tilde{p}_{A|Z,L}(a \mid Z,L)
			$
			independently of $U$; otherwise,
			$
			A \sim \tilde{p}_{A|U,L}(a \mid U,L)
			$
			independently of $Z$. By the law of total probability,
			\[
			p_{A|Z,U,L}(a\mid Z,U,L) 
			= p_0(U,L)\, \tilde{p}_{A|Z,L}(a\mid Z,L) 
			+ (1-p_0(U,L))\, \tilde{p}_{A|U,L}(a\mid U,L),
			\]
			which clearly satisfies the WAIV structure.
		\end{example}
		Next, we leverage the concepts of AIV and WAIV to identify the contrast $\mathbb{E}[Y(a_1)] - \mathbb{E}[Y(a_2)]$ for any $a_1, a_2 \in \mathring{\mathcal{A}}$.
		\begin{theorem}[Contrast identification]\label{thm: contrast identification}
			Under Assumptions~\ref{as: consistency}--\ref{as: positivity}, for two fixed values $a_1, a_2 \in \mathring{\mathcal{A}}$, let $\pi_1(Z,L)$ and $\pi_2(Z,L)$ be two RWFs for $A=a_1$ and $A=a_2$, respectively. Assume that either (1) $Z$ is an AIV for $A$, or (2) $Z$ is a WAIV for $A$ and $\mathbb{E}[Y(a_1)-Y(a_2)\mid U,L]\indep U\mid L$. Then,
			\[
			\mathbb{E}[Y(a_1)] - \mathbb{E}[Y(a_2)] = \mathbb{E}\left[ \mu_{\pi_1}^o(a_1,L) - \mu_{\pi_2}^o(a_2,L) \right].
			\]
		\end{theorem}
		We can see that the conditions of Theorem~\ref{thm: contrast identification} are strictly weaker than those of Theorem~\ref{thm: identification}.
		In addition, if $Z$ is a WAIV for $A$, then $\nabla_a \omega_{a,\pi}(U,L)\equiv 0$ by Proposition~\ref{prop: equivalent characterization of WAIV}. 
		\sy{Figure~\ref{fig: identification results} summarizes the identification results established in Theorem~\ref{thm: contrast identification}.}
		This property can be leveraged to identify the derivative $\nabla_a \mathbb{E}[Y(a)]$, as shown below.
		\begin{figure}[t]
			\centering
			\begin{tikzpicture}[
				node distance=2.0cm,
				box/.style={rectangle, draw, rounded corners, minimum width=3.5cm, minimum height=1cm, align=center, fill=gray!5},
				arrow/.style={-{Stealth[scale=1.2]}, thick}
				]
				% Scenario 1 (top)
				\node[box] (mixture) {Mixture model};
				\node[box, below=1cm of mixture] (aiv) {AIV};
				\node[box, below=2cm of aiv] (adrf) {Identify $\mathbb{E}[Y(a)]$};
				% Arrows Scenario 1
				\draw[arrow] (mixture) -- (aiv);
				\draw[arrow] (aiv) -- (adrf);
				% Scenario 2 (bottom)
				\node[box, right=of mixture] (confounded) {Confounded mixture model};
				\node[box, below=1cm of confounded] (waiv) {WAIV};
				\node[box, below=2cm of waiv] (contrast) {Identify $\mathbb{E}[Y(a_1)-Y(a_2)]$};
				% Arrows Scenario 2
				\draw[arrow] (confounded) -- (waiv);
				\draw[arrow] (waiv) -- (contrast) node[midway, right,align=left] {$\mathbb{E}[Y(a_1)-Y(a_2)| U,L]$\\$\indep U\mid L$};
				\draw[arrow, color=black,dashed] (mixture) -- (confounded);
				\draw[arrow, color=black,dashed] (aiv) -- (waiv);
				\draw[arrow, color=black,dashed] (adrf) -- (contrast);
			\end{tikzpicture}
			\caption{Summary of identification results}
			\label{fig: identification results}
		\end{figure}
		% In addition, if $Z$ is a WAIV for $A$, and $\nabla_a\mathbb{E}[Y(a)\mid U,L]\indep U\mid L$, 
		% then $\nabla_a \mathbb{E}[Y(a)]$ is identifiable if it exists.
		\begin{corollary}[Derivative identification]\label{cor: derivative identification}
			Under Assumptions~\ref{as: consistency}--\ref{as: positivity}, for any $a\in \mathring{\mathcal{A}}$, let $\pi(Z,L)$ be an RWF for $A=a$. Assume that either (1) $Z$ is an AIV for $A$, or (2) $\nabla_a\omega_{a,\pi}(U,L)\indep U\mid L$ and $\nabla_a\mathbb{E}[Y(a)\mid U,L]\indep U\mid L$. Then,
			$\nabla_a\mathbb{E}[Y(a)]=\nabla_a\mathbb{E}\left[\mu_{\pi}^o(a,L)\right]$ if they exist.
		\end{corollary}
		\begin{proof}[Proof of Corollary~\ref{cor: derivative identification}]
			By Equation~\eqref{eq: identification E[Y(a)] bias} in Theorem~\ref{thm: identification},
			\begin{align*}
				&\nabla_a\mathbb{E}[\mu_{\pi}^o(a,L)] - \nabla_a\mathbb{E}[Y(a)]
				=\nabla_a\mathbb{E}\left[\mathrm{Cov}\left\{\mathbb{E}[Y(a)\mid U,L], \omega_{a,\pi}(U,L) \mid L\right\}\right]\\
				=&\mathbb{E}\left[\mathrm{Cov}\left\{\nabla_a\mathbb{E}[Y(a)\mid U,L], \omega_{a,\pi}(U,L) \mid L\right\}\right] 
				+\mathbb{E}\left[\mathrm{Cov}\left\{\mathbb{E}[Y(a)\mid U,L], \nabla_a\omega_{a,\pi}(U,L) \mid L\right\}\right]\equiv 0.
			\end{align*}
			The final step holds since:
			\begin{itemize}
				\item If $Z$ is an AIV for $A$, then $\omega_{a,\pi}(U,L)\equiv 1$ by Proposition~\ref{prop: AIV characterization} and $\nabla_a\omega_{a,\pi}(U,L)\equiv 0.$
				\item If $\nabla_a\omega_{a,\pi}(U, L)\indep U\mid L$ and $\nabla_a\mathbb{E}[Y(a)\mid U, L]\indep U\mid L$, it is obvious that both covariance terms are zero.
			\end{itemize}
		\end{proof}
		% Finally, we emphasize a subtle distinction. While there exists a clear equivalent characterization showing that $\omega_{a,\pi}(U, L)$ does not depend on $(a,\pi)$, which corresponds exactly to the WAIV condition, no such equivalent characterization is currently known for the condition $\nabla_a \omega_{a,\pi}(U,L)\equiv 0$. Clarifying this gap is an interesting direction for future research.
		% }
	\sy{
		\section{Connections with existing IV methods}
		\subsection{NPIV and proximal interpretation of the AIV condition}\label{append: NPIV}
		In this section, we provide another interpretation of the AIV condition. 
		The identification can be motivated by nonparametric IV (NPIV) \citep{ai2003efficient,newey2003instrumental} and proximal causal inference \citep{miao2018identifying,cui2024semiparametric,tchetgen2024introduction}.
		In particular, for any treatment level $a\in\mathring{\mathcal{A}}$, consider the following nonparametric IV problem:
		\begin{align}\label{eq: NPIV}
			p_{A\mid Z,L}(a\mid Z,L)\mathbb{E}[Y \mid A=a,Z,L]= p_{A\mid Z,L}(a\mid Z,L)\{f_a^o(1,L)-f_a^o(0,L)\}+ f_a^o(0,L).
		\end{align}
		When the instrumental variable $Z$ takes more than two values, Equation~\eqref{eq: NPIV} becomes an overidentified system. That is, there may not exist a pair of functions $f_a^o(1,L)$ and $f_a^o(0,L)$  that satisfies the equation simultaneously for all values of $Z$. The next proposition shows that the AIV condition provides exactly such a compatibility restriction: under AIV, the overidentified system in  Equation~\eqref{eq: NPIV} admits a solution, and the resulting contrast recovers 
		the ADRF.
		
		\begin{proposition}[Existence]\label{prop: existence NPIV}
			Under Assumptions~\ref{as: consistency}--\ref{as: positivity}, for all $a\in\mathring{\mathcal{A}}$, consider
			\begin{align*}
				f_a^o(0,L) &= \mathrm{Cov}\!\{\mathbb{E}[Y(a)\mid U,L], p_{A|U,L}(a\mid U,L)\mid L\},\\
				f_a^o(1,L) &= \mathrm{Cov}\!\{\mathbb{E}[Y(a)\mid U,L], p_{A|U,L}(a\mid U,L)\mid L\} + \mathbb{E}[Y(a)\mid L].\notag
			\end{align*}
			If $Z$ is an AIV for $A$, then $f_a^o(c,L)$ ($c\in\{0,1\}$) is a solution to Equation~\eqref{eq: NPIV}. In particular, $\mathbb{E}[Y(a)]=\mathbb{E}[f_a^o(1,L)-f_a^o(0,L)].$
		\end{proposition}
		The previous proposition establishes that the AIV condition guarantees the existence of a solution to the nonparametric IV equation. We next show that, once existence is granted, the solution is in fact pinned down by the observed-data distribution. Specifically, under the IV relevance condition, the solution to Equation~\eqref{eq: NPIV} is unique. In addition, it must admit a closed-form representation, which directly leads to identification of the ADRF.
		\begin{proposition}[Uniqueness and closed form]\label{prop: uniqueness NPIV}
			Under Assumption~\ref{as: IV relevance}, for all $a\in\mathring{\mathcal{A}}$, if a solution $f_a^o(c,L)$ ($c\in\{0,1\}$) to Equation~\eqref{eq: NPIV} exists, then this solution is unique. Furthermore, if $\pi$ is an RWF for $A=a$, then $f_a^o(c,L)$ has the following closed form:
			\begin{align*}
				f_a^o(0,L) &= -\mu_{\pi}^o(a,L) p_{A|Z,L}(a|Z,L) + p_{A|Z,L}(a|Z,L)\mathbb{E}[Y\mid A=a,Z,L],\\
				f_a^o(1,L) &= \mu_{\pi}^o(a,L) \{1-p_{A|Z,L}(a|Z,L)\} + p_{A|Z,L}(a|Z,L)\mathbb{E}[Y\mid A=a,Z,L].
			\end{align*}
		\end{proposition}
		Finally, if we combine the results of Propositions~\ref{prop: existence NPIV} and \ref{prop: uniqueness NPIV}, we can identify the ADRF by noticing that
		$\mathbb{E}[Y(a)]=\mathbb{E}[f_a^o(1,L) - f_a^o(0,L)]=\mathbb{E}[\mu_{\pi}^o(a,L)].$

		%   \sy{
			% 	There is also a growing literature on identifying the distribution function of potential outcomes with continuous treatments. In particular, \citet{imbens2009identification, torgovitsky2015identification} study triangular structural models for treatment assignment and construct control variables to identify distributional features of potential outcomes. \citet{chernozhukov2005iv} propose IV quantile regression to identify structural quantile functions of potential outcomes. 
			% 	More recently, \citet{chernozhukov2024estimatingcausaleffectsdiscrete} use a copula invariance condition to identify potential outcome distributions with a binary instrument, while \citet{holovchak2025distributional} estimate such distributions under triangular structural models using an energy-loss-based criterion.
			% }
		\subsection{Comparison with existing IV methods}\label{append: comparison to existing literature}
		In this section, we briefly compare our framework with existing distributional IV approaches. These approaches and our framework rely on different structural assumptions and target different causal estimands. First, we acknowledge that, although our framework allows for a flexible outcome model, the AIV assumption imposes a restriction on the treatment assignment mechanism and may not hold in some commonly used parametric models, such as simple linear treatment models. Thus, the proposed framework is most suitable for applications where the additive structure in the treatment density is scientifically plausible.

		Second, \citet{imbens2009identification,torgovitsky2015identification, holovchak2025distributional} consider triangular models for treatment assignment and outcome generation. These approaches introduce latent variables to capture treatment-selection heterogeneity and often impose a monotone relationship between treatment and a scalar latent rank. 
		Such representations provide a powerful framework for identifying distributional features of potential outcomes, but may not naturally accommodate settings where treatment assignment arises from multiple latent mechanisms, including the mixture model considered in our article.
		
		Third, \citet{chernozhukov2005iv, chernozhukov2020instrumentalvariablequantileregression} develop the IV quantile regression (IVQR) framework to identify structural quantile functions. IVQR allows for flexible treatment mechanisms and characterizes heterogeneous effects across the outcome distribution through rank similarity requirements. In contrast, our framework focuses on the ADRF and allows for a flexible outcome-generating mechanism.
		
		Fourth, \citet{chernozhukov2024estimatingcausaleffectsdiscrete} study the identification of counterfactual distributions with binary instruments using a copula invariance condition. Their approach models the dependence structure between potential outcomes and latent treatment ranks, enabling identification of distributional features under the invariance restriction. Our framework instead directly targets the ADRF and does not require modeling the full counterfactual distribution when the ADRF is the primary target.
		
		Finally, \citet{newey2003instrumental, chen2012estimation} develop NPIV methods based on conditional moment restrictions. These methods provide a flexible framework for recovering structural functions without specifying parametric treatment or outcome models. Our framework is motivated by this line of work and extends the NPIV perspective to continuous treatment effects by exploiting an additive structure in the treatment density (See Appendix~\ref{append: NPIV} for more details). 
		Unlike standard NPIV formulations based on conditional moment restrictions, our approach uses the AIV structure to directly characterize the ADRF through an explicit identification functional.
		
		\begin{table}[t]
			\centering
			\caption{Comparison between different models. Notably, 
				% \citet{holovchak2025distributional} can accommodate discrete IVs only under additional structural restrictions on the treatment assignment. 
				% Furthermore, 
				in Remark~\ref{re: GRDF}, we have shown that the distributional function $F_{Y(a)}(y)$ can be identified under the AIV assumption.}
			\label{tbl: comparison between different models main}
			\resizebox{\textwidth}{!}{
				\begin{tabular}{C{3cm}|C{2cm}C{2cm}C{2cm}C{2cm}C{2cm}C{5cm}}
					\toprule 
					Model & Reduces to NUC
					& Flexible treatment assignment
					& Flexible outcome model & 
					Type of IV
					& Identifiability of $F_{Y(a)}(y)$
					& References
					\\\midrule
					
					Triangular model 
					& $\times$ 
					& $\times$ 
					& $\checkmark$
					& C
					& $\checkmark$
					& \citet{imbens2009identification}
					\\\midrule
					
					IVQR model
					& $\times$
					& $\checkmark$
					& $\times$
					& D/C
					& $\checkmark$
					&\citet{chernozhukov2005iv, chernozhukov2020instrumentalvariablequantileregression}
					\\\midrule
					
					NPIV model
					& $\checkmark$
					& $\checkmark$
					& $\times$
					& D/C
					& $\times$
					& \citet{newey2003instrumental,
						chen2012estimation}
					% bennett2019deep,
					% hartford2017deep,
					% saengkyongam2022exploiting
					\\\midrule
					
					Copula model 
					& $\checkmark$ 
					& \multicolumn{2}{c}{Copula invariance}
					& D
					& $\checkmark$
					& \citet{chernozhukov2024estimatingcausaleffectsdiscrete}
					\\\midrule
					
					AIV
					& $\checkmark$
					& $\times$
					& $\checkmark$
					& D/C
					& $\checkmark$
					& The proposed method
					\\\bottomrule
				\end{tabular}
			}
		\end{table}
		
		Overall, these approaches represent different strategies for addressing unmeasured confounding in IV settings. Triangular models, IVQR, copula-based approaches, and the proposed AIV framework impose different restrictions on the relationship between latent variables, treatment assignment, and potential outcomes. 
		We summarize the key features of different models in Table~\ref{tbl: comparison between different models main}. In the second column, ``Reduces to NUC'' indicates whether, in the absence of unmeasured confounding, the method reduces to the standard no-unmeasured-confounding setting without imposing assumptions beyond the basic IV conditions, namely IV relevance, IV independence, and the exclusion restriction. In the third and fourth columns, ``Flexible treatment assignment'' and ``Flexible outcome model'' indicate whether the model allows flexible treatment and outcome mechanisms, respectively. In the fifth column, ``Type of IV'' indicates whether the model accommodates a general class of IVs, including discrete (D) and continuous (C) IVs.}
	
	\section{Additional data analyses}\label{append: additional data analysis}
	\subsection{Data-generating process for Figure~\ref{fig: critical_problem}}
	\label{append: generating figures in section 2}
	In this subsection, we describe the data-generating process used to produce Figure~\ref{fig: critical_problem}. First, $L$, $\varepsilon_U$, $\varepsilon_Z$, and $\varepsilon_A$ are drawn independently from standard normal distributions, and we set \[U = -0.5 L + 2 \varepsilon_U, \quad Z = -0.5 L + 2 \varepsilon_Z.\] Next, let $e$ be a random variable drawn from a uniform distribution $\mathrm{Unif}(0,1)$. The treatment $A$ is then generated as $A = I(e \leq 0.5) (-0.5 L + Z) + I(e > 0.5) (-0.5 L + U) + \varepsilon_A.$ In this DGP, $Z$ is an AIV for $A$.
	\begin{figure}[htbp]
		\centering
		\includegraphics[width = \textwidth]{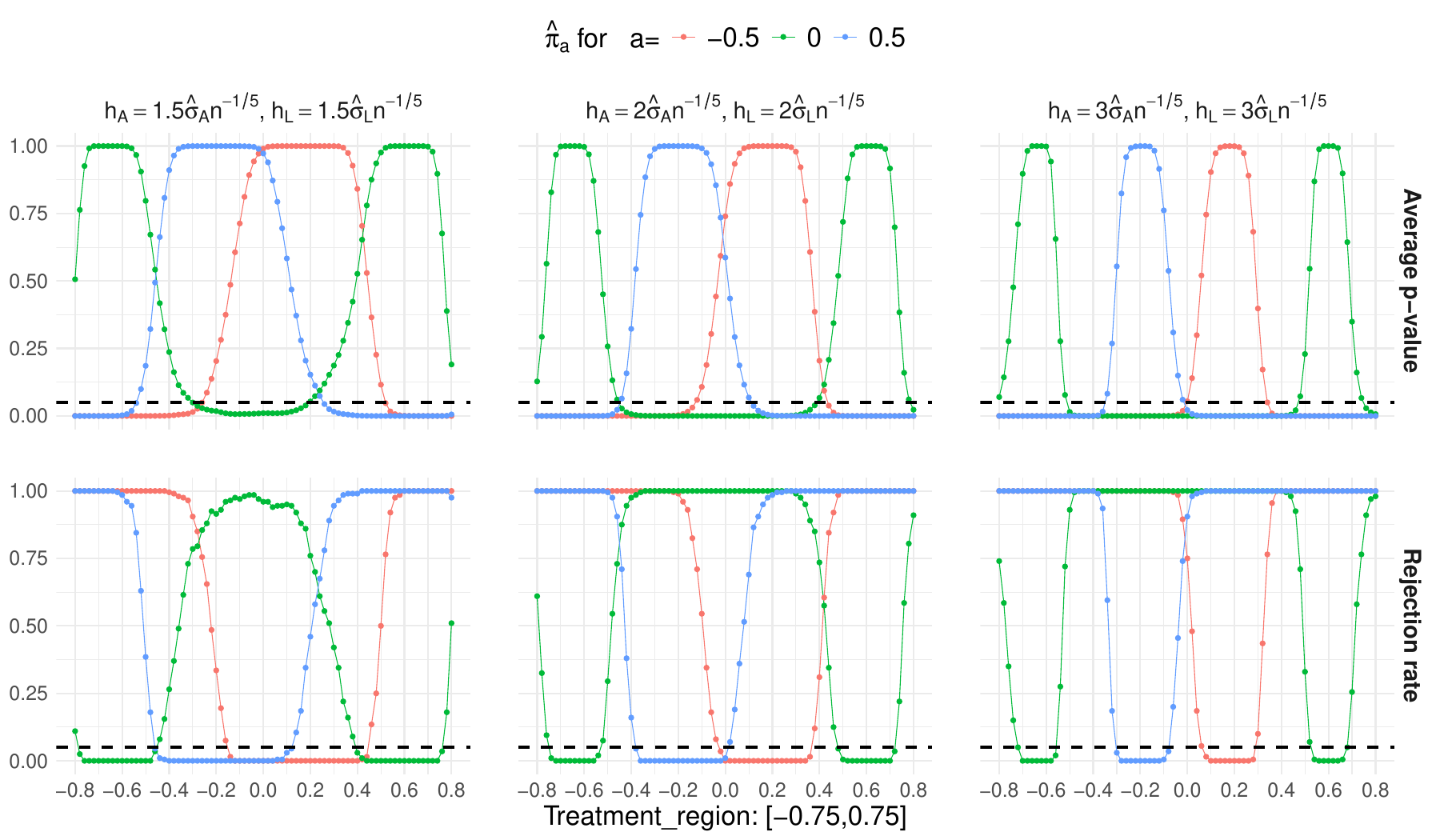}
		\caption{\sy{The upper panels display the average estimated $p$-value curves,
				$\{\frac{1}{S}\sum_{s=1}^S\hat{p}_{a_j^*,\pi_{[m]}}^{(n,s)}\}_{j=1}^J$,
				for $m=1,2,3$, under three bandwidth scaling factors, with different colors representing RWFs. The curves are obtained based on $200$ times repetitions. The lower panels show the corresponding rejection rate curves,
				$\{\frac{1}{S}\sum_{s=1}^S I(\hat{p}_{a_j^*,\pi_{[m]}}^{(n,s)}\leq \alpha)\}_{j=1}^J$,
				for $m=1,2,3$. The black dashed line indicates the significance level $\alpha=0.05$.}}
		\label{fig: RWF testing}
	\end{figure}
	\subsection{Simulations}\label{append: additional simulation results}
	In this subsection, we provide the details of the simulation results presented in Figure~\ref{fig: simulation results}. First, we implement Algorithm~\ref{alg: test RWF continuous} to obtain $p$-values for a sequence of treatment values $a_j^*=-0.77+0.02j$ with $j=1,\ldots, 76$, and three different weighting functions $\{\pi_{[m]}\}_{m=1}^3$. 
	\sy{The upper three panels in Figure~\ref{fig: RWF testing} display the average estimated $p$-value curves, $\{\frac{1}{S}\sum_{s=1}^S\hat{p}_{a_j^*,\pi_{[m]}}^{(n,s)}\}_{j=1}^J,$ obtained from Algorithm~\ref{alg: test RWF continuous} using 10,000 samples under three different bandwidth selections $h:=[h_A,h_L]=Cn^{-1/5}[\hat{\sigma}_A,\hat{\sigma}_L]$ for $C=1.5,2,3$, with $200$ repetitions. The lower three panels present the corresponding rejection rate curves, $\{\frac{1}{S}\sum_{s=1}^S I(\hat{p}_{a_j^*,\pi_{[m]}}^{(n,s)}\leq \alpha)\}_{j=1}^J,$ for $m=1,2,3$.} 
	Each color corresponds to a specific $\pi_{[m]}(Z, L)$, and each point represents the $p$-value at $a_j^*$. The dashed line represents the significance level $\alpha = 0.05$.
	In conclusion, the rejection rate curves remain below the significance level $\alpha=0.05$ for $a_j^*\in\mathcal{N}_m^c$ and approach one for $a_j^*\in\mathcal{N}_m$ across all $m=1,2,3$. These results provide empirical evidence that $\pi_{[m]}$ serves as a suitable candidate URWF for the corresponding neighborhood $\mathcal{N}_m$.

	Next, for $S=400$ Monte Carlo replications, we compute the localized integrated absolute mean bias (BIAS) and the RMSE as evaluation criteria as follows:
	\begin{align*}
		\begin{array}{ll}
			\widehat{\mathrm{BIAS}}(a):=\left|
			\frac{1}{S}\sum\limits_{s=1}^S \hat{\theta}_{\pi,\mathcal{N},s}^{(n)}(a)-\theta(a)
			\right|,&
			\widehat{\mathrm{BIAS}}(\mathcal{N}):=\int_{\mathcal{N}}\widehat{\mathrm{BIAS}}(a)\frac{p_A(a)}{\Pr(A\in\mathcal{N})}\mathrm{d}a,
			\\
			\widehat{\mathrm{RMSE}}(a)^2:=
			\frac{1}{S}\sum\limits_{s=1}^S \left\{\hat{\theta}_{\pi,\mathcal{N},s}^{(n)}(a)-\theta(a)\right\}^2,&
			\widehat{\mathrm{RMSE}}(\mathcal{N})^2:=\int_{\mathcal{N}}
			\widehat{\mathrm{RMSE}}(a)^2
			\frac{p_A(a)}{\Pr(A\in\mathcal{N})}\mathrm{d}a.
		\end{array}
	\end{align*}
	
	\begin{table}[htbp]
		\centering
		\footnotesize
		\caption{Comparison of six different estimators on $\mathcal{N}_3=[0.25,0.75]$.}
		\label{tbl: simulation results}
		\resizebox{0.7\textwidth}{!}{
			\begin{tabular}{ccc|cccc|cccc}
				\toprule
				\multirow{2}{*}{$n$}  & \multicolumn{2}{c|}{\multirow{2}{*}{Method}}
				& \multicolumn{4}{c|}{BIAS} & \multicolumn{4}{c}{RMSE} \\
				&& & 0.4 & 0.5 & 0.6 & $\mathcal{N}$ & 0.4 & 0.5 & 0.6 & $\mathcal{N}$\\ 
				\midrule
				\multirow{6}{*}{5,000} 
				& \multirow{3}{*}{IV}
				& AIPW 
				& .0026 & .0016 & .0011 & .0050 & .0871 & .0645 & .0682 & .0973\\
				&& IPW  
				& .0022 & .0093 & .0319 & .0314 & .0852 & .0714 & .0812 & .1129 \\
				&& OR   
				& .0539 & .0339 & .0025 & .0445 & .0897 & .0787 & .0774 & .0920 \\[0.5em]
				& \multirow{3}{*}{NUC}
				& AIPW 
				& .1216 & .1575 & .1919 & .1670 & .1280 & .1621 & .1958 & .1810 \\
				&& IPW  
				& .1273 & .1613 & .1888 & .1645 & .1333 & .1657 & .1927 & .1764 \\
				&& OR   
				& .1337 & .1690 & .2052 & .1794 & .1349 & .1699 & .2061 & .1878 \\[0.5em]
				\multirow{6}{*}{10,000} 
				& \multirow{3}{*}{IV}
				& AIPW 
				& .0008 & .0002 & .0001 & .0028 & .0597 & .0538 & .0564 & .0701\\
				&& IPW  
				& .0005 & .0026 & .0079 & .0130 & .0550 & .0533 & .0531 & .0741 \\
				&& OR   
				& .0423 & .0384 & .0174 & .0340 & .0685 & .0634 & .0544 & .0655 \\[0.5em]
				& \multirow{3}{*}{NUC}
				& AIPW 
				& .1221 & .1578 & .1909 & .1667 & .1258 & .1606 & .1931 & .1782 \\
				&& IPW  
				& .1271 & .1629 & .1937 & .1673 & .1309 & .1654 & .1958 & .1771 \\
				&& OR   
				& .1299 & .1654 & .2027 & .1768 & .1307 & .1660 & .2031 & .1852 \\[0.5em]
				\bottomrule
		\end{tabular}}
		\vspace{1em}
		\centering
		\footnotesize
		\caption{Comparison of six different estimators on $\mathcal{N}_1=[-0.75,-0.25]$.}
		\label{tbl: simulation results2}
		\resizebox{0.7\textwidth}{!}{
			\begin{tabular}{ccc|cccc|cccc}
				\toprule
				\multirow{2}{*}{$n$}  & \multicolumn{2}{c|}{\multirow{2}{*}{Method}}
				& \multicolumn{4}{c|}{BIAS} & \multicolumn{4}{c}{RMSE} \\
				&& & -0.4 & -0.5 & -0.6 & $\mathcal{N}$ & -0.4 & -0.5 & -0.6 & $\mathcal{N}$\\ 
				\midrule
				\multirow{6}{*}{5,000} 
				& \multirow{3}{*}{IV}
				& AIPW 
				& .0008 & .0003 & .0038 & .0020 & .0913 & .0766 & .0728 & .1012\\
				&& IPW  
				& .0029 & .0111 & .0284 & .0288 & .0833 & .0825 & .0735 & .1106 \\
				&& OR   
				& .0536 & .0344 & .0023 & .0406 & .0919 & .0774 & .0711 & .0873 \\[0.5em]
				& \multirow{3}{*}{NUC}
				& AIPW 
				& .1201 & .1556 & .1921 & .1655 & .1268 & .1611 & .1962 & .1808 \\
				&& IPW  
				& .1268 & .1605 & .1922 & .1636 & .1333 & .1655 & .1965 & .1769 \\
				&& OR   
				& .1331 & .1686 & .2051 & .1771 & .1342 & .1695 & .2059 & .1857 \\[0.5em]
				
				\multirow{6}{*}{10,000} 
				& \multirow{3}{*}{IV}
				& AIPW 
				& .0001 & .0054 & .0014 & .0032 & .0628 & .0512 & .0482 & .0689\\
				&& IPW  
				& .0021 & .0054 & .0037 & .0121 & .0568 & .0497 & .0485 & .0713 \\
				&& OR   
				& .0371 & .0329 & .0160 & .0284 & .0671 & .0604 & .0539 & .0637 \\[0.5em]
				& \multirow{3}{*}{NUC}
				& AIPW 
				& .1218 & .1587 & .1945 & .1672 & .1253 & .1615 & .1969 & .1792 \\
				&& IPW  
				& .1272 & .1643 & .1980 & .1682 & .1305 & .1671 & .2005 & .1786 \\
				&& OR   
				& .1315 & .1670 & .2040 & .1760 & .1323 & .1676 & .2044 & .1845 \\
				\bottomrule
		\end{tabular}}
		\vspace{1em}
		\centering
		\footnotesize
		\caption{Comparison of six different estimators on $\mathcal{N}_2=[-0.25,0.25]$.}
		\label{tbl: simulation results3}
		\resizebox{0.7\textwidth}{!}{
			\begin{tabular}{ccc|cccc|cccc}
				\toprule
				\multirow{2}{*}{$n$}  & \multicolumn{2}{c|}{\multirow{2}{*}{Method}}
				& \multicolumn{4}{c|}{BIAS} & \multicolumn{4}{c}{RMSE} \\
				&& & -0.1 & 0.0 & 0.1 & $\mathcal{N}$ & -0.1 & 0.0 & 0.1 & $\mathcal{N}$\\ 
				\midrule
				\multirow{6}{*}{5,000} 
				& \multirow{3}{*}{IV}
				& AIPW 
				& .0053 & .0019 & .0038 & .0040 & .0876 & .0765 & .0899 & .1135\\
				&& IPW  
				& .0009 & .0025 & .0028 & .0037 & .0870 & .0773 & .0866 & .1132 \\
				&& OR   
				& .0151 & .0014 & .0115 & .0149 & .0880 & .0892 & .0929 & .0917 \\[0.5em]
				& \multirow{3}{*}{NUC}
				& AIPW 
				& .0258 & .0014 & .0282 & .0352 & .0533 & .0397 & .0515 & .0655 \\
				&& IPW  
				& .0273 & .0008 & .0295 & .0369 & .0539 & .0402 & .0529 & .0680 \\
				&& OR   
				& .0324 & .0003 & .0330 & .0421 & .0362 & .0162 & .0369 & .0512 \\[0.5em]
				
				\multirow{6}{*}{10,000} 
				& \multirow{3}{*}{IV}
				& AIPW 
				& .0044 & .0004 & .0022 & .0037 & .0668 & .0542 & .0665 & .0763\\
				&& IPW  
				& .0048 & .0024 & .0022 & .0039 & .0651 & .0567 & .0684 & .0773 \\
				&& OR   
				& .0136 & .0050 & .0056 & .0106 & .0676 & .0634 & .0600 & .0646 \\[0.5em]
				& \multirow{3}{*}{NUC}
				& AIPW 
				& .0291 & .0004 & .0295 & .0389 & .0454 & .0331 & .0467 & .0586 \\
				&& IPW  
				& .0298 & .0014 & .0311 & .0401 & .0459 & .0336 & .0485 & .0599 \\
				&& OR   
				& .0319 & .0001 & .0315 & .0409 & .0348 & .0135 & .0342 & .0490 \\
				\bottomrule
		\end{tabular}}
	\end{table}
	The simulation results are shown in Tables~\ref{tbl: simulation results}--\ref{tbl: simulation results3} for $400$ repetitions of the above process. We demonstrate the BIAS and RMSE at three specific points, including $A=-0.4,\,-0.5,\,-0.6\in\mathcal{N}_1$, $A=-0.1,\,0,\,0.1\in\mathcal{N}_2$, and $A=0.4,\,0.5,\,0.6\in\mathcal{N}_3$, and report the integrated BIAS and RMSE for each $\mathcal{N}_m$. 
	As we can see, AIPW under the IV framework exhibits stable performance. 
	
	For instance, in Table~\ref{tbl: simulation results} with sample size $n=5000$, the bias of AIPW with IV at points $0.4$, $0.5$, $0.6$, and overall $\mathcal{N}_3$ remains around $0.0050$, substantially lower than other methods, especially OR, indicating near-unbiased estimation. The corresponding RMSE values further demonstrate that the estimates are precise and stable. In addition, AIPW under the IV framework provides robust estimates with low BIAS. The RMSE of the AIPW method occasionally exceeds that of the IPW method, which is likely because the nuisance functions $\hat{\mu}_{\pi}^{(n,-k)}$ are not accurately estimated. This can also be inferred from the biased estimates produced by the OR method, indicating that the misspecification of $\hat{\mu}_{\pi}^{(n,-k)}$ may undermine the efficiency of the AIPW estimator.

	In contrast, in Table~\ref{tbl: simulation results}, the NUC framework may lead to notable bias, since there indeed exists an unmeasured confounder in our simulation study. For instance, AIPW bias ranges from $0.1216$ to $0.1919$, with the overall $\mathcal{N}_3$ bias at $0.1670$, much higher than $0.0050$ (IV bias). Correspondingly, RMSE values also increase, indicating that estimates deviate further from the truth and exhibit greater variability. This suggests that in the absence of proper IV control, even AIPW can be affected by unmeasured confounding and produce biased estimates.
	
	\begin{figure}[t]
		\centering
		\includegraphics[width = \textwidth]{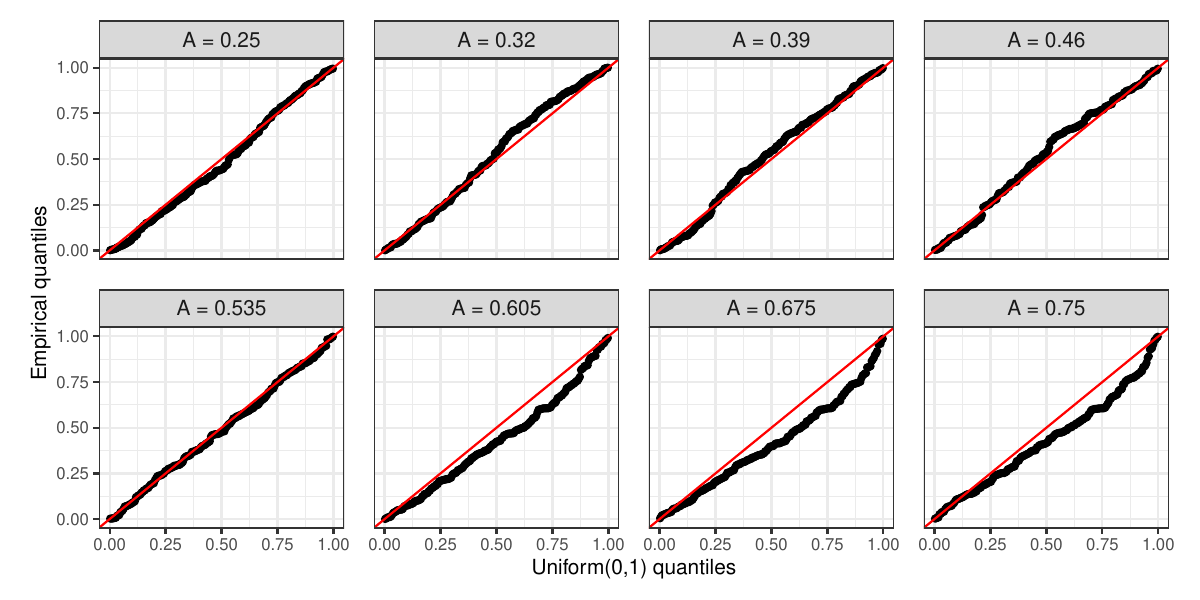}
		\caption{\sy{
				Q-Q plots of the simulated $p$-values across different treatment levels. The red diagonal line represents the reference line corresponding to a perfectly calibrated $p$-value distribution.
		}}
		\label{fig: falsification test}
	\end{figure}
	\sy{
		%    \subsection{A falsification test}
		Finally, to assess the finite-sample validity and performance of the proposed falsification test in Section~\ref{subsec: falsification test}, we conduct a Monte Carlo experiment with a sample size of $n=5000$ and examine the calibration of the resulting $p$-values under the null hypothesis. Concretely, we first construct $\hat{\pi}_{0.35}$ and $\hat{\pi}_{0.65}$ from one realization and use them as fixed URWFs over the interval $[0.25, 0.75]$. 
		We then repeat the simulation procedure 200 times, and the empirical distributions of the resulting $p$-values are presented in Figure~\ref{fig: falsification test}.
		The empirical quantiles closely align with the theoretical quantiles across most treatment levels, providing evidence that the proposed testing procedure yields well-calibrated $p$-values in finite samples. 
		Consequently, under this simulation setting, there is insufficient evidence against the necessary implication of the AIV condition examined by our falsification test.
	}

	\subsection{JTPA}\label{append: JTPA}
	For the real data analysis in Section~\ref{sec: empirical illustration}, we provide additional summary statistics of the AIPW scores before merging them. Concretely, Table~\ref{tbl: estimated pre-program total annual earning} presents the estimated pre-program total annual earnings at different levels of education. 
	The results report mean point estimates (EST) and standard deviations (SD) from both the IV approach, using the estimated density functions $\hat{\pi}_{9}$, $\hat{\pi}_{12}$, and $\hat{\pi}_{15}$, and the degenerate NUC score, with 400 bootstrap replications. 
	The estimates suggest a positive relationship between educational attainment and pre-program earnings. 
	A notable difference between the IV and NUC methods is observed at lower education levels ($A=8,9,10$), where IV-based earnings are substantially lower than the NUC estimates. 
	In addition, estimates when $L=\emptyset$ (the left panel of Table~\ref{tbl: estimated pre-program total annual earning}) are more stable than those obtained when sex is included as a confounder (the right panel of Table~\ref{tbl: estimated pre-program total annual earning}), likely because including a covariate introduces additional variability of $\kappa_{\pi}^o(A,L)$.  
	
	\begin{table}
		\centering
		\footnotesize
		\caption{Estimated pre-program total annual earnings (EST), measured in dollars, and the corresponding standard deviation (SD). A backslash ``$\backslash$'' indicates that the estimated standard deviation is excessively large (exceeding 10,000).}
		\label{tbl: estimated pre-program total annual earning}
		\setlength{\tabcolsep}{4pt}
		\resizebox{\textwidth}{!}{
			\begin{tabular}{c | ccccccc| ccccccc}
				\toprule
				\multirow{2}{*}{Years of education}&
				\multicolumn{7}{c|}{No confounder}&
				\multicolumn{7}{c}{Confounder: sex}\\
				& 8 & 9 & 10 & 12 & 13 & 14 & 15 
				& 8 & 9 & 10 & 12 & 13 & 14 & 15\\
				\midrule
				EST (IV, $\hat{\pi}_{9}$) 
				& 580 & 1687 & 940 & 6412 & 4582 & 5143 & 4090 
				& 506 & 1649 & 449 & 5577 & 4065 & 4774 & 3960 \\ 
				SD (IV, $\hat{\pi}_{9}$) 
				& 380 & 552 & 895 & 2812 & 1480 & 998 & 643 
				& 421 & 597 & 5842 & 5547 & 6556 & 2327 & 731 \\ [0.5em]
				EST (IV, $\hat{\pi}_{12}$)   
				& 231 & -156 & -1476 & 8010 & 4444 & $\backslash$ & 4358 
				& 247 & $\backslash$ & $\backslash$ & 6960 & $\backslash$ & $\backslash$ & $\backslash$ \\ 
				SD (IV, $\hat{\pi}_{12}$)    
				& 443 & 947 & 1796 & 2448 & 2047 & $\backslash$ & 3828 
				& 479 & $\backslash$ & $\backslash$ & 3806 & $\backslash$ & $\backslash$ & $\backslash$ \\ [0.5em]
				EST (IV, $\hat{\pi}_{15}$)   
				& 1016 & 2132 & 2221 & 5297 & 4364 & 4367 & 3833 
				& 892 & 1984 & $\backslash$ & $\backslash$ & 3829 & 4443 & 3600 \\ 
				SD (IV, $\hat{\pi}_{15}$)   
				& 450 & 564 & 825 & 8411 & 1557 & 847 & 629 
				& 553 & 726 & $\backslash$ & $\backslash$ & 8741 & 847 & 732 \\ [0.5em]
				EST (IV, merge)   
				& 580 & 1687 & 940 & 8010 & 4364 & 4367 & 3833 & 506 & 1649 & 449 & 6960 & 3829 & 4443 & 3600 \\ 
				SD (IV, merge)
				& 380 & 552 & 895 & 2448 & 1557 & 847 & 629 & 421 & 597 & 5842 & 3806 & 8741 &  847 & 732 \\ [0.5em]
				EST (NUC)  
				& 1977 & 2555 & 2778 & 3456 & 3914 & 3907 & 4074 
				& 1983 & 2613 & 2740 & 3497 & 3919 & 3782 & 3941 \\ 
				SD (NUC) 
				& 134 & 140 & 112 & 68 & 168 & 183 & 217 
				& 128 & 147 & 111 & 69 & 180 & 184 & 202 \\ 
				\bottomrule
		\end{tabular}}
	\end{table}
	
	\magenta{Meanwhile, the IV-based standard deviations are also larger, especially at $A=12$, likely due to the potential violation of the RWF condition.}
	In contrast, the NUC estimates are more stable, providing a conservative depiction of the overall trend. At higher education levels, the estimated earnings with IV slightly decrease, suggesting that additional education beyond a certain point may not further increase pre-program earnings. 
	
	\sy{
		Next, we present an AIV testing procedure that differs slightly from the one discussed in Section~\ref{subsec: falsification test}, as the treatment region considered here is discrete. 
		In particular, for a series of RWFs $\Pi:=\{\pi_j\}_{j=1}^J$ for $A=a$, we define 
		\begin{align*}
			&\mathcal{L}_{\Pi,h,a}^{(n,\mathrm{dis})}(\beta,Q):= 
			\left\|\frac{1}{K}\sum_{k=1}^{K}\mathbb{E}_{nk}\left[\Psi_{\Pi,h,a}^{(\mathrm{dis})}(O;\beta,\hat\alpha_{\Pi}^{(n,-k)}, \hat{P}_O^{(n,-k)})\right]
			\right\|_{Q}^2,\\
			&\Psi_{\Pi,h,a}^{(\mathrm{dis})}(O;\beta,\alpha_{\Pi},P_O):= I(A = a)\left\{\varphi_{\Pi}(O; \alpha_{\Pi}, P_O)-\mathbf{1}_J \beta\right\}\in\mathbb{R}^J.
		\end{align*}
		We define an initial estimator $\hat{\beta}_{\Pi,h}^{(n,\mathrm{dis},\mathrm{init})}(a):=\underset{\beta\in\mathbb{R}}{\arg\min}\{\mathcal{L}_{\Pi,h,a}^{(n,\mathrm{dis})}(\beta,\mathbf{I}_{J})\}$, and a debiased GMM estimator is defined as 
		$\hat{\beta}_{\Pi,h}^{(n,\mathrm{dis})}(a) :=\underset{\beta\in\mathbb{R}}{\arg\min}\{\mathcal{L}_{\Pi,h,a}^{(n,\mathrm{dis})}(\beta,\hat{Q}_{\Pi,h}^{(n,\mathrm{dis})}(a)^{-1})\},$ where
		\begin{align*}
			\hat{Q}_{\Pi,h}^{(n,\mathrm{dis})}(a)=&\frac{1}{K}\sum_{k=1}^{K}\mathbb{E}_{nk}\left[\begin{array}{l}
				\Psi_{\Pi,h,a}^{(\mathrm{dis})}(O;\hat{\beta}_{\Pi,h}^{(n,\mathrm{dis},\mathrm{init})}(a),\hat\alpha_{\Pi}^{(n,-k)}, \hat{P}_O^{(n,-k)})\\\Psi_{\Pi,h,a}^{(\mathrm{dis})}(O;\hat{\beta}_{\Pi,h}^{(n,\mathrm{dis},\mathrm{init})}(a),\hat\alpha_{\Pi}^{(n,-k)}, \hat{P}_O^{(n,-k)})^{\top}
			\end{array}\right]\in\mathbb{R}^{J\times J}.
		\end{align*}
		Similar to Theorem~\ref{thm: falsify AIV condition}, under Assumptions~\ref{as: consistency}--\ref{as: positivity}, if $Z$ is an AIV for $A=a$, and Conditions~(f), (g) and (h) in Assumption~\ref{as: regular condition} are satisfied with $c(n)=n^{-1/2}$, then $$n\mathcal{L}_{\Pi,h,a}^{(n,\mathrm{dis})}\left(\hat{\beta}_{\Pi,h}^{(n,\mathrm{dis})}(a),\hat{Q}_{\Pi,h}^{(n,\mathrm{dis})}(a)^{-1}\right)\overset{d}{\to}\chi_{J-1}^2.$$ In this setting, the null hypothesis is given by $\widetilde{\mathcal{H}}_0(a):$ for all $j_1,j_2=1,\ldots,J,m_{j_1}(a)=m_{j_2}(a)$, and the alternative hypothesis is given by $\widetilde{\mathcal{H}}_1(a):$ there exists $j_1\neq j_2,\text{ such that }m_{j_1}(a)\neq m_{j_2}(a)$. The corresponding $p$-value is defined as $\hat{p}_{a,\Pi}^{(n,\mathrm{dis})}:=1-F_{\chi_{J-1}^2}(n\mathcal{L}_{\Pi,h,a}^{(n,\mathrm{dis})}(\hat{\beta}_{\Pi,h}^{(n,\mathrm{dis})}(a),\hat{Q}_{\Pi,h}^{(n,\mathrm{dis})}(a)^{-1}))$.

		Concretely, we use $\hat{\pi}_{8}$, $\hat{\pi}_{9}$, and $\hat{\pi}_{10}$ as the RWFs for $A=8,9,10$; $\hat{\pi}_{12}$ and $\hat{\pi}_{13}$ as the RWFs for $A=12,13$; and $\hat{\pi}_{14}$, $\hat{\pi}_{15}$ as the RWFs for $A=14,15$. We then report the $p$-values proposed in Section~\ref{subsec: falsification test} to assess the falsifiable implication of the AIV condition. 
		The results are summarized in Table~\ref{tbl: falsification test}. The reported $p$-values range from $0.198$ to $0.799$ across the treatment levels considered. None of these $p$-values is particularly small, and thus the proposed test does not provide clear evidence against the corresponding falsifiable implication of the AIV condition at any of these treatment levels. Overall, the results suggest that the AIV condition is not obviously contradicted by the observed data in this application.}
	
	\begin{table}[htbp]
		\centering
		\sy{
			\begin{tabular}{c|ccc|cc|cc}
				\toprule
				Years of education & 8 & 9 & 10 & 12 & 13 & 14 & 15\\
				RWFs 
				&\multicolumn{3}{c|}{$\{\hat{\pi}_{8},\hat{\pi}_{9},\hat{\pi}_{10}\}$}
				&\multicolumn{2}{c|}{$\{\hat{\pi}_{12},\hat{\pi}_{13}\}$}
				&\multicolumn{2}{c}{$\{\hat{\pi}_{14},\hat{\pi}_{15}\}$}\\\midrule
				$\hat{p}_{a,\Pi}^{(n,\mathrm{dis})}$& 
				0.278 &0.198 &0.791 &0.402 &0.315 &0.383 &0.799\\
				\bottomrule
		\end{tabular}}
		\caption{\sy{Falsification test results for the AIV condition across years of education.}}
		\label{tbl: falsification test}
	\end{table}
	
	\section{Additional discussion}\label{append: additional discussion}
	\subsection{Alternative score functions}\label{append: naive alternative scores}
	
	In this subsection, we summarize the \textit{inverse probability weighting} (IPW) and \textit{outcome regression} (OR) scores as naive alternatives to the AIPW score defined in Section~\ref{subsec: semiparametric theory} in the IV setting. 
	Recall that for $Z_{\pi}:=\pi(Z,L)$, we define the nuisance functions
	\begin{align*}
		\begin{array}{ll}
			\rho_{\pi}^o(L):=\mathbb{E}[Z_{\pi}\mid L], &\kappa_{\pi}^o(A,L) := \mathbb{E}[Z_{\pi}\mid A,L] - \mathbb{E}[Z_{\pi}\mid L], \\
			\eta^o(A,L):= \mathbb{E}[Y\mid A,L], &\delta^o(A,L) :=  p_{A}(A)/p_{A|L}(A\mid L).
		\end{array}
	\end{align*}
	For the IPW score, we define 
	$
	\varphi_{\pi,\mathrm{IPW}}(O;\alpha_{\pi},P_O):=
	\delta(A,L)
	\frac{(Z_{\pi}-\rho_{\pi}(L))Y}{\kappa_{\pi}(A,L)}.
	$
	For the OR score, we set 
	$
	\varphi_{\pi,\mathrm{OR}}(O;\alpha_{\pi},P_O):=
	\int \mu_{\pi}(A,l)\,\mathrm{d}P_O(o).
	$
	It is straightforward to verify that when $Z$ is an AIV for $A=a$, and $\pi$ is an RWF for $A=a$,
	\[
	\theta(a) = \mathbb{E}\big[\varphi_{\pi,\mathrm{IPW}}(O;\alpha_{\pi}^o,P_O)\mid A=a\big]
	=\mathbb{E}\big[\varphi_{\pi,\mathrm{OR}}(O;\alpha_{\pi}^o,P_O)\mid A=a\big].
	\]
	Hence, the corresponding estimation procedure can be implemented simply by replacing the AIPW score function $\varphi_{\pi}(O;\alpha_{\pi},P_O)$ in Equation~\eqref{eq: varphi defn} with the IPW or OR score function. All subsequent analyses proceed analogously. 
	
	Notably, although these scores serve as alternatives to the AIPW score, they do not inherit the mixed bias property discussed in Lemma~\ref{lem: mixed bias}. 
	For reference, in \citet{kennedy2017non}, the corresponding AIPW, IPW, and OR score functions are listed as: 
	\begin{align*}
		\varphi_{\mathrm{AIPW}}^{\mathrm{nuc}}(O;\alpha_{\pi},P_O)
		&:= 
		\delta(A,L)(Y-\eta(A,L))
		+\int \eta(A,l)\,\mathrm{d}P_O(o),\\
		\varphi_{\mathrm{IPW}}^{\mathrm{nuc}}(O;\alpha_{\pi},P_O)
		&:=
		\delta(A,L)Y,\quad
		\varphi_{\mathrm{OR}}^{\mathrm{nuc}}(O;\alpha_{\pi},P_O):=
		\int \eta(A,l)\,\mathrm{d}P_O(o).
	\end{align*}
	All three of these score functions can lead to non-negligible bias in the presence of unmeasured confounders~$U$.
	
	\subsection{Multi-categorical treatments}\label{append: multi-categorical treatments}
	In this subsection, we consider a degenerate case where $A$ is multi-categorical. We briefly compare our methods to those of \citet{chen2025identificationdebiasedlearningcausal}.
	We assume that $\pi(Z,L)$ is an RWF for $A=a$.
	First, we define nuisance functions as
	\begin{align*}
		&\begin{array}{ll}
			\rho_{\pi}^o(L):=\mathbb{E}[Z_{\pi}\mid L], &\kappa_{\pi}^o(A,L) := \mathbb{E}[Z_{\pi}\mid A,L] - \mathbb{E}[Z_{\pi}\mid L], \\
			\eta^o(A,L):= \mathbb{E}[Y\mid A,L], &\Delta^o(a,L) := \mathbb{E}[I(A=a)\mid L],
		\end{array}\\
		&\mu_{\pi}^o(A,L):=\dfrac{\mathbb{E}[YZ_{\pi}\mid A,L] - \mathbb{E}[Y\mid A,L]\mathbb{E}[Z_{\pi}\mid L]}{\mathbb{E}[Z_{\pi}\mid A,L] - \mathbb{E}[Z_{\pi}\mid L]}.
	\end{align*}
	We then unify them into a nuisance vector
	\begin{align*}
		\alpha_{\pi}^o(A,L) = [\mu_{\pi}^o(A,L),\rho_{\pi}^o(L), \kappa_{\pi}^o(A,L), \eta^o(A,L), \Delta^o(A,L)].
	\end{align*}
	Then we can simply verify that for any $a\mathcal{A}$,
	\begin{align*}
		\mathbb{E}[Y(a)] = \mathbb{E}[\mu_{\pi}^o(a,L)] = 
		\mathbb{E}\left[\dfrac{I(A=a)}{\Delta^o(A,L)}\dfrac{Z_{\pi}-\rho_{\pi}^o(L)}{\kappa_{\pi}^o(A,L)}Y\right].
	\end{align*}
	If we define $\theta_{\pi,a}^o:=\mathbb{E}\left[\dfrac{I(A=a)}{\Delta^o(A,L)}\dfrac{Z_{\pi}-\rho_{\pi}^o(L)}{\kappa_{\pi}^o(A,L)}Y\right]$, then the EIF of $\theta_{\pi,a}^o$ in \citet{chen2025identificationdebiasedlearningcausal} can be derived as $\varphi_{\pi,alt}(O;\theta_{\pi,a}^o,\alpha_{\pi}^o)$, where 
	
	\begin{equation}\label{exp: 23}
		\begin{aligned}
			\varphi_{\pi,alt}(O;\theta_{\pi,a},\alpha_{\pi}) :=& \dfrac{I(A=a)}{\Delta(A,L)}
			\dfrac{\{Z_{\pi}- \rho_{\pi}(L)\}}{\kappa_{\pi}(A,L)}(Y-\mu_{\pi}(A,L))-\theta_{\pi,a}\\
			&+ \mu_{\pi}(a,L) - \dfrac{\{Z_{\pi}-\rho_{\pi}(L)\}}{\kappa_{\pi}(a,L)}(\eta(a,L)-\mu_{\pi}(a,L)).
		\end{aligned}
	\end{equation}
	Notably, this EIF no longer depends on the nuisance distribution $P_O(o)$. As shown in \citet{chen2025identificationdebiasedlearningcausal}, this EIF also exhibits a mixed-bias property when $\alpha_{\pi}^o$ is plugged into estimation. In fact, the EIF proposed in Equation~\eqref{eq: varphi defn} can be regarded as a generalization of Equation~\eqref{exp: 23}. The main difference between these two approaches is that the framework developed in this article accommodates both multi-categorical and continuous treatments. 
	
	In fact, when $A$ is an ordinal treatment, \citet{chen2025identificationdebiasedlearningcausal} can only identify $\mathbb{E}[Y(a)]$. In contrast, our method identifies $\mathbb{E}[Y(a)]$ as a continuous curve, treating multi-categorical treatments $A$ as observations sampled from a continuous space.
	Therefore, under additional smoothness assumptions, the fitted curve can be interpolated between observed treatment levels.

	\subsection{Degenerate AIPW scores when $L$ is empty}\label{append: degenerate AIPW score}
	When $L = \emptyset$, the AIPW score function reduces to a simpler form compared to the discussion in Section~\ref{subsec: semiparametric theory}.
	Choose a function $\pi(Z)$ that serves as a URWF for $\mathcal{N}$ and set $Z_{\pi} := \pi(Z)$ as before.
	We then define the corresponding degenerate nuisance function as
	\begin{align*}
		&\rho_{\pi,\emptyset}^o:=\mathbb{E}[Z_{\pi}],
		\;
		\mu_{\pi,\emptyset}^o(A):=\dfrac{\mathbb{E}[Y (Z_{\pi}- \mathbb{E}[Z_{\pi}])\mid A]}{\mathbb{E}[Z_{\pi}\mid A]-\mathbb{E}[Z_{\pi}]},
		\;\\
		&\kappa_{\pi,\emptyset}^o(A):= \mathbb{E}[Z_{\pi}\mid A]-\mathbb{E}[Z_{\pi}],
		\;
		\eta_{\emptyset}^o(A):=\mathbb{E}[Y\mid A].
	\end{align*}
	Then, we summarize them into a unified nuisance vector
	$\alpha_{\pi,\emptyset}^o = [\rho_{\pi,\emptyset}^o, \mu_{\pi,\emptyset}^o, \kappa_{\pi,\emptyset}^o,\eta_{\emptyset}^o].$
	In this setting, if $Z$ serves as an AIV for $A$, one can easily show that $\theta(a)=\mu_{\pi,\emptyset}^o(a)$. 
	Next, we define the degenerate AIPW score as
	\begin{align*}
		\varphi_{\pi,\emptyset}(O;\alpha_{\pi,\emptyset},P_O) := &
		\dfrac{(Z_{\pi}- \rho_{\pi,\emptyset})(Y-\mu_{\pi,\emptyset}(A))}{\kappa_{\pi,\emptyset}(A)}+\mu_{\pi,\emptyset}(A)
		\\&-\left(\int z_{\pi}\mathrm{d} P_O(o) - \rho_{\pi,\emptyset}\right)\dfrac{\eta_{\emptyset}(A) - \mu_{\pi,\emptyset}(A)}{\kappa_{\pi,\emptyset}(A)}.
	\end{align*}
	As in Lemma~\ref{lem: mixed bias}, the mixed bias property still holds:
	\begin{equation}
		\begin{aligned}\label{exp: 26}
			&\mathbb{E}\left[\varphi_{\pi,\emptyset}(O;\alpha_{\pi,\emptyset},P_O)\middle| A=a\right]-\mu_{\pi,\emptyset}^o(a)\\
			=&\{\mu_{\pi,\emptyset}(a)-\mu_{\pi,\emptyset}^o(a)\}\left(1-\dfrac{\kappa_{\pi,\emptyset}^o(a)}{\kappa_{\pi,\emptyset}(a)}\right)+\left(\rho_{\pi,\emptyset} - \rho_{\pi,\emptyset}^o\right)
			\dfrac{\eta_{\emptyset}(a) - \eta_{\emptyset}^o(a)}{\kappa_{\pi,\emptyset}(a)}.
		\end{aligned}
	\end{equation}
	Next, the estimating procedure proceeds similarly to the discussion in Section~\ref{subsec: estimation procedure}.
	\begin{proof}[Proof of Equation~\eqref{exp: 26}:]
		We can directly calculate that
		\begin{align*}
			&\mathbb{E}\left[\varphi_{\pi,\emptyset}(O;\alpha_{\pi,\emptyset},P_O)\middle| A=a\right]-\mu_{\pi,\emptyset}^o(a)\\
			=&\mathbb{E}\left[
			\dfrac{(Z_{\pi}- \rho_{\pi,\emptyset})(Y-\mu_{\pi,\emptyset}(A))}{\kappa_{\pi,\emptyset}(A)}+\mu_{\pi,\emptyset}(A)
			\middle| A=a
			\right]-\mu_{\pi,\emptyset}^o(a)\\
			&-\left(\int z_{\pi}\mathrm{d} P_O(o) - \rho_{\pi,\emptyset}\right)
			\dfrac{\eta_{\emptyset}(a) - \mu_{\pi,\emptyset}(a)}{\kappa_{\pi,\emptyset}(a)}\\
			=&\dfrac{1}{\kappa_{\pi,\emptyset}(a)}\mathbb{E}\left[
			(Z_{\pi}- \rho_{\pi,\emptyset})(Y-\mu_{\pi,\emptyset}(a))\middle| A=a
			\right]+\mu_{\pi,\emptyset}(a)
			-\mu_{\pi,\emptyset}^o(a)
			\\&-\left(\rho_{\pi,\emptyset}^o - \rho_{\pi,\emptyset}\right)
			\dfrac{\eta_{\emptyset}(a) - \mu_{\pi,\emptyset}(a)}{\kappa_{\pi,\emptyset}(a)}\\
			=&\dfrac{1}{\kappa_{\pi,\emptyset}(a)}\mathbb{E}\left[
			Z_{\pi}Y+\rho_{\pi,\emptyset}\mu_{\pi,\emptyset}(a)-\mathbb{E}[Z_{\pi}|A=a]\mu_{\pi,\emptyset}(a)-\rho_{\pi,\emptyset}\mathbb{E}[Y|A=a]
			\middle| A=a\right]\\
			&+\mu_{\pi,\emptyset}(a)-\mu_{\pi,\emptyset}^o(a)-\left(\rho_{\pi,\emptyset}^o - \rho_{\pi,\emptyset}\right)
			\dfrac{\eta_{\emptyset}(a) - \mu_{\pi,\emptyset}(a)}{\kappa_{\pi,\emptyset}(a)}\\
			=&\dfrac{1}{\kappa_{\pi,\emptyset}(a)}\left\{\begin{array}{l}
				\mu_{\pi,\emptyset}^o(a)\kappa_{\pi,\emptyset}^o(a)+\rho_{\pi,\emptyset}^o\mathbb{E}[Y|A=a]\\
				+\rho_{\pi,\emptyset}\mu_{\pi,\emptyset}(a)
				-(\kappa_{\pi,\emptyset}^o(a)+\rho_{\pi,\emptyset}^o)\mu_{\pi,\emptyset}(a)-\rho_{\pi,\emptyset}\mathbb{E}[Y|A=a]
			\end{array}\right\}\\
			&+\mu_{\pi,\emptyset}(a)-\mu_{\pi,\emptyset}^o(a)-\left(\rho_{\pi,\emptyset}^o - \rho_{\pi,\emptyset}\right)
			\dfrac{\eta_{\emptyset}(a) - \mu_{\pi,\emptyset}(a)}{\kappa_{\pi,\emptyset}(a)}\\
			=&\dfrac{1}{\kappa_{\pi,\emptyset}(a)}\left\{\begin{array}{l}
				\{\mu_{\pi,\emptyset}(a)-\mu_{\pi,\emptyset}^o(a)\}\kappa_{\pi,\emptyset}^o(a)
				+\{\rho_{\pi,\emptyset}^o-\rho_{\pi,\emptyset}\}\eta_{\emptyset}^o(a)
				+\{\rho_{\pi,\emptyset}-\rho_{\pi,\emptyset}^o\}\mu_{\pi,\emptyset}(a)
			\end{array}\right\}\\
			&+\mu_{\pi,\emptyset}(a)-\mu_{\pi,\emptyset}^o(a)+\left(\rho_{\pi,\emptyset} - \rho_{\pi,\emptyset}^o\right)
			\dfrac{\eta_{\emptyset}(a) - \mu_{\pi,\emptyset}(a)}{\kappa_{\pi,\emptyset}(a)}\\
			=&\{\mu_{\pi,\emptyset}(a)-\mu_{\pi,\emptyset}^o(a)\}\left(1-\dfrac{\kappa_{\pi,\emptyset}^o(a)}{\kappa_{\pi,\emptyset}(a)}\right)+\left(\rho_{\pi,\emptyset} - \rho_{\pi,\emptyset}^o\right)
			\dfrac{\eta_{\emptyset}(a) - \eta_{\emptyset}^o(a)}{\kappa_{\pi,\emptyset}(a)}.
		\end{align*}
		Thus, Equation~\eqref{exp: 26} is verified.
	\end{proof}

	\subsection{Alternative DML framework for learning ADRFs}\label{append: alternative DML framework for learning ADRF}
	As mentioned in \citet{Colangelo16072025,bonvini2022fastconvergenceratesdoseresponse}, there exists an alternative debiased learning framework that does not rely on  $\varphi_{\pi}(O;\alpha_{\pi}^o, P_O)$. Instead, they simply used the kernel function $K_h(A-a)$ to substitute $I(A=a)$ in Equation~\eqref{exp: 23}:
	\begin{align*}
		&\varphi_{\pi,alt}(O;\theta_{\pi,a}^o,\alpha_{\pi}^o,h) := \frac{K_h(A-a)}{\Delta^o(a,L)} \frac{(Z_{\pi}- \rho_{\pi}^o(L))(Y - \mu_{\pi}^o(a,L))}{\kappa_{\pi}^o(a,L)}\\
		&+ \mu_{\pi}^o(a,L) - (Z_{\pi}- \rho_{\pi}^o(L)) \frac{\eta^o(a,L) - \mu_{\pi}^o(a,L)}{\kappa_{\pi}^o(a,L)}-\theta_{\pi,a}^o,
	\end{align*}
	where $\Delta^o(a,L) := p_{A|L}(a|L)$.
	The main distinction of this type of estimator from the methodology proposed in the main text is that it simultaneously defines the efficient influence function and incorporates the kernel function directly into the estimation procedure. To some extent, this constrains the learning process, as it must rely on kernel regression. \citet{bonvini2022fastconvergenceratesdoseresponse} further refine this approach by introducing the higher-order influence function corrections \citep{RobinsLiTchetgenVanderVaart2008}, which are beyond the scope of this paper.

	\begin{algorithm}[t]
		\caption{Nested cross-fitting procedure}
		\label{alg: cross-fitting alt}
		\begin{algorithmic}[1]
			\State \textbf{Input:} Pair of integers $[K,J]$, samples $\{O_i\}_{i=1}^n$.
			\State Randomly split the sample $\{O_i\}_{i=1}^n$ into $K$ folds $\{I_k\}_{k=1}^K$.
			\For{$k = 1, \ldots, K$}
			\State Train nuisance estimators on $I_{-k}$:
			\begin{align*}
				\hat{\alpha}_{\pi}^{(n,-k)} := 
				\big[\hat\mu_{\pi}^{(n,-k)}, \hat\rho_{\pi}^{(n,-k)}, 
				\hat\kappa_{\pi}^{(n,-k)}, \hat\eta^{(n,-k)}, 
				\hat\delta^{(n,-k)}\big].
			\end{align*}
			\State Randomly split the sample $I_k$ into $J$ folds $\{I_{k,j}\}_{j=1}^J$.
			Denote $I_{k,-j}:=\bigcup_{s\neq j} I_{k,s}$.
			\For{$j = 1, \ldots, J$} 
			\State Estimate the empirical distribution $\hat{P}_O^{(n,k,-j)}$ using $I_{k,-j}$: for any $f(O)$, 
			\begin{align*}
				\int f(o)\, d\hat{P}_O^{(n,k,-j)}(o)
				:= \frac{1}{|I_{k,-j}|}\sum_{i\in I_{k,-j}} f(O_i).
			\end{align*}
			\EndFor
			\EndFor
			\State \textbf{Output:} A vector $[\varphi_{\pi}(O_i; \hat\alpha_{\pi}^{(n,-k_i)}, \hat{P}_O^{(n,k_i,-j_i)}) ]_{i=1}^n$, where for all $i=1,\ldots, n$, $O_i\in I_{k_i,j_i}$. 
		\end{algorithmic}
	\end{algorithm}
	\subsection{Alternative cross-fitting procedure}\label{append: alternative cross-fitting procedure}
	We first give an alternative cross-fitting procedure in Algorithm~\ref{alg: cross-fitting alt}, named nested cross-fitting. Then, Equation~\eqref{eq: LLKR} can be modified to
	\begin{align*}
		e_1^{\top} \underset{\beta \in \mathbb{R}^2}{\arg\min}
		\frac{1}{n}\sum_{k=1}^{K} \sum_{j=1}^{J}
		\sum_{i\in I_{k,j}}
		K_h(A_i - a)
		\left\{\varphi_{\pi}(O_i; \hat\alpha_{\pi}^{(n,-k)}, \hat{P}_O^{(n,k,-j)})
		- g\left(\dfrac{A_i-a}{h}\right)^{\top}\beta\right\}^2.
	\end{align*}  
	Notably, the subsamples $I_{k,-j}$ used to estimate $\hat{P}_O^{(n,k,-j)}$ are independent of $I_{-k}$, the subsamples used to train $\hat{\alpha}_{\pi}^{(n,-k)}$. 
	Moreover, for all $i = 1, \ldots, n$, no observation $O_i\in I_{k_i,j_i}$ belongs to $I_{k_i,-j_i}$, ensuring mutual independence across folds. 
	Equivalently, for each pair $(k,j)$, the samples in $I_{k,j}$ are independent of those in $I_{k,-j}$ and $I_{-k}$. 
	Therefore, Theorems~\ref{thm: convergence rate} and~\ref{thm: clt}, as well as the estimation procedure described in Section~\ref{subsec: estimation procedure}, remain valid when Algorithm~\ref{alg: cross-fitting} is replaced by Algorithm~\ref{alg: cross-fitting alt}. 
	In practice, we recommend choosing $K \in \{2, 5, 10\}$ and $J = 2$.

	\section{Additional lemma and notation}\label{append: additional lemma and notation}
	% \subsection{Lemmas in local linear kernel regression}
	For completeness, we first give a brief review of an explicit solution to LLKR and define some necessary quantities. Recall that $e_1:=[1,0]^{\top}$ and $g(a):=[1,a]^{\top}$. We define 
	\begin{align*}
		&Q_{h}(a):=\left\{\sum_{i=1}^n K_h(A_i-a)g\left(\dfrac{A_i-a}{h}\right)g\left(\dfrac{A_i-a}{h}\right)^{\top}\right\}^{-1},\\
		&w_{ni}(a,h) := e_1^{\top} Q_{h}(a)g\left(\dfrac{A_i-a}{h}\right)K_h(A_i-a).
	\end{align*}
	In particular, $w_{ni}(a,h)$ is a random weight that relies only on $\{A_i\}_{i=1}^n$ and the selection of the bandwidth $h$.
	Then by solving Equation~\eqref{eq: LLKR} explicitly, we know that
	\begin{align*}
		&\hat{\theta}_{\pi,h}^{(n)}(a):=\sum_{k=1}^K\sum_{i\in I_k}w_{ni}(a,h)\varphi_{\pi}(O_i;\hat\alpha_{\pi}^{(n,-k)},\hat{P}_O^{(n,-k)}).
	\end{align*}

	\begin{lemma}\label{lem: lower bounded eigenvalue}
		Under Assumptions~\ref{as: continuity}~and~\ref{as: regular condition},
		\begin{align*}
			Q_{h}(a)=&\dfrac{1}{n}\left[\begin{array}{cc}
				\dfrac{1}{p_A(a)} & 0\\
				0 & \dfrac{1}{\int K(s)s^2\mathrm{d}s\times p_A(a)}
			\end{array}\right]+o_p(\dfrac{1}{n})
		\end{align*}
	\end{lemma}
	
	\begin{lemma}\label{lem: LLKR}
		For any $a\in\mathring{\mathcal{A}}$ and a bandwidth $h$, 
		$\sum_{i=1}^nw_{ni}(a,h)=1,$ $\sum_{i=1}^nw_{ni}(a,h)(A_i-a)=0.$
	\end{lemma}

	\subsection{Mixed bias property}
	Since the true distribution $P_{O}(o)$ and the nuisance function $\alpha_{\pi}(A, L)$ are unknown in practice, they must be properly estimated. 
	Consequently, it is necessary to demonstrate the mixed bias property when $\alpha_{\pi}(A, L)$ is plugged into the estimation.
	\begin{lemma}[Mixed bias property]\label{lem: mixed bias}
		Under Assumption~\ref{as: positivity}, suppose that $\pi(Z,L)$ is an RWF for $A=a$. Then, for any regular nuisance function $\alpha_{\pi}(A,L)$, $\mathbb{E}[\varphi_{\pi}(O;\alpha_{\pi}, P_O)\mid A=a] - \mathbb{E}[\mu_{\pi}^o(a,L)]$ equals
		\begin{align*}
			&\mathbb{E}\left[
			\dfrac{\delta(a,L)}{\kappa_{\pi}(a,L)} 
			(\rho_{\pi}(L)-\rho_{\pi}^o(L))(\eta(a,L)-\eta^o(a,L))
			\middle|  A=a\right]\\
			&+\mathbb{E}\left[ \left\{\delta^o(a,L)-
			\delta(a,L)\dfrac{\kappa_{\pi}^o(a,L)}{\kappa_{\pi}(a,L)}\right\}
			\{\mu_{\pi}(a,L)-\mu_{\pi}^o(a,L)\}\middle| A=a\right]\\
			&+\mathbb{E}\left[
			(\rho_{\pi}(L)-\rho_{\pi}^o(L))\left(\delta^o(a,L)-\delta(a,L)\right)\dfrac{\eta(a,L)-\mu_{\pi}(a,L)}{\kappa_{\pi}(a,L)}
			\middle| A=a\right].
		\end{align*}	
	\end{lemma}
	\begin{lemma}\label{lem: mixed bias property}
		Under Assumptions~\ref{as: positivity} and \ref{as: regular condition}, suppose that $\pi(Z,L)$ is a URWF for the support of $K_h(A-a)$. For any bounded and measurable function $f(a)$ with $|f(a)|$ monotone increasing with $a$, 
		\begin{align*}
			&\mathbb{E}_{k}\left[K_h(A-a)f(\dfrac{A-a}{h})\left\{\varphi_{\pi}(O;\hat\alpha_{\pi}^{(n,-k)},P_O) - \varphi_{\pi}(O;\alpha_{\pi}^o, P_O)\right\}\right]\\
			\lesssim &\left\|\hat\rho_{\pi}^{(n,-k)}-\rho_{\pi}^o\right\|_2\times
			\left(\left\|\hat\eta^{(n,-k)}-\eta^o\right\|_{\mathcal{B}(a,h),2}
			+\left\|\hat\delta^{(n,-k)}-\delta^o\right\|_{\mathcal{B}(a,h),2}\right)\\
			&+\left\|\hat\mu_{\pi}^{(n,-k)}-\mu_{\pi}^o\right\|_{\mathcal{B}(a,h),2}\times
			\left\|\hat\kappa_{\pi}^{(n,-k)}-\kappa_{\pi}^o\right\|_{\mathcal{B}(a,h),2}\\
			&+\left\|\hat\mu_{\pi}^{(n,-k)}-\mu_{\pi}^o\right\|_{\mathcal{B}(a,h),2}\times
			\left\|\hat\delta^{(n,-k)}-\delta^o\right\|_{\mathcal{B}(a,h),2} =o_p(c(n)).
		\end{align*}
	\end{lemma}
	Lemma~\ref{lem: mixed bias} shows that when $\alpha_{\pi}(A,L)$ is not well estimated, the bias can be decomposed into three components: (i) the interaction between the estimation errors of $\rho_{\pi}(L)$ and $\eta(a,L)$; (ii) the discrepancy in $\mu_{\pi}(a,L)$ together with the mismatch between $\delta(a,L)$ and $\kappa_{\pi}(a,L)$; and (iii) the joint misspecification of $\rho_{\pi}(L)$ and $\delta(a,L)$. Consequently, the bias remains small provided that, within each term, at least one of the involved nuisance components is estimated sufficiently well.
	
	Notably, Assumptions~\ref{as: consistency}--\ref{as: continuity} are not required for Lemma~\ref{lem: mixed bias}, as the definitions of $\psi_{\pi,q}^o$ and $\varphi_{\pi}(O;\alpha_{\pi},P_O)$ do not depend on them. Moreover, Assumption~\ref{as: IV relevance} is implicitly satisfied whenever the existence of an RWF is ensured. Finally, $Z$ is not required to be an AIV for $A$, and the semiparametric analysis is always conducted within the nonparametric model space.

	\subsection{Lemmas in mathematical analysis}
	The lemmas listed in this subsection are all fundamental and classical in mathematical analysis. Thus, we omit their proofs. For clarity, the finite open cover Lemma~\ref{lem: finite open cover theorem} is required for the proof of Proposition~\ref{prop: finite covering number}.
	The implicit function Lemma~\ref{lem: implicit function theorem} is needed in the proof of Proposition~\ref{prop: non-existence of a URWF}.
	The mean value Lemma~\ref{lem: mean value theorem for integrals} for integrals is adopted in the proof of Propositions~\ref{prop: weakness of binary IV}~and~\ref{prop: non-existence of a URWF}.
	
	\begin{lemma}[Finite open cover theorem]\label{lem: finite open cover theorem}
		Let $K \subseteq \mathbb{R}^n$ be compact, and let $\{U_\alpha\}$ be an open cover of $K$, i.e.,
		$K \subseteq \bigcup_{\alpha} U_\alpha,\; U_\alpha$ is open.
		Then there exists a finite collection $\{\alpha_1, \dots, \alpha_m\}$ such that
		$
		K \subseteq \bigcup_{i=1}^m U_{\alpha_i}.
		$
	\end{lemma}
	
	% \begin{lemma}[Implicit function theorem]\label{lem: implicit function theorem}
	% 	Let $F:\mathbb{R} \times \mathbb{R}^m \to \mathbb{R}$ be continuous near $(a_0, l_0)$, and assume there exists a direction $d_0 \in \mathbb{R}^m$ such that the directional derivative along $d_0$ 
	% 	\[
	% 	\nabla_l^{d_0} F(a_0, l_0) := \lim_{t \to 0} \frac{F(a_0, l_0 + t d_0) - F(a_0, l_0)}{t}
	% 	\] 
	% 	exists and is nonzero.
	% 	Then there exists a neighborhood $\mathcal{B}(a_0,h)$ and a continuous function $l = l(a)$ defined on $\mathcal{B}(a_0,h)$ such that
	% 	$F(a, l(a)) = F(a_0, l_0),$ for all $a \in \mathcal{B}(a_0,h).$
	% \end{lemma}

    \begin{lemma}[Implicit function theorem]\label{lem: implicit function theorem}
    	Let $F:\mathbb{R} \times \mathbb{R}^m \to \mathbb{R}$ be continuously differentiable near $(a_0,l_0)$, and assume that there exists a direction $d_0\in\mathbb{R}^m$ such that
    	\[
    	\nabla_l^{d_0}F(a_0,l_0)
    	:=
    	\lim_{t\to0}
    	\frac{F(a_0,l_0+td_0)-F(a_0,l_0)}{t}
    	\]
    	exists and is nonzero.
    	Then there exists a neighborhood $\mathcal{B}(a_0,h)$ and a continuous function $l=l(a)$ defined on $\mathcal{B}(a_0,h)$ such that $F(a,l(a))=F(a_0,l_0),$ for all $a\in\mathcal{B}(a_0,h)$.
    \end{lemma}
	
	\begin{lemma}[Mean value theorem for integrals]\label{lem: mean value theorem for integrals}
		Let $f:[a,b]\to \mathbb{R}$ and $g:[a,b]\to \mathbb{R}$ be continuous functions, with $g(x) >0$ on $[a,b]$. Then there exists $\varepsilon \in (a,b)$ such that
		\[\int_a^b f(x) g(x) \, dx = f(\varepsilon) \int_a^b g(x) \, dx.\]
	\end{lemma}

	\section{Proof of propositions}\label{append: proof of propositions}
	\subsection{Proof of Proposition~\ref{prop: positivity equivalence}}
	First, we observe that
	\begin{align*}
		\chi^2\left[ p_{Z|A,L}(\cdot|a,L)\| p_{Z|L}(\cdot|L)\right]=
		\mathrm{Var}\left[\dfrac{p_{Z|A,L}(Z\mid a,L)}{p_{Z|L}(Z\mid L)} \middle| L\right]=
		\mathrm{Var}\left[\dfrac{p_{A|Z,L}(a\mid Z,L)}{p_{A|L}(a\mid L)} \middle| L\right].
	\end{align*}
	Then, if there exists an $\epsilon_3(a)>0$ such that 
	$\mathrm{Var}[p_{A|Z,L}(a\mid Z,L) \mid L] \geq \epsilon_3(a) \text{ almost surely,}$
	\begin{align*}
		&|\epsilon_3(a)|\lesssim p_{A|L}(a| l)^2\mathrm{Var}\left[\dfrac{p_{A|Z,L}(a\mid Z,L)}{p_{A|L}(a\mid L)} \middle| L=l\right]\\
		\lesssim&
		\left|\mathrm{Var}\left[\dfrac{p_{A|Z,L}(a\mid Z,L)}{p_{A|L}(a\mid L)} \middle| L=l\right]\right|\sup_{l\in\mathcal{L}}\{p_{A|L}(a|l)^2\}\\
		\lesssim&
		\mathrm{Var}\left[\dfrac{p_{A|Z,L}(a\mid Z,L)}{p_{A|L}(a\mid L)} \middle| L=l\right]
		=\chi^2\left[ p_{Z|A,L}(\cdot|a,l)\| p_{Z|L}(\cdot|l)\right].
	\end{align*}
	It follows from Assumption~\ref{as: continuity} that $p_{A|L}(a|L)$ and $p_{A|Z,L}(a\mid Z,L)$ are uniformly bounded almost surely. Then Assumption~\ref{as: IV relevance} is satisfied.
	
	Conversely, under Assumptions~\ref{as: positivity}~and~\ref{as: IV relevance}, there exists $\epsilon_0(a)$ such that
	\begin{align*}
		&\mathrm{Var}[p_{A|Z,L}(a\mid Z,L) \mid L=l]\\
		=&\mathrm{Var}[p_{A|Z,L}(a\mid Z,L)/p_{A|L}(a\mid L) \mid L=l]\times p_{A|L}(a|l)^2
		\gtrsim\epsilon_0(a)^3.
	\end{align*}
	This completes the proof of Proposition~\ref{prop: positivity equivalence}.
	
	\subsection{Proof of Proposition~\ref{prop: weakness of binary IV}}
	To prove Proposition~\ref{prop: weakness of binary IV}, we first present a necessary condition for Assumption~\ref{as: IV relevance}. This lemma states that the IV relevance condition in Assumption~\ref{as: IV relevance} is violated if it leaves the density $p_{A|Z,L}(a|z,l)$ unchanged for all $z$.
	\begin{lemma}\label{lem: necessary condition for IV relevance}
		Under Assumption~\ref{as: IV relevance}, for any $a\in\mathring{\mathcal{A}}$ and $l\in\mathcal{L}$ such that $p_{A\mid Z,L}(a\mid z,l)$ is continuous in $a$, there exist two distinct values $z_0,z_1\in\mathcal{Z}$ such that 
		$p_{A\mid Z,L}(a\mid z_0,l) \neq p_{A\mid Z,L}(a\mid z_1,l).$
	\end{lemma}
	\begin{proof}[Proof of Lemma~\ref{lem: necessary condition for IV relevance}]
		We prove this statement by contradiction. If there exist $a_0\in\mathring{\mathcal{A}}$ and $l_0\in\mathcal{L}$, such that
		$p_{A|Z,L}(a_0\mid z_0,l_0)= p_{A|Z,L}(a_0\mid z_1,l_0)$ for any $z_0,z_1\in\mathcal{Z}$, we have
		\begin{align*}
			p_{A|L}(a_0\mid l_0)\equiv p_{A|Z,L}(a_0\mid z,l_0) 
		\end{align*}
		for any $z\in\mathcal{Z}$. Thus,
		$\mathrm{Var}\left[\dfrac{p_{A|Z,L}(a_0\mid Z,L)}{p_{A|L}(a_0\mid L)} \middle| L=l_0\right]
		=\mathrm{Var}\left[\dfrac{p_{A|L}(a_0\mid L)}{p_{A|L}(a_0\mid L)} \middle| L=l_0\right]=0.$
		Thus Assumption~\ref{as: IV relevance} is violated at $A=a_0$ and $L=l_0$. This finishes the proof of Lemma~\ref{lem: necessary condition for IV relevance}.     
	\end{proof}
	Next, when $Z$ is a binary IV, we define $R(a,l)=\dfrac{p_{A|Z,L}(a\mid 1,l)}{p_{A|Z,L}(a\mid 0,l)}.$ Since the support of $p_{A|Z,L}(a\mid z,L)$ is just $\mathcal{A}$, it is obvious that
	\begin{align*}
		&\mathbb{E}\left[R(A,L)\mid Z=0,L=l\right]=\int_{\mathcal{A}} R(a,l) p_{A|Z,L}(a\mid 0,l)\mathrm{d}  a\\
		=&\int_{\mathcal{A}} \dfrac{p_{A|Z,L}(a\mid 1,l)}{p_{A|Z,L}(a\mid 0,l)} p_{A|Z,L}(a\mid 0,l)\mathrm{d} a
		=\int_{\mathcal{A}}  p_{A|Z,L}(a\mid 1,l)\mathrm{d} a\equiv 1.
	\end{align*}
	That is, for any $l\in\mathcal{L}$, $\mathbb{E}\left[R(A,L)\mid Z=0,L=l\right]=1$. 
	Since $\mathcal{A}$ is a closed interval in $\mathbb{R}$, and $R(a,l)$ is a continuous function with respect to $a$ at $l$, from the mean value theorem for integrals (Lemma~\ref{lem: mean value theorem for integrals}), there must exist some $a_0\in\mathring{\mathcal{A}}$, such that $$1=R(a_0,l)\int_{\mathcal{A}} p_{A|Z,L}(a\mid 0,l)\mathrm{d}  a=R(a_0,l).$$ Thus, for any $l\in\mathcal{L}$, we can find an interior point $a_0\in\mathring{\mathcal{A}}$ such that $R(a_0,l)=1$, i.e., $p_{A|Z,L}(a_0\mid 1,l)=p_{A|Z,L}(a_0\mid 0,l)$. By Lemma~\ref{lem: necessary condition for IV relevance}, we know that Assumption~\ref{as: IV relevance} is violated, finishing the proof of Proposition~\ref{prop: weakness of binary IV}.

	\subsection{Proof of Proposition~\ref{prop: existence of RWFs}}

	We omit the subscript for the conditional density function for brevity.
	First, we note that $\mathbb{E}[p(a\mid Z,L)\mid L]=\int p(a\mid z,L)p(z\mid L)\mathrm{d} z=p(a\mid L)$, then
	\begin{align*}
		&\chi^2\left[p_{Z|A,L}(\cdot|a,L)\| p_{Z|L}(\cdot|L)\right]
		=\mathrm{Var}\left[\dfrac{p(Z|a,L)}{p(Z|L)}\middle| L\right]
		=\mathrm{Var}\left[\dfrac{p(a|Z,L)}{p(a|L)}\middle| L\right]\\
		=&\dfrac{1}{p(a|L)^2}\mathrm{Var}\left[p(a\mid Z,L)\mid L\right]
		=\dfrac{1}{p(a|L)^2}\mathbb{E}[(p(a\mid Z,L)-p(a\mid L))^2\mid L]\\
		=&\dfrac{1}{p(a|L)^2}\int \left\{
		p(a\mid z,L)-p(a\mid L)
		\right\}^2 p(z\mid L)\mathrm{d} z.
	\end{align*}
	From the Cauchy-Schwarz inequality, we know that
	\begin{align}
		&\mathbb{E}[\pi(Z,L)\mid A=a,L]-\mathbb{E}[\pi(Z,L)\mid L]\notag\\
		=&\int \pi(z,L)\left\{
		p(z\mid a,L)-
		p(z\mid L)
		\right\} \mathrm{d} z\notag\\
		=&\int \pi(z,L)\left\{
		\dfrac{p(z\mid L)p(a\mid z,L)}{p(a\mid L)}-
		p(z\mid L)
		\right\} \mathrm{d} z\notag\\
		=&\int \pi(z,L)\left\{
		\dfrac{p(a\mid z,L)-p(a\mid L)}{p(a\mid L)}
		\right\} p(z\mid L)\mathrm{d} z\notag\\
		=&\dfrac{1}{p(a\mid L)}\int \pi(z,L)\left\{
		p(a\mid z,L)-p(a\mid L)
		\right\} p(z\mid L)\mathrm{d} z\notag\\
		=&\dfrac{1}{p(a\mid L)}\int \left\{\pi(z,L) - \mathbb{E}[\pi(Z,L)\mid L]\right\}\left\{
		p(a\mid z,L)-p(a\mid L)
		\right\} p(z\mid L)\mathrm{d} z\label{exp: 5}\\
		\leq&
		\dfrac{1}{p(a\mid L)}\left[\int \left\{\pi(z,L) - \mathbb{E}[\pi(Z,L)\mid L]\right\}^2 p(z\mid L)\mathrm{d} z\right]^{1/2}\notag\\&\times
		\left[\int \left\{
		p(a\mid z,L)-p(a\mid L)
		\right\}^2 p(z\mid L)\mathrm{d} z\right]^{1/2}\notag\\
		=&
		\left[\int \left\{\pi(z,L) - \mathbb{E}[\pi(Z,L)\mid L]\right\}^2 p(z\mid L)\mathrm{d} z\right]^{1/2}\notag\\
		&\times
		\left[\int \left\{
		\dfrac{p(a\mid z,L)}{p(a\mid L)}-1
		\right\}^2 p(z\mid L)\mathrm{d} z\right]^{1/2}
		\lesssim
		(\mathrm{Var}\left[\dfrac{p(a\mid Z,L)}{p(a\mid L)}\middle| L\right])^{1/2}.\label{exp: 6}
	\end{align}
	The final inequality follows from the uniform boundedness of $\pi(Z,L)$.
	If $\pi(Z,L)$ is an RWF for $A=a$, then we know from Equation~\eqref{exp: 6} that there exists $\epsilon_0(a)$ such that
	\begin{align*}
		\mathrm{Var}\left[\dfrac{p(a\mid Z,L)}{p(a\mid L)}\middle| L\right]\geq \left|\mathbb{E}[\pi(Z,L)\mid A=a,L]-\mathbb{E}[\pi(Z,L)\mid L]\right|^2\geq \epsilon_0(a)^2.
	\end{align*}
	Conversely, if there exists $\epsilon_0(a)$, such that $\mathrm{Var}\left[\dfrac{p(a\mid Z,L)}{p(a\mid L)}\middle| L\right]\geq \epsilon_0(a)$, then from Equation~\eqref{exp: 5}, we can just take $\pi(Z,L) = p(a\mid Z,L)/p(a\mid L)$ to demonstrate that
	\begin{align}
		&\mathbb{E}[p(a\mid Z,L)/p(a\mid L)\mid A=a,L]-\mathbb{E}[p(a\mid Z,L)/p(a\mid L)\mid L]\notag\\
		=&\dfrac{1}{p(a\mid L)^2}\int \left\{
		p(a\mid z,L)-p(a\mid L)
		\right\}^2 p(z\mid L)\mathrm{d} z\notag\\
		=&\mathrm{Var}\left[\dfrac{p(a\mid Z,L)}{p(a\mid L)}\middle| L\right]\gtrsim \epsilon_0(a).\label{exp: 7}
	\end{align}
	The final inequality follows from the boundedness of $p(a\mid L)$. This completes the proof of the first statement.
	
	Next, if an RWF $\pi(Z,L)$ exists for $A=a$, then we can use Equations~\eqref{exp: 7} and \eqref{exp: 6} sequentially to deduce that 
	\begin{align*}
		&\mathbb{E}[p(a\mid Z,L)/p(a\mid L)\mid A=a,L]-\mathbb{E}[p(a\mid Z,L)/p(a\mid L)\mid L]\\
		=&\dfrac{1}{p(a\mid L)^2}\mathrm{Var}\left[p(a\mid Z,L)\mid L\right]\\
		\gtrsim& |\mathbb{E}[\pi(Z,L)\mid A=a,L]-\mathbb{E}[\pi(Z,L)\mid L]|^2\geq \epsilon_0(a).
	\end{align*}
	Thus, $\pi_a(Z,L) := p(a\mid Z,L)/p(a\mid L)$ is an RWF for $A=a$. Now we finish the proof of Proposition~\ref{prop: existence of RWFs}.

	\subsection{Proof of Proposition~\ref{prop: local stability of RWF}}
	Since $\pi(Z,L)$ is an RWF for $A=a_0\in\mathring{\mathcal{A}}$, there exists $\epsilon_0(a_0)$, such that
	$|\kappa_{\pi}^o(a_0,l)|\geq \epsilon_0(a_0)$. Then since $\kappa_{\pi}^o(a,l)$ is equicontinuous  at $a_0$ over $L$ almost surely, there exists an $r_0$, such that for any $a\in\mathcal{B}(a_0, r_0)$ and $l\in\mathcal{L}$, $|\kappa_{\pi}^o(a,l)-\kappa_{\pi}^o(a_0,l)|\leq\epsilon_0(a_0)/2$.
	\begin{align*}
		|\kappa_{\pi}^o(a,l)|\geq& 	|\kappa_{\pi}^o(a_0,l)| - |\kappa_{\pi}^o(a,l)-\kappa_{\pi}^o(a_0,l)|
		\geq\epsilon_0(a_0) - \epsilon_0(a_0)/2\geq \epsilon_0(a_0)/2.
	\end{align*}
	This asserts that $\pi(Z,L)$ is a URWF for $\mathcal{B}(a_0,h)$.

	\subsection{Proof of Proposition~\ref{prop: non-existence of a URWF}}
	For any $l_0\in\mathcal{L}$, since $\mathbb{E}[\kappa_{\pi}^o(A,L)\mid L=l_0]=\mathbb{E}[\mathbb{E}[\pi(Z,L)\mid A,L]\mid L=l_0]-\mathbb{E}[\pi(Z,L)\mid L=l_0]=0$, 
	\begin{align*}
		\int_{\mathcal{A}} \kappa_{\pi}^o(a,l_0)p_{A|L}(a|l_0)\mathrm{d} a=0.
	\end{align*}
	Since $p_{A|L}(a|l_0)$ is positive and continuous for all $a\in\mathcal{A}$ almost surely (we assumed that $p_{A|L}(a|L)$ and $p_{A}(a)$ share the same support $\mathcal{A}$ throughout the paper), and $\kappa_{\pi}^o(a,l_0)$ is continuous in $\mathcal{A}$, by the mean value theorem of integrals, there must exist some $a(l_0)\in\mathring{\mathcal{A}}$, such that $\kappa_{\pi}^o(a(l_0),l_0)=0$. This means that $\pi(Z, L)$ is not an RWF for $A=a(l_0)$.
	
	In addition, if for  fixed $l_0\in\mathcal{L}$, the directional derivative of $\kappa_{\pi}^o(a,l)$ with respect to $l$ is nonzero at $(a(l_0),l_0)$, then by the implicit function theorem, there exists a small $h$, such that for any $a\in\mathcal{B}(a(l_0),h)$, there exists $l(a)\in\mathcal{L}$, such that $\kappa_{\pi}^o(a,l(a))\equiv 0.$ 
	This indicates that $\pi(Z,L)$ is not an RWF for all $A=a\in\mathcal{B}(a(l_0),h)$.
	This completes the proof of Proposition~\ref{prop: non-existence of a URWF}.

	\subsection{Proof of Proposition~\ref{prop: finite covering number}}
	For any fixed $a\in\mathring{\mathcal{A}}$, by Proposition~\ref{prop: existence of RWFs}, under Assumptions~\ref{as: continuity},~\ref{as: IV relevance}, $\pi_a(Z,L)$ is already an RWF for $A=a$. Moreover, by Proposition~\ref{prop: local stability of RWF}, since $\kappa_{\pi_a}^o(A,L)$ is equicontinuous with respect to all $a\in\mathring{\mathcal{A}}$, there exists a small constant $h(a)$, such that $\pi_a(Z,L)$ is a URWF for $\mathcal{B}(a,h(a))$. 
	
	Thus, for any compact subset $\mathcal{A}_c\subseteq \mathring{\mathcal{A}}$,
	\begin{align*}
		\mathcal{A}_c\subseteq \bigcup_{a\in\mathring{\mathcal{A}}}\mathcal{B}(a,h(a)).
	\end{align*}
	By the finite open cover theorem, we know that there exists a finite collection $\{a_{m}\}_{m=1}^M$, such that
	\begin{align*}
		\mathcal{A}_c\subseteq \bigcup_{m=1}^M \mathcal{B}(a_{m},h(a_{m})).
	\end{align*}
	Finally, due to our selection of $h(a_{m})$,
	$\pi_{a_{m}}(Z,L)$ is a URWF for $\mathcal{B}(a_{m},h(a_{m}))$, 
	finishing the proof of Proposition~\ref{prop: finite covering number}.

	\subsection{Proof of Proposition~\ref{prop: AIV characterization}}
    Without loss of generality, we define $0/0=1$.
	We omit the subscript for the conditional density function for brevity. 
	\begin{align}
		&\mathbb{E}\left[p(a\mid U,L)\mathbb{E}[\pi(Z,L)\mid A=a,U,L]\middle| L\right]\notag\\
		=&\mathbb{E}\left[p(a\mid U,L)\int \pi(z,L)p(z\mid a,U,L)\mathrm{d} z\middle| L\right]\notag\\
		=&\mathbb{E}\left[\int \pi(z,L)p(z,a\mid U,L)\mathrm{d} z\middle| L\right]\notag\\
		=&\mathbb{E}\left[\int \pi(z,L)p(a\mid z,U,L)p(z\mid U,L)\mathrm{d} z\middle| L\right]\notag\\
		=&\mathbb{E}\left[\int \pi(z,L)p(a\mid z,U,L)p(z\mid L)\mathrm{d} z\middle| L\right]\label{exp: 8}\\
		=&\int \pi(z,L)\mathbb{E}\left[p(a\mid z,U,L)\middle| L\right]p(z\mid L)\mathrm{d} z\notag\\
		=&\int \pi(z,L)\mathbb{E}\left[p(a\mid z,U,L)\middle| Z=z,L\right]p(z\mid L)\mathrm{d} z\label{exp: 9}\\
		=&\int \pi(z,L)p(a\mid z,L)p(z\mid L)\mathrm{d} z
		=\int \pi(z,L)p(z\mid a,L)\mathrm{d} z\times p(a\mid L)\notag\\
		=&\mathbb{E}[\pi(Z,L)\mid A=a,L]p(a\mid L)\notag.
	\end{align}
	Here Equations~\eqref{exp: 8},~\eqref{exp: 9} follow from Assumption~\ref{as: IV independence} that $Z\indep U\mid L$.
	Thus, we can verify that
	\begin{align*}
		\mathbb{E}[\omega_{a,\pi}(U,L)\mid L]
		=&\dfrac{1}{p(a\mid L)}
		\dfrac{\mathbb{E}[p(a\mid U,L)\mathbb{E}[\pi(Z,L)\mid A=a,U,L]\mid L] - p(a\mid L)\mathbb{E}[\pi(Z,L)\mid L]}{\mathbb{E}[\pi(Z,L)\mid A=a,L] - \mathbb{E}[\pi(Z,L)\mid L]}\\
		=&\dfrac{1}{p(a\mid L)}
		\dfrac{p(a\mid L)\mathbb{E}[\pi(Z,L)\mid A=a,L] - p(a\mid L)\mathbb{E}[\pi(Z,L)\mid L]}{\mathbb{E}[\pi(Z,L)\mid A=a,L] - \mathbb{E}[\pi(Z,L)\mid L]}\\
		=&\dfrac{\mathbb{E}[\pi(Z,L)\mid A=a,L] - \mathbb{E}[\pi(Z,L)\mid L]}{\mathbb{E}[\pi(Z,L)\mid A=a,L] - \mathbb{E}[\pi(Z,L)\mid L]}\equiv 1.
	\end{align*}
	This concludes the proof of the first statement. 
	Next, we prove the second result. 
	If $Z$ is an AIV for $A=a$, there exist $b_a(U,L)$ and $c_a(Z,L)$, such that 
	\begin{align*}
		p(a|Z,U,L)= b_a(U,L) + c_a(Z,L).
	\end{align*}
	Next, we can calculate that
	\begin{align*}
		&\mathbb{E}[p(a\mid U,L)\left\{\pi(Z,L)-\mathbb{E}[\pi(Z,L)|L]\right\}\mid A=a,U,L]\\
		=&p(a\mid U,L)\int \left\{\pi(z,L)-\mathbb{E}[\pi(Z,L)|L]\right\}p(z\mid a,U,L)\mathrm{d} z\\
		=&\int \left\{\pi(z,L)-\mathbb{E}[\pi(Z,L)|L]\right\}p(a\mid z,U,L)p(z\mid U,L)\mathrm{d} z\\
		=&\int \left\{\pi(z,L)-\mathbb{E}[\pi(Z,L)|L]\right\}\{b_a(U,L) + c_a(z,L)\}p(z\mid L)\mathrm{d} z\\
		=&\int \left\{\pi(z,L)-\mathbb{E}[\pi(Z,L)|L]\right\}c_a(z,L)p(z\mid L)\mathrm{d} z\\
		&+b_a(U,L)\int \left\{\pi(z,L)-\mathbb{E}[\pi(Z,L)|L]\right\}p(z\mid L)\mathrm{d} z\\
		=&\int \left\{\pi(z,L)-\mathbb{E}[\pi(Z,L)|L]\right\}c_a(z,L)p(z\mid L)\mathrm{d} z
	\end{align*}
	The last equality holds since $\mathbb{E}[\pi(Z,L)|L]=\int \pi(z,L)p(z\mid L)\mathrm{d} z$.
	Now we see that the numerator part of $\omega_{a,\pi}(U,L)$ no longer depends on $U$, so 
	\begin{align*}
		\omega_{a,\pi}(U,L)\equiv\mathbb{E}[\omega_{a,\pi}(U,L)\mid L]\equiv 1.
	\end{align*}
	Conversely, if $\omega_{a,\pi}(U,L)\equiv 1$, i.e., for any $\pi(Z,L)$,
	\begin{align*}
		&p(a\mid U,L)\left\{\mathbb{E}[\pi(Z,L)\mid A=a,U,L] - \mathbb{E}[\pi(Z,L)\mid L]\right\}\\
		=&p(a\mid L)\left\{\mathbb{E}[\pi(Z,L)\mid A=a,L] - \mathbb{E}[\pi(Z,L)\mid L]\right\}.
	\end{align*}
	In particular, if we take $\pi(Z,L):= I\{Z\in \mathcal{S}\}$, 
	\begin{align*}
		&p(a\mid U,L)\left\{\mathbb{E}[I\{Z\in \mathcal{S}\}\mid A=a,U,L] - \mathbb{E}[I\{Z\in \mathcal{S}\}\mid L]\right\}\\
		=&p(a\mid U,L)\mathbb{E}[\left\{I\{Z\in \mathcal{S}\}-\mathbb{E}[I\{Z\in \mathcal{S}\}\mid L]\right\}\mid A=a,U,L]\\
		=&p(a\mid U,L)\mathbb{E}[\left\{I\{Z\in \mathcal{S}\}-\mathbb{E}[I\{Z\in \mathcal{S}\}\mid L]\right\}\dfrac{p(a\mid Z,U,L)}{p(a\mid U,L)}\mid U,L]\\
		=&\mathbb{E}[\left\{I\{Z\in \mathcal{S}\}-\mathbb{E}[I\{Z\in \mathcal{S}\}\mid L]\right\}p(a\mid Z,U,L)\mid U,L]\\
		=&\mathbb{E}[\left\{I\{Z\in \mathcal{S}\}-\mathbb{E}[I\{Z\in \mathcal{S}\}\mid U,L]\right\}p(a\mid Z,U,L)\mid U,L]\\
		=&\mathbb{E}[I\{Z\in \mathcal{S}\}\left\{p(a\mid Z,U,L)-p(a\mid U,L)\right\}\mid U,L]\\
		=&\mathbb{E}[\left\{p(a\mid Z,U,L)-p(a\mid U,L)\right\}\mid Z\in\mathcal{S},U,L]\Pr(Z\in\mathcal{S}\mid L).
	\end{align*}
	Analogously, one can compute that
	\begin{align*}
		&p(a\mid L)\left\{\mathbb{E}[\pi(Z,L)\mid A=a,L] - \mathbb{E}[\pi(Z,L)\mid L]\right\}\\
		=&\mathbb{E}[\left\{I\{Z\in \mathcal{S}\}-\mathbb{E}[I\{Z\in \mathcal{S}\}\mid L]\right\}p(a\mid Z,L)\mid L]\\
		=&\mathbb{E}[I\{Z\in \mathcal{S}\}\left\{p(a\mid Z,L)-p(a\mid L)\right\}\mid L]\\
		=&\mathbb{E}[\left\{p(a\mid Z,L)-p(a\mid L)\right\}\mid Z\in \mathcal{S},L]\Pr(Z\in\mathcal{S}\mid L).
	\end{align*}
	These two equalities hold for every measurable set $\mathcal{S}\subseteq\mathcal{Z}$. Thus, we know that
	\begin{align*}
		&\mathbb{E}[\left\{p(a\mid Z,U,L)-p(a\mid U,L)\right\}\mid Z,U,L]
		=\mathbb{E}[\left\{p(a\mid Z,L)-p(a\mid L)\right\}\mid Z,L]\\
		\Rightarrow&\quad
		p(a\mid Z,U,L)-p(a\mid U,L)\equiv p(a\mid Z,L)-p(a\mid L)\\
		\Rightarrow&\quad
		p(a\mid Z,U,L)\equiv p(a\mid U,L)+\left\{p(a\mid Z,L)-p(a\mid L)\right\}.
	\end{align*}
	Thus, we can just take $b_a(U,L)=p(a\mid U,L)$ and $c_a(Z,L)=p(a\mid Z,L)-p(a\mid L)$ to formulate $p(a\mid Z,U,L)=b_a(U,L)+c_a(Z,L).$
	This completes the proof of Proposition~\ref{prop: AIV characterization}.

	\subsection{Proof of Proposition~\ref{prop: equivalent characterization of WAIV}}
	When $Z$ is a WAIV for $A$, there exist functions $b_a(U,L)$, $c_a(Z,L)$, and $d(U,L)$ such that
	\[
	p_{A\mid Z,U,L}(a \mid Z, U, L) = b_a(U, L) + d(U,L)c_a(Z, L).
	\]
	Then, we can verify that
	\begin{align*}
		\omega_{a,\pi}(U,L)
		=&\frac{\mathrm{Cov}\{ b_a(U, L) + d(U,L)c_a(Z, L),\, \pi(Z,L) \mid U,L \}}
		{\mathrm{Cov}\{ b_a(U, L) + d(U,L)c_a(Z, L),\, \pi(Z,L) \mid L \}}\\
		=&\frac{d(U,L)\mathrm{Cov}\{ c_a(Z, L),\, \pi(Z,L) \mid U,L \}}
		{\mathrm{Cov}\{ d(U,L)c_a(Z, L),\, \pi(Z,L) \mid L \}}\\
		=&\frac{d(U,L)\mathrm{Cov}\{ c_a(Z, L),\, \pi(Z,L) \mid L \}}
		{\mathbb{E}[d(U,L)\mid L]\mathbb{E}[c_a(Z, L)\pi(Z,L)\mid L] 
			- \mathbb{E}[d(U,L)\mid L]\mathbb{E}[c_a(Z, L)\mid L]\mathbb{E}[\pi(Z,L)\mid L]}\\
		=&\frac{d(U,L)}{\mathbb{E}[d(U,L)\mid L]}.
	\end{align*}
	Here the third equality adopts Assumption~\ref{as: IV independence}.
	Hence, $\omega_{a,\pi}(U,L)$ does not depend on $a$ and $\pi$ apparently.
	
	Conversely, suppose that $\nabla_a \omega_{a,\pi}(U,L)\equiv 0$ for all $\pi(Z,L)$. Taking $\pi(Z,L)=I\{Z\in\mathcal{S}\}$, there exists $\omega_0(U,L)$ such that
	\begin{align*}
		\omega_0(U,L)
		=&\frac{\mathrm{Cov}\{ p_{A\mid Z,U,L}(a\mid Z,U,L),\, I\{Z\in\mathcal{S}\} \mid U,L \}}
		{\mathrm{Cov}\{ p_{A\mid Z,L}(a\mid Z,L),\, I\{Z\in\mathcal{S}\} \mid L \}}\\
		=&\frac{p_{A\mid Z,U,L}(a\mid \mathcal{S},U,L)-p_{A\mid U,L}(a\mid U,L)}
		{p_{A\mid Z,L}(a\mid \mathcal{S},L)-p_{A\mid L}(a\mid L)}.
	\end{align*}
	Letting $\mathcal{S}=[z-h,z+h]$ and taking $h\to 0$, we obtain
	\begin{align*}
		\omega_0(U,L)
		=\frac{p_{A\mid Z,U,L}(a\mid z,U,L)-p_{A\mid U,L}(a\mid U,L)}
		{p_{A\mid Z,L}(a\mid z,L)-p_{A\mid L}(a\mid L)},
	\end{align*}
	which implies
	\[
	p_{A\mid Z,U,L}(a\mid Z,U,L)
	= p_{A\mid U,L}(a\mid U,L)
	+ \omega_0(U,L)\Big\{p_{A\mid Z,L}(a\mid Z,L)-p_{A\mid L}(a\mid L)\Big\}.
	\]
	This recovers the WAIV structure and completes the proof of Proposition~\ref{prop: equivalent characterization of WAIV}.

	\subsection{Proof of Proposition~\ref{prop: AIV transformation}}
	We provide the proof of the transformation invariance property for the AIV condition only. The corresponding result for the WAIV condition can be established by analogous arguments and is therefore omitted for brevity.
	If $Z$ is an AIV for $A$, then there exist functions $b_a(U,L)$ and $c_a(Z,L)$ such that
	\begin{align*}
		p_{A|Z,U,L}(a|Z,U,L) = b_a(U,L) + c_a(Z,L).
	\end{align*}
	
	Denote $\tilde{A} = q(A)$. By the change-of-variable formula, we have
	\begin{align*}
		p_{\tilde{A}|Z,U,L}(\tilde{a}|Z,U,L) 
		&=  \frac{p_{A|Z,U,L}\big(q^{-1}(\tilde{a})\,\big|\,Z,U,L\big)}{|q'\big(q^{-1}(\tilde{a})\big)|}
		= \frac{b_{q^{-1}(\tilde{a})}(U,L) + c_{q^{-1}(\tilde{a})}(Z,L)}{|q'\big(q^{-1}(\tilde{a})\big)|}.
	\end{align*}
	
	Setting 
	\[
	\tilde{b}_{\tilde{a}}(U,L) = \frac{b_{q^{-1}(\tilde{a})}(U,L)}{|q'\big(q^{-1}(\tilde{a})\big)|},\quad 
	\tilde{c}_{\tilde{a}}(Z,L) = \frac{c_{q^{-1}(\tilde{a})}(Z,L)}{|q'\big(q^{-1}(\tilde{a})\big)|},
	\]
	we obtain $p_{\tilde{A}|Z,U,L}(\tilde{a}|Z,U,L) = \tilde{b}_{\tilde{a}}(U,L) + \tilde{c}_{\tilde{a}}(Z,L).$ This completes the proof of Proposition~\ref{prop: AIV transformation}.

	\subsection{Proof of Proposition~\ref{prop: existence NPIV}}
	We can directly calculate that
	\begin{align*}
		&p_{A|Z,L}(a|Z,L)\mathbb{E}[Y|A=a,Z,L]=
		p_{A|Z,L}(a|Z,L)\mathbb{E}[\mathbb{E}[Y|A=a,U,Z,L]|A=a,Z,L]\\
		=&p_{A|Z,L}(a|Z,L)\mathbb{E}[\mathbb{E}[Y(a)|U,L]|A=a,Z,L]\\
		=&\mathbb{E}[p_{A|Z,U,L}(a|Z,U,L)\mathbb{E}[Y(a)|U,L]|Z,L]
		=\mathbb{E}[\{b_a(U,L) + c_a(Z,L)\}\mathbb{E}[Y(a)|U,L]|Z,L]\tag{By the AIV condition}\\
		=&\mathbb{E}[b_a(U,L)\mathbb{E}[Y(a)|U,L]\mid L]
		+c_a(Z,L)\mathbb{E}[\mathbb{E}[Y(a)|U,L]\mid L]\tag{By Assumption~\ref{as: IV independence} that $Z\indep U\mid L$}\\
		=&\mathbb{E}[b_a(U,L)Y(a)\mid L]
		+c_a(Z,L)\mathbb{E}[Y(a)\mid L]\\
		=&\mathbb{E}[b_a(U,L)Y(a)\mid L]- \mathbb{E}[b_a(U,L)\mid L]\mathbb{E}[Y(a)\mid L]
		+p_{A|Z,L}(a|Z,L)\mathbb{E}[Y(a)\mid L]\\
		=&\mathrm{Cov}\{b_a(U,L),\;Y(a)\mid L\}
		+p_{A|Z,L}(a|Z,L)\mathbb{E}[Y(a)\mid L]\\
		=&\mathrm{Cov}\{b_a(U,L)+c_a(Z,L),\;\mathbb{E}[Y(a)\mid U,L]\mid L\}
		+p_{A|Z,L}(a|Z,L)\mathbb{E}[Y(a)\mid L]\\
		=&\mathrm{Cov}\!\{\mathbb{E}[Y(a)\mid U,L],\;p_{A|U,L}(a\mid U,L)\mid L\}
		+p_{A|Z,L}(a|Z,L)\mathbb{E}[Y(a)\mid L].
	\end{align*}
	Substituting the explicit forms of $f_a^o(0,L)$ and $f_a^o(1,L)$ into 
	Equation~\eqref{eq: NPIV} verifies the desired identity, thereby completing 
	the proof of Proposition~\ref{prop: existence NPIV}.
	
	\subsection{Proof of Proposition~\ref{prop: uniqueness NPIV}}
	If there exist two different solutions $f_a^o(c,L),g_a^o(c,L)$ ($c=0,1$) to Equation~\eqref{eq: NPIV}, then
	\begin{align*}
		p_{A\mid Z,L}(a\mid Z,L)\mathbb{E}[Y \mid A=a,Z,L]
		= p_{A\mid Z,L}(a\mid Z,L)
		\{f_a^o(1,L)-f_a^o(0,L)\}
		+ f_a^o(0,L);\\
		p_{A\mid Z,L}(a\mid Z,L)\mathbb{E}[Y \mid A=a,Z,L]
		= p_{A\mid Z,L}(a\mid Z,L)
		\{g_a^o(1,L)-g_a^o(0,L)\}
		+ g_a^o(0,L).
	\end{align*}
	This implies that
	\begin{align*}
		&p_{A\mid Z,L}(a\mid Z,L)[\{f_a^o(1,L)-f_a^o(0,L)\}-\{g_a^o(1,L)-g_a^o(0,L)\}]
		= g_a^o(0,L)-f_a^o(0,L)\\
		=&p_{A\mid L}(a\mid L)[\{f_a^o(1,L)-f_a^o(0,L)\}-\{g_a^o(1,L)-g_a^o(0,L)\}]\\
		\Rightarrow\quad&
		\left\{\dfrac{p_{A\mid Z,L}(a\mid Z,L)}{p_{A\mid L}(a\mid L)}-1\right\}
		[\{f_a^o(1,L)-f_a^o(0,L)\}-\{g_a^o(1,L)-g_a^o(0,L)\}]\equiv 0.\\
		\Rightarrow\quad&
		\mathbb{E}\left[\left\{\dfrac{p_{A\mid Z,L}(a\mid Z,L)}{p_{A\mid L}(a\mid L)}-1\right\}^2\middle| L\right]
		[\{f_a^o(1,L)-f_a^o(0,L)\}-\{g_a^o(1,L)-g_a^o(0,L)\}]\equiv 0.
	\end{align*}
	By Assumption~\ref{as: IV relevance}, $\mathbb{E}[\{\frac{p_{A\mid Z,L}(a\mid Z,L)}{p_{A\mid L}(a\mid L)}-1\}^2| L]$ is always bounded away from zero, implying that
	\begin{align*}
		\{f_a^o(1,L)-f_a^o(0,L)\}-\{g_a^o(1,L)-g_a^o(0,L)\}\equiv 0,\text{ and }
		g_a^o(0,L)-f_a^o(0,L)\equiv 0.
	\end{align*}
	This implies that the solution to Equation~\eqref{eq: NPIV} is unique. 
	Next, if we substitute the solution to Equation~\eqref{eq: NPIV} into the definition of $\mu_{\pi}^o(a,l)$,
	\begin{align*}
		\mu_{\pi}^o(a,l) :=& \frac{\mathbb{E}[Y Z_{\pi} \mid A=a,L=l] - \mathbb{E}[Y \mid A=a,L=l]\,\mathbb{E}[Z_{\pi} \mid L=l]}{\mathbb{E}[Z_{\pi} \mid A=a,L=l] - \mathbb{E}[Z_{\pi} \mid L=l]}\\
		=&\frac{\mathbb{E}[p_{A|Z,L}(a|Z,L)\{f_a^o(1,L)-f_a^o(0,L)\} \{Z_{\pi} -\mathbb{E}[Z_{\pi} \mid L=l]\}\mid A=a,L=l]}{p_{A|Z,L}(a|Z,L)\{\mathbb{E}[Z_{\pi} \mid A=a,L=l] - \mathbb{E}[Z_{\pi} \mid L=l]\}}\\
		=&f_a^o(1,l)-f_a^o(0,l).
	\end{align*}
	Then, we know from Equation~\eqref{eq: NPIV} that
	\begin{align*}
		p_{A\mid Z,L}(a\mid Z,L)\mathbb{E}[Y \mid A=a,Z,L]
		= p_{A\mid Z,L}(a\mid Z,L)\mu_{\pi}^o(a,L) + f_a^o(0,L).
	\end{align*}
	This indicates that
	\begin{align*}
		f_a^o(0,L) =& p_{A\mid Z,L}(a\mid Z,L)\mathbb{E}[Y \mid A=a,Z,L] - p_{A\mid Z,L}(a\mid Z,L)\mu_{\pi}^o(a,L);\\
		f_a^o(1,L)
		=& f_a^o(0,L) + \mu_{\pi}^o(a,L)\\
		=& p_{A\mid Z,L}(a\mid Z,L)\mathbb{E}[Y \mid A=a,Z,L] + \{1 - p_{A\mid Z,L}(a\mid Z,L)\}\mu_{\pi}^o(a,L).
	\end{align*}
	This finishes the proof of Proposition~\ref{prop: uniqueness NPIV}.

	\section{Proof of theorems}
	\label{append: proof of theorems}
	\subsection{Proof of Theorem~\ref{thm: identification}}
	For any fixed $a\in\mathring{\mathcal{A}}$, denote $g(u,L):=\mathbb{E}[Y(a)\mid U,L]$. We first observe that
	\begin{align*}
		&\mathbb{E}[Y\pi(Z,L)\mid A=a,L]=\mathbb{E}[\mathbb{E}[Y\mid Z,A=a,U,L]\pi(Z,L)\mid A=a,L]\\
		=&\mathbb{E}[\mathbb{E}[Y(a)\mid U,L]\pi(Z,L)\mid A=a,L]\\
		=&\mathbb{E}[\mathbb{E}[Y(a)\mid U,L]\pi(Z,L)\dfrac{p(a|U,Z,L)}{p(a|L)}\mid L].
	\end{align*}
	Notably, the second equality uses Assumption~\ref{as: latent ignorability} that $Y(a)\indep Z, A\mid U, L$.
	Similarly, we can calculate that
	\begin{align*}
		&\mathbb{E}[Y\left\{\pi(Z,L)-\mathbb{E}[\pi(Z,L)\mid L]\right\}\mid A=a,L]\\
		=&\mathbb{E}[\mathbb{E}[Y(a)\mid U,L]\left\{\pi(Z,L)-\mathbb{E}[\pi(Z,L)\mid  L]\right\}\dfrac{p(a|U,Z,L)}{p(a|L)}\mid L]\\
		=&\mathbb{E}[\mathbb{E}[Y(a)\mid U,L]\left\{\pi(Z,L)-\mathbb{E}[\pi(Z,L)\mid U, L]\right\}\dfrac{p(a|U,Z,L)}{p(a|L)}\mid L]\\
		=&\mathbb{E}[\mathbb{E}[Y(a)\mid U,L]\pi(Z,L)\left\{
		\dfrac{p(a|U,Z,L)}{p(a|L)}-\mathbb{E}\left[\dfrac{p(a|U,Z,L)}{p(a|L)}\middle| U,L\right]
		\right\}\mid L]\\
		=&\mathbb{E}[\mathbb{E}[Y(a)\mid U,L]\pi(Z,L)\left\{
		\dfrac{p(a|U,Z,L)}{p(a|L)}-\dfrac{p(a|U,L)}{p(a|L)}
		\right\}\mid L]\\
		=&\dfrac{1}{p(a|L)}\mathbb{E}[\mathbb{E}[Y(a)\mid U,L]\pi(Z,L)\left\{
		p(a|U,Z,L)-p(a|U,L)
		\right\}\mid L]\\
		=&\dfrac{1}{p(a|L)}\mathbb{E}[\mathbb{E}[Y(a)\mid U,L]\left\{\begin{array}{l}
			\mathbb{E}[\pi(Z,L)p(a|U,Z,L)\mid U,L]\\
			-\mathbb{E}[\pi(Z,L)\mid U,L]p(a|U,L)
		\end{array}\right\}\mid L]\\
		=&\dfrac{1}{p(a|L)}\mathbb{E}[\mathbb{E}[Y(a)\mid U,L]\left\{\begin{array}{l}
			\mathbb{E}[\pi(Z,L)\mid A=a,U,L]p(a|U,L)\\
			-\mathbb{E}[\pi(Z,L)\mid L]p(a|U,L)
		\end{array}\right\}\mid L]\\
		=&\mathbb{E}[\mathbb{E}[Y(a)\mid U,L]\dfrac{p(a|U,L)}{p(a|L)}\left\{\begin{array}{l}
			\mathbb{E}[\pi(Z,L)\mid A=a,U,L]
			-\mathbb{E}[\pi(Z,L)\mid L]
		\end{array}\right\}\mid L].
	\end{align*}
	These equalities just repeatedly adopt the Bayes formula and the fact that $Z\indep U\mid L$ from Assumption~\ref{as: IV independence}.
	Therefore,
	\begin{align*}
		&\dfrac{\mathbb{E}[Y\left\{\pi(Z,L)-\mathbb{E}[\pi(Z,L)\mid L]\right\}\mid A=a,L]}
		{\mathbb{E}[\pi(Z,L)\mid A=a,L] - \mathbb{E}[\pi(Z,L)\mid L]}\\
		=&\mathbb{E}[\mathbb{E}[Y(a)\mid U,L]\dfrac{p(a|U,L)}{p(a|L)}\dfrac{\left\{
			\mathbb{E}[\pi(Z,L)\mid A=a,U,L]
			-\mathbb{E}[\pi(Z,L)\mid L]\right\}}
		{\mathbb{E}[\pi(Z,L)\mid A=a,L] - \mathbb{E}[\pi(Z,L)\mid L]}\mid L]\\
		=&\mathbb{E}[\mathbb{E}[Y(a)\mid U,L]\omega_{a,\pi}(U,L)\mid L].
	\end{align*}
	The final step follows the definition of $\omega_{a,\pi}(U,L)$. Finally, we can integrate $L$ out to deduce that
	\begin{align*}
		\mathbb{E}\left[\dfrac{\mathbb{E}[Y\left\{\pi(Z,L)-\mathbb{E}[\pi(Z,L)\mid L]\right\}\mid A=a,L]}
		{\mathbb{E}[\pi(Z,L)\mid A=a,L] - \mathbb{E}[\pi(Z,L)\mid L]}\right]
		=\mathbb{E}[\mathbb{E}[Y(a)\mid U,L]\omega_{a,\pi}(U,L)].
	\end{align*}
	The left side of the Equation above equals $\mathbb{E}[\mu_{\pi}^o(a,L)]$, thus finishing the proof of Equation~\ref{eq: identification E[Y(a)]}.
	\begin{align*}
		&\mathbb{E}[\mu_{\pi}^o(a,L)] - \mathbb{E}[Y(a)] =
		\mathbb{E}[\mathbb{E}[Y(a)\mid U,L]\omega_{a,\pi}(U,L)]-
		\mathbb{E}[\mathbb{E}[Y(a)\mid U,L]]\\
		=&\mathbb{E}[\mathbb{E}[Y(a)\mid U,L]\{\omega_{a,\pi}(U,L)-1\}]\\
		=&\mathbb{E}[\mathbb{E}[Y(a)\mid U,L]\{\omega_{a,\pi}(U,L)-\mathbb{E}[\omega_{a,\pi}(U,L)|L]\}]\\
		=&\mathbb{E}\left[\mathrm{Cov}\{\mathbb{E}\{Y(a)\mid U,L\},\,\omega_{a,\pi}(U,L)\mid L\}\right],
	\end{align*}
	where the third equality follows from $\mathbb{E}[\omega_{a,\pi}(U,L)|L]\equiv1$ from Proposition~\ref{prop: AIV characterization}. This finishes the proof of Equation~\eqref{eq: identification E[Y(a)] bias}.
	If $Z$ is an AIV for $A=a$, then $\omega_{a,\pi}(U,L)\equiv 1$ by Proposition~\ref{prop: AIV characterization}. Thus,
	\begin{align*}
		\mathbb{E}\left[\mu_{\pi}^o(a,L)\right]-\mathbb{E}[Y(a)]=\mathbb{E}\left[\mathrm{Cov}\{\mathbb{E}\{Y(a)\mid U,L\},\,1\mid L\}\right]\equiv 0.
	\end{align*}
	This finishes the proof of Theorem~\ref{thm: identification}.

	\subsection{Proof of Theorem~\ref{thm: psi_q}}
	Notably, for $Z_{\pi}:=\pi(Z,L)$,
	\begin{align*}
		\mu_{\pi}^o(A,L) - \mathbb{E}[Y\mid A,L]
		=\dfrac{\mathbb{E}[YZ_{\pi} \mid A,L] - \mathbb{E}[Y \mid A,L]\mathbb{E}[Z_{\pi} \mid A,L]}
		{\mathbb{E}[Z_{\pi} \mid A,L] - \mathbb{E}[Z_{\pi} \mid L]}. 
	\end{align*}
	For any pathwise differentiable parameterization $p_\theta(O)$, we denote $s(O):=\nabla_{\theta}\log p_\theta(O)\eval_{\theta=0}$ as the score function, $\mathbb{E}_\theta$ as taking expectation with respect to $p_\theta(O)$.
	We calculate the EIF for $\psi$ as follows.
	\begin{align*}
		&\nabla_{\theta}\int \int \mu_{\pi,\theta}(a, l) q(a)\mathrm{d} P_L(l;\theta) \mathrm{d} P_A(a;\theta)\eval_{\theta=0}\\
		=&\nabla_{\theta} \int \int
		\dfrac{\left(\begin{array}{c}
				\mathbb{E}_\theta[YZ_{\pi} \mid A=a,L=l]\\
				-\mathbb{E}_\theta[Y \mid A=a,L=l]\mathbb{E}_\theta[Z_{\pi} \mid L=l]
			\end{array}\right)}
		{\mathbb{E}_\theta[Z_{\pi} \mid A=a,L=l] - \mathbb{E}_\theta[Z_{\pi} \mid L=l]}
		q(a)\mathrm{d} P_L(l;\theta) \mathrm{d} P_A(a;\theta)\eval_{\theta=0}\\
		=&\int \int
		\nabla_{\theta}\left(
		\dfrac{\left\{\begin{array}{c}
				\mathbb{E}_\theta[YZ_{\pi} \mid A=a,L=l]\\
				- \mathbb{E}_\theta[Y \mid A=a,L=l]\mathbb{E}_\theta[Z_{\pi} \mid L=l]
			\end{array}\right\}}
		{\mathbb{E}_\theta[Z_{\pi} \mid A=a,L=l] - \mathbb{E}_\theta[Z_{\pi} \mid L=l]}\right)\eval_{\theta=0}
		q(a)\mathrm{d} P_L(l;\theta) \mathrm{d} P_A(a;\theta)\\
		&+\nabla_{\theta} \int \int \mu_{\pi}(a,l)
		\mathrm{d} P_L(l;\theta) q(a)\mathrm{d} P_A(a;\theta)\eval_{\theta=0}\\
		=&\int \int
		\dfrac{\mathbb{E}[(YZ_{\pi}-\mathbb{E}[YZ_{\pi}\mid A,L])s(O)\mid A=a,L=l]}
		{\mathbb{E}[Z_{\pi} \mid A=a,L=l] - \mathbb{E}[Z_{\pi} \mid L=l]}
		q(a)\mathrm{d} P_L(l) \mathrm{d} P_A(a)\\
		&-\int \int
		\dfrac{\mathbb{E}[(Y - \mathbb{E}[Y\mid A,L])s(O)\mid A=a,L=l]}
		{\mathbb{E}[Z_{\pi} \mid A=a,L=l] - \mathbb{E}[Z_{\pi} \mid L=l]} \mathbb{E}[Z_{\pi}\mid L=l]
		q(a)\mathrm{d} P_L(l) \mathrm{d} P_A(a)\\
		&-\int \int
		\dfrac{\mathbb{E}[(Z_{\pi}- \mathbb{E}[Z_{\pi}\mid L])s(O)\mid L=l]}
		{\mathbb{E}[Z_{\pi} \mid A=a,L=l] - \mathbb{E}[Z_{\pi} \mid L=l]} \mathbb{E}[Y\mid A=a,L=l]
		q(a)\mathrm{d} P_L(l) \mathrm{d} P_A(a)\\
		&-\int \int \mu_{\pi}(a,l)
		\dfrac{\left\{\begin{array}{l}
				\mathbb{E}[(Z_{\pi}-\mathbb{E}[Z_{\pi}\mid A,L])s(O) \mid A=a,L=l]\\
				- \mathbb{E}[(Z_{\pi}- \mathbb{E}[Z_{\pi}\mid L])s(O) \mid L=l]
			\end{array}\right\}}
		{\mathbb{E}[Z_{\pi} \mid A=a,L=l] - \mathbb{E}[Z_{\pi} \mid L=l]}	
		q(a)\mathrm{d} P_L(l) \mathrm{d} P_A(a)\\
		&+\nabla_{\theta} \mathbb{E}_\theta \left[q(A)\int \mu_{\pi}(A,l)
		\mathrm{d} P_L(l) \right]\eval_{\theta=0}
		+\nabla_{\theta} \mathbb{E} \left[\int \mu_{\pi}(a,L) q(a)\mathrm{d} P_A(a)\right]\eval_{\theta=0}\\
		=&\int \int
		\mathbb{E}\left[
		\dfrac{(YZ_{\pi}-\mathbb{E}[YZ_{\pi}\mid A,L])}
		{\mathbb{E}[Z_{\pi} \mid A,L] - \mathbb{E}[Z_{\pi} \mid L]}s(O)\middle| A=a,L=l\right]
		q(a)\mathrm{d} P_L(l) \mathrm{d} P_A(a)\\
		&-\int \int
		\mathbb{E}\left[\dfrac{(Y - \mathbb{E}[Y\mid A,L])}
		{\mathbb{E}[Z_{\pi} \mid A,L] - \mathbb{E}[Z_{\pi} \mid L]}\mathbb{E}[Z_{\pi}\mid L]s(O)
		\middle| A=a,L=l\right]
		q(a)\mathrm{d} P_L(l) \mathrm{d} P_A(a)\\
		&-\int \int
		\mathbb{E}\left[
		\dfrac{\mathbb{E}[(Z_{\pi}- \mathbb{E}[Z_{\pi}\mid L])s(O)\mid L]}
		{\mathbb{E}[Z_{\pi} \mid A,L] - \mathbb{E}[Z_{\pi} \mid L]}
		Y\middle| A=a,L=l\right]
		q(a)\mathrm{d} P_L(l) \mathrm{d} P_A(a)\\
		&-\int \int 
		\mathbb{E}\left[\mu_{\pi}(A,L)\dfrac{(Z_{\pi}-\mathbb{E}[Z_{\pi}\mid A,L])s(O) }
		{\mathbb{E}[Z_{\pi} \mid A,L] - \mathbb{E}[Z_{\pi} \mid L]}\middle| A=a,L=l\right]	
		q(a)\mathrm{d} P_L(l) \mathrm{d} P_A(a)\\
		&+\int \int \mu_{\pi}(a,l)
		\dfrac{\mathbb{E}[(Z_{\pi}- \mathbb{E}[Z_{\pi}\mid L])s(O) \mid L=l]}
		{\mathbb{E}[Z_{\pi} \mid A=a,L=l] - \mathbb{E}[Z_{\pi} \mid L=l]}	
		q(a)\mathrm{d} P_L(l) \mathrm{d} P_A(a)\\
		&+\mathbb{E} \left[\left(q(A)\int \mu_{\pi}(A,l)\mathrm{d} P_L(l) - \psi_{\pi,q}^o\right)s(O) \right]\\
		&+\mathbb{E} \left[\left(\int \mu_{\pi}(a,L) q(a)\mathrm{d} P_A(a) - \psi_{\pi,q}^o\right)s(O)\right]\\
		=&\mathbb{E}\left[
		\dfrac{p_A(A)q(A)}{p_{A|L}(A\mid L)}
		\dfrac{(YZ_{\pi}-\mathbb{E}[YZ_{\pi}\mid A,L])}
		{\mathbb{E}[Z_{\pi} \mid A,L] - \mathbb{E}[Z_{\pi} \mid L]}s(O)\right]\\
		&-\mathbb{E}\left[
		\dfrac{p_A(A)q(A)}{p_{A|L}(A\mid L)}
		\dfrac{(Y - \mathbb{E}[Y\mid A,L])}
		{\mathbb{E}[Z_{\pi} \mid A,L] - \mathbb{E}[Z_{\pi} \mid L]}\mathbb{E}[Z_{\pi}\mid L]s(O)\right]\\
		&-\mathbb{E}\left[
		\dfrac{p_A(A)q(A)}{p_{A|L}(A\mid L)}Y
		\dfrac{\mathbb{E}[(Z_{\pi}- \mathbb{E}[Z_{\pi}\mid L])s(O)\mid L]}
		{\mathbb{E}[Z_{\pi} \mid A,L] - \mathbb{E}[Z_{\pi} \mid L]}
		\right]\\
		&-\mathbb{E}\left[
		\dfrac{p_A(A)q(A)}{p_{A|L}(A\mid L)}   
		\mu_{\pi}(A,L)\dfrac{(Z_{\pi}-\mathbb{E}[Z_{\pi}\mid A,L])s(O) }
		{\mathbb{E}[Z_{\pi} \mid A,L] - \mathbb{E}[Z_{\pi} \mid L]}\right]\\
		&+\mathbb{E}\left[
		\dfrac{p_A(A)q(A)}{p_{A|L}(A\mid L)}\mu_{\pi}(A,L)
		\dfrac{\mathbb{E}[(Z_{\pi}- \mathbb{E}[Z_{\pi}\mid L])s(O) \mid L]}
		{\mathbb{E}[Z_{\pi} \mid A,L] - \mathbb{E}[Z_{\pi} \mid L]}
		\right]\\
		&+\mathbb{E} \left[\left(q(A)\int \mu_{\pi}(A,l)\mathrm{d} P_L(l) - \psi_{\pi,q}^o\right)s(O) \right]\\
		&+\mathbb{E} \left[\left(\int \mu_{\pi}(a,L) q(a)\mathrm{d} P_A(a) - \psi_{\pi,q}^o\right)s(O)\right]\\
		=&\mathbb{E}\left[
		\dfrac{p_A(A)q(A)}{p_{A|L}(A\mid L)}
		\left\{\dfrac{Y (Z_{\pi}- \mathbb{E}[Z_{\pi}\mid L])}
		{\mathbb{E}[Z_{\pi} \mid A,L] - \mathbb{E}[Z_{\pi} \mid L]} - \mu_{\pi}(A,L)\right\} s(O)\right]\\
		&-\mathbb{E}\left[
		\dfrac{p_A(A)q(A)}{p_{A|L}(A\mid L)}   
		\mu_{\pi}(A,L)\dfrac{(Z_{\pi}-\mathbb{E}[Z_{\pi}\mid A,L])s(O) }
		{\mathbb{E}[Z_{\pi} \mid A,L] - \mathbb{E}[Z_{\pi} \mid L]}\right]\\
		&+\mathbb{E}\left[
		\dfrac{p_A(A)q(A)}{p_{A|L}(A\mid L)}
		\dfrac{\mathrm{Cov}\{Y,Z_{\pi}\mid A,L\}}
		{(\mathbb{E}[Z_{\pi} \mid A,L] - \mathbb{E}[Z_{\pi} \mid L])^2}\mathbb{E}[(Z_{\pi}- \mathbb{E}[Z_{\pi}\mid L])s(O) \mid L]
		\right]\\
		&+\mathbb{E} \left[\left(q(A)\int \mu_{\pi}(A,l)\mathrm{d} P_L(l) - \psi_{\pi,q}^o\right)s(O) \right]\\
		&+\mathbb{E} \left[\left(\int \mu_{\pi}(a,L) q(a)\mathrm{d} P_A(a) - \psi_{\pi,q}^o\right)s(O)\right]\\
		=&\mathbb{E}\left[
		\dfrac{p_A(A)q(A)}{p_{A|L}(A\mid L)}
		\dfrac{ (Y-\mu_{\pi}(A,L)) (Z_{\pi}- \mathbb{E}[Z_{\pi}\mid L])}
		{\mathbb{E}[Z_{\pi} \mid A,L] - \mathbb{E}[Z_{\pi} \mid L]} s(O)\right]\\
		&+\mathbb{E}\left[\mathbb{E}\left[
		\dfrac{p_A(A)q(A)}{p_{A|L}(A\mid L)}
		\dfrac{\mathrm{Cov}\{Y,Z_{\pi}\mid A,L\}}
		{(\mathbb{E}[Z_{\pi} \mid A,L] - \mathbb{E}[Z_{\pi} \mid L])^2}\middle| L\right]
		(Z_{\pi}- \mathbb{E}[Z_{\pi}\mid L])s(O) 
		\right]\\
		&+\mathbb{E} \left[\left(q(A)\int \mu_{\pi}(A,l)\mathrm{d} P_L(l) - \psi_{\pi,q}^o\right)s(O) \right]\\
		&+\mathbb{E} \left[\left(\int \mu_{\pi}(a,L) q(a)\mathrm{d} P_A(a) - \psi_{\pi,q}^o\right)s(O)\right].
	\end{align*}
	Then the influence function for $\psi_{\pi,q}^o$ is
	\begin{align*}
		&\dfrac{p_{A}(A)}{p_{A|L}(A\mid L)}
		\left\{\begin{array}{l}
			\dfrac{ (Y-\mu_{\pi}(A,L)) (Z_{\pi}- \mathbb{E}[Z_{\pi}\mid L])}
			{\mathbb{E}[Z_{\pi} \mid A,L] - \mathbb{E}[Z_{\pi} \mid L]}q(A)\\
			-\displaystyle\int \dfrac{\mathbb{E}[Y\mid A=a,L] - \mu_{\pi}(a,L)}
			{\mathbb{E}[Z_{\pi} \mid A=a,L] - \mathbb{E}[Z_{\pi} \mid L]} q(a)\mathrm{d} P_A(a)\\
			\times(Z_{\pi}- \mathbb{E}[Z_{\pi}\mid L])
		\end{array}\right\}\\
		&+q(A)\int \mu_{\pi}(A,l)\mathrm{d} P_L(l) - \psi_{\pi,q}^o
		+\int \mu_{\pi}(a,L) q(a)\mathrm{d} P_A(a) - \psi_{\pi,q}^o.
	\end{align*}
	Since the model is nonparametric, the influence function for $\psi_{\pi,q}^o$ is just the EIF for $\psi_{\pi,q}^o$, finishing the proof of Theorem~\ref{thm: psi_q}.

	\subsection{Proof of Theorem~\ref{thm: convergence rate}}
	Adopting Lemma~\ref{lem: LLKR}, we can decompose $\hat{\theta}_{\pi,h}^{(n)}(a)-\theta(a)$ as
	\begin{align*}
		&\hat{\theta}_{\pi,h}^{(n)}(a)-\theta(a)\\
		=&\sum_{k=1}^K\sum_{i\in I_k}w_{ni}(a,h)\left\{\varphi_{\pi}(O_i;\hat\alpha_{\pi}^{(n,-k)},\hat{P}_O^{(n,-k)})-\varphi_{\pi}(O_i;\hat\alpha_{\pi}^{(n,-k)}, P_O)\right\}\\
		&+\sum_{k=1}^K\sum_{i\in I_k}w_{ni}(a,h)\left\{\varphi_{\pi}(O_i;\hat\alpha_{\pi}^{(n,-k)}, P_O)-\varphi_{\pi}(O_i;\alpha_{\pi}^o, P_O)\right\}\\
		&+\sum_{k=1}^K\sum_{i\in I_k}w_{ni}(a,h)\left\{\varphi_{\pi}(O_i;\alpha_{\pi}^o, P_O)-\theta(a)\right\}\\
		=&e_1^{\top} Q_{h}(a)\times\left\{\begin{array}{l}
			\displaystyle\sum_{k=1}^K\sum_{i\in I_k}g\left(\dfrac{A_i-a}{h}\right)
			\left\{\varphi_{\pi}(O_i;\hat\alpha_{\pi}^{(n,-k)},\hat{P}_O^{(n,-k)})-\varphi_{\pi}(O_i;\hat\alpha_{\pi}^{(n,-k)}, P_O)\right\}\\
			+\displaystyle\sum_{k=1}^K\sum_{i\in I_k}g\left(\dfrac{A_i-a}{h}\right)
			\left\{\varphi_{\pi}(O_i;\hat\alpha_{\pi}^{(n,-k)}, P_O)-\varphi_{\pi}(O_i;\alpha_{\pi}^o, P_O)\right\}\\
			+\displaystyle\sum_{k=1}^K\sum_{i\in I_k}g\left(\dfrac{A_i-a}{h}\right)
			\left\{\varphi_{\pi}(O_i;\alpha_{\pi}^o, P_O)-\theta(a)\right\}\\
		\end{array}\right\}\\
		=&e_1^{\top} Q_{h}(a)\sum_{k=1}^K\sum_{i\in I_k}\left[\begin{array}{l}
			K_h(A_i-a),\\
			K_h(A_i-a)\dfrac{A_i-a}{h}
		\end{array}\right]\\&\times\displaystyle\left\{\begin{array}{l}
			\left\{\varphi_{\pi}(O_i;\hat\alpha_{\pi}^{(n,-k)},\hat{P}_O^{(n,-k)})-\varphi_{\pi}(O_i;\hat\alpha_{\pi}^{(n,-k)}, P_O)\right\}\\
			+\left\{\varphi_{\pi}(O_i;\hat\alpha_{\pi}^{(n,-k)}, P_O)-\varphi_{\pi}(O_i;\alpha_{\pi}^o, P_O)\right\}\\
			+\left\{\varphi_{\pi}(O_i;\alpha_{\pi}^o, P_O)-\mathbb{E}[\varphi_{\pi}(O_i;\alpha_{\pi}^o, P_O)\mid A_i]\right\}\\
			+\left\{\theta(A_i)-\theta(a)-\theta'(a)(A_i-a)\right\}\\
		\end{array}\right\}
	\end{align*}
	First, by Lemma~\ref{lem: lower bounded eigenvalue} and Assumptions~\ref{as: regular condition}\ref{cd: bandwidth rate},~\ref{cd: kernel condition} in Assumption~\ref{as: regular condition}, $e_1^{\top}  Q_{h}(a)=O_p(1/n)$.
	Thus, we need only to control the following quantity for $j=0,1$ and $k=1,\ldots, K$:
	\begin{align}
		R_{1kj}:=&\sum_{i\in I_k} K_h(A_i-a)\left(\dfrac{A_i-a}{h}\right)^j
		\left\{\varphi_{\pi}(O_i;\hat\alpha_{\pi}^{(n,-k)},\hat{P}_O^{(n,-k)})-\varphi_{\pi}(O_i;\hat\alpha_{\pi}^{(n,-k)}, P_O)\right\}\notag\\
		=& o_p(\sqrt{n/h}),\label{exp: 1}\\
		R_{2kj}:=&\sum_{i\in I_k} K_h(A_i-a)\left(\dfrac{A_i-a}{h}\right)^j
		\left\{\varphi_{\pi}(O_i;\hat\alpha_{\pi}^{(n,-k)}, P_O)-\varphi_{\pi}(O_i;\alpha_{\pi}^o, P_O)\right\}\notag\\ 
		=& o_p(\sqrt{n/h}) + o_p(nc(n)),\label{exp: 2}\\
		R_{3kj}:=&\sum_{i\in I_k} K_h(A_i-a)\left(\dfrac{A_i-a}{h}\right)^j
		\left\{\varphi_{\pi}(O_i;\alpha_{\pi}^o, P_O)-\mathbb{E}[\varphi_{\pi}(O_i;\alpha_{\pi}^o, P_O)\mid A_i]\right\}\notag\\
		=& O_p(\sqrt{n/h}),\label{exp: 3}\\
		R_{4kj}:=&\sum_{i\in I_k} K_h(A_i-a)\left(\dfrac{A_i-a}{h}\right)^j
		\left\{\theta(A_i)-\theta(a)-\theta'(a)(A_i-a)\right\}\notag\\
		=& O_p(nh^2).\label{exp: 4}
	\end{align}
	The final step is to combine them as
	\begin{align*}
		|\hat{\theta}_{\pi,h}^{(n)}(a)-\theta(a)|=&
		\left|e_1^{\top} Q_{h}(a)\sum_{k=1}^K\left[\begin{array}{l}
			R_{1k0} +R_{2k0} + R_{3k0} + R_{4k0}\\
			R_{1k1} +R_{2k1} + R_{3k1} + R_{4k1}
		\end{array}\right]\right|\\
		\lesssim&\dfrac{1}{n}\sum_{k=1}^K\sum_{j=0,1}\{|R_{1kj}| +|R_{2kj}| + |R_{3kj}| + R_{4kj}|\}\\
		=&O_p(\dfrac{1}{\sqrt{nh}}) + O_p(h^2) +  o_p(c(n)).
	\end{align*}
	This finishes the proof of Equation~\eqref{eq: convergence rate}.
	Next, we show the proof of Equation~\eqref{exp: 1}--\eqref{exp: 4}, respectively.
	
	For clarity, we require the sample independence between $I_{-k}^{(2)}$ and $I_{-k}^{(1)}\cup I_k$ in the proof of Equation~\eqref{exp: 1}, and the independence between $I_{-k}^{(1)}$ and $I_k$ in the proof of Equation~\eqref{exp: 2}. 
	The mixed-bias property (Lemma~\ref{lem: mixed bias}) is only invoked in the proof of Equation~\eqref{exp: 2}. 
	For the proofs of Equations~\eqref{exp: 3} and~\eqref{exp: 4}, one can apply the standard techniques developed in the LLKR literature \citep{fan1996local, tsybakov2009introduction}.

	\begin{proof}[Proof of Equation~\eqref{exp: 1}:]
		First, we control $R_{1kj}$. Recall that
		\begin{align*}
			&\varphi_{\pi}(O;\alpha_{\pi}, \hat{P}_O^{(n,-k)}) - \varphi_{\pi}(O;\alpha_{\pi}, P_O) \\
			=& \int \mu_{\pi}(A,l) - (z_{\pi}-\rho_{\pi}(l))\dfrac{\eta(A,l)-\mu_{\pi}(A,l)}{\kappa_{\pi}(A,l)}\mathrm{d} (\hat{P}_O^{(n,-k)}(o) - P_O(o)).
		\end{align*}
		It is easy to observe that for any $\alpha_{\pi}$ trained from $I_{-k}^{(1)}$,
		\begin{align*}
			&\mathbb{E}_{-k}^{(2)}\left[\varphi_{\pi}(O;\alpha_{\pi}, \hat{P}_O^{(n,-k)}) - \varphi_{\pi}(O;\alpha_{\pi}, P_O)\right]\\
			=&\mathbb{E}_{-k}^{(2)}\left[\int\mu_{\pi}(A,l)
			-(z_{\pi}-\rho_{\pi}(l))\dfrac{\eta(A,l)-\mu_{\pi}(A,l)}{\kappa_{\pi}(A,l)}
			\mathrm{d} (\hat{P}_O^{(n,-k)}(o) - P_O(o))\right]\\
			=&\mathbb{E}_{-k}^{(2)}\left[\dfrac{1}{|I_{-k}^{(2)}|}\sum_{i\in I_{-k}^{(2)}} 
			\left\{\mu_{\pi}(A,L_i)- (Z_{\pi,i}-\rho_{\pi}(L_i))\dfrac{\eta(A,L_i)-\mu_{\pi}(A,L_i)}{\kappa_{\pi}(A,L_i)}\right\}\right]\\
			&-\int \mu_{\pi}(A,l) - (z_{\pi}-\rho_{\pi}(l))\dfrac{\eta(A,l)-\mu_{\pi}(A,l)}{\kappa_{\pi}(A,l)}\mathrm{d} P_O(o)=0
		\end{align*}
		This is because in Algorithm~\ref{alg: cross-fitting}, samples in $I_{-k}^{(2)}$ are independent of those in $I_{k}$ and $I_{-k}^{(1)}$.
		Thus, we know that $\mathbb{E}_{-k}^{(2)}[R_{1kj}]=0$, and
		\begin{align}
			&\mathbb{E}_k[\mathbb{E}_{-k}^{(2)}[R_{1kj}^2]]\notag\\
			=&\mathbb{E}_k\left\{\mathbb{E}_{-k}^{(2)}\left[
			\left(\begin{array}{c}
				\displaystyle\sum_{i\in I_k} K_h(A_i-a)\left(\dfrac{A_i-a}{h}\right)^j\\
				\times\left\{\varphi_{\pi}(O_i;\hat\alpha_{\pi}^{(n,-k)},\hat{P}_O^{(n,-k)})-\varphi_{\pi}(O_i;\hat\alpha_{\pi}^{(n,-k)}, P_O)\right\}
			\end{array}\right)^2\right]\right\}\notag\\
			=&\mathbb{E}_k\left\{
			\sum_{i\in I_k}\left(\begin{array}{c}
				K_h(A_i-a)^2\left(\dfrac{A_i-a}{h}\right)^{2j}\\
				\times\mathbb{E}_{-k}^{(2)}\left[
				\left\{\varphi_{\pi}(O_i;\hat\alpha_{\pi}^{(n,-k)},\hat{P}_O^{(n,-k)})-\varphi_{\pi}(O_i;\hat\alpha_{\pi}^{(n,-k)}, P_O)\right\}^2\right]
			\end{array}\right)
			\right\}\label{exp: 10}\\
			\lesssim&\dfrac{1}{h^2} \sum_{i\in I_k}\mathbb{E}_k\left\{
			I\{|A_i-a|\leq h\}\times\mathbb{E}_{-k}^{(2)}\left[ 
			\left\{\begin{array}{c}
				\varphi_{\pi}(O_i;\hat\alpha_{\pi}^{(n,-k)},\hat{P}_O^{(n,-k)})\\
				-\varphi_{\pi}(O_i;\hat\alpha_{\pi}^{(n,-k)}, P_O)
			\end{array}\right\}^2\right]
			\right\}\label{exp: 11}\\
			\lesssim&\dfrac{1}{h^2}\sum_{i\in I_k}\mathbb{E}_k\Bigg\{I\{|A_i-a|\leq h\}\label{exp: 12}\\
			&\times \mathbb{E}_{-k}^{(2)}\left[ 
			\left\{\int \left(\begin{array}{l}
				\hat\mu_{\pi}^{(n,-k)}(A_i,l)- (z_{\pi}-\hat\rho_{\pi}^{(n,-k)}(l))\\
				\times\dfrac{\hat\eta^{(n,-k)}(A_i,l)-\hat\mu_{\pi}^{(n,-k)}(A_i,l)}
				{\hat\kappa_{\pi}^{(n,-k)}(A_i,l)}
			\end{array}\right)
			\mathrm{d} (\hat{P}_O^{(n,-k)}(o) - P_O(o))\right\}^2\right]\Bigg\}\notag\\
			\lesssim&\dfrac{1}{h^2|I_{-k}^{(2)}|}\sum_{i\in I_k}\mathbb{E}_k\left\{I\{|A_i-a|\leq h\}
			\int \left(\begin{array}{l}
				\hat\mu_{\pi}^{(n,-k)}(A_i,l)- (z_{\pi}-\hat\rho_{\pi}^{(n,-k)}(l))\\
				\times\dfrac{\hat\eta^{(n,-k)}(A_i,l)-\hat\mu_{\pi}^{(n,-k)}(A_i,l)}
				{\hat\kappa_{\pi}^{(n,-k)}(A_i,l)}
			\end{array}\right)^2\mathrm{d} P_O(o)\right\}\label{exp: 13}\\
			\lesssim&\dfrac{1}{h^2|I_{-k}^{(2)}|}\sum_{i\in I_k}\mathbb{E}_k\left[I\{|A_i-a|\leq h\}\right].
			\label{exp: 14}
		\end{align}
		Inequalities in Equations~\eqref{exp: 10}~and~\eqref{exp: 13} use the fact that the samples in $I_{-k}^{(2)}$ are independent from the samples in $I_k$ and $I_{-k}^{(1)}$.
		Equation~\eqref{exp: 11} follows from $K_h(A_i-a)=K((A_i-a)/h)/h\lesssim I\{|A_i-a|\leq h\}/h$ by Assumption~\ref{as: regular condition}\ref{cd: kernel condition}.
		Equation~\eqref{exp: 12} follows from the basic inequality that $(a+b)^2\leq 2(a^2+b^2)$, and Assumption~\ref{as: regular condition}\ref{cd: alpha uniform boundedness} that $|\hat\kappa_{\pi}^{(n,-k)}(A,L)|$ are lower bounded away from zero uniformly for all $n$, and $A,L\in\mathcal{N}\times \mathcal{L}$.
		Equation~\eqref{exp: 14} follows from the uniform boundedness of $\hat\mu_{\pi}^{(n,-k)}(A,l)$ and $\hat\eta^{(n,-k)}(A,l)$ in $\mathcal{N}\times \mathcal{L}$ from Assumptions~\ref{as: regular condition}\ref{cd: ZY boundedness}~and~\ref{cd: alpha uniform boundedness}.
		Finally, we observe that
		\begin{align*}
			&\mathbb{E}\left[R_{1kj}^2\right]
			\lesssim\mathbb{E}\left[\dfrac{1}{h^2|I_{-k}^{(2)}|}\sum_{i\in I_k}I\{|A_i-a|\leq h\}\right]\\
			\lesssim& \dfrac{n}{h^2|I_{-k}^{(2)}|}\Pr\{|A-a|\leq h\}
			\lesssim\dfrac{n}{h|I_{-k}^{(2)}|}.
		\end{align*}
		The final inequality holds since $p_A(a):=\lim\limits_{h\to 0^+}\Pr\{|A-a|\leq h\}/(2h)>0$ at $a\in\mathring{\mathcal{A}}$.
		Thus, we know that $R_{1kj}^2=O_p(\frac{n}{h|I_{-k}^{(2)}|}) =o_p(\frac{n}{h})$ from Assumption~\ref{as: regular condition}\ref{cd: plnk rate}.
		
	\end{proof}
	
	\begin{proof}[Proof of Equation~\eqref{exp: 2}:]
		Next, we control $R_{2kj}$. First, using Lemma~\ref{lem: mixed bias property}, we have $\mathbb{E}_k[R_{2kj}]=o(nc(n))$ by Assumption~\ref{as: regular condition}\ref{cd: alpha rate}, it remains only to control $R_{2kj} - \mathbb{E}_k[R_{2kj}]$. Concretely,
		\begin{align}
			&\mathbb{E}_k[(R_{2kj} - \mathbb{E}_k[R_{2kj}])^2] = \mathrm{Var}_k[R_{2kj}]\notag\\
			= &\sum_{i\in I_k}\mathrm{Var}_k\left[K_h(A-a)\left(\dfrac{A-a}{h}\right)^j
			\left\{\varphi_{\pi}(O;\hat\alpha_{\pi}^{(n,-k)},P_O)-\varphi_{\pi}(O;\alpha_{\pi}^o, P_O)\right\}\right]\label{exp: 15}\\
			\lesssim &\sum_{i\in I_k}\mathbb{E}_k\left[K_h(A-a)^2\left(\dfrac{A-a}{h}\right)^{2j}
			\left\{\varphi_{\pi}(O;\hat\alpha_{\pi}^{(n,-k)},P_O)-\varphi_{\pi}(O;\alpha_{\pi}^o, P_O)\right\}^2\right]\notag\\
			\lesssim &\dfrac{1}{h^2}\sum_{i\in I_k}\mathbb{E}_k\left[I(|A-a|\leq h)
			\mathbb{E}\left[\left\{\varphi_{\pi}(O;\hat\alpha_{\pi}^{(n,-k)},P_O)-\varphi_{\pi}(O;\alpha_{\pi}^o, P_O)\right\}^2\middle| A\right]\right]\label{exp: 16}\\
			\lesssim &\dfrac{1}{h^2}\sum_{i\in I_k}\mathbb{E}_k\left[I(|A-a|\leq h)\right]
			\sup_{s\in\mathcal{B}(a,h)}\mathbb{E}\left[\left\{\varphi_{\pi}(O;\hat\alpha_{\pi}^{(n,-k)},P_O)-\varphi_{\pi}(O;\alpha_{\pi}^o, P_O)\right\}^2\middle| A=s\right]\notag\\
			\lesssim &\dfrac{n}{h}
			\sup_{s\in\mathcal{B}(a,h)}\mathbb{E}\left[\left\{\varphi_{\pi}(O;\hat\alpha_{\pi}^{(n,-k)},P_O)-\varphi_{\pi}(O;\alpha_{\pi}^o, P_O)\right\}^2\middle| A=s\right]\label{exp: 17}\\
			\lesssim &\dfrac{n}{h}
			\left\{\begin{array}{l}
				\left\|\hat\rho_{\pi}^{(n,-k)}-\rho_{\pi}^o\right\|_2^2+
				\left\|\hat\eta^{(n,-k)}-\eta^o\right\|_{\mathcal{B}(a,h),2}^2
				+\left\|\hat\mu_{\pi}^{(n,-k)}-\mu_{\pi}^o\right\|_{\mathcal{B}(a,h),2}^2\\+
				\left\|\hat\kappa_{\pi}^{(n,-k)}-\kappa_{\pi}^o\right\|_{\mathcal{B}(a,h),2}^2
				+\left\|\hat\delta^{(n,-k)}-\delta^o\right\|_{\mathcal{B}(a,h),2}^2\\
			\end{array}\right\}\notag\\
			\lesssim &\dfrac{n}{h} \left\|\hat\alpha_{\pi}^{(n,-k)}-\alpha_{\pi}^o\right\|_{\mathcal{B}(a,h),2}^2.\notag
		\end{align}
		Equation~\eqref{exp: 15} follows from the independent structure of the samples between $I_k$ and $I_{-k}$.
		Equation~\eqref{exp: 16} follows from $K_h(A_i-a)=K((A_i-a)/h)/h\lesssim I\{|A_i-a|\leq h\}/h$ by Assumption~\ref{as: regular condition}\ref{cd: kernel condition}.
		Equation~\eqref{exp: 17} follows from the fact that $p_A(a):=\lim_{h\to 0}\Pr\{|A-a|\leq h\}/(2h)>0$ at $a\in\mathring{\mathcal{A}}$.
		The final two inequalities can be seen from the definition of $\varphi$ in Equation~\eqref{eq: varphi defn} and Assumptions~\ref{as: regular condition}\ref{cd: ZY boundedness}~and~\ref{cd: alpha uniform boundedness}.

		Thus, from Assumption~\ref{as: regular condition}\ref{cd: alpha rate}, we know that 
		\begin{align*}
			\mathbb{E}\left[(R_{2kj} - \mathbb{E}_k[R_{2kj}])^2\right]
			\lesssim\dfrac{n}{h} \mathbb{E}\left[\left\|\hat\alpha_{\pi}^{(n,-k)}-\alpha_{\pi}^o\right\|_{\mathcal{B}(a,h),2}^2\right]=o(\dfrac{n}{h}).
		\end{align*}
		This is equivalent to saying $R_{2kj}=\left(R_{2kj} - \mathbb{E}_k[R_{2kj}]\right) +\mathbb{E}_k[R_{2kj}]= o_p(\sqrt{n/h}) + o_p(nc(n))$.	
	\end{proof}
	
	\begin{proof}[Proof of Equation~\eqref{exp: 3}:]
		We control $R_{3kj}$ by
		\begin{align}
			&\mathbb{E}[(R_{3kj})^2\mid A_i, i\in I_k]\notag\\
			= &\mathbb{E}\left[\left(\sum_{i\in I_k}K_h(A_i-a)\left(\dfrac{A_i-a}{h}\right)^j
			\left\{\varphi_{\pi}(O_i;\alpha_{\pi}^o, P_O)-\mathbb{E}[\varphi_{\pi}(O_i;\alpha_{\pi}^o, P_O)\mid A_i]\right\}\right)^2
			\middle| A_i, i\in I_k\right]\notag\\
			= &\sum_{i\in I_k}K_h(A_i-a)^2\left(\dfrac{A_i-a}{h}\right)^{2j}
			\mathbb{E}\left[\left\{\varphi_{\pi}(O_i;\alpha_{\pi}^o, P_O)-\mathbb{E}[\varphi_{\pi}(O_i;\alpha_{\pi}^o, P_O)\mid A_i]\right\}^2\middle| A_i\right]\label{exp: 18}\\
			\lesssim &\dfrac{1}{h^2}\sum_{i\in I_k} I\left\{|A_i-a|\leq h\right\}\mathrm{Var}\left[\varphi_{\pi}(O;\alpha_{\pi}^o, P_O)\middle| A_i\right]
			\lesssim \dfrac{n_K}{h}\dfrac{1}{n_K}\sum_{i\in I_k} \dfrac{I\left\{|A_i-a|\leq h\right\}}{h},\label{exp: 19}
		\end{align}
		where $n_K=n/K$ is the sample number in $|I_k|$.
		Equation~\eqref{exp: 18} follows from the fact that $\{A_i\}_{i\in I_k}$ are mutually independent from each other.
		Equation~\eqref{exp: 19} follows from $K_h(A_i-a)=K((A_i-a)/h)/h\lesssim I\{|A_i-a|\leq h\}/h$ by Assumption~\ref{as: regular condition}\ref{cd: kernel condition}, and $\mathrm{Var}\left[\varphi_{\pi}(O;\alpha_{\pi}^o, P_O)\middle| A\right]$ can be uniformly bounded by some constant from Assumption~\ref{as: regular condition}\ref{cd: varphi}.
		Thus we have
		\begin{align*}
			\mathbb{E}[(R_{3kj})^2]
			\lesssim \dfrac{n}{h}\mathbb{E}\left[\dfrac{1}{n_K}\sum_{i\in I_k} \dfrac{I\left\{|A_i-a|\leq h\right\}}{h}\right]=O(\dfrac{n}{h}).
		\end{align*}
		Thus we end with $R_{3kj}=O_p(\sqrt{n/h})$.	
	\end{proof}
	
	\begin{proof}[Proof of Equation~\eqref{exp: 4}:]
		Next, we control $R_{4kj}$ by
		\begin{align}
			R_{4kj}
			\lesssim&\left|\sum_{i\in I_k}K_h(A_i-a)\left(\dfrac{A_i-a}{h}\right)^{j}\left\{\theta(A_i)-\theta(a)-\theta'(a)(A_i-a)\right\}\right|\notag\\
			\lesssim&\sum_{i\in I_k}\left|K_h(A_i-a)\right| \left(\dfrac{A_i-a}{h}\right)^{j}(A_i-a)^2\label{exp: 20}\\
			\lesssim&\dfrac{1}{h}\sum_{i\in I_k} I\{|A_i-a|\leq h\}(A_i-a)^2
			\lesssim\sum_{i\in I_k}I\{|A_i-a|\leq h\} h=O_p(nh^2).\label{exp: 21}
		\end{align}
		Equation~\eqref{exp: 20} follows from Assumption~\ref{as: regular condition}\ref{cd: theta condition} that $\theta(a)$ is twice continuously differentiable in a neighborhood of $a$. 
		The final equality holds by the positive and finite probability density at $p_A(a)$.	
	\end{proof}

	\subsection{Proof of Theorem~\ref{thm: clt}}
	Since we have shown that $R_{1kj},R_{2kj}=o_p(\sqrt{n/h}),$ which can be omitted, we can deduce that
	\begin{align}
		&\sqrt{nh}\left\{\hat{\theta}_{\pi,h}^{(n)}(a)-\theta(a)-\dfrac{1}{np_A(a)}\sum_{k=1}^KR_{4k0}\right\}\notag\\
		=&\sqrt{nh}e_1^{\top} Q_{h}(a)\left[\begin{array}{l}
			\sum_{k=1}^K R_{3k0}+R_{4k0}\\
			\sum_{k=1}^K R_{3k1}+R_{4k1}
		\end{array}\right] - \sqrt{\dfrac{h}{n}} \dfrac{1}{p_A(a)}\sum_{k=1}^K R_{4k0}+ o_p(1)\notag\\
		=&\sqrt{\dfrac{h}{n}}[1,0]\left\{\begin{array}{cc}
			\dfrac{1}{p_A(a)}&0\\
			0&\dfrac{1}{\int K(a)a^2\mathrm{d}a p_A(a)}
		\end{array}\right\}\left[\begin{array}{l}
			\sum_{k=1}^K \{R_{3k0}+R_{4k0}\}\\
			\sum_{k=1}^K \{R_{3k1}+R_{4k1}\}
		\end{array}\right] \label{exp: 22}\\
		&- \sqrt{\dfrac{h}{n}} \dfrac{1}{p_A(a)}\sum_{k=1}^K R_{4k0}+ o_p(1)\notag\\
		=&\dfrac{1}{p_A(a)}\sqrt{\dfrac{h}{n}}\sum_{k=1}^{K}R_{3k0}+o_p(1)\notag\\
		=&\dfrac{1}{p_A(a)}\sqrt{\dfrac{h}{n}}\sum_{i=1}^n K_{h}(A_i-a)
		\left\{\varphi_{\pi}(O_i;\alpha_{\pi}^o, P_O)-\mathbb{E}[\varphi_{\pi}(O_i;\alpha_{\pi}^o, P_O)\mid A_i]\right\}+o_p(1).\notag
	\end{align}
	Equation~\eqref{exp: 22} follows from Lemma~\ref{lem: lower bounded eigenvalue}.
	
	For any $i=1,\ldots, n$, define 
	$
	X_{ni} := \sqrt{\frac{h}{n}} K_{h}(A_i-a)
	\left\{\varphi_{\pi}(O_i;\alpha_{\pi}^o, P_O)-\mathbb{E}[\varphi_{\pi}(O_i;\alpha_{\pi}^o, P_O)\mid A_i]\right\}.
	$
	We can verify that
	\begin{align*}
		\mathrm{Var}\left[\sum_{i=1}^n X_{ni}\right]
		=&\mathrm{Var}\left[\sqrt{\dfrac{h}{n}}\sum_{i=1}^n K_{h}(A_i-a)
		\left\{\varphi_{\pi}(O_i;\alpha_{\pi}^o, P_O)-\mathbb{E}[\varphi_{\pi}(O_i;\alpha_{\pi}^o, P_O)\mid A_i]\right\}\right]\\
		=&h\mathrm{Var}\left[ K_{h}(A-a)
		\left\{\varphi_{\pi}(O;\alpha_{\pi}^o, P_O)-\mathbb{E}[\varphi_{\pi}(O;\alpha_{\pi}^o, P_O)\mid A]\right\}\right]\\
		=&h\mathbb{E}\left[ K_{h}(A-a)^2
		\left\{\varphi_{\pi}(O;\alpha_{\pi}^o, P_O)-\mathbb{E}[\varphi_{\pi}(O;\alpha_{\pi}^o, P_O)\mid A]\right\}^2\right]\\
		=&\dfrac{1}{h}\mathbb{E}\left[ K(\dfrac{A-a}{h})^2\mathrm{Var}\left[\varphi_{\pi}(O;\alpha_{\pi}^o, P_O)\mid A\right]\right]\\
		=&\dfrac{1}{h}\int K(\dfrac{a'-a}{h})^2\mathrm{Var}\left[\varphi_{\pi}(O;\alpha_{\pi}^o, P_O)\mid A=a'\right]p_A(a')\mathrm{d} a'\\
		=&\int K(s)^2\mathrm{Var}\left[\varphi_{\pi}(O;\alpha_{\pi}^o, P_O)\mid A=a+sh\right]p_A(a+sh)\mathrm{d} s\\
		=&\int K(s)^2\mathrm{d} s \times p_A(a)\mathrm{Var}\left[\varphi_{\pi}(O;\alpha_{\pi}^o, P_O)\mid A=a\right]+o(1).
	\end{align*}
	For any $\epsilon>0$, by Assumption~\ref{as: regular condition}\ref{cd: bandwidth rate}, for any $i$,
	\begin{align*}
		&\mathbb{E}\left[X_{ni}^2 I\{|X_{ni}|\geq \epsilon\}\right]\\
		\lesssim&\mathbb{E}\left[\dfrac{h}{n} K_{h}(A-a)^2 
		\left\{\varphi_{\pi}(O;\alpha_{\pi}^o, P_O)-\mathbb{E}[\varphi_{\pi}(O;\alpha_{\pi}^o, P_O)\mid A]\right\}^2 I\{\frac{1}{\sqrt{nh}}\geq \epsilon\}\right]\\
		=&\mathbb{E}\left[\dfrac{h}{n} K_{h}(A-a)^2\mathrm{Var}[\varphi_{\pi}(O;\alpha_{\pi}^o, P_O)\mid A]I\{nh\leq \frac{1}{\epsilon^2}\}\right]\\
		\lesssim& 
		\mathbb{E}\left[\dfrac{1}{nh} I\{|A-a|\leq h\}\mathrm{Var}[\varphi_{\pi}(O;\alpha_{\pi}^o, P_O)\mid A]
		I\{nh\leq \frac{1}{\epsilon^2}\}\right]=o(\dfrac{1}{n}).
	\end{align*}
	%	Here, the last inequality follows from Assumption~\ref{as: regular condition}\ref{cd: kernel condition}.
	Thus, by Lindeberg's central limit theorem, we know that
	\begin{align*}
		&\sqrt{nh}\left\{\hat{\theta}_{\pi,h}^{(n)}(a)-\theta(a)-\dfrac{1}{np_A(a)}\sum_{k=1}^K R_{4k0} \right\}
		\overset{d}{\to}N(0, \sigma_{\pi,\theta}^o(a)^2),\quad\text{ where }\\
		&\sigma_{\pi,\theta}^o(a)^2:=\dfrac{\int K(s)^2\mathrm{d} s}{p_A(a)}\mathrm{Var}\left[\varphi_{\pi}(O;\alpha_{\pi}^o, P_O)\mid A=a\right].
	\end{align*}
	Now we define $\mathrm{bias}(a):=\dfrac{1}{np_A(a)}\sum_{k=1}^K R_{4k0}$. From proof of Equation~\eqref{exp: 4}, we have seen that
	\begin{align*}
		\mathrm{bias}(a)=&\dfrac{1}{np_A(a)}\sum_{i=1}^n K_{h}(A_i-a)
		\left\{\theta(A_i)-\theta(a)-\theta'(a)(A_i-a)\right\}\\
		=&\dfrac{\theta''(a)}{2np_A(a)}\sum_{i=1}^n K_{h}(A_i-a)(A_i-a)^2 + o(h^2)\\
		=&\dfrac{\theta''(a)}{2p_A(a)}\mathbb{E} [K_{h}(A-a)(A-a)^2] + o(h^2)\\
		=&\dfrac{\theta''(a)}{2p_A(a)}\int K_{h}(s-a)(s-a)^2 p_A(s)\mathrm{d} s + o(h^2)\\
		=&\dfrac{h^2\theta''(a)}{2p_A(a)}\int K(s)s^2p_A(a+sh)\mathrm{d} s + o(h^2)\\
		=&\dfrac{h^2}{2}\theta''(a)\int K(s)s^2\mathrm{d} s + o(h^2).
	\end{align*}
	This finishes the proof of Equation~\eqref{eq: CLT}.
	\sy{
		\subsection{Proof of Theorem~\ref{thm: falsify AIV condition}}
		First, under Assumption~\ref{as: consistency}--\ref{as: positivity}, if $Z$ is an AIV for $A$, then $\mathbb{E}\left[\varphi_{\Pi}(O; \alpha_{\Pi}^o, P_O)\middle| A=a\right]=\mathbf{1}_{J}\theta(a)$. Therefore, we can define the oracle optimal weighting matrix
		\begin{align*}
			&Q_{\Pi,h}^o(a):=\lim_{h\to 0}h\mathbb{E}\left[\Psi_{\Pi,h,a}(O;\beta_{h}^o(a),\alpha_{\Pi}^o, P_O)\Psi_{\Pi,h,a}(O;\beta_{h}^o(a),\alpha_{\Pi}^o, P_O)^{\top}\right]\\
			=&\lim_{h\to 0}\frac{1}{h}\mathbb{E}\left[\begin{array}{l}
				K\left(\dfrac{A - a}{h}\right)^2 g\left(\dfrac{A-a}{h}\right)g\left(\dfrac{A-a}{h}\right)^{\top}\\\otimes \bigg(
				\left\{\varphi_{\Pi}(O; \alpha_{\Pi}^o, P_O)
				- \mathbf{1}_J g\left(\dfrac{A-a}{h}\right)^{\top} \beta_{h}^o(a)\right\}\\\times
				\left\{\varphi_{\Pi}(O; \alpha_{\Pi}^o, P_O)
				- \mathbf{1}_J g\left(\dfrac{A-a}{h}\right)^{\top} \beta_{h}^o(a)\right\}^{\top}\bigg)
			\end{array}\right]\\
			=&\lim_{h\to 0}\frac{1}{h}\mathbb{E}\left[
			K\left(\dfrac{A - a}{h}\right)^2 g\left(\dfrac{A-a}{h}\right)g\left(\dfrac{A-a}{h}\right)^{\top}\otimes 
			\mathrm{Var}\left\{\varphi_{\Pi}(O; \alpha_{\Pi}^o, P_O)\mid A\right\}\right]\\
			&+\lim_{h\to 0}\frac{1}{h}\mathbb{E}\left[
			\begin{array}{l}
				K\left(\dfrac{A - a}{h}\right)^2 g\left(\dfrac{A-a}{h}\right)g\left(\dfrac{A-a}{h}\right)^{\top}\\
				\otimes\bigg(\left\{\mathbb{E}\left[\varphi_{\Pi}(O; \alpha_{\Pi}^o, P_O)\mid A\right] - \mathbf{1}_Jg\left(\dfrac{A-a}{h}\right)^{\top} \beta_{h}^o(a)\right\}\\
				\times \left\{\mathbb{E}\left[\varphi_{\Pi}(O; \alpha_{\Pi}^o, P_O)\mid A\right] - \mathbf{1}_Jg\left(\dfrac{A-a}{h}\right)^{\top} \beta_{h}^o(a)\right\}^{\top}\bigg)
			\end{array}\right]\\
			=&\lim_{h\to 0}\frac{1}{h}\int 
			K\left(\dfrac{s-a}{h}\right)^2 g\left(\dfrac{s-a}{h}\right)g\left(\dfrac{s-a}{h}\right)^{\top}\otimes 
			\mathrm{Var}\left\{\varphi_{\Pi}(O; \alpha_{\Pi}^o, P_O)\mid A=s\right\}p_A(s)\mathrm{d}s\\
			&+\lim_{h\to 0}\frac{1}{h}\mathbb{E}\left[
			\begin{array}{l}
				K\left(\dfrac{A - a}{h}\right)^2 g\left(\dfrac{A-a}{h}\right)g\left(\dfrac{A-a}{h}\right)^{\top}\\
				\left\{\theta(A) - g\left(\dfrac{A-a}{h}\right)^{\top} \beta_{h}^o(a)\right\}^2
			\end{array}\right]\otimes 
			\mathbf{1}_J\mathbf{1}_J^{\top}\\
			=&\lim_{h\to 0}\int 
			K\left(s\right)^2 g\left(s\right)g\left(s\right)^{\top}\otimes 
			\mathrm{Var}\left\{\varphi_{\Pi}(O; \alpha_{\Pi}^o, P_O)\mid A=a+hs\right\}p_A(a+hs)\mathrm{d}s\\
			&+\lim_{h\to 0}\int 
			\left(\begin{array}{l}
				K\left(s\right)^2 g\left(s\right)g\left(s\right)^{\top}\\
				\left\{\theta(a+sh) - \theta(a) - \theta'(a)sh\right\}^2
			\end{array}\right)p_{A}(a+sh)\mathrm{d}s\otimes 
			\mathbf{1}_J\mathbf{1}_J^{\top}\\
			=&p_A(a)\int 
			K\left(s\right)^2 g\left(s\right)g\left(s\right)^{\top} \mathrm{d}s \otimes 
			\mathrm{Var}\left\{\varphi_{\Pi}(O; \alpha_{\Pi}^o, P_O)\mid A=a\right\}\in\mathbb{R}^{2J\times 2J}.
		\end{align*}
		The final step follows from the fact that $\theta(a+sh) - \theta(a) - \theta'(a)sh=O(h^2)$.
		Define the oracle AIPW estimator 
		\begin{align*}
			\tilde{\beta}_{\Pi,h}^{(n)}(a):=\underset{\beta \in \mathbb{R}^2}{\arg\min}\left\{
			\mathbb{E}_{n}\left[\Psi_{\Pi,h,a}(O;\beta,\alpha_{\Pi}^o, P_O)\right]^{\top} Q_{\Pi,h}^o(a)^{-1}
			\mathbb{E}_{n}\left[\Psi_{\Pi,h,a}(O;\beta,\alpha_{\Pi}^o, P_O)\right]\right\},
		\end{align*}
		We can observe that 
		\begin{align*}
			0=&2
			\mathbb{E}_{n}\left[
			\nabla_{\beta}\Psi_{\Pi,h,a}(O;\tilde{\beta}_{\Pi,h}^{(n)}(a),\alpha_{\Pi}^o, P_O)\right]^{\top} Q_{\Pi,h}^o(a)^{-1}
			\mathbb{E}_{n}\left[\Psi_{\Pi,h,a}(O;\tilde{\beta}_{\Pi,h}^{(n)}(a),\alpha_{\Pi}^o, P_O)\right]\\
			=&2\mathbb{E}_{n}\left[
			\nabla_{\beta}\Psi_{\Pi,h,a}(O;\tilde{\beta}_{\Pi,h}^{(n)}(a),\alpha_{\Pi}^o, P_O)\right]^{\top} Q_{\Pi,h}^o(a)^{-1}\\
			&\left\{\mathbb{E}_{n}\left[\Psi_{\Pi,h,a}(O;\beta_{h}^o(a),\alpha_{\Pi}^o, P_O)\right] +
			\left(\mathbb{E}_{n}\left[\nabla_{\beta}\Psi_{\Pi,h,a}(O;\beta_{h}^o(a),\alpha_{\Pi}^o, P_O)\right] + o_p(1)\right) (\tilde{\beta}_{\Pi,h}^{(n)}(a)-\beta_{h}^o(a))\right\}\\
			=&2\left(\mathbb{E}\left[\nabla_{\beta}\Psi_{\Pi,h,a}(O;\tilde{\beta}_{\Pi,h}^{(n)}(a),\alpha_{\Pi}^o, P_O)\right]+o_p(1)\right)^{\top} Q_{\Pi,h}^o(a)^{-1}\\
			&\left\{\mathbb{E}_{n}\left[\Psi_{\Pi,h,a}(O;\beta_{h}^o(a),\alpha_{\Pi}^o, P_O)\right] +
			\left\{\mathbb{E}\left[\nabla_{\beta}\Psi_{\Pi,h,a}(O;\beta_{h}^o(a),\alpha_{\Pi}^o, P_O)\right] + o_p(1)\right\} (\tilde{\beta}_{\Pi,h}^{(n)}(a)-\beta_{h}^o(a))\right\}.
		\end{align*}
		Now, we know that 
		\begin{align*}
			&\sqrt{nh}(\tilde{\beta}_{\Pi,h}^{(n)}(a)-\beta_{h}^o(a))\\=&-\left\{\mathbb{E}\left[\nabla_{\beta}\Psi_{\Pi,h,a}(O;\beta_{h}^o(a),\alpha_{\Pi}^o, P_O)\right]^{\top} Q_{\Pi,h}^o(a)^{-1}\mathbb{E}\left[\nabla_{\beta}\Psi_{\Pi,h,a}(O;\beta_{h}^o(a),\alpha_{\Pi}^o, P_O)\right]\right\}^{-1}\\
			&\times\mathbb{E}\left[\nabla_{\beta}\Psi_{\Pi,h,a}(O;\beta_{h}^o(a),\alpha_{\Pi}^o, P_O)\right]^{\top} Q_{\Pi,h}^o(a)^{-1}\sqrt{\dfrac{h}{n}}\sum_{i=1}^n\Psi_{\Pi,h,a}(O;\beta_{h}^o(a),\alpha_{\Pi}^o, P_O) + o_p(1).
		\end{align*}
		Now, we can calculate
		\begin{align*}
			&\lim_{h\to 0}\mathbb{E}\left[\nabla_{\beta}\Psi_{\Pi,h,a}\right]
			=-\lim_{h\to 0}\mathbb{E}\left[K_h(A-a)g\left(\dfrac{A-a}{h}\right)\otimes \left\{\mathbf{1}_J g\left(\dfrac{A-a}{h}\right)^{\top}\right\}\right]\\
			=&-\mathbf{I}_2\otimes\mathbf{1}_J\lim_{h\to 0}\mathbb{E}\left[K_h(A-a)g\left(\dfrac{A-a}{h}\right) g\left(\dfrac{A-a}{h}\right)^{\top}\right]
			\\
			=&-p_A(a) \mathbf{I}_2\otimes\mathbf{1}_J\int K(s)g\left(s\right) g\left(s\right)^{\top}\mathrm{d}s
			=-p_A(a)(\mathbf{I}_2\otimes\mathbf{1}_J)\left[\begin{array}{ll}
				1 & 0 \\ 0 & \int K(s)s^2\mathrm{d}s
			\end{array}\right]\\
			=&-p_A(a)\left(\left[\begin{array}{ll}
				1 & 0 \\ 0 & \int K(s)s^2\mathrm{d}s
			\end{array}\right]\otimes\mathbf{1}_J\right).
		\end{align*}
		Therefore, the asymptotic distribution of $\tilde{\beta}_{\Pi,h}^{(n)}(a)$ consists of 
		\begin{align*}
			&\lim_{h\to 0}\left\{\mathbb{E}\left[\nabla_{\beta}\Psi_{\Pi,h,a}\right]^{\top} Q_{\Pi,h}^o(a)^{-1}\mathbb{E}\left[\nabla_{\beta}\Psi_{\Pi,h,a}\right]\right\}^{-1}\\
			=&\dfrac{1}{p_A(a)}\left\{\left[\begin{array}{ll}
				1 & 0 \\ 0 & \int K(s)s^2\mathrm{d}s
			\end{array}\right]\otimes\mathbf{1}_J^{\top}\left(\begin{array}{l}
				(\int 
				K\left(s\right)^2 g\left(s\right)g\left(s\right)^{\top} \mathrm{d}s)^{-1}\;\otimes\\
				\mathrm{Var}\left\{\varphi_{\Pi}(O; \alpha_{\Pi}^o, P_O)| A=a\right\}^{-1}
			\end{array}\right)\left[\begin{array}{cc}
				1 & 0 \\ 0 & \int K(s)s^2\mathrm{d}s
			\end{array}\right]\otimes
			\mathbf{1}_J
			\right\}^{-1}\\
			=&\dfrac{1}{p_A(a)}\left(\left[\begin{array}{cc}
				\frac{1}{\int K\left(s\right)^2 \mathrm{d}s} & 0 \\
				0 & \dfrac{(\int K(s)s^2\mathrm{d}s)^2}{\int K\left(s\right)^2 s^2\mathrm{d}s}
			\end{array}\right]\otimes 
			(\mathbf{1}_J^{\top}
			\mathrm{Var}\left\{\varphi_{\Pi}(O; \alpha_{\Pi}^o, P_O)| A=a\right\}^{-1}\mathbf{1}_J)
			\right)^{-1}\\
			=&\mathrm{diag}\left\{
			\int K\left(s\right)^2 \mathrm{d}s, \dfrac{\int K\left(s\right)^2 s^2\mathrm{d}s}{(\int K(s)s^2\mathrm{d}s)^2}\right\}\left(p_A(a)\mathbf{1}_J^{\top}
			\mathrm{Var}\left\{\varphi_{\Pi}(O; \alpha_{\Pi}^o, P_O)| A=a\right\}^{-1}\mathbf{1}_J\right)^{-1}.
		\end{align*}
		This finishes the proof of the first conclusion of Theorem~\ref{thm: falsify AIV condition}.
		
		Next, to compute the limiting distribution of $(nh)\cdot\mathcal{L}_{\Pi,h,a}^{(n)}(\hat{\beta}_{\Pi,h}^{(n)}(a),\hat{Q}_{\Pi,h}^{(n)}(a)^{-1})$, we first calculate the limiting distribution of $(nh)\cdot\mathcal{L}_{\Pi,h,a}^{(n)}(\tilde{\beta}_{\Pi,h}^{(n)}(a),Q_{\Pi,h}^o(a)^{-1})$. In particular, we can decompose
		\begin{align*}
			&(nh)\cdot\mathcal{L}_{\Pi,h,a}^{(n)}(\tilde{\beta}_{\Pi,h}^{(n)}(a),Q_{\Pi,h}^o(a)^{-1})\\
			=&nh
			\mathbb{E}_{n}\left[\Psi_{\Pi,h,a}(O;\tilde{\beta}_{\Pi,h}^{(n)}(a),\alpha_{\Pi}^o, P_O)\right]^{\top} Q_{\Pi,h}^o(a)^{-1}
			\mathbb{E}_{n}\left[\Psi_{\Pi,h,a}(O;\tilde{\beta}_{\Pi,h}^{(n)}(a),\alpha_{\Pi}^o, P_O)\right] + o_p(1)\\
			=&\left\|\begin{array}{l}
				\sqrt{nh}\mathbb{E}_{n}\left[\Psi_{\Pi,h,a}(O;\beta_{h}^o(a),\alpha_{\Pi}^o, P_O)\right] \\+ \sqrt{nh}\mathbb{E}\left[\nabla_{\beta}\Psi_{\Pi,h,a}(O;\beta_{h}^o(a),\alpha_{\Pi}^o, P_O)\right](\tilde{\beta}_{\Pi,h}^{(n)}(a)-\beta_{h}^o(a))
			\end{array}\right\|_{Q_{\Pi,h}^o(a)^{-1}}^2 + o_p(1)\\
			=&\left\|\left(\begin{array}{l}
				I_{2J}- \mathbb{E}\left[\nabla_{\beta}\Psi_{\Pi,h,a}\right]
				\left\{\mathbb{E}\left[\nabla_{\beta}\Psi_{\Pi,h,a}\right]^{\top} Q_{\Pi,h}^o(a)^{-1}\mathbb{E}\left[\nabla_{\beta}\Psi_{\Pi,h,a}\right]\right\}^{-1}\\
				\times\mathbb{E}\left[\nabla_{\beta}\Psi_{\Pi,h,a}\right]^TQ_{\Pi,h}^o(a)^{-1}
			\end{array}\right)
			\sqrt{nh}\mathbb{E}_{n}\left[\Psi_{\Pi,h,a}\right]\right\|_{Q_{\Pi,h}^o(a)^{-1}}^2 + o_p(1)
		\end{align*}
		Since $\sqrt{nh}\mathbb{E}_{n}\left[\Psi_{\Pi,h,a}(O;\beta_{h}^o(a),\alpha_{\Pi}^o,P_O)\right]$ converges to a normal distribution, we know that the limiting distribution of $(nh)\cdot\mathcal{L}_{\Pi,h,a}^{(n)}(\tilde{\beta}_{\Pi,h}^{(n)}(a),Q_{\Pi,h}^o(a)^{-1})$ can be approximated by a $\chi^2$ distribution, where its degree of freedom equals
		\begin{align*}
			&\lim_{n\to\infty}\mathbb{E}\left[\mathrm{tr}\left(
			\left\|\begin{array}{l}
				\Big(I_{2J}- \mathbb{E}\left[\nabla_{\beta}\Psi_{\Pi,h,a}\right]
				\left\{\mathbb{E}\left[\nabla_{\beta}\Psi_{\Pi,h,a}\right]^{\top} Q_{\Pi,h}^o(a)^{-1}\mathbb{E}\left[\nabla_{\beta}\Psi_{\Pi,h,a}\right]\right\}^{-1}\\
				\times\mathbb{E}\left[\nabla_{\beta}\Psi_{\Pi,h,a}\right]^TQ_{\Pi,h}^o(a)^{-1}\Big)\sqrt{nh}\mathbb{E}_{n}\left[\Psi_{\Pi,h,a}\right]
			\end{array}
			\right\|_{Q_{\Pi,h}^o(a)^{-1}}^2
			\right)\right]\\
			=&\mathrm{tr}\left(\begin{array}{l}
				\lim\limits_{n\to\infty}\mathbb{E}\left[nh\mathbb{E}_{n}\left[\Psi_{\Pi,h,a}\right]\mathbb{E}_{n}\left[\Psi_{\Pi,h,a}\right]^{\top}\right]
				Q_{\Pi,h}^o(a)^{-1}\\
				\left\{\begin{array}{l}
					I_{2J}- \mathbb{E}\left[\nabla_{\beta}\Psi_{\Pi,h,a}\right]
					\left\{\mathbb{E}\left[\nabla_{\beta}\Psi_{\Pi,h,a}\right]^{\top} Q_{\Pi,h}^o(a)^{-1}\mathbb{E}\left[\nabla_{\beta}\Psi_{\Pi,h,a}\right]\right\}^{-1}\\
					\times\mathbb{E}\left[\nabla_{\beta}\Psi_{\Pi,h,a}\right]^TQ_{\Pi,h}^o(a)^{-1}
				\end{array}\right\}^2
			\end{array}\right)\\
			=&\mathrm{tr}\left(
			I_{2J}- \mathbb{E}\left[\nabla_{\beta}\Psi_{\Pi,h,a}\right]
			\left\{\mathbb{E}\left[\nabla_{\beta}\Psi_{\Pi,h,a}\right]^{\top} Q_{\Pi,h}^o(a)^{-1}\mathbb{E}\left[\nabla_{\beta}\Psi_{\Pi,h,a}\right]\right\}^{-1}
			\mathbb{E}\left[\nabla_{\beta}\Psi_{\Pi,h,a}\right]^TQ_{\Pi,h}^o(a)^{-1}\right)\\
			=&2J-\mathrm{tr}\left(
			\left\{\mathbb{E}\left[\nabla_{\beta}\Psi_{\Pi,h,a}\right]^{\top} Q_{\Pi,h}^o(a)^{-1}\mathbb{E}\left[\nabla_{\beta}\Psi_{\Pi,h,a}\right]\right\}^{-1}
			\mathbb{E}\left[\nabla_{\beta}\Psi_{\Pi,h,a}\right]^TQ_{\Pi,h}^o(a)^{-1}\mathbb{E}\left[\nabla_{\beta}\Psi_{\Pi,h,a}\right]\right)
			=2J-2.
		\end{align*}
		The second equality follows from the fact that $\lim\limits_{n\to\infty}\mathbb{E}\left[nh\mathbb{E}_{n}\left[\Psi_{\Pi,h,a}\right]\mathbb{E}_{n}\left[\Psi_{\Pi,h,a}\right]^{\top}\right]=Q_{\Pi,h}^o(a)$.
		Now, we finish the proof of 
		$(nh)\cdot\mathcal{L}_{\Pi,h,a}^{(n)}(\tilde{\beta}_{\Pi,h}^{(n)}(a),Q_{\Pi,h}^o(a)^{-1})\overset{d}{\to}\chi^2_{2J-2}.$
		Therefore, the remaining things is to show that  $\hat{Q}_{\Pi,h}^{(n)}(a)-Q_{\Pi,h}^o(a)=o_p(1)$, and 
		\begin{align}\label{expadd: 10}
			\hat{\beta}_{\Pi,h}^{(n)}(a) - \tilde{\beta}_{\Pi,h}^{(n)}(a)= o_p(1/\sqrt{nh}).
		\end{align}
		First, we show the proof of
		\begin{align}
			&\hat{\beta}_{\Pi,h}^{(n,\mathrm{init})}(a)-\tilde{\beta}_{\Pi,h}^{(n,\mathrm{init})}(a)=o_p(1/\sqrt{nh}),\text{ where }\label{expadd: 11}\\
			&\tilde{\beta}_{\Pi,h}^{(n,\mathrm{init})}(a):=
			\underset{\beta\in\mathbb{R}^2}{\arg\min}\left\{
			\mathbb{E}_{n}\left[\Psi_{\Pi,h,a}(O;\beta,\alpha_{\Pi}^o, P_O)\right]^{\top} 
			\mathbb{E}_{n}\left[\Psi_{\Pi,h,a}(O;\beta,\alpha_{\Pi}^o, P_O)\right]\right\},\notag\\
			&\hat{\beta}_{\Pi,h}^{(n,\mathrm{init})}(a):=
			\underset{\beta\in\mathbb{R}^2}{\arg\min}\left\{
			\begin{array}{l}
				\left(\dfrac{1}{K}\sum\limits_{k=1}^K\mathbb{E}_{nk}\left[\Psi_{\Pi,h,a}(O;\beta,\hat\alpha_{\Pi}^{(n,-k)}, \hat{P}_O^{(n,-k)})\right]\right)^{\top}\\
				\times\dfrac{1}{K}\sum\limits_{k=1}^K\mathbb{E}_{nk}\left[\Psi_{\Pi,h,a}(O;\beta,\hat\alpha_{\Pi}^{(n,-k)}, \hat{P}_O^{(n,-k)})\right]
			\end{array}\right\}.\notag
		\end{align}
		
		\begin{proof}[Proof of Equation~\eqref{expadd: 11}:]
			We assume that $\tilde{\beta}_{\Pi,h}^{(n,\mathrm{init})}(a)$ and $\hat{\beta}_{\Pi,h}^{(n,\mathrm{init})}(a)$ belong to a compact subset $\mathcal{S}\subseteq\mathbb{R}^2$. Then, we control 
			\begin{align*}
				&\nabla_{\beta}\left(\begin{array}{l}
					\mathbb{E}_{n}\left[\Psi_{\Pi,h,a}(O;\beta,\alpha_{\Pi}^o, P_O)\right]^{\top} 
					\mathbb{E}_{n}\left[\Psi_{\Pi,h,a}(O;\beta,\alpha_{\Pi}^o, P_O)\right]\\
					-\left(\dfrac{1}{K}\sum\limits_{k=1}^K\mathbb{E}_{nk}\left[\Psi_{\Pi,h,a}(O;\beta,\hat\alpha_{\Pi}^{(n,-k)}, \hat{P}_O^{(n,-k)})\right]\right)^{\top}\\
					\times\dfrac{1}{K}\sum\limits_{k=1}^K\mathbb{E}_{nk}\left[\Psi_{\Pi,h,a}(O;\beta,\hat\alpha_{\Pi}^{(n,-k)}, \hat{P}_O^{(n,-k)})\right]
				\end{array}\right)\\
				=&2\mathbb{E}_{n}\left[\nabla_{\beta}\Psi_{\Pi,h,a}(O;\beta,\alpha_{\Pi}^o, P_O)\right]^{\top} 
				\mathbb{E}_{n}\left[\Psi_{\Pi,h,a}(O;\beta,\alpha_{\Pi}^o, P_O)\right]\\
				&-2\left(\dfrac{1}{K}\sum\limits_{k=1}^K\mathbb{E}_{nk}\left[\nabla_{\beta}\Psi_{\Pi,h,a}(O;\beta,\hat\alpha_{\Pi}^{(n,-k)}, \hat{P}_O^{(n,-k)})\right]\right)^{\top}
				\times\dfrac{1}{K}\sum\limits_{k=1}^K\mathbb{E}_{nk}\left[\Psi_{\Pi,h,a}(O;\beta,\hat\alpha_{\Pi}^{(n,-k)}, \hat{P}_O^{(n,-k)})\right]\\
				=&2\mathbb{E}_{n}\left[\nabla_{\beta}\Psi_{\Pi,h,a}(O;\beta,\alpha_{\Pi}^o, P_O)\right]^{\top} 
				\left\{\mathbb{E}_{n}\left[\Psi_{\Pi,h,a}(O;\beta,\alpha_{\Pi}^o, P_O)\right]-\dfrac{1}{K}\sum\limits_{k=1}^K\mathbb{E}_{nk}\left[\Psi_{\Pi,h,a}(O;\beta,\hat\alpha_{\Pi}^{(n,-k)}, \hat{P}_O^{(n,-k)})\right]\right\}\\
				=&2\mathbb{E}_{n}\left[\nabla_{\beta}\Psi_{\Pi,h,a}(O;\beta,\alpha_{\Pi}^o, P_O)\right]^{\top} 
				\left\{\begin{array}{l}
					\mathbb{E}_{n}\left[K_h(A - a) g\left(\dfrac{A-a}{h}\right)\otimes\varphi_{\Pi}(O; \alpha_{\Pi}^o, P_O)\right]\\
					-\dfrac{1}{K}\sum\limits_{k=1}^K\mathbb{E}_{nk}\left[K_h(A - a)g\left(\dfrac{A-a}{h}\right)\otimes
					\varphi_{\Pi}(O;\hat\alpha_{\Pi}^{(n,-k)}, \hat{P}_O^{(n,-k)})\right]
				\end{array}\right\}.
			\end{align*}
			To control this quantity, for $n_K=n/K$, we calculate
			\begin{align*}
				&\sqrt{n_Kh}\mathbb{E}_{nk}\left[K_h(A - a)g\left(\dfrac{A-a}{h}\right)\otimes
				\left\{\varphi_{\Pi}(O;\hat\alpha_{\Pi}^{(n,-k)}, \hat{P}_O^{(n,-k)})-\varphi_{\Pi}(O; \alpha_{\Pi}^o, P_O)\right\}\right]\\
				=&\sqrt{n_Kh}\mathbb{E}_{nk}\left[K_h(A - a)g\left(\dfrac{A-a}{h}\right)\otimes
				\left\{\varphi_{\Pi}(O;\hat\alpha_{\Pi}^{(n,-k)}, \hat{P}_O^{(n,-k)})-\varphi_{\Pi}(O; \hat\alpha_{\Pi}^{(n,-k)}, P_O)\right\}\right]\\
				&+\mathbb{G}_{n_K}\left[\sqrt{h}K_h(A - a)g\left(\dfrac{A-a}{h}\right)\otimes
				\left\{\varphi_{\Pi}(O;\hat\alpha_{\Pi}^{(n,-k)}, P_O)-\varphi_{\Pi}(O; \alpha_{\Pi}^o, P_O)\right\}\right]\\
				&+\sqrt{n_Kh}\mathbb{E}_k\left[K_h(A - a)g\left(\dfrac{A-a}{h}\right)\otimes
				\left\{\mathbb{E}\left[\varphi_{\Pi}(O;\hat\alpha_{\Pi}^{(n,-k)}, P_O)-\varphi_{\Pi}(O; \alpha_{\Pi}^o, P_O)\middle| A\right]\right\}\right]=o_p(1).
			\end{align*}
			The final step follows from the fact stated in Equations~\eqref{exp: 1}--\eqref{exp: 3}. This indicates that
			\begin{align*}
				\sup_{\beta\in\mathcal{S}}\left|\nabla_{\beta}\left(\begin{array}{l}
					\mathbb{E}_{n}\left[\Psi_{\Pi,h,a}(O;\beta,\alpha_{\Pi}^o, P_O)\right]^{\top} 
					\mathbb{E}_{n}\left[\Psi_{\Pi,h,a}(O;\beta,\alpha_{\Pi}^o, P_O)\right]\\
					-\left(\dfrac{1}{K}\sum\limits_{k=1}^K\mathbb{E}_{nk}\left[\Psi_{\Pi,h,a}(O;\beta,\hat\alpha_{\Pi}^{(n,-k)}, \hat{P}_O^{(n,-k)})\right]\right)^{\top}\\
					\times\dfrac{1}{K}\sum\limits_{k=1}^K\mathbb{E}_{nk}\left[\Psi_{\Pi,h,a}(O;\beta,\hat\alpha_{\Pi}^{(n,-k)}, \hat{P}_O^{(n,-k)})\right]
				\end{array}\right)\right| = o_p\left(\frac{1}{\sqrt{nh}}\right).
			\end{align*}
			This indicates that $\tilde{\beta}_{\Pi,h}^{(n,\mathrm{init})}(a) - \hat{\beta}_{\Pi,h}^{(n,\mathrm{init})}(a) = o_p(1/\sqrt{nh})$.
		\end{proof}
		%		\begin{align*}
			%			\Psi_{\Pi,h,a}(O;\beta,\alpha_{\Pi},P_O):= K_h(A - a)g\left(\dfrac{A-a}{h}\right)\otimes
			%			\left\{\varphi_{\Pi}(O; \alpha_{\Pi}, P_O)
			%			- \mathbf{1}_J g\left(\dfrac{A-a}{h}\right)^{\top} \beta\right\}
			%		\end{align*}
		\begin{proof}[Proof of $\hat{Q}_{\Pi,h}^{(n)}(a) - Q_{\Pi,h}^o(a)=o_p(1)$:]
			It is direct to verify that $\hat{Q}_{\Pi,h}^{(n)}(a) - Q_{\Pi,h}^o(a)$ equals
			\begin{align*}
				&\frac{h}{K}\sum_{k=1}^{K}\mathbb{E}_{nk}\left[\Psi_{\Pi,h,a}(O;\hat{\beta}_{\Pi,h}^{(n,\mathrm{init})}(a),\hat\alpha_{\Pi}^{(n,-k)}, \hat{P}_O^{(n,-k)})\Psi_{\Pi,h,a}(O;\hat{\beta}_{\Pi,h}^{(n,\mathrm{init})}(a),\hat\alpha_{\Pi}^{(n,-k)}, \hat{P}_O^{(n,-k)})^{\top}\right]\\
				&-h\mathbb{E}\left[\Psi_{\Pi,h,a}(O;\beta_{h}^o(a),\alpha_{\Pi}^o, P_O)\Psi_{\Pi,h,a}(O;\beta_{h}^o(a),\alpha_{\Pi}^o, P_O)^{\top}\right]+o(1)\\
				=&\frac{h}{K}\sum_{k=1}^{K}\mathbb{E}_{k}\left[\begin{array}{l}
					\Psi_{\Pi,h,a}(O;\hat{\beta}_{\Pi,h}^{(n,\mathrm{init})}(a),\hat\alpha_{\Pi}^{(n,-k)}, \hat{P}_O^{(n,-k)})\Psi_{\Pi,h,a}(O;\hat{\beta}_{\Pi,h}^{(n,\mathrm{init})}(a),\hat\alpha_{\Pi}^{(n,-k)}, \hat{P}_O^{(n,-k)})^{\top}\\
					-\Psi_{\Pi,h,a}(O;\beta_{h}^o(a),\alpha_{\Pi}^o, P_O)\Psi_{\Pi,h,a}(O;\beta_{h}^o(a),\alpha_{\Pi}^o, P_O)^{\top}
				\end{array}\right]+o_p(1)\\
				=&O_p(\|\hat{\beta}_{\Pi,h}^{(n,\mathrm{init})}(a)-\beta_{h}^o(a)\|_2)
				+O_p(\|\hat\alpha_{\Pi}^{(n,-k)}-\alpha_{\Pi}^o\|_{\mathcal{B}(a,h),2})+o_p(1)=o_p(1).
			\end{align*}
			The final step follows from Equation~\eqref{expadd: 11} and the fact that $\mathbb{E}[\|\hat\alpha_{\Pi}^{(n,-k)}-\alpha_{\Pi}^o\|_{\mathcal{B}(a,h),2}^2]=o(1)$.
		\end{proof}
		Finally, the proof of Equation~\eqref{expadd: 10} is similar to that of Equation~\eqref{expadd: 11}, thus omitted.
		\begin{align*}
			&nh\mathcal{L}_{\Pi,h,a}^{(n)}(\hat{\beta}_{\Pi,h}^{(n)}(a),\hat{Q}_{\Pi,h}^{(n)}(a)^{-1})-
			nh\mathcal{L}_{\Pi,h,a}^{(n)}(\tilde{\beta}_{\Pi,h}^{(n)}(a),Q_{\Pi,h}^o(a)^{-1})\\
			=&nh\left\|\frac{1}{K}\sum_{k=1}^K\mathbb{E}_{nk}\left[\Psi_{\Pi,h,a}(O;\hat{\beta}_{\Pi,h}^{(n)}(a),\hat{\alpha}_{\Pi}^{(n,-k)}, \hat{P}_O^{(n,-k)})\right]\right\|_{\hat{Q}_{\Pi,h}^{(n)}(a)^{-1}}^2
			\\&-nh\left\|\frac{1}{K}\sum_{k=1}^K
			\mathbb{E}_{nk}\left[\Psi_{\Pi,h,a}(O;\tilde{\beta}_{\Pi,h}^{(n)}(a),\hat{\alpha}_{\Pi}^{(n,-k)}, \hat{P}_O^{(n,-k)})\right]\right\|_{Q_{\Pi,h}^o(a)^{-1}}^2\\
			=&nh\left\|\frac{1}{K}\sum_{k=1}^K\mathbb{E}_{nk}\left[\Psi_{\Pi,h,a}(O;\tilde{\beta}_{\Pi,h}^{(n)}(a),\hat{\alpha}_{\Pi}^{(n,-k)}, \hat{P}_O^{(n,-k)})\right]
			+\mathbb{E}_{n}\left[\nabla_{\beta}\Psi_{\Pi,h,a}\right](\hat{\beta}_{\Pi,h}^{(n)}(a)-\tilde{\beta}_{\Pi,h}^{(n)}(a))\right\|_{Q_{\Pi,h}^o(a)^{-1}}^2\\
			&-nh\left\|\frac{1}{K}\sum_{k=1}^K
			\mathbb{E}_{nk}\left[\Psi_{\Pi,h,a}(O;\tilde{\beta}_{\Pi,h}^{(n)}(a),\hat{\alpha}_{\Pi}^{(n,-k)}, \hat{P}_O^{(n,-k)})\right]\right\|_{Q_{\Pi,h}^o(a)^{-1}}^2+o_p(1)\\
			=&o_p(nh\|\hat{\beta}_{\Pi,h}^{(n)}(a)-\tilde{\beta}_{\Pi,h}^{(n)}(a)\|_2^2)+o_p(1)=o_p(1).
		\end{align*}
		Combining all these results together, we know that $nh\mathcal{L}_{\Pi,h,a}^{(n)}(\hat{\beta}_{\Pi,h}^{(n)}(a),\hat{Q}_{\Pi,h}^{(n)}(a)^{-1})\overset{d}{\to} \chi^2_{2J-2}$.
		This finishes the proof of the latter part of Theorem~\ref{thm: falsify AIV condition}.
	}
	\sy{
		\subsection{Proof of Theorem~\ref{thm: RWF testing 1}}
		We first prove the result for $L=\emptyset$. The case of discrete $L$ follows by applying the same argument within each stratum $L=l$. From Algorithm~\ref{alg: test RWF discrete}, we can see that
		\begin{align*}
			&\tilde{\zeta}_{\pi,h_A}^{(n)}(a,l) =\frac{1}{n}\sum_{k=1}^K \sum_{i\in I_k}(K_{h_A}(A_i-a) - \hat{\gamma}_{h_A,a}^{(n,-k)}(L_i))(\pi(Z_i,L_i) - \hat{\rho}_{\pi}^{(n,-k)}(L_i))\\
			=&\frac{1}{K}\sum_{k=1}^K\frac{1}{n_K}\sum_{i\in I_k} \left\{\begin{array}{l}
				(K_{h_A}(A_i-a) - \hat{\gamma}_{h_A,a}^{(n,-k)}(L))(\pi(Z_i,L_i) - \hat{\rho}_{\pi}^{(n,-k)}(L_i))\\
				-\mathbb{E}_k\left[
				(K_{h_A}(A-a) - \hat{\gamma}_{h_A,a}^{(n,-k)}(L))(\pi(Z,L) - \hat{\rho}_{\pi}^{(n,-k)}(L))
				\right]
			\end{array}
			\right\}\\
			&+\dfrac{1}{K}\sum_{k=1}^K
			\mathbb{E}_k\left[
			(K_{h_A}(A-a) - \hat{\gamma}_{h_A,a}^{(n,-k)}(L))(\pi(Z,L) - \hat{\rho}_{\pi}^{(n,-k)}(L))
			\right]\\
			=&\frac{1}{K}\sum_{k=1}^K\{\mathbb{E}_{nk}-\mathbb{E}_k\} \left[
			(K_{h_A}(A-a) - \gamma_{h_A,a}^o(L))(\pi(Z,L) - \rho_{\pi}^o(L))\right]\\
			&+\frac{1}{K}\sum_{k=1}^K\{\mathbb{E}_{nk}-\mathbb{E}_k\} \left[\begin{array}{l}
				(K_{h_A}(A_i-a) - \hat{\gamma}_{h_A,a}^{(n,-k)}(L_i))(\pi(Z_i,L_i) - \hat{\rho}_{\pi}^{(n,-k)}(L_i))\\
				-(K_{h_A}(A_i-a) - \gamma_{h_A,a}^o(L_i))(\pi(Z_i,L_i) - \rho_{\pi}^o(L_i))
			\end{array}\right]\\
			&+\dfrac{1}{K}\sum_{k=1}^K
			\mathbb{E}_k\left[
			(K_{h_A}(A-a) - \hat{\gamma}_{h_A,a}^{(n,-k)}(L))(\pi(Z,L) - \hat{\rho}_{\pi}^{(n,-k)}(L))
			\right]\\
			=&\frac{1}{K}\sum_{k=1}^K\{\mathbb{E}_{nk}-\mathbb{E}_k\} \left[
			(K_{h_A}(A-a) - \gamma_{h_A,a}^o(L))(\pi(Z,L) - \rho_{\pi}^o(L))\right]+
			\frac{1}{K}\sum_{k=1}^K\{R_{8k} + R_{9k}\}.
		\end{align*}
		First, we control $R_{8k}$ as follows. We verify that
		\begin{align*}
			&\mathrm{Var}_k\left[\begin{array}{l}
				(K_{h_A}(A_i-a) - \hat{\gamma}_{h_A,a}^{(n,-k)}(L))(\pi(Z_i,L_i) - \hat{\rho}_{\pi}^{(n,-k)}(L_i))\\
				-(K_{h_A}(A_i-a) - \gamma_{h_A,a}^o(L))(\pi(Z_i,L_i) - \rho_{\pi}^o(L_i))
			\end{array}\right]\\
			\lesssim& \mathbb{E}_k\left[(\hat{\gamma}_{h_A,a}^{(n,-k)}(L) - \gamma_{h_A,a}^o(L))^2\right]
			+ \frac{1}{h_A}\mathbb{E}_k\left[(\hat{\rho}_{\pi}^{(n,-k)}(L) - \rho_{\pi}^o(L))^2\right].
			% \\=&\|\hat{\gamma}_{h_A,a}^{(n,-k)}-\gamma_{h_A,a}\|_k^2 + \|\hat{\rho}_{\pi}^{(n,-k)}(L) - \rho_{\pi}\|_k^2 
		\end{align*}
		Then, use the Markov inequality to deduce that for any $\epsilon$,
		\begin{align*}
			&\Pr(\sqrt{nh_A}|R_{8k}|> \epsilon)\leq \dfrac{nh_A}{\epsilon^2}\mathbb{E}[R_{8k}^2]\\
			\leq& \dfrac{1}{\epsilon^2}\left\{h_A\mathbb{E}\left[(\hat{\gamma}_{h_A,a}^{(n,-k)}(L) - \gamma_{h_A,a}^o(L))^2\right]
			+ \mathbb{E}\left[(\hat{\rho}_{\pi}^{(n,-k)}(L) - \rho_{\pi}^o(L))^2\right]\right\}=o(1).
		\end{align*}
		This is equivalent to saying that $R_{8k}=o(\frac{1}{\sqrt{nh_A}})$.
		Next, we control $R_{9k}$ by 
		\begin{align*}
			R_{9k}=&\mathbb{E}_k\left[
			(K_{h_A}(A-a) - \hat{\gamma}_{h_A,a}^{(n,-k)}(L))(\pi(Z,L) - \hat{\rho}_{\pi}^{(n,-k)}(L))
			\right]\\
			=& \mathbb{E}_k\left[
			(K_{h_A}(A-a) -\gamma_{h_A,a}^o(L))(\pi(Z,L) - \hat{\rho}_{\pi}^{(n,-k)}(L))\right]\\
			&+\mathbb{E}_k\left[
			(\gamma_{h_A,a}^o(L) - \hat{\gamma}_{h_A,a}^{(n,-k)}(L))(\pi(Z,L) - \hat{\rho}_{\pi}^{(n,-k)}(L))\right]\\
			=& \mathbb{E}_k\left[
			(K_{h_A}(A-a) -\gamma_{h_A,a}^o(L))(\pi(Z,L) -\rho_{\pi}^o(L))\right]\\
			&+\mathbb{E}_k\left[
			(K_{h_A}(A-a) -\gamma_{h_A,a}^o(L))(\rho_{\pi}^o(L)-\hat{\rho}_{\pi}^{(n,-k)}(L))\right]\\
			&+\mathbb{E}_k\left[
			(\gamma_{h_A,a}^o(L) - \hat{\gamma}_{h_A,a}^{(n,-k)}(L))(\rho_{\pi}^o(L) - \hat{\rho}_{\pi}^{(n,-k)}(L))\right]\\
			=& \mathbb{E}_k\left[
			(K_{h_A}(A-a) -\gamma_{h_A,a}^o(L))(\pi(Z,L) -\rho_{\pi}^o(L))\right]\\
			&+\mathbb{E}_k\left[
			(\gamma_{h_A,a}^o(L) - \hat{\gamma}_{h_A,a}^{(n,-k)}(L))(\rho_{\pi}^o(L) - \hat{\rho}_{\pi}^{(n,-k)}(L))\right]\\
			=&\mathbb{E}_k\left[ K_{h_A}(A-a)(\pi(Z,L) -\rho_{\pi}^o(L))\right] + o_p(\dfrac{1}{\sqrt{nh_A}}).
		\end{align*}
		These follow by the fact that $\gamma_{h_A,a}^o(L) = \mathbb{E}[K_{h_A}(A-a)|L]$, $\rho_{\pi}^o(L)=\mathbb{E}[\pi(Z,L)|L]$, and the Cauchy--Schwarz inequality that
		\begin{align*}
			&\mathbb{E}_k\left[
			(\gamma_{h_A,a}^o(L) - \hat{\gamma}_{h_A,a}^{(n,-k)}(L))(\rho_{\pi}^o(L) - \hat{\rho}_{\pi}^{(n,-k)}(L))\right]\\
			\leq& \|\hat{\gamma}_{h_A,a}^{(n,-k)}-\gamma_{h_A,a}^o\|_k\times 
			\|\hat{\rho}_{\pi}^{(n,-k)} - \rho_{\pi}^o\|_k = o_p(\dfrac{1}{\sqrt{nh_A}}).
		\end{align*}
		Next, we observe that
		\begin{align*}
			&\mathbb{E}\left[ K_{h_A}(A-a)(\pi(Z,L) -\rho_{\pi}^o(L))\right]-
			\zeta_{\pi}^o(a,l)\\
			=&\mathbb{E}\left[ K_{h_A}(A-a)(\pi(Z,L) -\rho_{\pi}^o(L))\right]-
			\mathbb{E}\left[ p(a|Z,L)(\pi(Z,L) -\rho_{\pi}^o(L))\right]\\
			=&\mathbb{E}\left[ \left\{\int \dfrac{1}{h_A}K(\dfrac{x-a}{h_A})p(x|Z,L)\mathrm{d}x - p(a|Z,L)\right\} (\pi(Z,L) -\rho_{\pi}^o(L))\right]\\
			=&\mathbb{E}\left[ \int K(x)\{p(a+h_Ax|Z,L)- p(a|Z,L)\}\mathrm{d}x  (\pi(Z,L) -\rho_{\pi}^o(L))\right]\\
			=&\mathbb{E}\left[ \int K(x)\frac{1}{2}p''(a_x|Z,L)h_A^2x^2\mathrm{d}x  (\pi(Z,L) -\rho_{\pi}^o(L))\right]\\
			=&\frac{h_A^2}{2}\mathbb{E}\left[ p''(a|Z,L)  (\pi(Z,L) -\rho_{\pi}^o(L))\right]\int K(x)x^2\mathrm{d}x + o(h_A^2)
			=\mathrm{bias}_{\zeta}(a,l) + o(h_A^2).
		\end{align*}
		Here $a_x$ is an intermediate value between $a$ and $a+xh_A$. Thus, we have verified that $$R_{9k}=\zeta_{\pi}^o(a,l)+\mathrm{bias}_{\zeta}(a,l)+o(h_A^2) + o_p(\frac{1}{\sqrt{nh_A}}).$$
		
		Now, we have shown that
		\begin{align*}
			&\tilde{\zeta}_{\pi,h_A}^{(n)}(a,l) - \zeta_{\pi}^o(a,l) -\mathrm{bias}_{\zeta}(a,l)\\
			=&\{\mathbb{E}_{n}-\mathbb{E}\} \left[
			(K_{h_A}(A-a) - \gamma_{h_A,a}^o(L))(\pi(Z,L) - \rho_{\pi}^o(L))\right]+
			o_p(\frac{1}{\sqrt{nh_A}}) +o(h_A^2).
		\end{align*}
		Now it is easy to apply the Lindeberg central limit theorem again to deduce that 
		\begin{align*}
			&\sqrt{nh_A}\{\mathbb{E}_{n}-\mathbb{E}\} \left[
			(K_{h_A}(A-a) - \gamma_{h_A,a}^o(L))(\pi(Z,L) - \rho_{\pi}^o(L))\right]
			\overset{d}{\to} N(0, \sigma_{\pi,\zeta}^*(a,l)^2),\text{ where }\\
			&\sigma_{\pi,\zeta}^*(a,l)^2:=\lim_{h_A\to 0^+} h_A\mathrm{Var}\{(K_{h_A}(A-a) - \gamma_{h_A,a}^o(L))(\pi(Z,L) - \rho_{\pi}^o(L))\}\\
			=&\lim_{h_A\to 0^+}
			h_A\mathbb{E}[K_{h_A}(A-a)^2 (\pi(Z,L) - \rho_{\pi}^o(L))^2] - h_A\mathbb{E}[K_{h_A}(A-a) (\pi(Z,L) - \rho_{\pi}^o(L))]^2\\
			=&\lim_{h_A\to 0^+}
			\mathbb{E}[\int K(x)^2p(a+xh_A|Z,L)\mathrm{d}x (\pi(Z,L) - \rho_{\pi}^o(L))^2] - h_A\mathbb{E}\left[ p(a|L)\kappa_\pi^o(a,L)\right]^2\\
			=&\int K(x)^2\mathrm{d}x\mathbb{E}[p(a|Z,L) (\pi(Z,L) - \rho_{\pi}^o(L))^2].
		\end{align*}
		As a result, when $h_A=O(n^{-1/5})$, we can show that
		\begin{align*}
			\sqrt{nh_A}\left\{\tilde{\zeta}_{\pi,h_A}^{(n)}(a,l) - \zeta_{\pi}^o(a,l) -\mathrm{bias}_{\zeta}(a,l)\right\}\overset{d}{\to} N(0,\sigma_{\pi,\zeta}^*(a,l)^2).
		\end{align*}
		This completes the proof of Theorem~\ref{thm: RWF testing 1} for the case $L=\emptyset$. When $L$ is discrete, the result follows from an analogous argument.

		\subsection{Proof of Theorem~\ref{thm: RWF testing 2}}
		We complete the proof in the following steps.
		
		\paragraph{Bias decomposition.}
		We first prove the pointwise convergence and asymptotic normality for a fixed
		$l\in\mathcal{L}$. Recall that
		$g(l)=\left[1,l\right]^{\top}$
		and
		$
		S_{\pi,h_A}^o(a,O_i)=\{K_{h_A}(A_i-a)-\gamma_{h_A,a}^o(L_i)\}
		\{\pi(Z_i,L_i)-\rho_\pi^o(L_i)\}.
		$
		The local linear estimator can be written as
		\begin{align*}
			\hat{\zeta}_{\pi,h}^{(n)}(a,l)
			={}&
			e_1^{\top}
			\left\{
			\frac{1}{n}
			\sum_{j=1}^n
			K_{h_L}(L_j-l)g\left(\dfrac{L_j-l}{h_L}\right)g\left(\dfrac{L_j-l}{h_L}\right)^{\top}
			\right\}^{-1}
			\frac{1}{n}
			\sum_{i=1}^n
			K_{h_L}(L_i-l)g\left(\dfrac{L_i-l}{h_L}\right)
			\hat{S}_{\pi,h_A}(a,O_i).
		\end{align*}
		Since
		$
		g\left(\dfrac{L_i-l}{h_L}\right)^{\top}
		\begin{bmatrix}
			\zeta_{\pi}^o(a,l)\\
			h_L\nabla_l\zeta_{\pi}^o(a,l)
		\end{bmatrix}
		=
		\zeta_{\pi}^o(a,l)
		+
		(L_i-l)\nabla_l\zeta_{\pi}^o(a,l),
		$
		we have
		\begin{align*}
			&\hat{\zeta}_{\pi,h}^{(n)}(a,l)-\zeta_{\pi}^o(a,l)
			=e_1^{\top}
			\left\{
			\frac{1}{n}
			\sum_{j=1}^n
			K_{h_L}(L_j-l)g\left(\dfrac{L_j-l}{h_L}\right)g\left(\dfrac{L_j-l}{h_L}\right)^{\top}
			\right\}^{-1}\\
			&\frac{1}{n}
			\sum_{i=1}^n
			K_{h_L}(L_i-l)g\left(\dfrac{L_i-l}{h_L}\right)\times
			\left\{\begin{array}{l}
				\hat{S}_{\pi,h_A}(a,O_i)-S_{\pi,h_A}^o(a,O_i)\\
				+S_{\pi,h_A}^o(a,O_i)
				-\mathbb{E}\{S_{\pi,h_A}^o(a,O_i)\mid L_i\}\\
				+\mathbb{E}\{S_{\pi,h_A}^o(a,O_i)\mid L_i\}
				-\zeta_{\pi}^o(a,L_i)\\
				+\zeta_{\pi}^o(a,L_i)-
				\zeta_{\pi}^o(a,l)-
				(L_i-l)\nabla_l\zeta_{\pi}^o(a,l)
			\end{array}\right\}.
		\end{align*}
		By standard kernel arguments and the symmetry of $K$,
		\begin{align*}
			\frac{1}{n}
			\sum_{j=1}^n
			K_{h_L}(L_j-l)g\left(\dfrac{L_j-l}{h_L}\right)g\left(\dfrac{L_j-l}{h_L}\right)^{\top}
			=
			p_L(l)
			\begin{bmatrix}
				1&0\\
				0&\int u^2K(u)\,\mathrm{d}u
			\end{bmatrix}
			+o_p(1).
		\end{align*}
		Therefore, we have the decomposition that
		\begin{align*}
			\hat{\zeta}_{\pi,h}^{(n)}(a,l)-\zeta_{\pi}^o(a,l)
			=&
			\frac{1+o_p(1)}{np_L(l)}
			\sum_{i=1}^n
			K_{h_L}(L_i-l)
			\left\{\begin{array}{l}
				\hat{S}_{\pi,h_A}(a,O_i)-S_{\pi,h_A}^o(a,O_i)\\
				+S_{\pi,h_A}^o(a,O_i)
				-\mathbb{E}\{S_{\pi,h_A}^o(a,O_i)\mid L_i\}\\
				+\mathbb{E}\{S_{\pi,h_A}^o(a,O_i)\mid L_i\}
				-\zeta_{\pi}^o(a,L_i)\\
				+\zeta_{\pi}^o(a,L_i)-
				\zeta_{\pi}^o(a,l)-
				(L_i-l)\nabla_l\zeta_{\pi}^o(a,l)
			\end{array}\right\}\\
			=&\frac{1+o_p(1)}{p_L(l)}\left\{
			E_{11}(l) + E_{12}(l) + E_{13}(l) + E_{14}(l)
			\right\}.
		\end{align*}
		We next consider the effect of nuisance estimation. Specifically,
		\begin{align*}
			&E_{11}(l)=\dfrac{1}{n}\displaystyle\sum_{i=1}^n K_{h_L}(L_i - l) \left\{\hat{S}_{\pi,h_A}(a,O_i)-
			S_{\pi,h_A}^o(a,O_i)\right\}\\
			=&\dfrac{1}{n}\displaystyle\sum_{i=1}^n K_{h_L}(L_i - l)
			\left\{\hat{S}_{\pi,h_A}(a,O_i)-S_{\pi,h_A}^o(a,O_i)
			-\mathbb{E}_{-k}[\hat{S}_{\pi,h_A}(a,O_i)-S_{\pi,h_A}^o(a,O_i)\mid L_i]
			\right\}\\
			&+\dfrac{1}{n}\displaystyle\sum_{i=1}^n K_{h_L}(L_i - l)
			\left\{\mathbb{E}_{-k}[\hat{S}_{\pi,h_A}(a,O_i)\mid L_i]
			-\mathbb{E}_{-k}[S_{\pi,h_A}^o(a,O_i)\mid L_i]\right\}\\
			=&O_p\left(\dfrac{\sum_{k=1}^K\mathbb{E}[\|\hat\rho_{\pi}^{(n,-k)} - \rho_{\pi}^o\|_{\infty}^2]}{\sqrt{nh_Ah_L}}
			+\dfrac{\sum_{k=1}^K\mathbb{E}[\|\hat\gamma_{h_A,a}^{(n,-k)}- \gamma_{h_A,a}^o\|_{\infty}^2]}{\sqrt{nh_L}}
			\right)\\&+ 
			\dfrac{1}{n}\displaystyle\sum_{i=1}^n K_{h_L}(L_i - l)
			\{\hat\gamma_{h_A,a}^{(n,-k_i)}(L_i) - \gamma_{h_A,a}^o(L_i)\}\{\hat\rho_{\pi}^{(n,-k_i)}(L_i) - \rho_{\pi}^o(L_i)\}\\
			\lesssim&o_p(\dfrac{1}{\sqrt{nh_Ah_L}})+
			\dfrac{1}{n}\displaystyle\sum_{i=1}^n K_{h_L}(L_i - l)
			\|\hat\gamma_{h_A,a}^{(n,-k_i)}- \gamma_{h_A,a}^o\|_{\infty}
			\times\|\hat\rho_{\pi}^{(n,-k_i)} - \rho_{\pi}^o\|_{\infty}\\
			=&o_p(\dfrac{1}{\sqrt{nh_Ah_L}}) + O_p(\|\hat\gamma_{h_A,a}^{(n,-k_i)}- \gamma_{h_A,a}^o\|_{\infty}
			\times\|\hat\rho_{\pi}^{(n,-k_i)} - \rho_{\pi}^o\|_{\infty})
			=o_p(\dfrac{1}{\sqrt{nh_Ah_L}}).
		\end{align*}
		The third equality holds since for each $k$,
		\begin{align*}
			&\mathrm{Var}\left\{
			\dfrac{1}{n_K}\displaystyle\sum_{i\in I_{k}} K_{h_L}(L_i - l)
			\left(\hat{S}_{\pi,h_A}(a,O_i)-S_{\pi,h_A}^o(a,O_i)\right)\middle| I_{-k}
			\right\}\\
			\lesssim&
			\dfrac{1}{n_K}\mathbb{E}\left[K_{h_L}(L_i - l)^2
			\left(\hat{S}_{\pi,h_A}(a,O_i)-S_{\pi,h_A}^o(a,O_i)\right)^2
			\middle| I_{-k}\right]\\
			\lesssim&
			\dfrac{1}{n_K}\mathbb{E}\left[K_{h_L}(L_i - l)^2
			K_{h_A}(A_i - a)^2\{\hat\rho_{\pi}^{(n,-k)}(L_i) - \rho_{\pi}^o(L_i)\}^2
			\middle| I_{-k}\right]\\
			&+\dfrac{1}{n_K}\mathbb{E}\left[K_{h_L}(L_i - l)^2
			\{\hat\gamma_{h_A,a}^{(n,-k)}(L_i)- \gamma_{h_A,a}^o(L_i)\}^2
			\{\pi(Z_i,L_i) - \rho_{\pi}^o(L_i)\}^2
			\middle| I_{-k}\right]\\
			\lesssim&\dfrac{\|\hat\rho_{\pi}^{(n,-k)} - \rho_{\pi}^o\|_{\infty}^2}{nh_Ah_L} + 
			\dfrac{\|\hat\gamma_{h_A,a}^{(n,-k)}- \gamma_{h_A,a}^o\|_{\infty}^2}{nh_L}.
		\end{align*}
		Then we can apply the Chebyshev inequality to deduce that
		\begin{align*}
			&\dfrac{1}{n}\displaystyle\sum_{i\in I_{k}} K_{h_L}(L_i - l)
			\left(\hat{S}_{\pi,h_A}(a,O_i)-S_{\pi,h_A}^o(a,O_i)-\mathbb{E}[\hat{S}_{\pi,h_A}(a,O_i)-S_{\pi,h_A}^o(a,O_i)|L_i]\right)\\
			=&O_p\left(\sqrt{\dfrac{\|\hat\rho_{\pi}^{(n,-k)} - \rho_{\pi}^o\|_{\infty}^2}{nh_Ah_L} + 
				\dfrac{\|\hat\gamma_{h_A,a}^{(n,-k)}- \gamma_{h_A,a}^o\|_{\infty}^2}{nh_L}}\right).
		\end{align*}

		We now examine the smoothing bias with respect to $A$. Notice that
		\begin{align*}
			\mathbb{E}\{S_{\pi,h_A}^o(a,O_i)\mid L=l\}
			=
			\mathbb{E}\left[
			K_{h_A}(A-a)
			\{\pi(Z,L)-\rho_\pi^o(L)\}
			\mid L=l
			\right],
		\end{align*}
		because $\mathbb{E}\{\pi(Z,L)-\rho_\pi^o(L)\mid L\}=0.$
		Furthermore,
		\begin{align*}
			&\mathbb{E}\left[
			K_{h_A}(A-a)
			\{\pi(Z,L)-\rho_\pi^o(L)\}
			\mid L=l
			\right]\\
			={}&
			\iint
			\frac{1}{h_A}
			K\left(\frac{s-a}{h_A}\right)
			\{\pi(z,l)-\rho_\pi^o(l)\}
			p_{A,Z\mid L}(s,z\mid l)
			\,\mathrm{d}s\,\mathrm{d}z\\
			={}&
			\iint
			K(u)
			\{\pi(z,l)-\rho_\pi^o(l)\}
			p_{A,Z\mid L}(a+h_Au,z\mid l)
			\,\mathrm{d}u\,\mathrm{d}z\\
			={}&
			\iint
			K(u)
			\{\pi(z,l)-\rho_\pi^o(l)\}
			\{p_{A,Z\mid L}(a+h_Au,z\mid l)-p_{A,Z\mid L}(a,z\mid l)\}
			\,\mathrm{d}u\,\mathrm{d}z + \zeta_{\pi}^o(a,l).
		\end{align*}
		By a second-order Taylor expansion and the symmetry of $K$,
		$\mathbb{E}\{S_{\pi,h_A}^o(a,O_i)\mid L=l\}-\zeta_{\pi}^o(a,l)=O(h_A^2).$
		It follows that
		\begin{align*}
			E_{13}(l)=&\frac{1}{n}
			\sum_{i=1}^n
			K_{h_L}(L_i-l)
			\left\{\mathbb{E}\{S_{\pi,h_A}^o(a,O_i)\mid L_i\}
			-\zeta_{\pi}^o(a,L_i)\right\}
			=O_p(h_A^2).
		\end{align*}
		
		Similarly, a second-order Taylor expansion with respect to $l$ gives
		\begin{align*}
			\zeta_{\pi}^o(a,L_i)
			={}&
			\zeta_{\pi}^o(a,l)
			+
			(L_i-l)\nabla_l\zeta_{\pi}^o(a,l)+
			\frac{1}{2}(L_i-l)^2
			\nabla_l^2\zeta_{\pi}^o(a,\tilde{L}_i),
		\end{align*}
		where $\tilde{L}_i$ lies between $L_i$ and $l$. Therefore,
		\begin{align*}
			E_{14}(l)=&\frac{1}{n}
			\sum_{i=1}^n
			K_{h_L}(L_i-l)
			\left[
			\zeta_{\pi}^o(a,L_i)
			-
			\zeta_{\pi}^o(a,l)
			-
			(L_i-l)\nabla_l\zeta_{\pi}^o(a,l)
			\right]
			=
			O_p(h_L^2).
		\end{align*}
		Combining these facts together,	we obtain for all fixed $a\in\mathring{\mathcal{A}}$ and $l\in\mathcal{L}$,
		\begin{align}
			&\hat{\zeta}_{\pi,h}^{(n)}(a,l)-\zeta_{\pi}^o(a,l)
			=\dfrac{1+o_p(1)}{p_L(l)}E_{12}(l)\notag\\
			=&\frac{1}{np_L(l)}
			\sum_{i=1}^n
			K_{h_L}(L_i-l)
			\left\{S_{\pi,h_A}^o(a,O_i)-\mathbb{E}[S_{\pi,h_A}^o(a,O_i)\mid L_i]\right\}\notag\\&+
			O_p(h_A^2+h_L^2)
			+
			o_p\left(
			\frac{1}{\sqrt{nh_Ah_L}}
			\right).
			\label{eq: asymptotic representation}
		\end{align}
		
		\paragraph{Calculate the asymptotic distribution.}
		We next calculate the variance of the leading term. 
		Since the summands are independent and have mean zero,
		\begin{align*}
			&\operatorname{Var}
			\left[
			\frac{1}{n}
			\sum_{i=1}^n
			K_{h_L}(L_i-l)
			\left\{
			S_{\pi,h_A}^o(a,O_i)
			-
			\mathbb{E}[S_{\pi,h_A}^o(a,O_i)\mid L_i]
			\right\}
			\right]
			\nonumber\\
			={}&
			\frac{1}{n}
			\mathbb{E}
			\left[
			K_{h_L}(L-l)^2
			\operatorname{Var}
			\{S_{\pi,h_A}^o(a,O_i)\mid L\}
			\right].
		\end{align*}
		Moreover,
		\begin{align*}
			\operatorname{Var}
			\{S_{\pi,h_A}^o(a,O_i)\mid L\}
			={}&
			\mathbb{E}
			\left[
			K_{h_A}(A-a)^2
			\{\pi(Z,L)-\rho_\pi^o(L)\}^2
			\mid L
			\right]
			+
			o(1).
		\end{align*}
		The $o(1)$ terms are negligible compared with the leading term, which
		is of order $h_A^{-1}$. Hence,
		\begin{align*}
			&\mathbb{E}
			\left[
			K_{h_L}(L-l)^2
			\operatorname{Var}
			\{S_{\pi,h_A}^o(a,O_i)\mid L\}
			\right]\\
			={}&
			\{1+o(1)\}
			\mathbb{E}
			\left[
			K_{h_L}(L-l)^2
			K_{h_A}(A-a)^2
			\{\pi(Z,L)-\rho_\pi^o(L)\}^2
			\right].
		\end{align*}
		By the changes of variables $L=l+h_Lt,\;A=a+h_As,$ we obtain
		\begin{align*}
			&\mathbb{E}
			\left[
			K_{h_L}(L-l)^2
			K_{h_A}(A-a)^2
			\{\pi(Z,L)-\rho_\pi^o(L)\}^2
			\right]
			\nonumber\\
			={}&
			\frac{1+o(1)}{h_Ah_L}
			\left\{\int K(u)^2\,\mathrm{d}u\right\}^2
			\int
			\{\pi(z,l)-\rho_\pi^o(l)\}^2
			p_{A,Z,L}(a,z,l)\,\mathrm{d}z
			\nonumber\\
			={}&
			\frac{1+o(1)}{h_Ah_L}
			\left\{\int K(u)^2\,\mathrm{d}u\right\}^2
			p_{A,L}(a,l)
			\mathbb{E}
			\left[
			\{\pi(Z,L)-\rho_\pi^o(L)\}^2
			\mid A=a,L=l
			\right].
		\end{align*}
		Consequently, $\mathrm{Var}
		\left[
		\sqrt{nh_Ah_L}
		\left\{
		\hat{\zeta}_{\pi,h}^{(n)}(a,l)
		-\zeta_{\pi}^o(a,l)
		\right\}
		\right]\to
		\sigma_{\pi,\zeta}^o(a,l)^2.$
		This also shows that
		\[
		\hat{\zeta}_{\pi,h}^{(n)}(a,l)-\zeta_{\pi}^o(a,l)
		=
		O_p\left(
		\frac{1}{\sqrt{nh_Ah_L}}
		+
		h_A^2+h_L^2
		\right).
		\]
		
		To establish asymptotic normality, define
		\begin{align*}
			X_{n,i}:=
			\frac{\sqrt{h_Ah_L}}{\sqrt{n} p_L(l)}
			K_{h_L}(L_i-l)
			\left[
			S_{\pi,h_A}^o(a,O_i)
			-
			\mathbb{E}\{S_{\pi,h_A}^o(a,O_i)\mid L_i\}
			\right].
		\end{align*}
		Then $
		%	\sqrt{\frac{1}{n}}
		\sum_{i=1}^n X_{n,i}$
		is the normalized leading term in
		\eqref{eq: asymptotic representation}. Suppose that, for some
		$\delta>0$,
		\[
		\sup_{a,l}
		\mathbb{E}
		\left[
		|\pi(Z,L)-\rho_\pi^o(L)|^{2+\delta}
		\mid A=a,L=l
		\right]
		<\infty.
		\]
		Standard kernel moment calculations yield
		\begin{align*}
			\sum_{i=1}^n\mathbb{E}|X_{n,i}|^{2+\delta}
			=O\left\{(nh_Ah_L)^{-\delta/2}\right\}=o(1).
		\end{align*}
		Thus, the Lyapunov condition, and hence the Lindeberg condition, holds.
		By the Lindeberg--Feller central limit theorem,
		$\sum_{i=1}^n X_{n,i}
		\overset{d}{\to}
		N\{0,\sigma_{\pi,\zeta}^o(a,l)^2\}.$
		If $h_A,h_L=o(n^{-1/6}),$ then $\sqrt{nh_Ah_L}(h_A^2+h_L^2)=o(1)$.
		It follows from
		\eqref{eq: asymptotic representation} and Slutsky's theorem that
		\begin{align*}
			\sqrt{nh_Ah_L}
			\left\{
			\hat{\zeta}_{\pi,h}^{(n)}(a,l)
			-
			\zeta_{\pi}^o(a,l)
			\right\}
			\overset{d}{\to}
			N\{0,\sigma_{\pi,\zeta}^o(a,l)^2\}.
		\end{align*}
		
		\subsection{Proof of Theorem~\ref{thm: RWF testing 3}}
		We next establish the simultaneous result. First, we need to show that 
		$$\sup_{l\in\mathcal{L}}\{|E_{11}(l)|\},\;
		\sup_{l\in\mathcal{L}}\{|E_{13}(l)|\},\;
		\sup_{l\in\mathcal{L}}\{|E_{14}(l)|\}
		=o_p\left(\dfrac{1}{\sqrt{nh_Ah_L\log(n)}}\right).$$
		To prove these three facts, we need the following two lemmas.
		\begin{lemma}[Talagrand Inequality]\label{lem: talagrand inequality}
			Assume that $\{O_i\}_{i=1}^n$ are independent random variables, and let $\mathcal{F}$ be a uniformly bounded subspace of $L_2(O)$.
			With probability $1-\epsilon_n$,
			\begin{align*}
				&\sup_{f\in\mathcal{F}}\dfrac{1}{n}\sum_{i=1}^n f(O_i) \lesssim \mathbb{E}\left[\sup_{f\in\mathcal{F}}\dfrac{1}{n}\sum_{i=1}^n f(O_i)\right] + 
				\sqrt{\sup_{f\in\mathcal{F}}\mathbb{E}\left[\frac{1}{n}\sum_{i=1}^nf(O_i)^2\right]\dfrac{\log(1/\epsilon_n)}{n}} +
				\|\mathcal{F}\|_{\infty}	\dfrac{\log(1/\epsilon_n)}{n}.
			\end{align*}
		\end{lemma}
		\begin{lemma}[Maximal inequality for VC-type classes]
			\label{lem: VC maximal inequality}
			Let $O_1,\ldots,O_n$ be i.i.d. with distribution $P$, and let
			$\mathcal{F}$ be a pointwise measurable functional class with envelope $F$ and finite VC dimension. Suppose that
			$
			\mathbb{E}[f(O)]=0,\;
			\|F\|_\infty\leq b,\;
			\sup_{f\in\mathcal{F}}\mathbb{E}[f(O)^2]\leq\sigma^2\leq b^2.
			$
			Then
			\[
			\mathbb{E}\left[
			\sup_{f\in\mathcal{F}}
			\left|
			\frac{1}{\sqrt n}\sum_{i=1}^n f(O_i)
			\right|
			\right]
			\lesssim
			\sigma\sqrt{\log\left(\frac{b}{\sigma}\right)}
			+\frac{b}{\sqrt n}\log\left(\frac{b}{\sigma}\right).
			\]
		\end{lemma}
		Applying the Lemma~\ref{lem: talagrand inequality}, we know that with probability $1-\epsilon_n$,
		\begin{align*}
			&\sup_{l\in\mathcal{L}}\{|E_{11}(l)|\}=\sup_{l\in\mathcal{L}}\left\{
			\dfrac{1}{n}\displaystyle\sum_{i=1}^n K_{h_L}(L_i - l) \mathbb{E}_{-k}\left[\hat{S}_{\pi,h_A}(a,O_i)-S_{\pi,h_A}^o(a,O_i)\middle| L_i\right]\right\}\\
			&+\sup_{l\in\mathcal{L}}\left\{
			\dfrac{1}{n}\displaystyle\sum_{i=1}^n K_{h_L}(L_i - l) \left(\hat{S}_{\pi,h_A}(a,O_i)-S_{\pi,h_A}^o(a,O_i)-\mathbb{E}_{-k}\left[\hat{S}_{\pi,h_A}(a,O_i)-S_{\pi,h_A}^o(a,O_i)\middle| L_i\right]\right)\right\}\\
			\lesssim& \sup_{l\in\mathcal{L}}\left\{
			\dfrac{1}{n}\displaystyle\sum_{i=1}^n K_{h_L}(L_i - l) \right\}
			\|\hat{\gamma}_{h_A,a}^{(n,-k)}-\gamma_{h_A,a}^o\|_\infty \, \|\hat{\rho}_{\pi}^{(n,-k)}-\rho_{\pi}^o\|_\infty\\
			&+\sum_{k=1}^K\mathbb{E}\left[
			\sup_{l\in\mathcal{L}}\left\{
			\dfrac{1}{n}\displaystyle\sum_{i\in I_k} K_{h_L}(L_i - l) \left(\begin{array}{l}
				\hat{S}_{\pi,h_A}(a,O_i)-S_{\pi,h_A}^o(a,O_i)\\
				-\mathbb{E}_{-k}\left[\hat{S}_{\pi,h_A}(a,O_i)-S_{\pi,h_A}^o(a,O_i)\middle| L_i\right]
			\end{array}\right)\right\}
			\right]\\
			&+\sqrt{\sup_{l\in\mathcal{L}}\mathbb{E}\left[K_{h_L}(L_i - l)^2 \left(\hat{S}_{\pi,h_A}(a,O_i)-S_{\pi,h_A}^o(a,O_i)\right)^2\right]
				\dfrac{\log(1/\epsilon_n)}{n}} + \dfrac{\log(1/\epsilon_n)}{nh_Ah_L}\\
			\lesssim&
			O_p(\|\hat{\gamma}_{h_A,a}^{(n,-k)}-\gamma_{h_A,a}^o\|_\infty \, \|\hat{\rho}_{\pi}^{(n,-k)}-\rho_{\pi}^o\|_\infty)\\
			&+\sum_{k=1}^K\mathbb{E}\left[
			\sup_{l\in\mathcal{L}}\left\{
			\dfrac{1}{n}\displaystyle\sum_{i\in I_k} K_{h_L}(L_i - l) \left(\begin{array}{l}
				\hat{S}_{\pi,h_A}(a,O_i)-S_{\pi,h_A}^o(a,O_i)\\
				-\mathbb{E}_{-k}\left[\hat{S}_{\pi,h_A}(a,O_i)-S_{\pi,h_A}^o(a,O_i)\middle| L_i\right]
			\end{array}\right)\right\}
			\right]\\
			&+\sqrt{\left\{\dfrac{\|\hat\rho_{\pi}^{(n,-k)} - \rho_{\pi}^o\|_{\infty}^2}{h_Ah_L} + 
				\dfrac{\|\hat\gamma_{h_A,a}^{(n,-k)}- \gamma_{h_A,a}^o\|_{\infty}^2}{h_L}\right\}
				\dfrac{\log(1/\epsilon_n)}{n}} 
			+\dfrac{\log(1/\epsilon_n)}{nh_Ah_L}\\
			\lesssim&o_p(\dfrac{1}{\sqrt{nh_Ah_L\log(n)}})
			+\sqrt{\dfrac{o_p(1)}{h_Ah_L(\log n)^2}
				\dfrac{\log(1/\epsilon_n)}{n}} 
			+\dfrac{\log(1/\epsilon_n)}{nh_Ah_L}.
		\end{align*}
		The final step follows from Lemma~\ref{lem: VC maximal inequality} by setting $b=\frac{1}{h_Ah_L}$ and $$\sigma^2 = \dfrac{\log(1/h_L)}{h_Ah_L} 
		\left\{\|\hat\rho_{\pi}^{(n,-k)}-\rho_{\pi}^o\|_{\infty}^2 + 
		h_A\|\hat\gamma_{h_A,a}^{(n,-k)}-\gamma_{h_A,a}^o\|_{\infty}^2\right\},$$
		deriving that  
		\begin{align*}
			&\mathbb{E}\left[
			\sup_{l\in\mathcal{L}}\left\{
			\dfrac{1}{n}\displaystyle\sum_{i\in I_k} K_{h_L}(L_i - l) \left(\begin{array}{l}
				\hat{S}_{\pi,h_A}(a,O_i)-S_{\pi,h_A}^o(a,O_i)\\
				-\mathbb{E}_{-k}\left[\hat{S}_{\pi,h_A}(a,O_i)-S_{\pi,h_A}^o(a,O_i)\middle| L_i\right]
			\end{array}\right)\right\}
			\right]\\
			\lesssim&\mathbb{E}\left[\dfrac{1}{\sqrt{n}}
			\left\{\sigma\sqrt{\log\left(\frac{b}{\sigma}\right)}
			+\frac{b}{\sqrt n}\log\left(\frac{b}{\sigma}\right)\right\}\right]\\
			\lesssim&O\left(\sqrt{\dfrac{\log(n)\log(1/h_L)}{n}
				\left\{\dfrac{\mathbb{E}[\|\hat\rho_{\pi}^{(n,-k)} - \rho_{\pi}^o\|_{\infty}^2]}{h_Ah_L} + 
				\dfrac{\mathbb{E}[\|\hat\gamma_{h_A,a}^{(n,-k)}- \gamma_{h_A,a}^o\|_{\infty}^2]}{h_L}\right\}
			}\right)\\
			=&O\left(\sqrt{\dfrac{1}{nh_Ah_L\log(n)}}\right),
		\end{align*}
		since $\mathbb{E}[\|\hat\rho_{\pi}^{(n,-k)} - \rho_{\pi}^o\|_{\infty}^2]=o(h_A/\log(n)^3)$, and 
		$\mathbb{E}[\|\hat\gamma_{h_A,a}^{(n,-k)}- \gamma_{h_A,a}^o\|_{\infty}^2] = o(\log(n)^{-3})$.
		This implies that $\sup_{l\in\mathcal{L}}|E_{11}(l)| = o_p(\dfrac{1}{\sqrt{nh_Ah_L\log(n)}})$.
		Furthermore, by the fact that $h_A,h_L=O(n^{-c})$ for some $c>1/6$, one can analogously prove that
		\begin{align*}
			&\sup_{l\in\mathcal{L}}|E_{13}(l)| = O_p(h_A^2) =o_p(\dfrac{1}{\sqrt{nh_Ah_L\log(n)}}),\\
			&\sup_{l\in\mathcal{L}}|E_{14}(l)| = O_p(h_L^2) =o_p(\dfrac{1}{\sqrt{nh_Ah_L\log(n)}}).
		\end{align*}
		Now, by Equation~\eqref{eq: asymptotic representation}, we have proved that
		\begin{align*}
			\sup_{l\in\mathcal{L}}
			\left|\hat{\zeta}_{\pi,h}^{(n)}(a,l)-\zeta_{\pi}^o(a,l) - \dfrac{E_{12}(l)}{p_L(l)}\right|
			\lesssim \sum_{j=1,3,4} \sup_{l\in\mathcal{L}} |E_{1j}(l)|=
			o_p(\dfrac{1}{\sqrt{nh_Ah_L\log(n)}}).
		\end{align*}
		Following the same arguments as in the proof of Theorem~3 of \citet{takatsu2023debiasedinferencecovariateadjustedregression}, one  can verify that
		\begin{align*}
			&\sup_{l\in\mathcal{L}}\left|\dfrac{E_{12}(l)}{p_L(l)}\right| = \sup_{l\in\mathcal{L}}\left|\frac{1}{np_L(l)}
			\sum_{i=1}^n
			K_{h_L}(L_i-l)
			\left\{S_{\pi,h_A}^o(a,O_i)-\mathbb{E}[S_{\pi,h_A}^o(a,O_i)\mid L_i]\right\}\right| 
			= O_p(\sqrt{\dfrac{\log (n)}{nh_Ah_L}}).
		\end{align*}
		This indicates that $\sup_{l\in\mathcal{L}}\left|\hat{\zeta}_{\pi,h}^{(n)}(a,l)-\zeta_{\pi}^o(a,l)\right|=O_p(\sqrt{\dfrac{\log (n)}{nh_Ah_L}})$.
		Next, we decompose the discrepancy between the distribution of the estimated maximum of empirical process and that of its Gaussian bootstrap analogue as
		\begin{align*}
			&\sup_{x\in\mathbb{R}}\left|
			\Pr\left(\sup_{l\in\mathcal{L}}\sqrt{nh_Ah_L}
			\left|\frac{\hat{\zeta}_{\pi,h}^{(n)}(a,l)-\zeta_{\pi}^o(a,l)}
			{\hat{\sigma}_{\pi,\zeta,h}^{(n)}(a,l)}
			\right|\leq x\right)
			-\Pr\left(\max_{t=1,\ldots,T}|\hat{d}_{b,l_t}|\leq x
			\middle|\{O_i\}_{i=1}^n\right)\right|\\
			\leq&\sup_{x\in\mathbb{R}}\left|
			\Pr\left(\sup_{l\in\mathcal{L}}\sqrt{nh_Ah_L}
			\left|\frac{\hat{\zeta}_{\pi,h}^{(n)}(a,l)-\zeta_{\pi}^o(a,l)}
			{\hat{\sigma}_{\pi,\zeta,h}^{(n)}(a,l)}
			\right|\leq x\right)
			-\Pr\left(\sup_{l\in\mathcal{L}}|d_{b,l}^o|\leq x\right)\right|\\
			&+\sup_{x\in\mathbb{R}}\left|\Pr\left(\sup_{l\in\mathcal{L}}|d_{b,l}^o|\leq x\right) - \Pr\left(\sup_{t=1,\ldots,T}|d_{b,l_t}^o|\leq x\right)\right|\\
			&+\sup_{x\in\mathbb{R}}\left|\Pr\left(\sup_{t=1,\ldots,T}|d_{b,l_t}^o|\leq x\right) - \Pr\left(\max_{t=1,\ldots,T}|\hat{d}_{b,l_t}|\leq x
			\middle|\{O_i\}_{i=1}^n\right) \right|\\
			=&\sup_{x\in\mathbb{R}}|E_{21}(x)| + 
			\sup_{x\in\mathbb{R}}|E_{22}(x)| + 
			\sup_{x\in\mathbb{R}}|E_{23}(x)|.
			%		=o_p(1).
		\end{align*}
		Let $d_{b,l}^o$ be a zero-mean Gaussian process with covariance structure 
		\begin{align*}
			&\mathrm{Cov}\left\{d_{b,l}^o,d_{b,l'}^o\right\}=
			\dfrac{\Gamma_{\pi,h}^o(a,l,l')}
			{\sqrt{\Gamma_{\pi,h}^o(a,l,l)\Gamma_{\pi,h}^o(a,l',l')}},\text{ where}\\
			&\Gamma_{\pi,h}^o(a,l,l'):=\dfrac{h_Ah_L}{p_L(l)p_L(l')}\mathbb{E}\left[K_{h_L}(L-l)K_{h_L}(L-l')\left\{S_{\pi,h_A}^o(a,O)-\mathbb{E}[S_{\pi,h_A}^o(a,O)\mid L]\right\}^2\right].
		\end{align*}
		Next, we show that $\sup_{x\in\mathbb{R}}|E_{2j}(x)|=o(1)$ for $j=1,2,3$, finishing the proof of Theorem~\ref{thm: RWF testing 3}.
		
		\begin{proof}[Proof of $\sup_{x\in\mathbb{R}}|E_{21}(x)|=o(1)$:]
			Following the same argument as in the proof of Theorem~3 of \citet{takatsu2023debiasedinferencecovariateadjustedregression}, we know that for any $\epsilon>0$,
			\begin{align*}
				&\sup_{x\in\mathbb{R}}|E_{21}(x)| = 
				\sup_{x\in\mathbb{R}}\left|
				\Pr\left(\sup_{l\in\mathcal{L}}\sqrt{nh_Ah_L}
				\left|\frac{\hat{\zeta}_{\pi,h}^{(n)}(a,l)-\zeta_{\pi}^o(a,l)}
				{\hat{\sigma}_{\pi,\zeta,h}^{(n)}(a,l)}
				\right|\leq x\right)
				-\Pr\left(\sup_{l\in\mathcal{L}}|d_{b,l}^o|\leq x\right)\right|\\
				=&\sup_{x\in\mathbb{R}}\left|
				\Pr\left(\sup_{l\in\mathcal{L}}\sqrt{nh_Ah_L}
				\left|\frac{\sum_{j=1}^4 E_{1j}(l)/p_L(l)}
				{\hat{\sigma}_{\pi,\zeta,h}^{(n)}(a,l)}
				\right|\leq x\right)
				-\Pr\left(\sup_{l\in\mathcal{L}}|d_{b,l}^o|\leq x\right)\right|\\
				\lesssim& 
				\sup_{x\in\mathbb{R}}\left|
				\Pr\left(\sup_{l\in\mathcal{L}}\sqrt{nh_Ah_L}
				\left|\frac{E_{12}(l)/p_L(l)}
				{\sqrt{\Gamma_{\pi,h}^o(a,l,l)}}
				\right|\leq x\right)
				-\Pr\left(\sup_{l\in\mathcal{L}}|d_{b,l}^o|\leq x\right)\right|\\
				&+\sup_{x\in\mathbb{R}}\Pr\left(\left|\sup_{l\in\mathcal{L}}|d_{b,l}^o|-x\right|\leq
				\sum_{j=1,3,4} \left|\sqrt{nh_Ah_L}\dfrac{E_{1j}(l)}{p_L(l)\sqrt{\Gamma_{\pi,h}^o(a,l,l)}}\right|
				\right)\\
				\lesssim& 
				\sup_{x\in\mathbb{R}}\left|
				\Pr\left(\sup_{l\in\mathcal{L}}\sqrt{nh_Ah_L}
				\left|\frac{E_{12}(l)}
				{p_L(l)\sqrt{\Gamma_{\pi,h}^o(a,l,l)}}
				\right|\leq x\right)
				-\Pr\left(\sup_{l\in\mathcal{L}}|d_{b,l}^o|\leq x\right)\right|\\
				&+\sup_{x\in\mathbb{R}}\Pr\left(\left|\sup_{l\in\mathcal{L}}|d_{b,l}^o|-x\right|\leq
				\dfrac{\epsilon}{\sqrt{\log(n)}}
				\right) \\&+ 
				\Pr\left(\sum_{j=1,3,4} \sup_{l\in\mathcal{L}} \left|\sqrt{nh_Ah_L}\dfrac{E_{1j}(l)}{p_L(l)\sqrt{\Gamma_{\pi,h}^o(a,l,l)}}\right|
				\geq\dfrac{\epsilon}{\sqrt{\log(n)}}\right)\\
				\lesssim& 
				\Pr\left(\left|\sup_{l\in\mathcal{L}}\sqrt{nh_Ah_L}
				\left|\frac{E_{12}(l)}
				{p_L(l)\sqrt{\Gamma_{\pi,h}^o(a,l,l)}}\right| - 
				\sup_{l\in\mathcal{L}}\left|d_{b,l}^o\right|\right|\geq \dfrac{\epsilon}{\sqrt{\log(n)}}\right)\\
				&+\sup_{x\in\mathbb{R}}\Pr\left(\left|\sup_{l\in\mathcal{L}}|d_{b,l}^o|-x\right|\leq
				\dfrac{\epsilon}{\sqrt{\log(n)}}
				\right) \\&+ 
				\Pr\left(\sum_{j=1,3,4} \sup_{l\in\mathcal{L}} \left|\sqrt{nh_Ah_L}\dfrac{E_{1j}(l)}{p_L(l)\sqrt{\Gamma_{\pi,h}^o(a,l,l)}}\right|
				\geq\dfrac{\epsilon}{\sqrt{\log(n)}}\right)\\
				=&E_{211} + E_{212} + E_{213}.
			\end{align*}
			For any $\epsilon>0$, if we take $b=(h_Ah_L)^{-1/2}$ in Corollary~2.2 of \citet{chernozhukov2014gaussian}, 
			%	$E_{211}$ converges to zero by 
			\begin{align*}
				&\Pr\left(\left|\sup_{l\in\mathcal{L}}
				\left|\frac{\sqrt{nh_Ah_L} E_{12}(l)}
				{p_L(l)\sqrt{\Gamma_{\pi,h}^o(a,l,l)}}\right| - 
				\sup_{l\in\mathcal{L}}\left|d_{b,l}^o\right|\right|\geq 
				\dfrac{b \log(n)}{\gamma^{1/2}n^{1/2}} + 
				\dfrac{b^{1/2} \log(n)^{3/4}}{\gamma^{1/2}n^{1/4}} + 
				\dfrac{b^{1/3}\log(n)^{2/3}}{\gamma^{1/3}n^{1/6}}
				\right)\\&\lesssim \gamma + \dfrac{\log(n)}{n}.
			\end{align*}
			Thus, we need 
			$\dfrac{(h_Ah_L)^{-1/2} \log(n)}{\gamma^{1/2}n^{1/2}}
			+\dfrac{(h_Ah_L)^{-1/4} \log(n)^{3/4}}{\gamma^{1/2}n^{1/4}}
			+\dfrac{(h_Ah_L)^{-1/6}\log(n)^{2/3}}{\gamma^{1/3}n^{1/6}}\leq \dfrac{\epsilon}{3\sqrt{\log(n)}}$.
			This is equivalent to requiring that
			\begin{align*}
				\gamma\gtrsim \dfrac{\log(n)^3}{\epsilon^2nh_Ah_L}+
				\dfrac{\log(n)^{5/2}}{\epsilon^2(nh_Ah_L)^{1/2}}
				+\dfrac{\log(n)^{7/2}}{\epsilon^3(nh_Ah_L)^{1/2}}.
			\end{align*}
			Since $\log(n)^7=o(nh_Ah_L)$, we know that
			\begin{align*}
				E_{211}:=&\Pr\left(\left|\sup_{l\in\mathcal{L}}\sqrt{nh_Ah_L}
				\left|\frac{E_{12}(l)}
				{p_L(l)\sqrt{\Gamma_{\pi,h}^o(a,l,l)}}\right| - 
				\sup_{l\in\mathcal{L}}\left|d_{b,l}^o\right|\right|\geq \dfrac{\epsilon}{\sqrt{\log(n)}}\right)\\
				\lesssim&
				\dfrac{\log(n)^3}{\epsilon^2nh_Ah_L}+
				\dfrac{\log(n)^{5/2}}{\epsilon^2(nh_Ah_L)^{1/2}}
				+\dfrac{\log(n)^{7/2}}{\epsilon^3(nh_Ah_L)^{1/2}} + \dfrac{\log(n)}{n}=o(1).
			\end{align*}	
			The $E_{212}$ can be controlled through the anti-concentration inequality (Lemma~A.1 of \cite{chernozhukov2014gaussian}) as
			\begin{align}
				E_{212}=&\sup_{x\in\mathbb{R}}\Pr\left(\left|\sup_{l\in\mathcal{L}}|d_{b,l}^o|-x\right|\leq
				\dfrac{\epsilon}{\sqrt{\log(n)}}
				\right)\notag\\
				\lesssim &C\dfrac{\epsilon}{\sqrt{\log(n)}}\left\{
				\mathbb{E}\left[\sup_{l\in\mathcal{L}} d_{b,l}^o\right]
				+\sqrt{1 + \log(\dfrac{\sqrt{\log(n)}}{\epsilon})}
				\right\}\notag\\
				\lesssim &C\dfrac{\epsilon}{\sqrt{\log(n)}}\left\{
				\int_0^1 \sqrt{\log(\frac{1}{h_L\epsilon})}\mathrm{d}\epsilon
				+\sqrt{1 + \log(\dfrac{\sqrt{\log(n)}}{\epsilon})}
				\right\}\notag\\
				\lesssim &C\dfrac{\epsilon}{\sqrt{\log(n)}}\left\{
				1+\sqrt{\log(\frac{1}{h_L})}
				+ \log(\dfrac{\sqrt{\log(n)}}{\epsilon})
				\right\}\lesssim C\epsilon + o(1).\label{expadd: 1}
			\end{align}
			The third inequality uses the Dudley integral lemma.
			The final step follows from the fact that $h_L^{-(1+p)}\lesssim T\lesssim n^q$ for some $p,q>0$. Furthermore, $E_{213}$ converges to zero in probability since 
			$
			\sup_{l\in\mathcal{L}} \left|\dfrac{E_{1j}(l)}{p_L(l)\sqrt{\Gamma_{\pi,h}^o(a,l,l)}}\right|=o_p(\dfrac{1}{\sqrt{\log(n)nh_Ah_L}})
			$ for $j=1,3,4$.
		\end{proof}
		
		\begin{proof}[Proof of $\sup_{x\in\mathbb{R}}|E_{22}(x)|=o(1)$:]
			We directly verify that
			\begin{align*}
				|E_{22}(x)|=&\left|\Pr\left(\sup_{l\in\mathcal{L}}|d_{b,l}^o|\leq x\right) - \Pr\left(\sup_{t=1,\ldots,T}|d_{b,l_t}^o|\leq x\right)\right|\\
				=&\left|\begin{array}{l}
					\Pr\left(\sup\limits_{l\in\mathcal{L}}|d_{b,l}^o|\leq x,
					\sup\limits_{t=1,\ldots,T}|d_{b,l_t}^o|\leq x\right) 
					+ \Pr\left(\sup\limits_{l\in\mathcal{L}}|d_{b,l}^o|\leq x,
					\sup\limits_{t=1,\ldots,T}|d_{b,l_t}^o|\geq x\right) \\
					- \Pr\left(\sup\limits_{l\in\mathcal{L}}|d_{b,l}^o|\leq x,
					\sup\limits_{t=1,\ldots,T}|d_{b,l_t}^o|\leq x\right) 
					- \Pr\left(\sup\limits_{l\in\mathcal{L}}|d_{b,l}^o|\geq x,
					\sup\limits_{t=1,\ldots,T}|d_{b,l_t}^o|\leq x\right) 
				\end{array}\right|\\
				\leq&\Pr\left(\sup\limits_{l\in\mathcal{L}}|d_{b,l}^o|\leq x,
				\sup\limits_{t=1,\ldots,T}|d_{b,l_t}^o|\geq x\right)
				+\Pr\left(\sup\limits_{l\in\mathcal{L}}|d_{b,l}^o|\geq x,
				\sup\limits_{t=1,\ldots,T}|d_{b,l_t}^o|\leq x\right).
			\end{align*}
			We control the first term as
			\begin{align*}
				&\Pr\left(\sup\limits_{l\in\mathcal{L}}|d_{b,l}^o|\leq x,
				\sup\limits_{t=1,\ldots,T}|d_{b,l_t}^o|\geq x\right)\\
				\leq& \Pr\left(\sup\limits_{l\in\mathcal{L}}|d_{b,l}^o|\leq x-\dfrac{\epsilon}{\sqrt{\log(n)}},
				\sup\limits_{t=1,\ldots,T}|d_{b,l_t}^o|\geq x\right)
				+\Pr\left(x-\dfrac{\epsilon}{\sqrt{\log(n)}}\leq\sup\limits_{l\in\mathcal{L}}|d_{b,l}^o|\leq x\right)\\
				\leq&\Pr\left(\left|\sup\limits_{l\in\mathcal{L}}|d_{b,l}^o|-
				\sup\limits_{t=1,\ldots,T}|d_{b,l_t}^o|\right|\geq 
				\dfrac{\epsilon}{\sqrt{\log(n)}}\right)+
				\sup_{x\in\mathbb{R}}\left|
				\Pr\left(\left|\sup\limits_{l\in\mathcal{L}}|d_{b,l}^o|-x\right|\leq \dfrac{\epsilon}{\sqrt{\log(n)}}\right)
				\right|.
			\end{align*}
			The second term can be controlled as in Equation~\eqref{expadd: 1}.
			For the first term,
			\begin{align*}
				&\Pr\left(\left|\sup\limits_{l\in\mathcal{L}}|d_{b,l}^o|-
				\sup\limits_{t=1,\ldots,T}|d_{b,l_t}^o|\right|\geq 
				\dfrac{\epsilon}{\sqrt{\log(n)}}\right)\leq 
				\Pr\left(\sup\limits_{l,l'\in\mathcal{L},\;|l-l'|\leq 1/T}|d_{b,l}^o-d_{b,l'}^o|\geq 
				\dfrac{\epsilon}{\sqrt{\log(n)}}\right)\\
				\leq&\dfrac{\sqrt{\log(n)}}{\epsilon}\mathbb{E}\left[
				\sup\limits_{l,l'\in\mathcal{L},\;|l-l'|\leq 1/T}|d_{b,l}^o-d_{b,l'}^o|
				\right]
			\end{align*}
			Define a norm vector $\rho_L(l,l')^2:=\mathbb{E}[(d_{b,l}^o-d_{b,l'}^o)^2]$,
			\begin{align*}
				&\Gamma_{\pi,h}^o(a,l,l'):=\dfrac{h_Ah_L}{p_L(l)p_L(l')}\mathbb{E}\left[K_{h_L}(L-l)K_{h_L}(L-l')\left\{S_{\pi,h_A}^o(a,O)-\mathbb{E}[S_{\pi,h_A}^o(a,O)\mid L]\right\}^2\right];\\
				&\rho_L(l,l')^2:=2-2\mathrm{Cov}\{d_{b,l}^o,d_{b,l'}^o\}
				=2-2\dfrac{\Gamma_{\pi,h}^o(a,l,l')}
				{\sqrt{\Gamma_{\pi,h}^o(a,l,l)\Gamma_{\pi,h}^o(a,l',l')}}\\
				\lesssim& 
				\dfrac{\sqrt{\Gamma_{\pi,h}^o(a,l',l')}-\sqrt{\Gamma_{\pi,h}^o(a,l,l')}}{\sqrt{\Gamma_{\pi,h}^o(a,l',l')}}
				-\dfrac{\sqrt{\Gamma_{\pi,h}^o(a,l,l')}-\sqrt{\Gamma_{\pi,h}^o(a,l,l)}}{\sqrt{\Gamma_{\pi,h}^o(a,l,l)}}\\
				\lesssim&
				\dfrac{\left|\Gamma_{\pi,h}^o(a,l',l')-\Gamma_{\pi,h}^o(a,l,l')\right|}{\sqrt{\Gamma_{\pi,h}^o(a,l',l')}\{\sqrt{\Gamma_{\pi,h}^o(a,l',l')}+\sqrt{\Gamma_{\pi,h}^o(a,l,l')}\}}
				+\dfrac{\left|\Gamma_{\pi,h}^o(a,l,l')-\Gamma_{\pi,h}^o(a,l,l)\right|}{\sqrt{\Gamma_{\pi,h}^o(a,l,l)}\{\sqrt{\Gamma_{\pi,h}^o(a,l,l')}+\sqrt{\Gamma_{\pi,h}^o(a,l,l)}\}}\\
				\lesssim&
				\dfrac{h_Ah_L}{p_L(l)p_L(l')}\mathbb{E}\left[
				|K_{h_L}(L-l')-K_{h_L}(L-l)|\cdot K_{h_L}(L-l')
				\left\{S_{\pi,h_A}^o(a,O)-\mathbb{E}[S_{\pi,h_A}^o(a,O)\mid L]\right\}^2\right]\\
				&+\dfrac{h_Ah_L}{p_L(l)p_L(l')}\mathbb{E}\left[
				|K_{h_L}(L-l')-K_{h_L}(L-l)|\cdot K_{h_L}(L-l)
				\left\{S_{\pi,h_A}^o(a,O)-\mathbb{E}[S_{\pi,h_A}^o(a,O)\mid L]\right\}^2\right]\\
				\lesssim&
				\dfrac{h_Ah_L}{p_L(l)p_L(l')}\dfrac{|l-l'|}{h_L^2}\mathbb{E}\left[
				\{K_{h_L}(L-l')+K_{h_L}(L-l)\}
				\left\{S_{\pi,h_A}^o(a,O)-\mathbb{E}[S_{\pi,h_A}^o(a,O)\mid L]\right\}^2\right]\\
				\lesssim&
				\dfrac{h_A |l-l'|}{h_L}\int 
				\{K_{h_L}(t-l')+K_{h_L}(t-l)\}
				K_{h_A}(r-a)^2\{\pi(z,t)-\rho_{\pi}^o(t)\}^2 
				p_{Z,A,L}(z,r,t)\mathrm{d} z\mathrm{d} r\mathrm{d} t\\
				\lesssim&
				\dfrac{|l-l'|}{h_L}\int 
				K(t)
				K(r)^2\{\pi(z,t)-\rho_{\pi}^o(t)\}^2 
				\left\{\begin{array}{l}
					p_{Z,A,L}(z,a+rh_A,l+rh_L) \\ + p_{Z,A,L}(z,a+rh_A,l'+rh_L)
				\end{array}\right\}
				\mathrm{d} z\mathrm{d} r\mathrm{d} t
				\lesssim \dfrac{|l-l'|}{h_L}.
			\end{align*}
			Therefore, we know that $\rho_L(l,l')\lesssim \sqrt{\frac{|l-l'|}{h_L}}$, and
			\begin{align*}
				&\Pr\left(\left|\sup\limits_{l\in\mathcal{L}}|d_{b,l}^o|-
				\sup\limits_{t=1,\ldots,T}|d_{b,l_t}^o|\right|\geq 
				\dfrac{\epsilon}{\sqrt{\log(n)}}\right)
				\lesssim
				\dfrac{\sqrt{\log(n)}}{\epsilon}\int_0^{\sqrt{\frac{1}{Th_L}}} \sqrt{\log\{N(\epsilon, \mathcal{L}, \rho_L)\}}\mathrm{d} \epsilon\\
				\lesssim& 
				\dfrac{\sqrt{\log(n)}}{\epsilon} \int_0^{\sqrt{\frac{1}{Th_L}}} \sqrt{\log(\dfrac{1}{h_L}\dfrac{1}{\epsilon})}\mathrm{d} \epsilon
				\lesssim 
				\dfrac{1}{\epsilon} 
				\sqrt{\dfrac{\log(n)\log(1/h_L)}{Th_L}}
				\sqrt{\log(Th_L)}=o(1).
			\end{align*}
			The final step follows by $T^{-1}=o(h_L^{1+p})$ and $T=O(n^q)$ for some $p,q>0$.
			Consequently, we conclude that $\sup_{x\in\mathbb{R}}|E_{22}(x)|=o(1)$.
		\end{proof}
		
		\begin{proof}[Proof of $\sup_{x\in\mathbb{R}}|E_{23}(x)|=o(1)$:]
			First, recall that
			\begin{align*}
				&\mathrm{Cov}\left\{\hat{d}_{b,l_{t_1}},\hat{d}_{b,l_{t_2}}\right\}=
				\frac{\hat{\Gamma}_{\pi,h}^{(n)}(a,l_{t_1},l_{t_2})}{\sqrt{\hat{\Gamma}_{\pi,h}^{(n)}(a,l_{t_1},l_{t_1})\hat{\Gamma}_{\pi,h}^{(n)}(a,l_{t_2},l_{t_2})}},\\
				&\mathrm{Cov}\left\{d_{b,l_{t_1}}^o,d_{b,l_{t_2}}^o\right\}=
				\dfrac{\Gamma_{\pi,h}^o(a,l_{t_1},l_{t_2})}
				{\sqrt{\Gamma_{\pi,h}^o(a,l_{t_1},l_{t_1})\Gamma_{\pi,h}^o(a,l_{t_2},l_{t_2})}},\text{ where}\\
				&\hat{\Gamma}_{\pi,h}^{(n)}(a,l,l') 
				:=h_Ah_L
				e_1^{\top}
				\left\{
				\mathbb{E}_n\left[
				K_{h_L}(L-l)g\left(\frac{L-l}{h_L}\right)
				g\left(\frac{L-l'}{h_L}\right)^{\top}\right]
				\right\}^{-1}\\
				&\mathbb{E}_n\left[\begin{array}{l}
					K_{h_L}(L-l)
					K_{h_L}(L-l')
					g\left(\frac{L-l}{h_L}\right)
					g\left(\frac{L-l'}{h_L}\right)^{\top}\\
					\left\{\hat{S}_{\pi,h_A}(a,O)
					-g\left(\frac{L-l}{h_L}\right)^{\top}
					\hat\beta_{\pi,h}^{(n)}(a,l)\right\}
					\left\{\hat{S}_{\pi,h_A}(a,O)
					-g\left(\frac{L-l'}{h_L}\right)^{\top}
					\hat\beta_{\pi,h}^{(n)}(a,l')\right\}
				\end{array}\right]\\
				&\times
				\left\{
				\mathbb{E}_n\left[
				K_{h_L}(L-l')g\left(\frac{L-l'}{h_L}\right)
				g\left(\frac{L-l'}{h_L}\right)^{\top}\right]
				\right\}^{-1}e_1,\\
				&\Gamma_{\pi,h}^o(a,l,l'):=\dfrac{h_Ah_L}{p_L(l)p_L(l')}\mathbb{E}\left[K_{h_L}(L-l)K_{h_L}(L-l')\left\{S_{\pi,h_A}^o(a,O)-\mathbb{E}[S_{\pi,h_A}^o(a,O)\mid L]\right\}^2\right].
			\end{align*}
			Define $\Delta_n:=\sup\limits_{t,t'=1,\ldots, T}\left|
			\hat{\Gamma}_{\pi,h}^{(n)}(a,l_{t},l_{t'})
			-\Gamma_{\pi,h}^o(a,l_{t},l_{t'})
			\right|$. 
			Since $\Delta_n=o_p(1/\log(n)^2)$, 
			%	it is direct to verify that $h_Ah_L\sqrt{n}\Delta_n = O_p(\sqrt{\log(1/h_L)}) = O_p(\sqrt{\log(n)})$. Next, 
			we adopt Theorem~2 of \citet{chernozhukov2015comparison} to deduce that
			\begin{align*}
				\sup_{x\in\mathbb{R}}|E_{23}(x)|=&\sup_{x\in\mathbb{R}}\left|\Pr\left(\sup_{t=1,\ldots,T}|d_{b,l_t}^o|\leq x\right) - \Pr\left(\max_{t=1,\ldots,T}|\hat{d}_{b,l_t}|\leq x
				\middle|\{O_i\}_{i=1}^n\right) \right|\\
				\lesssim&
				\left[\Delta_n \log(T) \left\{
				\log(\dfrac{1}{h_Ah_L})
				+\log(\dfrac{1}{\Delta_n})
				\right\}\right]^{1/3}\\
				\lesssim&
				o_p\left(\left[ \frac{\log(T/h_Ah_L)}{\log(n)^2} + \sqrt{\Delta_n}\log(\dfrac{1}{\Delta_n}) \right]^{1/3}\right)=o_p(1).
			\end{align*}
			The final step follows from the fact that $T=O(n^q)$ for some $q>0$.
			This completes the proof of $\sup_{x\in\mathbb{R}}|E_{23}(x)|=o(1)$.
			%	Theorem~\ref{thm: RWF testing 3}.
		\end{proof}

		\subsection{Proof of Theorem~\ref{thm: RWF testing 4}}
		For any $P_0\in \mathcal{H}_0(a)$, there exists $l_0\in\mathcal{L}$, such that $\zeta_{\pi}^o(a,l_0)=0$. Therefore, 
		\begin{align*}
			&\Pr\nolimits_{P_0}\left(\min_{t=1,\ldots, T}\sqrt{nh_Ah_L}\left|
			\frac{\hat{\zeta}_{\pi,h}^{(n)}(a,l_t)}{\hat{\sigma}_{\pi,\zeta,h}^{(n)}(a,l_t)}
			\right| > x\right)=
			\Pr\nolimits_{P_0}\left(\inf_{l\in\mathcal{L}}\sqrt{nh_Ah_L}\left|
			\frac{\hat{\zeta}_{\pi,h}^{(n)}(a,l)}{\hat{\sigma}_{\pi,\zeta,h}^{(n)}(a,l)}
			\right| > x\right) + o(1)\\
			=&\Pr\nolimits_{P_0}\left(\inf_{l\in\mathcal{L}}\sqrt{nh_Ah_L}\left|
			\frac{\hat{\zeta}_{\pi,h}^{(n)}(a,l)}{\hat{\sigma}_{\pi,\zeta,h}^{(n)}(a,l)}
			\right| > x\middle| \;\exists l_0\in\mathcal{L},\text{ such that } \zeta_{\pi}^o(a,l_0)=0\right)+o(1)\\
			\leq&\Pr\nolimits_{P_0}\left(\sqrt{nh_Ah_L}\left|
			\frac{\hat{\zeta}_{\pi,h}^{(n)}(a,l_0)- \zeta_{\pi}^o(a,l_0)}{\hat{\sigma}_{\pi,\zeta,h}^{(n)}(a,l_0)}
			\right| > x\middle| \;\exists l_0\in\mathcal{L},\text{ such that } \zeta_{\pi}^o(a,l_0)=0\right)+o(1)\\
			\leq&\Pr\nolimits_{P_0}\left(\sup_{l\in\mathcal{L}}\sqrt{nh_Ah_L}\left|
			\frac{\hat{\zeta}_{\pi,h}^{(n)}(a,l)- \zeta_{\pi}^o(a,l)}{\hat{\sigma}_{\pi,\zeta,h}^{(n)}(a,l)}
			\right| > x\right)+o(1)\\
			=&\Pr\nolimits_{P_0}\left(\max\limits_{t=1,\ldots,T}|\hat{d}_{b,l_t}|> x\middle| \{O_i\}_{i=1}^n\right) + o_p(1)
			=\dfrac{1}{B}\sum_{b=1}^B I\{\hat{D}_{b}^*> x\}+ o_p(1).
		\end{align*} 
		If we define $\hat{c}_{1-\alpha}:=\inf\{x: \frac{1}{B}\sum_{b=1}^B I\{\hat{D}_{b}^*\leq x\}\geq 1-\alpha\}$,
		this indicates that
		\begin{align*}
			\underset{n\to+\infty}{\overline{\lim}}\Pr\nolimits_{P_0}\left(
			\hat{p}_{a,\pi}^{(n)}\leq \alpha\right)
			=&\underset{n\to+\infty}{\overline{\lim}}\Pr\nolimits_{P_0}\left(\dfrac{1}{B}\sum\limits_{b=1}^B I\left\{\hat{D}_b^*\geq \min\limits_{t=1,\ldots,T} 
			\left|\dfrac{\sqrt{nh_Ah_L}\hat{\zeta}_{\pi,h}^{(n)}(a,l_t)}{\hat{\sigma}_{\pi,\zeta,h}^{(n)}(a,l_t)}\right|\right\}\leq \alpha\right)\\
			=&\underset{n\to+\infty}{\overline{\lim}}\Pr\nolimits_{P_0}\left(\min_{t=1,\ldots, T}\sqrt{nh_Ah_L}\left|
			\frac{\hat{\zeta}_{\pi,h}^{(n)}(a,l_t)}{\hat{\sigma}_{\pi,\zeta,h}^{(n)}(a,l_t)}
			\right| > \hat{c}_{1-\alpha}\right)\\
			\leq&\underset{n\to+\infty}{\overline{\lim}}
			\underset{B\to+\infty}{\overline{\lim}}
			\dfrac{1}{B}\sum_{b=1}^B I\{\max\limits_{t=1,\ldots,T}|\hat{d}_{b,l_t}|> \hat{c}_{1-\alpha}\}\leq \alpha.
		\end{align*}
		In contrast, for any $P_1\in \mathcal{H}_1(a)$, since $\zeta_{\pi}^o(a,l)$ is continuous on the compact set $\mathcal{L}\subseteq \mathbb{R}$, there exists $\epsilon_0>0$, such that $|\zeta_{\pi}^o(a,l)| \geq \epsilon_0$. Therefore,
		\begin{align*}
			&\Pr\nolimits_{P_1}\left(\min_{t=1,\ldots, T}\sqrt{nh_Ah_L}\left|
			\frac{\hat{\zeta}_{\pi,h}^{(n)}(a,l_t)}{\hat{\sigma}_{\pi,\zeta,h}^{(n)}(a,l_t)}
			\right| > x\right)=
			\Pr\nolimits_{P_1}\left(\inf_{l\in\mathcal{L}}\sqrt{nh_Ah_L}\left|
			\frac{\hat{\zeta}_{\pi,h}^{(n)}(a,l)}{\hat{\sigma}_{\pi,\zeta,h}^{(n)}(a,l)}
			\right| > x\right) + o(1)\\
			=&\Pr\nolimits_{P_1}\left(\inf_{l\in\mathcal{L}}\sqrt{nh_Ah_L}\left|
			\frac{\hat{\zeta}_{\pi,h}^{(n)}(a,l)}{\hat{\sigma}_{\pi,\zeta,h}^{(n)}(a,l)}
			\right| > x\middle| \;\forall l\in\mathcal{L}, |\zeta_{\pi}^o(a,l)|\geq\epsilon_0\right)+o(1)\\
			\geq&\Pr\nolimits_{P_1}\left(
			\inf\limits_{l\in\mathcal{L}}\left|
			\frac{\sqrt{nh_Ah_L}\zeta_{\pi}^o(a,l) }{\hat{\sigma}_{\pi,\zeta,h}^{(n)}(a,l)}
			\right|>
			\sup\limits_{l\in\mathcal{L}}\sqrt{nh_Ah_L}\left|
			\frac{\hat{\zeta}_{\pi,h}^{(n)}(a,l) -\zeta_{\pi}^o(a,l) }{\hat{\sigma}_{\pi,\zeta,h}^{(n)}(a,l)}
			\right| + x
			\middle| \;\forall l\in\mathcal{L}, |\zeta_{\pi}^o(a,l)|\geq\epsilon_0\right)+o(1)\\
			\geq&\Pr\nolimits_{P_1}\left(
			\sup\limits_{l\in\mathcal{L}}\sqrt{nh_Ah_L}\left|
			\frac{\hat{\zeta}_{\pi,h}^{(n)}(a,l) -\zeta_{\pi}^o(a,l) }{\hat{\sigma}_{\pi,\zeta,h}^{(n)}(a,l)}
			\right|\leq
			\inf\limits_{l\in\mathcal{L}}\left|
			\frac{\sqrt{nh_Ah_L} }{\hat{\sigma}_{\pi,\zeta,h}^{(n)}(a,l)}
			\right|\epsilon_0
			-x\right)+o(1)\\
			=&\Pr\nolimits_{P_1}\left(
			\sup\limits_{l\in\mathcal{L}}\left|
			d_{b,l}^o
			\right|\leq
			\inf\limits_{l\in\mathcal{L}}\left|
			\frac{\sqrt{nh_Ah_L} }{\hat{\sigma}_{\pi,\zeta,h}^{(n)}(a,l)}
			\right|\epsilon_0
			-x\right)+o(1).
		\end{align*}
		Since $\sup\limits_{l\in\mathcal{L}}\left|
		d_{b,l}^o\right| = O_p(\log(1/h_Ah_L))=O_p(\log(n)) = o_p(\sqrt{nh_Ah_L})$, and $\hat{c}_{1-\alpha} =O_p(\log(n))= o_p(\sqrt{nh_Ah_L})$ we know that for any $P_1\in \mathcal{H}_1(a)$,
		\begin{align*}
			&\underset{n\to+\infty}{\underline{\lim}}\Pr\nolimits_{P_1}\left(
			\hat{p}_{a,\pi}^{(n)}\leq \alpha\right)=\underset{n\to+\infty}{\underline{\lim}}\Pr\nolimits_{P_1}\left(\min_{t=1,\ldots, T}\sqrt{nh_Ah_L}\left|
			\frac{\hat{\zeta}_{\pi,h}^{(n)}(a,l_t)}{\hat{\sigma}_{\pi,\zeta,h}^{(n)}(a,l_t)}
			\right| > \hat{c}_{1-\alpha}\right)\\
			\geq&\underset{n\to+\infty}{\underline{\lim}}
			\Pr\nolimits_{P_1}\left(
			\sup\limits_{l\in\mathcal{L}}\left|
			d_{b,l}^o
			\right| \leq
			\inf\limits_{l\in\mathcal{L}}\left|
			\frac{\sqrt{nh_Ah_L}}{\hat{\sigma}_{\pi,\zeta,h}^{(n)}(a,l)}
			\right|\epsilon_0- \hat{c}_{1-\alpha}\right) \\
			=&\underset{n\to+\infty}{\underline{\lim}}
			\Pr\nolimits_{P_1}\left(
			\sup\limits_{l\in\mathcal{L}}\left|
			d_{b,l}^o
			\right|  + \hat{c}_{1-\alpha}\leq
			\frac{\sqrt{nh_Ah_L}\epsilon_0 }{\sup\limits_{l\in\mathcal{L}}\hat{\sigma}_{\pi,\zeta,h}^{(n)}(a,l)}
			\right)=1.
		\end{align*}
		This finishes the proof of Theorem~\ref{thm: RWF testing 4}.}

	\subsection{Proof of Theorem~\ref{thm: oracle bounds for ERM}}
	We use $\mathbb{E}_{\mathrm{out}}$ to denote the expectation taken with respect to an independent out-of-sample observation drawn from the same distribution as the original data. 
	In addition, we use $\widetilde{O}_p$ to denote quantities that are bounded in expectation, and $\widetilde{o}_p$ to denote quantities that converge to zero in expectation.
	\begin{align*}
		&\Pr(A\in\mathcal{N})\left\{\|\hat{f}_{\pi,\mathcal{N}}^{(n)}-\theta\|_{\mathcal{N}}^2-
		\|f_{\mathcal{N}}^*-\theta\|_{\mathcal{N}}^2\right\}\\
		=&
		\left\{\mathbb{E}_{\mathrm{out}}\left[|\hat{f}_{\pi,\mathcal{N}}^{(n)}(A)-\theta(A)|^2
		I\left\{A\in\mathcal{N}\right\}\right]
		-\mathbb{E}_{\mathrm{out}}\left[| f_{\mathcal{N}}^{*}(A)-\theta(A)|^2
		I\left\{A\in\mathcal{N}\right\}\right]\right\}\\
		=&
		\left\{\begin{array}{l}
			\mathbb{E}_{\mathrm{out}}\left[|\hat{f}_{\pi,\mathcal{N}}^{(n)}(A)-\varphi_{\pi}(O;\alpha_{\pi}^o,P_O)|^2
			I\left\{A\in\mathcal{N}\right\}\right]\\
			-\mathbb{E}_{\mathrm{out}}\left[| f_{\mathcal{N}}^{*}(A)-\varphi_{\pi}(O;\alpha_{\pi}^o,P_O)|^2
			I\left\{A\in\mathcal{N}\right\}\right]
		\end{array}\right\}\\
		=&
		\left\{\begin{array}{l}
			\dfrac{1}{K}\displaystyle\sum_{k=1}^K\mathbb{E}_{\mathrm{out}}\left[
			\left\{\begin{array}{l}
				|\hat{f}_{\pi,\mathcal{N}}^{(n)}(A)-\varphi_{\pi}(O;\alpha_{\pi}^o,P_O)|^2\\
				-|\hat{f}_{\pi,\mathcal{N}}^{(n)}(A)-\varphi_{\pi}(O;\hat\alpha_{\pi}^{(n,-k)},\hat{P}_O^{(n,-k)})|^2
			\end{array}\right\}
			I\left\{A\in\mathcal{N}\right\}\right]\\
			-\dfrac{1}{K}\displaystyle\sum_{k=1}^K\mathbb{E}_{\mathrm{out}}\left[
			\left\{\begin{array}{l}
				|f_{\mathcal{N}}^{*}(A)-\varphi_{\pi}(O;\alpha_{\pi}^o,P_O)|^2\\
				-|f_{\mathcal{N}}^{*}(A)-\varphi_{\pi}(O;\hat\alpha_{\pi}^{(n,-k)},\hat{P}_O^{(n,-k)})|^2
			\end{array}\right\}
			I\left\{A\in\mathcal{N}\right\}\right]\\
		\end{array}\right\}\\
		&+\left\{\begin{array}{l}
			\dfrac{1}{K}\displaystyle\sum_{k=1}^K\{\mathbb{E}_{\mathrm{out}}-\mathbb{E}_{nk}\}\left[|\hat{f}_{\pi,\mathcal{N}}^{(n)}(A)-\varphi_{\pi}(O;\hat\alpha_{\pi}^{(n,-k)},\hat{P}_O^{(n,-k)})|^2
			I\left\{A\in\mathcal{N}\right\}\right]\\
			-\dfrac{1}{K}\displaystyle\sum_{k=1}^K\{\mathbb{E}_{\mathrm{out}}-\mathbb{E}_{nk}\}\left[| f_{\mathcal{N}}^{*}(A)-\varphi_{\pi}(O;\hat\alpha_{\pi}^{(n,-k)},\hat{P}_O^{(n,-k)})|^2
			I\left\{A\in\mathcal{N}\right\}\right]
		\end{array}\right\}\\
		&+\left\{\begin{array}{l}
			\dfrac{1}{K}\displaystyle\sum_{k=1}^K\mathbb{E}_{nk}\left[|\hat{f}_{\pi,\mathcal{N}}^{(n)}(A)-\varphi_{\pi}(O;\hat\alpha_{\pi}^{(n,-k)},\hat{P}_O^{(n,-k)})|^2
			I\left\{A\in\mathcal{N}\right\}\right]\\
			-\dfrac{1}{K}\displaystyle\sum_{k=1}^K\mathbb{E}_{nk}\left[| f_{\mathcal{N}}^{*}(A)-\varphi_{\pi}(O;\hat\alpha_{\pi}^{(n,-k)},\hat{P}_O^{(n,-k)})|^2
			I\left\{A\in\mathcal{N}\right\}\right]
		\end{array}\right\}\\
		:=&\dfrac{1}{K}\sum_{k=1}^K\{R_{5k}+R_{6k}+R_{7k}\}.
	\end{align*}
	First, it is easy to see that
	$\frac{1}{K}\sum_{k=1}^KR_{7k}\leq 0$ by the definition of $\hat{f}_{\pi,\mathcal{N}}^{(n)}(A)$ in Equation~\eqref{exp: 30}. 
	Second, from Lemma~12 in \citet{foster2023orthogonalstatisticallearning} (or Theorem~14.20 in \citet{wainwright2019high}), we know that
	for any $k=1, \ldots, K$, there exists a constant $c$, such that with probability more than $1-c\exp\{-cn\xi_n^2\}$, 
	\begin{align}
		&\mathbb{E}_{\mathrm{out}}\left[|\hat{f}_{\pi,\mathcal{N}}^{(n)}(A)-f_{\mathcal{N}}^*(A)|^2I(A\in\mathcal{N})\right]\leq \xi_n^2, \text{ or }\label{exp: 31}\\
		&\left|\{\mathbb{E}_{nk}-\mathbb{E}_{\mathrm{out}}\}\left(\begin{array}{l}
			|\hat{f}_{\pi,\mathcal{N}}^{(n)}(A)-\varphi_{\pi}(O;\hat\alpha_{\pi}^{(n,-k)},\hat{P}_O^{(n,-k)})|^2I(A\in\mathcal{N})\\
			-| f_{\mathcal{N}}^{*}(A)-\varphi_{\pi}(O;\hat\alpha_{\pi}^{(n,-k)},\hat{P}_O^{(n,-k)})|^2I(A\in\mathcal{N})
		\end{array}\right)\right|\notag\\
		&\lesssim \xi_n 
		\left(\mathbb{E}_{\mathrm{out}}\left[\left|\hat{f}_{\pi,\mathcal{N}}^{(n)}(A)-f_{\mathcal{N}}^{*}(A)\right|^2I(A\in\mathcal{N})\right]\right)^{1/2}\label{exp: 32}
	\end{align}
	We separately discuss these two settings. If Equation~\eqref{exp: 31} holds, then we can directly get that 
	\begin{align*}
		&\mathbb{E}_{\mathrm{out}}\left[|\hat{f}_{\pi,\mathcal{N}}^{(n)}(A)-\theta(A)|^2I(A\in\mathcal{N})\right]\\
		\lesssim&\mathbb{E}_{\mathrm{out}}\left[|f_{\mathcal{N}}^*(A)-\theta(A)|^2I(A\in\mathcal{N})\right]
		+\mathbb{E}_{\mathrm{out}}\left[|\hat{f}_{\pi,\mathcal{N}}^{(n)}(A)-f_{\mathcal{N}}^*(A)|^2I(A\in\mathcal{N})\right]\\
		\leq &\mathbb{E}_{\mathrm{out}}\left[|f_{\mathcal{N}}^*(A)-\theta(A)|^2I(A\in\mathcal{N})\right]+\xi_n^2.
	\end{align*}
	Otherwise, if Equation~\eqref{exp: 32} holds, we know that 
	\begin{align}
		|R_{6k}|=&\left|\{\mathbb{E}_{nk}-\mathbb{E}_{\mathrm{out}}\}\left(\begin{array}{l}
			|\hat{f}_{\pi,\mathcal{N}}^{(n)}(A)-\varphi_{\pi}(O;\hat\alpha_{\pi}^{(n,-k)},\hat{P}_O^{(n,-k)})|^2I(A\in\mathcal{N})\notag\\
			-| f_{\mathcal{N}}^{*}(A)-\varphi_{\pi}(O;\hat\alpha_{\pi}^{(n,-k)},\hat{P}_O^{(n,-k)})|^2I(A\in\mathcal{N})
		\end{array}\right)\right|\\
		\lesssim&\mathbb{E}_{\mathrm{out}}\left[\left|\hat{f}_{\pi,\mathcal{N}}^{(n)}(A)-f_{\mathcal{N}}^{*}(A)\right|^2I(A\in\mathcal{N})\right]^{1/2}\xi_n.\label{exp: 33}
	\end{align}
	Finally, for controlling $R_{5k}$,
	\begin{align*}
		|R_{5k}|=&2\left|\mathbb{E}_{\mathrm{out}}\left[\begin{array}{l}
			\left\{
			\varphi_{\pi}(O;\hat\alpha_{\pi}^{(n,-k)}, \hat{P}_O^{(n,-k)})
			-\varphi_{\pi}(O;\alpha_{\pi}^o, P_O)
			\right\}\\\times\left\{
			\hat{f}_{\pi,\mathcal{N}}^{(n)}(A) - f_{\mathcal{N}}^*(A)
			\right\}I\{A\in\mathcal{N}\}
		\end{array}\right]\right|\\
		\lesssim&\left|\mathbb{E}_{\mathrm{out}}\left[\begin{array}{l}
			\left\{
			\varphi_{\pi}(O;\hat\alpha_{\pi}^{(n,-k)}, P_O)
			-\varphi_{\pi}(O;\alpha_{\pi}^o, P_O)
			\right\}\\\times\left\{
			\hat{f}_{\pi,\mathcal{N}}^{(n)}(A) - f_{\mathcal{N}}^*(A)
			\right\}I\{A\in\mathcal{N}\}
		\end{array}\right]\right|\\
		&+\left|\mathbb{E}_{\mathrm{out}}\left[\begin{array}{l}
			\left\{
			\varphi_{\pi}(O;\hat\alpha_{\pi}^{(n,-k)}, \hat{P}_O^{(n,-k)})
			-\varphi_{\pi}(O;\hat\alpha_{\pi}^{(n,-k)}, P_O)
			\right\}\\\times\left\{
			\hat{f}_{\pi,\mathcal{N}}^{(n)}(A) - f_{\mathcal{N}}^*(A)
			\right\}I\{A\in\mathcal{N}\}
		\end{array}\right]\right|\\
		\lesssim& R_{5k1}+R_{5k2}
	\end{align*}	
	Thus, we only need to control for the final term $R_{5k1}$ and $R_{5k2}$ for $k=1,\ldots,K$.
	From Lemma~\ref{lem: mixed bias}, we know that
	\begin{align*}
		&\mathbb{E}_{\mathrm{out}}[\varphi_{\pi}(O;\hat\alpha_{\pi}^{(n,-k)}, P_O)\mid A] - \theta(A)\\
		=&\mathbb{E}_{\mathrm{out}}\left[\begin{array}{l}
			\dfrac{\hat\delta^{(n,-k)}(A,L)}{\hat\kappa_{\pi}^{(n,-k)}(A,L)}
			(\hat\rho_{\pi}^{(n,-k)}(L)-\rho_{\pi}^o(L))(\hat\eta^{(n,-k)}(A,L)-\eta^o(A,L))\\
			+ \left\{\delta^o(A,L)-
			\hat\delta^{(n,-k)}(A,L)\dfrac{\kappa_{\pi}^o(A,L)}{\hat\kappa_{\pi}^{(n,-k)}(A,L)}\right\}
			\{\hat\mu_{\pi}^{(n,-k)}(A,L)-\mu_{\pi}^o(A,L)\}\\	
			+(\hat\rho_{\pi}^{(n,-k)}(L)-\rho_{\pi}^o(L))\left(\delta^o(A,L)-\hat\delta^{(n,-k)}(A,L)\right)\\\times
			\dfrac{\hat\eta^{(n,-k)}(A,L)-\hat\mu_{\pi}^{(n,-k)}(A,L)}{\hat\kappa_{\pi}^{(n,-k)}(A,L)}
		\end{array}
		\middle|  A\right].
	\end{align*}
	Thus, we can use the formula to deduce that
	\begin{align}
		&R_{5k1}^2
		=\left|\mathbb{E}_{\mathrm{out}}\left[
		\left\{
		\mathbb{E}_{\mathrm{out}}\left[\varphi_{\pi}(O;\hat\alpha_{\pi}^{(n,-k)}, P_O)\middle|A\right]-\theta(A)
		\right\}
		\times\left\{\hat{f}_{\pi,\mathcal{N}}^{(n)}(A)-f_{\mathcal{N}}^*(A)\right\} I\{A\in\mathcal{N}\} \right]\right|^2\notag\\
		\leq&\left|\mathbb{E}_k\left[
		\left\{
		\mathbb{E}_{\mathrm{out}}\left[\varphi_{\pi}(O;\hat\alpha_{\pi}^{(n,-k)}, P_O)\middle|A\right]-\theta(A)
		\right\}^2 I\{A\in\mathcal{N}\} \right]\right|\notag\\
		&\times\mathbb{E}_{\mathrm{out}}\left[
		\left\{\hat{f}_{\pi,\mathcal{N}}^{(n)}(A)-f_{\mathcal{N}}^*(A)\right\}^2 I\{A\in\mathcal{N}\} \right]\notag\\
		\lesssim&
		\left|\mathbb{E}_{\mathrm{out}}\left[
		\mathbb{E}_{\mathrm{out}}\left[\begin{array}{l}
			(\hat\rho_{\pi}^{(n,-k)}(L)-\rho_{\pi}^o(L))^2(\hat\eta^{(n,-k)}(A,L)-\eta^o(A,L))^2\\
			+\left\{\delta^o(A,L)-
			\hat\delta^{(n,-k)}(A,L)\dfrac{\kappa_{\pi}^o(A,L)}{\hat\kappa_{\pi}^{(n,-k)}(A,L)}\right\}^2\\\times
			\{\hat\mu_{\pi}^{(n,-k)}(A,L)-\mu_{\pi}^o(A,L)\}^2\\
			+(\rho_{\pi}^{(n,-k)}(L)-\rho_{\pi}^o(L))^2\\
			\times\left(\delta^o(A,L)-\hat\delta^{(n,-k)}(A,L)\right)^2
		\end{array}\middle| A\right]I\{A\in\mathcal{N}\}
		\right]\right|\notag\\
		&\times\mathbb{E}_{\mathrm{out}}\left[
		\left\{\hat{f}_{\pi,\mathcal{N}}^{(n)}(A)-f_{\mathcal{N}}^*(A)\right\}^2 I\{A\in\mathcal{N}\} \right]\notag\\
		\lesssim&\left|
		\mathbb{E}_{\mathrm{out}}\left[\begin{array}{l}
			(\hat\rho_{\pi}^{(n,-k)}(L)-\rho_{\pi}^o(L))^2(\hat\eta^{(n,-k)}(A,L)-\eta^o(A,L))^2\\
			+\left\{\delta^o(A,L)-
			\hat\delta^{(n,-k)}(A,L)\dfrac{\kappa_{\pi}^o(A,L)}{\hat\kappa_{\pi}^{(n,-k)}(A,L)}\right\}^2\\\times
			\{\mu_{\pi}^{(n,-k)}(A,L)-\mu_{\pi}^o(A,L)\}^2\\
			+(\rho_{\pi}^{(n,-k)}(L)-\rho_{\pi}^o(L))^2\\
			\times\left(\delta^o(A,L)-\hat\delta^{(n,-k)}(A,L)\right)^2
		\end{array}\right]\right|\notag\\
		\lesssim&
		\left\{\begin{array}{l}
			\|(\hat\rho_{\pi}^{(n,-k)}-\rho_{\pi}^o)\times(\hat\eta^{(n,-k)}-\eta^o)\|_{\mathcal{N}}^2
			+\|(\hat\delta^{(n,-k)}-\delta^o)\times(\hat\mu_{\pi}^{(n,-k)}-\mu_{\pi}^o)\|_{\mathcal{N}}^2\notag\\
			+\|(\hat\kappa_{\pi}^{(n,-k)}-\kappa_{\pi}^o)\times(\hat\mu_{\pi}^{(n,-k)}-\mu_{\pi}^o)\|_{\mathcal{N}}^2
			+\|(\hat\rho_{\pi}^{(n,-k)}-\rho_{\pi}^o)\times(\hat\delta^{(n,-k)}-\delta^o)\|_{\mathcal{N}}^2
		\end{array}\right\}\\
		&\times\mathbb{E}_{\mathrm{out}}\left[
		\left\{\hat{f}_{\pi,\mathcal{N}}^{(n)}(A)-f_{\mathcal{N}}^*(A)\right\}^2 I\{A\in\mathcal{N}\} \right].\notag
	\end{align}
	Thus, we now deduce the result that
	\begin{align}
		&R_{5k1}\lesssim
		\mathbb{E}_{\mathrm{out}}\left[
		\left\{\hat{f}_{\pi,\mathcal{N}}^{(n)}(A)-f_{\mathcal{N}}^*(A)\right\}^2 I\{A\in\mathcal{N}\} \right]^{1/2}\notag\\
		&\times\left\{\begin{array}{l}
			\|(\hat\rho_{\pi}^{(n,-k)}-\rho_{\pi}^o)\times(\hat\eta^{(n,-k)}-\eta^o)\|_{\mathcal{N}}^2
			+\|(\hat\delta^{(n,-k)}-\delta^o)\times(\hat\mu_{\pi}^{(n,-k)}-\mu_{\pi}^o)\|_{\mathcal{N}}^2\\
			+\|(\hat\kappa_{\pi}^{(n,-k)}-\kappa_{\pi}^o)\times(\hat\mu_{\pi}^{(n,-k)}-\mu_{\pi}^o)\|_{\mathcal{N}}^2
			+\|(\hat\rho_{\pi}^{(n,-k)}-\rho_{\pi}^o)\times(\hat\delta^{(n,-k)}-\delta^o)\|_{\mathcal{N}}^2
		\end{array}\right\}^{1/2}\notag\\
		&\lesssim \mathbb{E}_{\mathrm{out}}\left[
		\left\{\hat{f}_{\pi,\mathcal{N}}^{(n)}(A)-f_{\mathcal{N}}^*(A)\right\}^2 I\{A\in\mathcal{N}\} \right]^{1/2}
		\times \widetilde{o}_p(d(n)^{1/2})\label{exp: 28}
	\end{align}
	
	For the term $R_{5k2}$, we follow the same logic as that in Equations~\eqref{exp: 10}--\eqref{exp: 14}:
	\begin{align}
		R_{5k2}^2\lesssim&
		\mathbb{E}_{\mathrm{out}}\left[
		\left(\mathbb{E}_{\mathrm{out}}\left[
		\varphi_{\pi}(O;\hat\alpha_{\pi}^{(n,-k)},\hat{P}_O^{(n,-k)})-\varphi_{\pi}(O;\hat\alpha_{\pi}^{(n,-k)}, P_O)
		\middle| A\right]\right)^2I\{A\in\mathcal{N}\}\right]\notag\\
		&\times\mathbb{E}_{\mathrm{out}}\left[
		\left\{\hat{f}_{\pi,\mathcal{N}}^{(n)}(A)-f_{\mathcal{N}}^*(A)\right\}^2 I\{A\in\mathcal{N}\} \right]\notag\\
		\lesssim&
		\mathbb{E}_{\mathrm{out}}\left[
		\left\{\varphi_{\pi}(O;\hat\alpha_{\pi}^{(n,-k)},\hat{P}_O^{(n,-k)})-\varphi_{\pi}(O;\hat\alpha_{\pi}^{(n,-k)}, P_O)\right\}^2
		I\{A\in\mathcal{N}\}\right]\notag\\
		&\times\mathbb{E}_{\mathrm{out}}\left[
		\left\{\hat{f}_{\pi,\mathcal{N}}^{(n)}(A)-f_{\mathcal{N}}^*(A)\right\}^2 I\{A\in\mathcal{N}\} \right].\notag
	\end{align}
	Use the independence structure,
	\begin{align*}
		&\mathbb{E}_{-k}^{(2)}\left[
		\left\{\varphi_{\pi}(O;\hat\alpha_{\pi}^{(n,-k)},\hat{P}_O^{(n,-k)})-\varphi_{\pi}(O;\hat\alpha_{\pi}^{(n,-k)}, P_O)\right\}^2\right]\\
		\lesssim&
		\dfrac{1}{|I_{-k}^{(2)}|}\left\{\int \left(\begin{array}{l}
			\hat\mu_{\pi}^{(n,-k)}(A,l)- (z_{\pi}-\hat\rho_{\pi}^{(n,-k)}(l))\\
			\times\dfrac{\hat\eta^{(n,-k)}(A,l)-\hat\mu_{\pi}^{(n,-k)}(A,l)}
			{\hat\kappa_{\pi}^{(n,-k)}(A,l)}
		\end{array}\right)\mathrm{d} P_O(o)\right\}^2.
	\end{align*}
	Then we use the law of large numbers to deduce that
	\begin{align}
		&\mathbb{E}_{\mathrm{out}}\left[
		\left\{\varphi_{\pi}(O;\hat\alpha_{\pi}^{(n,-k)},\hat{P}_O^{(n,-k)})-\varphi_{\pi}(O;\hat\alpha_{\pi}^{(n,-k)}, P_O)\right\}^2
		I\{A\in\mathcal{N}\}\right]=\widetilde{O}_p(\dfrac{1}{|I_{-k}^{(2)}|}),\text{ and }\notag\\
		&R_{5k2}=
		\mathbb{E}_{\mathrm{out}}\left[
		\left\{\hat{f}_{\pi,\mathcal{N}}^{(n)}(A)-f_{\mathcal{N}}^*(A)\right\}^2 I\{A\in\mathcal{N}\} \right]^{1/2}
		\times \widetilde{O}_p(\dfrac{1}{\sqrt{|I_{-k}^{(2)}|}})\label{exp: 29}
	\end{align}

	Finally, combining Equations~\eqref{exp: 28}--\eqref{exp: 29}, we know that
	\begin{align*}
		&R_{5k}\lesssim R_{5k1}+R_{5k2}\\
		\lesssim& \mathbb{E}_{\mathrm{out}}\left[
		\left\{\hat{f}_{\pi,\mathcal{N}}^{(n)}(A)-f_{\mathcal{N}}^*(A)\right\}^2 I\{A\in\mathcal{N}\} \right]^{1/2}
		\times\left\{ \widetilde{O}_p(\dfrac{1}{\sqrt{|I_{-k}^{(2)}|}})+\widetilde{o}_p(d(n)^{1/2})\right\}.
	\end{align*}
	Overall, from Equations~\eqref{exp: 33}, for any $M>0$, with probability more than $1-c\exp\{-cn\xi_n^2\}$
	\begin{align*}
		&\|\hat{f}_{\pi,\mathcal{N}}^{(n)}-\theta\|_{\mathcal{N}}^2-
		\|f_{\mathcal{N}}^*-\theta\|_{\mathcal{N}}^2
		\lesssim\dfrac{1}{K}\sum_{k=1}^KR_{5k}+R_{6k}\\
		\lesssim&
		\mathbb{E}_{\mathrm{out}}\left[
		\left\{\hat{f}_{\pi,\mathcal{N}}^{(n)}(A)-f_{\mathcal{N}}^*(A)\right\}^2 I\{A\in\mathcal{N}\} \right]^{1/2}
		\times\left\{ \widetilde{O}_p(\dfrac{1}{\sqrt{|I_{-k}^{(2)}|}})+\widetilde{o}_p(d(n)^{1/2})+\xi_n\right\}\\
		\lesssim& \|\hat{f}_{\pi,\mathcal{N}}^{(n)}-f_{\mathcal{N}}^*\|_{\mathcal{N}}
		\times\left\{ \widetilde{O}_p(\dfrac{1}{\sqrt{|I_{-k}^{(2)}|}})+\widetilde{o}_p(d(n)^{1/2})+\xi_n\right\}\\
		\lesssim& \dfrac{1}{M}\|\hat{f}_{\pi,\mathcal{N}}^{(n)}-f_{\mathcal{N}}^*\|_{\mathcal{N}}^2
		+M\left\{ \widetilde{O}_p(\dfrac{1}{|I_{-k}^{(2)}|})+\widetilde{o}_p(d(n))+\xi_n^2\right\}\\
		\lesssim& \dfrac{1}{M}\left\{\|\hat{f}_{\pi,\mathcal{N}}^{(n)}-\theta\|_{\mathcal{N}}^2+\| f_{\mathcal{N}}^*-\theta\|_{\mathcal{N}}^2\right\}
		+M\left\{ \widetilde{O}_p(\dfrac{1}{|I_{-k}^{(2)}|})+\widetilde{o}_p(d(n))+\xi_n^2\right\}.
	\end{align*}
	Equivalently, with probability more than $1-c\exp\{-cn\xi_n^2\}$,
	\begin{align*}
		\|\hat{f}_{\pi,\mathcal{N}}^{(n)}-\theta\|_{\mathcal{N}}^2
		\lesssim
		\|f_{\mathcal{N}}^*-\theta\|_{\mathcal{N}}^2
		+\widetilde{O}_p(\dfrac{1}{|I_{-k}^{(2)}|})+\widetilde{o}_p(d(n))+\xi_n^2.
	\end{align*}
	If we take expectation with respect to $(O_1,\ldots,O_n)$, we can show that
	\begin{align*}
		\mathbb{E}\left[\|\hat{f}_{\pi,\mathcal{N}}^{(n)}-\theta\|_{\mathcal{N}}^2\right]
		\lesssim
		\|f_{\mathcal{N}}^*-\theta\|_{\mathcal{N}}^2
		+O(\dfrac{1}{|I_{-k}^{(2)}|})+o(d(n))+\xi_n^2.
	\end{align*}
	This finishes the proof of Theorem~\ref{thm: oracle bounds for ERM}.

	\subsection{Proof of Theorem~\ref{thm: contrast identification}}
	From Theorem~\ref{thm: identification}, we know that
	\begin{align*}
		&\mathbb{E}[\mu_{\pi_1}^o(a_1,L)] -\mathbb{E}[\mu_{\pi_2}^o(a_2,L)] 
		- (\mathbb{E}[Y(a_1)]-\mathbb{E}[Y(a_2)])\tag*{(Adopting Equation~\eqref{eq: identification E[Y(a)] bias})}\\
		=&\mathbb{E}\left[\mathrm{Cov}\left\{\mathbb{E}[Y(a_1)|U,L], \omega_{a_1,\pi_1}(U,L) \mid L\right\}\right]
		-\mathbb{E}\left[\mathrm{Cov}\left\{\mathbb{E}[Y(a_2)|U,L], \omega_{a_2,\pi_2}(U,L) \mid L\right\}\right]
		\\
		=&\mathbb{E}\left[\mathrm{Cov}\left\{\mathbb{E}[Y(a_1)|U,L]-\mathbb{E}[Y(a_2)|U,L], \omega_{a_1,\pi_1}(U,L) \mid L\right\}\right]\\
		&-\mathbb{E}\left[\mathrm{Cov}\left\{\mathbb{E}[Y(a_2)|U,L], \omega_{a_2,\pi_2}(U,L) -\omega_{a_1,\pi_1}(U,L)\mid L\right\}\right]\\
		=&\mathbb{E}\left[\mathrm{Cov}\left\{\mathbb{E}[Y(a_1)-Y(a_2)|U,L], \omega_{a_1,\pi_1}(U,L) \mid L\right\}\right].
	\end{align*}
	The final equality follows from Proposition~\ref{prop: equivalent characterization of WAIV}, which implies that if $Z$ is a WAIV for $A=a$, then $\omega_{a,\pi}(U,L)$ does not depend on $a$ or $\pi$.  
	
	Moreover, if $Z$ is an AIV for $A$, Proposition~\ref{prop: AIV characterization} gives $\omega_{a_1,\pi_1}(U,L)\equiv 1$, so that
	\[
	\mathbb{E}\left[\mathrm{Cov}\left\{\mathbb{E}[Y(a_1)-Y(a_2)\mid U,L], \omega_{a_1,\pi_1}(U,L) \mid L\right\}\right] = 0.
	\]
	Alternatively, this equality also follows if $\mathbb{E}[Y(a_1)-Y(a_2)\mid U,L]$ does not depend on $U$. This completes the proof of Theorem~\ref{thm: contrast identification}.
	\section{Proof of lemmas}\label{append: proof of lemmas}

	\subsection{Proof of Lemma~\ref{lem: lower bounded eigenvalue}}
	First, since $K_h(A_i-a)=K((A_i-a)/h)/h\lesssim I\{|A_i-a|\leq h\}/h$ by Assumption~\ref{as: regular condition}\ref{cd: kernel condition}, we calculate that
	\begin{align*}
		&\mathrm{Var}\left[
		\dfrac{1}{n}\sum_{i=1}^nK_h(A_i-a)(\dfrac{A_i-a}{h})^j
		\right]
		=\dfrac{1}{n}\mathrm{Var}\left[
		K_h(A-a)(\dfrac{A-a}{h})^j
		\right]\\
		\lesssim&
		\dfrac{1}{nh^2}\mathbb{E}\left[I\{|A-a|\leq h\}(\dfrac{A-a}{h})^{2j}\right]\lesssim 
		\dfrac{1}{nh^2}\mathbb{E}\left[I\{|A-a|\leq h\}\right]=O(\dfrac{1}{nh}).
	\end{align*}
	The final step follows by $p_A(a)>0$ since $a\in\mathring{\mathcal{A}}$.
	From the law of large numbers, this indicates that
	\begin{align*}
		\dfrac{1}{n}\sum_{i=1}^nK_h(A_i-a)g\left(\dfrac{A_i-a}{h}\right) g\left(\dfrac{A_i-a}{h}\right)^{\top}-
		\mathbb{E}\left[K_h(A-a)g\left(\dfrac{A-a}{h}\right)g\left(\dfrac{A-a}{h}\right)^{\top}\right] = O_p(\dfrac{1}{\sqrt{nh}}).
	\end{align*}
	Next, we can deduce that
	\begin{align}
		\dfrac{1}{n}Q_{h}(a)^{-1}
		=&\dfrac{1}{n}\sum_{i=1}^nK_h(A_i-a)g\left(\dfrac{A_i-a}{h}\right) g\left(\dfrac{A_i-a}{h}\right)^{\top}\notag\\
		=&\mathbb{E}\left[K_h(A-a)g\left(\dfrac{A-a}{h}\right) g\left(\dfrac{A-a}{h}\right)^{\top}\right]+ \widetilde{o}_p(1)\notag\\
		=&\int K(s)g(s)g(s)^{\top} p_A(a+hs)\mathrm{d} s+ \widetilde{o}_p(1)\notag\\
		=&\int K(s)g(s)g(s)^{\top} \mathrm{d} s\times p_A(a)+ \widetilde{o}_p(1)\label{exp: 24}\\
		=&\left[\begin{array}{ll}
			\int K(s)\mathrm{d} s,&\int K(s)s \mathrm{d} s\\
			\int K(s)s \mathrm{d} s,&\int K(s)s^2 \mathrm{d} s
		\end{array}\right]\times p_A(a)+ \widetilde{o}_p(1)\notag\\
		=&\left[\begin{array}{ll}
			1,&0\\
			0,&\int K(s)s^2 \mathrm{d} s
		\end{array}\right]\times p_A(a)+ \widetilde{o}_p(1)\label{exp: 25}
	\end{align}	
	Equation~\eqref{exp: 24} follows from the continuity of the density $p_{A}(a)$ at $a$ by Assumption~\ref{as: continuity}.
	Equation~\eqref{exp: 25} follows from Assumption~\ref{as: regular condition}\ref{cd: kernel condition}.
	
	Therefore, since 
	$\int K(s)s^2 \mathrm{d} s>0$ and 
	$p_A(a)>0$, we know that
	\begin{align*}
		nQ_{h}(a)=\left[\begin{array}{cc}
			1,&0\\
			0,&\dfrac{1}{\int K(s)s^2 \mathrm{d} s}
		\end{array}\right]\times \dfrac{1}{p_A(a)}+ \widetilde{o}_p(1).
	\end{align*}
	This finishes the proof of Lemma~\ref{lem: lower bounded eigenvalue}.

	\subsection{Proof of Lemma~\ref{lem: LLKR}}
	The proof of this lemma is basic in LLKR. We refer to \citet{wasserman2006all,tsybakov2009introduction} for a proof of this lemma.

	\subsection{Proof of Lemma~\ref{lem: mixed bias}}
	Recall that
	\begin{align*}
		\varphi_{\pi}(O;\alpha_{\pi},P_O) := &
		\delta(A,L)\times
		\dfrac{(Z_{\pi}- \rho_{\pi}(L))(Y-\mu_{\pi}(A,L))}{\kappa_{\pi}(A,L)}\notag\\
		&+\int \mu_{\pi}(A,l) - (z_{\pi}-\rho_{\pi}(l))\dfrac{\eta(A,l)-\mu_{\pi}(A,l)}{\kappa_{\pi}(A,l)}\mathrm{d} P_{O}(o)\\
		=&\delta(A,L)\times
		\dfrac{(Z_{\pi}- \rho_{\pi}(L))(Y-\mu_{\pi}(A,L))}{\kappa_{\pi}(A,L)}\notag\\
		&+\int \mu_{\pi}(A,l) - (z_{\pi}-\rho_{\pi}(l))\dfrac{\eta(A,l)-\mu_{\pi}(A,l)}{\kappa_{\pi}(A,l)}\mathrm{d} P_{Z,L}(z,l).
	\end{align*}
	We calculate that
	\begin{align*}
		&\mathbb{E}[\varphi_{\pi}(O)\mid A=a] - \mathbb{E}[\mu_{\pi}^o(a,L)]\\
		=&\mathbb{E}\left[
		\dfrac{\delta(A,L)}{\kappa_{\pi}(A,L)}(Z_{\pi}- \rho_{\pi}(L))(Y-\mu_{\pi}(A,L)) \middle|  A=a
		\right]\\
		&+\mathbb{E}\left[\int \mu_{\pi}(A,l)\mathrm{d} P_O(o)-\int (z_{\pi}-\rho_{\pi}(l))\dfrac{\eta(A,l)-\mu_{\pi}(A,l)}{\kappa_{\pi}(A,l)}\mathrm{d} P_{Z,L}(z,l)\middle|A=a\right]\\&- \int \mu_{\pi}^o(a,l)\mathrm{d} P_O(o)\\
		=&\mathbb{E}\left[
		\dfrac{\delta(a,L)}{\kappa_{\pi}(a,L)}(Z_{\pi}- \rho_{\pi}(L))(Y-\mu_{\pi}(a,L)) \middle|  A=a
		\right]\\
		&+\int \left\{\mu_{\pi}(a,l)-\mu_{\pi}^o(a,l)\right\}\mathrm{d} P_O(o)
		-\int (z_{\pi}-\rho_{\pi}(l))\dfrac{\eta(a,l)-\mu_{\pi}(a,l)}{\kappa_{\pi}(a,l)}\mathrm{d} P_{Z,L}(z,l)\\
		=&\mathbb{E}\left[
		\dfrac{\delta(a,L)}{\kappa_{\pi}(a,L)} \left\{\begin{array}{l}
			\eta^o(A,L)\rho_{\pi}^o(L)+\kappa_{\pi}^o(A,L)\mu_{\pi}^o(A,L)\\
			- (\rho_{\pi}^o(L)+\kappa_{\pi}^o(A,L))\mu_{\pi}(A,L)\\
			-\rho_{\pi}(L)\eta^o(A,L)+\rho_{\pi}(L)\mu_{\pi}(a,L)
		\end{array}\right\}
		\middle|  A=a\right]\\
		&+\mathbb{E}\left[\mu_{\pi}(a,L)-\mu_{\pi}^o(a,L)\right]
		-\mathbb{E}\left[
		(Z_{\pi}-\rho_{\pi}(L))\dfrac{\eta(a,L)-\mu_{\pi}(a,L)}{\kappa_{\pi}(a ,L)}
		\right]\\
		=&\mathbb{E}\left[
		\dfrac{\delta(a,L)}{\kappa_{\pi}(a,L)} \left\{\begin{array}{l}
			\eta^o(a,L)(\rho_{\pi}^o(L)-\rho_{\pi}(L))\\
			+\kappa_{\pi}^o(a,L)\{\mu_{\pi}^o(a,L)-\mu_{\pi}(a,L)\}\\
			+(\rho_{\pi}(L)-\rho_{\pi}^o(L))\mu_{\pi}(a,L)
		\end{array}\right\}
		\middle|  A=a\right]\\
		&+\mathbb{E}\left[ \delta^o(a,L)\{\mu_{\pi}(a,L)-\mu_{\pi}^o(a,L)\}\middle| A=a\right]\\
		&-\mathbb{E}\left[
		\dfrac{p_{Z,L}(Z,L)}{p_{Z,L|A}(Z,L\mid A)}
		(Z_{\pi}-\rho_{\pi}(L))\dfrac{\eta(A,L)-\mu_{\pi}(A,L)}{\kappa_{\pi}(A,L)}
		\middle| A=a\right]\\
		=&\mathbb{E}\left[
		\dfrac{\delta(a,L)}{\kappa_{\pi}(a,L)} \left\{\begin{array}{l}
			\eta^o(a,L)(\rho_{\pi}^o(L)-\rho_{\pi}(L))\\
			+(\rho_{\pi}(L)-\rho_{\pi}^o(L))\mu_{\pi}(a,L)
		\end{array}\right\}
		\middle|  A=a\right]\\
		&+\mathbb{E}\left[ \left\{\delta^o(a,L)-
		\delta(a,L)\dfrac{\kappa_{\pi}^o(a,L)}{\kappa_{\pi}(a,L)}\right\}
		\{\mu_{\pi}(a,L)-\mu_{\pi}^o(a,L)\}\middle| A=a\right]\\
		&-\mathbb{E}\left[
		\mathbb{E}\left[\dfrac{p_{Z,L}(Z,L)}{p_{Z,L|A}(Z,L\mid A)}
		(Z_{\pi}-\rho_{\pi}(L))\middle|A,L\right]\dfrac{\eta(a,L)-\mu_{\pi}(a,L)}{\kappa_{\pi}(a,L)}
		\middle| A=a\right]\\
		=&\mathbb{E}\left[
		\dfrac{\delta(a,L)}{\kappa_{\pi}(a,L)} 
		(\rho_{\pi}(L)-\rho_{\pi}^o(L))(\mu_{\pi}(a,L)-\eta^o(a,L))
		\middle|  A=a\right]\\
		&+\mathbb{E}\left[ \left\{\delta^o(a,L)-
		\delta(a,L)\dfrac{\kappa_{\pi}^o(a,L)}{\kappa_{\pi}(a,L)}\right\}
		\{\mu_{\pi}(a,L)-\mu_{\pi}^o(a,L)\}\middle| A=a\right]\\
		&-\mathbb{E}\left[
		\mathbb{E}\left[
		\dfrac{p_{Z,L}(Z,L)}{p_{Z,L|A}(Z,L\mid A)}
		\dfrac{p_{A|L}(A|L)}{p_{A}(A)}(Z_{\pi}-\rho_{\pi}(L))
		\middle|A,L\right]
		\dfrac{p_{A}(a)}{p_{A|L}(a|L)}
		\dfrac{\eta(a,L)-\mu_{\pi}(a,L)}{\kappa_{\pi}(a,L)}
		\middle| A=a\right]\\
		=&\mathbb{E}\left[
		\dfrac{\delta(a,L)}{\kappa_{\pi}(a,L)} 
		(\rho_{\pi}(L)-\rho_{\pi}^o(L))(\mu_{\pi}(a,L)-\eta^o(a,L))
		\middle|  A=a\right]\\
		&+\mathbb{E}\left[ \left\{\delta^o(a,L)-
		\delta(a,L)\dfrac{\kappa_{\pi}^o(a,L)}{\kappa_{\pi}(a,L)}\right\}
		\{\mu_{\pi}(a,L)-\mu_{\pi}^o(a,L)\}\middle| A=a\right]\\
		&-\mathbb{E}\left[
		\begin{array}{c}
			\mathbb{E}\left[\dfrac{p_{Z,L}(Z,L)}{p_{Z,L|A}(Z,L\mid A)}
			\dfrac{p_{A|L}(A|L)}{p_{A}(A)}
			\dfrac{p_{A|Z,L}(A|Z,L)}{p_{A|L}(A|L)}(Z_{\pi}-\rho_{\pi}(L))
			\middle|L\right]\\
			\times\dfrac{p_{A}(a)}{p_{A|L}(a|L)}\dfrac{\eta(a,L)-\mu_{\pi}(a,L)}{\kappa_{\pi}(a,L)}
		\end{array}\middle| A=a\right]\\
		=&\mathbb{E}\left[
		\dfrac{\delta(a,L)}{\kappa_{\pi}(a,L)} 
		(\rho_{\pi}(L)-\rho_{\pi}^o(L))(\mu_{\pi}(a,L)-\eta^o(a,L))
		\middle|  A=a\right]\\
		&+\mathbb{E}\left[ \left\{\delta^o(a,L)-
		\delta(a,L)\dfrac{\kappa_{\pi}^o(a,L)}{\kappa_{\pi}(a,L)}\right\}
		\{\mu_{\pi}(a,L)-\mu_{\pi}^o(a,L)\}\middle| A=a\right]\\
		&+\mathbb{E}\left[
		(\rho_{\pi}(L)-\rho_{\pi}^o(L))
		\times\delta^o(a,L)\dfrac{\eta(a,L)-\mu_{\pi}(a,L)}{\kappa_{\pi}(a,L)}
		\middle| A=a\right]\\
		=&\mathbb{E}\left[
		\dfrac{\delta(a,L)}{\kappa_{\pi}(a,L)} 
		(\rho_{\pi}(L)-\rho_{\pi}^o(L))(\eta(a,L)-\eta^o(a,L))
		\middle|  A=a\right]\\
		&+\mathbb{E}\left[ \left\{\delta^o(a,L)-
		\delta(a,L)\dfrac{\kappa_{\pi}^o(a,L)}{\kappa_{\pi}(a,L)}\right\}
		\{\mu_{\pi}(a,L)-\mu_{\pi}^o(a,L)\}\middle| A=a\right]\\
		&+\mathbb{E}\left[
		(\rho_{\pi}(L)-\rho_{\pi}^o(L))
		\times\left(\delta^o(a,L)-\delta(a,L)\right)\dfrac{\eta(a,L)-\mu_{\pi}(a,L)}{\kappa_{\pi}(a,L)}
		\middle| A=a\right].
	\end{align*}

	\subsection{Proof of Lemma~\ref{lem: mixed bias property}}
	From Lemma~\ref{lem: mixed bias}, we can directly calculate
	\begin{align*}
		&\mathbb{E}_{k}\left[K_h(A-a)f(\dfrac{A-a}{h})\left\{\varphi_{\pi}(O;\hat\alpha_{\pi}^{(n,-k)},P_O) - \varphi_{\pi}(O;\alpha_{\pi}^o, P_O)\right\}\right]\\
		=&\mathbb{E}_{k}\left[K_h(A-a)f(\dfrac{A-a}{h})\mathbb{E}_k\left[\varphi_{\pi}(O;\hat\alpha_{\pi}^{(n,-k)},P_O) - \varphi_{\pi}(O;\alpha_{\pi}^o, P_O)\middle| A\right]\right]\\
		=&\mathbb{E}_{k}\left[\left(\begin{array}{l}
			K_h(A-a)f(\dfrac{A-a}{h})
			\hat\delta^{(n,-k)}(A,L)\\
			\times\dfrac{1}{\hat\kappa_{\pi}^{(n,-k)}(A,L)} 
			(\hat\rho_{\pi}^{(n,-k)}(L)-\rho_{\pi}^o(L))(\hat\eta^{(n,-k)}(A,L)-\eta^o(A,L))
		\end{array}\right)
		\right]\\
		&+\mathbb{E}_{k}\left[\left(\begin{array}{l}
			K_h(A-a)f(\dfrac{A-a}{h})\{\hat\mu_{\pi}^{(n,-k)}(A,L)-\mu_{\pi}^o(A,L)\}\\
			\times\left\{\delta^o(A,L)-
			\hat\delta^{(n,-k)}(A,L)\dfrac{\kappa_{\pi}^o(A,L)}{\hat\kappa_{\pi}^{(n,-k)}(A,L)}\right\}
		\end{array}\right)\right]\\
		&+\mathbb{E}_{k}\left[\begin{array}{l}
			K_h(A-a)f(\dfrac{A-a}{h})
			\left( \dfrac{\hat\eta^{(n,-k)}(A,L)-\hat\mu_{\pi}^{(n,-k)}(A,L)}{\hat\kappa_{\pi}^{(n,-k)}(A,L)}\right)\\
			(\hat\rho_{\pi}^{(n,-k)}(L)-\rho_{\pi}^o(L))
			\times\left(\delta^o(A,L)-\hat\delta^{(n,-k)}(A,L)\right)
		\end{array}\right]\\
		\lesssim&\mathbb{E}_{k}\left[K_h(A-a)|f(\dfrac{A-a}{h})|
		\left|(\hat\rho_{\pi}^{(n,-k)}(L)-\rho_{\pi}^o(L))(\hat\eta^{(n,-k)}(A,L)-\eta^o(A,L))\right|
		\right]\\
		&+\mathbb{E}_{k}\left[\begin{array}{l}
			K_h(A-a)f(\dfrac{A-a}{h})\{\hat\mu_{\pi}^{(n,-k)}(A,L)-\mu_{\pi}^o(A,L)\}\\
			\times\left\{\delta^o(A,L)-
			\hat\delta^{(n,-k)}(A,L)\right\}
		\end{array}\right]\\
		&+\mathbb{E}_{k}\left[\begin{array}{l}
			K_h(A-a)f(\dfrac{A-a}{h})\{\hat\mu_{\pi}^{(n,-k)}(A,L)-\mu_{\pi}^o(A,L)\}\\
			\times\hat\delta^{(n,-k)}(A,L)\left\{1-
			\dfrac{\kappa_{\pi}^o(A,L)}{\hat\kappa_{\pi}^{(n,-k)}(A,L)}\right\}
		\end{array}\right]\\
		&+\mathbb{E}_{k}\left[\begin{array}{l}
			K_h(A-a)f(\dfrac{A-a}{h})
			\left| \dfrac{\hat\eta^{(n,-k)}(A,L)-\hat\mu_{\pi}^{(n,-k)}(A,L)}{\hat\kappa_{\pi}^{(n,-k)}(A,L)}\right|\\
			(\hat\rho_{\pi}^{(n,-k)}(L)-\rho_{\pi}^o(L))
			\times\left|\delta^o(A,L)-\hat\delta^{(n,-k)}(A,L)\right|
		\end{array}\right]\\
		\lesssim&\dfrac{1}{h}\mathbb{E}_{k}\left[I\{|A-a|\leq h\}
		\left|(\hat\rho_{\pi}^{(n,-k)}(L)-\rho_{\pi}^o(L))(\hat\eta^{(n,-k)}(A,L)-\eta^o(A,L))\right|
		\right]\\
		&+\dfrac{1}{h}\mathbb{E}_{k}\left[\begin{array}{l}
			I\{|A-a|\leq h\}\left|\hat\mu_{\pi}^{(n,-k)}(A,L)-\mu_{\pi}^o(A,L)\right|\\
			\times\left|\left\{\delta^o(A,L)-
			\hat\delta^{(n,-k)}(A,L)\right\}
			\right|
		\end{array}\right]\\
		&+\dfrac{1}{h}\mathbb{E}_{k}\left[I\{|A-a|\leq h\}
		\left|\left\{\hat\kappa_{\pi}^{(n,-k)}(A,L)-\kappa_{\pi}^o(A,L)\right\}
		\{\hat\mu_{\pi}^{(n,-k)}(A,L)-\mu_{\pi}^o(A,L)\}\right|\right]\\
		&+\dfrac{1}{h}\mathbb{E}_{k}\left[\begin{array}{l}
			I\{|A-a|\leq h\}
			\left| \dfrac{\hat\eta^{(n,-k)}(A,L)-\hat\mu_{\pi}^{(n,-k)}(A,L)}{\hat\kappa_{\pi}^{(n,-k)}(A,L)}\right|\\
			\left|\hat\rho_{\pi}^{(n,-k)}(L)-\rho_{\pi}^o(L)\right|
			\times\left|\delta^o(A,L)-\hat\delta^{(n,-k)}(A,L)\right|
		\end{array}\right]\\
		\lesssim&\dfrac{1}{h}\mathbb{E}_{k}\left[I\{|A-a|\leq h\}\left|\hat\rho_{\pi}^{(n,-k)}(L)-\rho_{\pi}^o(L)\right|
		\sup_{A\in\mathcal{B}(a,h)}\left|(\hat\eta^{(n,-k)}(A,L)-\eta^o(A,L))\right|
		\right]\\
		&+\dfrac{1}{h}\mathbb{E}_{k}\left[I\{|A-a|\leq h\}\left\{\begin{array}{l}
			\sup_{A\in\mathcal{B}(a,h)}\left|\left\{\delta^o(A,L)-
			\hat\delta^{(n,-k)}(A,L)\right\}\right|\\
			\times\sup_{A\in\mathcal{B}(a,h)}\left|\{\hat\mu_{\pi}^{(n,-k)}(A,L)-\mu_{\pi}^o(A,L)\}\right|
		\end{array}\right\}\right]\\
		&+\dfrac{1}{h}\mathbb{E}_{k}\left[I\{|A-a|\leq h\}\left\{\begin{array}{l}
			\sup_{A\in\mathcal{B}(a,h)}\left|\left\{\hat\kappa_{\pi}^{(n,-k)}(A,L)-\kappa_{\pi}^o(A,L)\right\}\right|\\
			\times\sup_{A\in\mathcal{B}(a,h)}\left|\{\hat\mu_{\pi}^{(n,-k)}(A,L)-\mu_{\pi}^o(A,L)\}\right|
		\end{array}\right\}\right]\\
		&+\dfrac{1}{h}\mathbb{E}_{k}\left[I\{|A-a|\leq h\}\left\{\begin{array}{l}
			\sup_{A\in\mathcal{B}(a,h)}\left|\delta^o(A,L)-\hat\delta^{(n,-k)}(A,L)\right| \left|\hat\rho_{\pi}^{(n,-k)}(L)-\rho_{\pi}^o(L)\right|
		\end{array}\right\}\right]\\
		\lesssim&\mathbb{E}_{k}\left[\dfrac{\Pr\left(|A-a|\leq h| L\right)}{h}\left\{\begin{array}{l}
			\left|\hat\rho_{\pi}^{(n,-k)}(L)-\rho_{\pi}^o(L)\right|
			\sup_{A\in\mathcal{B}(a,h)}\left|(\hat\eta^{(n,-k)}(A,L)-\eta^o(A,L))\right|\\
			+\sup_{A\in\mathcal{B}(a,h)}\left|\left(\delta^o(A,L)-
			\hat\delta^{(n,-k)}(A,L)\right)\right|\\
			\times\sup_{A\in\mathcal{B}(a,h)}\left(\begin{array}{l}
				\left|\{\hat\mu_{\pi}^{(n,-k)}(A,L)-\mu_{\pi}^o(A,L)\}\right|\\
				+\left|\{\hat\rho_{\pi}^{(n,-k)}(L)-\rho_{\pi}^o(L)\}\right|
			\end{array}\right)\\
			+\sup_{A\in\mathcal{B}(a,h)}\left|\left\{\hat\kappa_{\pi}^{(n,-k)}(A,L)-\kappa_{\pi}^o(A,L)\right\}\right|\\
			\times\sup_{A\in\mathcal{B}(a,h)}\left|\{\hat\mu_{\pi}^{(n,-k)}(A,L)-\mu_{\pi}^o(A,L)\}\right|
		\end{array}\right\}\right]\\
		\lesssim&
		\left\|\hat\rho_{\pi}^{(n,-k)}-\rho_{\pi}^o\right\|_2\times
		\left\|\hat\eta^{(n,-k)}-\eta^o\right\|_{\mathcal{B}(a,h),2}
		+\left\|\hat\mu_{\pi}^{(n,-k)}-\mu_{\pi}^o\right\|_{\mathcal{B}(a,h),2}
		\times
		\left\|\hat\kappa_{\pi}^{(n,-k)}-\kappa_{\pi}^o\right\|_{\mathcal{B}(a,h),2}\\
		&+\left\{\left\|\hat\mu_{\pi}^{(n,-k)}-\mu_{\pi}^o\right\|_{\mathcal{B}(a,h),2}
		+\left\|\hat\rho_{\pi}^{(n,-k)}-\rho_{\pi}^o\right\|_{2}\right\}
		\left\|\hat\delta^{(n,-k)}
		-\delta^o\right\|_{\mathcal{B}(a,h),2}.
	\end{align*}
	% \bibliographystyle{plainnat}
	% \bibliography{bibfile_supp}
	{\setlength{\baselineskip}{0.75\baselineskip}
		\putbib[bibfile_main]
	}
\end{bibunit}

\end{document}